# INTUITIONISTIC FIRST-ORDER LOGIC
# CATEGORICAL SEMANTICS VIA THE CURRY-HOWARD ISOMORPHISM

MARCO BENINI

ABSTRACT. In this technical report, the achievements obtained by the author during his sabbatical period in the Department of Pure Mathematics, University of Leeds, UK, under the Marie Curie Intra-European Fellowship 'Predicative Theories and Grothendieck Toposes' are documented.

Synthetically, this reports introduces a novel sound and complete semantics for first-order intuitionistic logic, in the framework of category theory and by the computational interpretation of the logic based on the so-called Curry-Howard isomorphism. Aside, a sound and complete semantics for the corresponding $\lambda$-calculus is derived, too.

This semantics extends, in a way, the more traditional meanings given by Heyting categories, the topos-theoretic interpretation, and Kripke models, in the sense that a model in these interpretation can always be transformed in a model in the proposed semantics. The vice versa is still an open question, as it is certainly true that to every model in the proposed system corresponds at least a model in those semantics, but there is no evidence how to construct it.

The feature which justifies the introduction of this novel semantics is the fact that it is 'point-free', i.e., there is no universe whose elements are used to interpret the logical terms. In other words, terms do not denote individuals of some collection but, instead, they denote the 'glue' which keeps together the interpretations of statements, similarly to what happens in formal topology. This fact has a number of not yet fully understood consequences in the philosophy of mathematics.

Since the proposed semantics can be trivially extended to all the first-order logical theories based on the intuitionistic system (and, with some care, to minimal systems as well), the semantics covers also all the predicative theories, even if some peculiar aspects of these theories should be remarked. Moreover, since Grothendieck toposes are special cases of elementary toposes, and there is an embedding of topos-based models into the proposed semantics, the illustrated approach is more general than the initial objective of the research project and the obtained results have a wider range of applications. Nevertheless, some attention has been devoted to the framework of Grothendieck toposes and predicative theories, as this restricted environment allows to profitably use some deep mathematical tools which cannot be extended to the general setting.

*Date*: June 29, 2013.

2010 *Mathematics Subject Classification.* Primary: 03F55; Secondary: 03B20, 03B40, 03B55, 03C90, 03G30.

The author has been supported by a Marie Curie Intra European Fellowship, grant n. PIEF-GA-2010-271926, *Predicative Theories and Grothendieck Toposes*, within the 7th European Community Framework Programme.

## Contents

## 1. INTRODUCTION

This paper illustrates an alternative interpretation for first-order intuitionistic predicate logic in a multi-sorted environment. Also, this interpretation is naturally related with the corresponding $\lambda$-calculus and gives a meaning to the computational interpretation of logical proofs following the proposition-as-types paradigm expressed by the Curry-Howard isomorphism.

Differently from any other interpretation for the logical systems under examination, we will use category theory to provide a class of models where the interpretation of terms does not require to assume the existence of some universe of elements. Instead, similarly to points in formal topology, see [Sam03], terms are interpreted as the 'glue' keeping together the interpretation of formulae and proofs in a suitable categorical structure. We will refer to this peculiar aspect of the interpretation as being *point-free*.

Ideally, the paper is divided into four parts: the first part illustrates the proposed semantics, introducing the notion of 'logically distributive category' as an extension of distributive categories which takes care of handling quantification. The fundamental properties of the semantics, i.e., being sound and complete for the intended logical systems and the corresponding $\lambda$-calculi, are proved.

It is worth remarking that we consider the whole collection of theories which can be axiomatised in the language of first-order logic, thus, for example, the results for classical logic can be derived by considering the theory with the axiom coding the Excluded Middle principle.

A peculiar aspect of our presentation is that equality has no special role: this is done on purpose as it is possible to specialise the obtained results to equip equality with the standard interpretation as the diagonal along the lines illustrated, e.g., in [Gol06]; avoiding to treat equalities as special atomic formulae allows to use more general models which is well-known to be non-equivalent to standard models, see, for example, the difference between Kripke models and intensional Kripke models; at last, the proposed approach allows to explore different possibilities, for example coding in the $\lambda$-calculus an equality rule similar to the Univalence Axiom, thus linking this approach to the front line research in type theory.

The second part compares the semantics based on logically distributive categories with the standard ones, namely, Heyting categories, the internal logic of elementary toposes, and Kripke models. In all these cases, a functor mapping a model in these systems to a 'canonical' model in a logically distributive category is built. Some discussions on the possibility to 'invert' these functors is pursued, even if there is not yet a definitive result. As far as we know by now, these functors have no left or right adjoints but only a right Kan extension, which, of course, introduces a number of 'ideal' objects to form a canonical 'inverse' in the presheaves topos of the domain, a situation which has a number of not yet fully understood consequences.

The third part of the paper illustrates some specialisations that are worth considering when we limit ourselves to predicative theories and Grothendieck toposes, the initial framework of the research project. It will turn out that predicativity has a direct computational meaning in the framework of logically distributive categories which is not revealed in the more traditional approaches.

The fourth and last part discusses some aspects of the philosophical consequences of having a sound and complete point-free semantics. Although the outcomes of this section are still not conclusive, they show how a shift in the traditional view of mathematics is technically possible, and, in this light, an alternative interpretation of Hilbert's finitistic programme can be achieved, a development which is already planned as a future philosophical investigation born from these notes.

As the reader can see, the scope of this paper is wide and the achievements are spread across many aspects of contemporary logic. Evidently not every aspect is developed in full depth and many parts are still 'work-in-progress' the author will eventually pursue in the near future. Here and there, hints for links with related subjects are given, even if space and time has not permitted to draw them beyond the state of a sketchy intuition.

Section 2 illustrates the syntax of the logical systems we want to examine, i.e., intuitionistic first-order multi-sorted predicate logic eventually extended with arbitrary axioms in the same language. Similarly, section 3 introduces the $\lambda$-calculi constituting the computational counterpart of the the logical systems. It will be evident from the definitions that there is a natural bijective correspondence between the logical systems and the $\lambda$-calculi which is referred to as the Curry-Howard isomorphism.

Section 4 describes the categorical semantics we want to use to give a meaning to the $\lambda$-calculi and to the logical systems. The section defines such a semantics, proving its basic properties. The study of the categorical semantics continues in section 5 which proves a general soundness theorem for the $\lambda$-calculi with respect to the defined semantics, while section 6 proves a general completeness theorem for the $\lambda$-calculi. The section concludes with a soundness and completeness result for the logical part as well.

Section 7 recalls the semantics for the considered logical systems which is based on Heyting categories, and section 8 relates this interpretation with the semantics we have introduced in section 4.

Section 9 recalls the semantics of first-order intuitionistic logic inside a topos, and section 10 shows how this semantics is related to logically distributive categories.

Section 11 illustrates the standard Kripke semantics for the considered logic inside the framework of category theory, in particular considering the model in the topoi of presheaves on partially ordered sets. It illustrates the strict relation between these models and the ones based on logically distributive categories.

Section 12 specialises the previous results when considering Grothendieck toposes, i.e., the categories of sheaves on a site.

Section 13 discusses predicative and impredicative theories in the context of logically distributive categories.

Section 14 discusses some philosophical aspects of the obtained results about the 'point-free' nature of the semantics. In particular, this aspect is related to Hilbert's programme and it ends suggesting a different reading of the program which entails a vision of (first-order) mathematics as the formal science studying 'computationally feasible' objects, which is the correct definition of Hilbert's 'concrete objects' in the light of our outcomes, and, indirectly, 'ideal' objects, which have no 'existence' apart allowing interpretation in a imposed rich environment.

## 2. First-order intuitionistic logic

This section illustrates the syntax of the logical systems we will use. This is a variant of multi-sorted, first-order intuitionistic logic. It is not the purpose of this section to explain the properties of this logic, nor the reason why this system is interesting. The reader may refer to any good textbook, e.g., [TS00].

**Definition 2.1** (Logical signature)**.** A *logical signature* $\Sigma = \langle S, F, R\rangle$ is a structure where

(1) $S$ is a set of *sort symbols*;
(2) $F$ is a set of *function symbols*, each one decorated as $f : s_1 \times \cdots \times s_n \to s_0$ where, for every $0 \le i \le n$, $s_i \in S$;
(3) $R$ is a set of *relation symbols*, each one decorated as $r : s_1 \times \cdots \times s_n$ where, for every $1 \le i \le n$, $s_i \in S$.

We assume the standard hypotheses on the language, i.e., all the sets are disjoint and each symbol has a unique decoration. Moreover, when an effective language is required, we assume that all sets are recursive.

The purpose of the hypotheses on the language is to make it non-ambiguous. Essentially, every symbol must be distinct from any other and from any possible construction so that any expression can be understood in a unique way. An effective language is a language which can be recognised by a computer program. Following the standard practice, we will be a bit more relaxed, as the standard language of mathematics can always be reduced to our constraints, even if the usual notation does not strictly adhere to them.

**Definition 2.2** (Logical term)**.** Let $\Sigma = \langle S, F, R\rangle$ be a (logical) signature. For each sort $s \in S$, we assume there is a denumerable set $V_s$ of *variables* with the standard hypotheses on the language. A *logical term on the signature* $\Sigma$ is inductively defined, along with its *free variables*, as:

(1) if $v \in V_s$ then $v : s \in \mathrm{LTerms}(\Sigma)$ and $\mathrm{FV}(v : s) = \{v : s\}$;
(2) if $f : s_1 \times \cdots \times s_n \to s_0 \in F$ and $t_1 : s_1, \ldots, t_n : s_n \in \mathrm{LTerms}(\Sigma)$ then $f(t_1, \ldots, t_n) : s_0 \in \mathrm{LTerms}(\Sigma)$ and $\mathrm{FV}(f(t_1, \ldots, t_n) : s_0) = \bigcup_{i=1}^{n} \mathrm{FV}(t_i : s_i)$.

The special case where $n = 0$ in the second clause gives *logical constants*.

**Definition 2.3** (Substitution on logical terms)**.** Fixed a signature $\Sigma = \langle S, F, R\rangle$, any function $\sigma \colon \bigcup_{s\in S} V_s \to \mathrm{LTerms}(\Sigma)$ is called a *substitution of (logical) variables*.

It naturally extends to a map $\sigma' \colon \mathrm{LTerms}(\Sigma) \to \mathrm{LTerms}(\Sigma)$, called *substitution of (logical) terms* as follows:

(1) if $v \in V_s$ then $\sigma'(v : s) = \sigma(v)$;
(2) if $f : s_1 \times \cdots \times s_n \to s_0 \in F$ then $\sigma'(f(t_1, \ldots, t_n) : s_0) = f(\sigma'(t_1 : s_1), \ldots, \sigma'(t_n : s_n)) : s_0$.

When $\sigma$ is the identity map ($x \in V_s \mapsto x : s$) for every variable except a finite subset $\{x_1, \ldots, x_n\}$, we write $t : s[r_1/x_1, \ldots, r_n/x_n]$ for $\sigma'(t : s)$ where, for each $1 \le i \le n$, $\sigma(x_i) = r_i : s_i$ and $\sigma'$ is the map extending $\sigma$ as above.

**Definition 2.4** (Logical formula)**.** Given a signature $\Sigma = \langle S, F, R\rangle$, a *logical formula on* $\Sigma$ is inductively defined, along its free variables, as:

(1) $\top, \bot \in \mathrm{LForms}(\Sigma)$ and $\mathrm{FV}(\top) = \mathrm{FV}(\bot) = \emptyset$;
(2) if $p : s_1 \times \cdots \times s_n \in R$ and $t_1 : s_1, \ldots, t_n : s_n \in \mathrm{LTerms}(\Sigma)$ then $p(t_1, \ldots, t_n) \in \mathrm{LForms}(\Sigma)$ and $\mathrm{FV}(p(t_1, \ldots, t_n)) = \bigcup_{i=1}^{n} \mathrm{FV}(t_i : s_i)$;
(3) if $A, B \in \mathrm{LForms}(\Sigma)$ then $A \wedge B, A \vee B, A \supset B \in \mathrm{LForms}(\Sigma)$ and $\mathrm{FV}(A \wedge B) = \mathrm{FV}(A \vee B) = \mathrm{FV}(A \supset B) = \mathrm{FV}(A) \cup \mathrm{FV}(B)$;
(4) if $x \in V_s$ and $A \in \mathrm{LForms}(\Sigma)$ then $\forall x : s.\, A, \exists x : s.\, A \in \mathrm{LForms}(\Sigma)$ and $\mathrm{FV}(\forall x : s.\, A) = \mathrm{FV}(\exists x : s.\, A) = \mathrm{FV}(A) \setminus \{x : s\}$.

As usual, the standard hypotheses on the language apply.

**Definition 2.5** (Substition on logical formulae)**.** Fixed a signature $\Sigma = \langle S, F, R\rangle$, a substitution on variables $\sigma\colon \bigcup_{s\in S} \to \mathrm{LTerms}(\Sigma)$ naturally extends to a map $\sigma''\colon \mathrm{LForms}(\Sigma) \to \mathrm{LForms}(\Sigma)$, called *substitution on (logical) formulae* as follows:

(1) $\sigma''(\top) = \top$, $\sigma''(\bot) = \bot$;
(2) $\sigma''(p(t_1,\dots,t_n)) = p(\sigma'(t_1 : s_1),\dots,\sigma'(t_n : s_n))$ where $\sigma'\colon \mathrm{LTerms}(\Sigma) \to \mathrm{LTerms}(\Sigma)$ is the substitution on logical terms extending $\sigma$, see definition 2.3, and $p : s_1 \times \cdots \times s_n \in R$;
(3) $\sigma''(A \wedge B) = \sigma''(A) \wedge \sigma''(B)$, $\sigma''(A \vee B) = \sigma''(A) \vee \sigma''(B)$, and $\sigma''(A \supset B) = \sigma''(A) \supset \sigma''(B)$;
(4) $\sigma''(\forall x : s.\, A) = \forall x : s.\, \theta''(A)$, $\sigma''(\exists x : s.\, A) = \exists x : s.\, \theta''(A)$ where $\theta''$ is the substitution on logical formulae extending $\theta(y) = \sigma(y)$, for all $y \not\equiv x$ and $\theta(x) = x$.

**Definition 2.6** (Logical calculus)**.** Given a multi-set $\Gamma$ of logical formulae, called *assumptions*, on a fixed signature $\Sigma = \langle S, F, R\rangle$, and a set Ax of logical formulae, called *axioms*, a *derivation* or *proof* $\Gamma \vdash A$ of a logical formula $A$, the *conclusion*, from $\Gamma$ is inductively defined by the following rules:

(AI) $(\Gamma \vdash A) = A$ if $A \in \Gamma$ (assumption);
(Ax) $(\Gamma \vdash A) = A$ if $A \in \mathrm{Ax}$ (axiom);
($\wedge$I) $(\Gamma \vdash A) = ((\Gamma \vdash B)(\Gamma \vdash C) \vdash B \wedge C)$ where $A \equiv B \wedge C$ (and introduction);
($\wedge \mathrm{E_l}$) $(\Gamma \vdash A) = ((\Gamma \vdash A \wedge B) \vdash A)$ (and elimination left);
($\wedge \mathrm{E_r}$) $(\Gamma \vdash A) = ((\Gamma \vdash B \wedge A) \vdash A)$ (and elimination right);
($\vee \mathrm{I_l}$) $(\Gamma \vdash A) = (\Gamma \vdash B) \vdash B \vee C$ where $A \equiv B \vee C$ (or introduction left);
($\vee \mathrm{I_r}$) $(\Gamma \vdash A) = (\Gamma \vdash C) \vdash B \vee C$ where $A \equiv B \vee C$ (or introduction right);
($\vee$E) $(\Gamma \vdash A) = ((\Gamma \vdash B \vee C)(\Gamma \vdash B \supset A)(\Gamma \vdash C \supset A) \vdash A)$ (or elimination);
($\supset$ I) $(\Gamma \vdash A) = ((\Gamma, C \vdash B) \vdash C \supset B)$ where $A \equiv C \supset B$ and the assumption $C$ is *discharged*, i.e., it is not part of $\Gamma$ in the whole proof, but it is added to $\Gamma$ in the proof $\Gamma, C \vdash B$ (implication introduction);
($\supset$ E) $(\Gamma \vdash A) = ((\Gamma \vdash B \supset A)(\Gamma \vdash B) \vdash A)$ (implication elimination);
($\top$I) $(\Gamma \vdash A) = \top$ (truth introduction);
($\bot$E) $(\Gamma \vdash A) = \bot \supset B$ where $A \equiv \bot \supset B$ (falsity elimination);
($\forall$I) $(\Gamma \vdash A) = ((\Gamma \vdash B) \vdash \forall x : s.\, B)$ where $A \equiv \forall x : s.\, B$, provided that $x \in V_s$ is *free in* $\Gamma \vdash B$, i.e., $x : s \not\in \bigcup_{C\in\Gamma} \mathrm{FV}(C)$ (for all introduction);
($\forall$E) $(\Gamma \vdash A) = ((\Gamma \vdash \forall x : s.\, B) \vdash B[t/x])$ where $A \equiv B[t/x]$ and $t : s \in \mathrm{LTerms}(\Sigma)$ (for all elimination);
($\exists$I) $(\Gamma \vdash A) = ((\Gamma \vdash B[t/x]) \vdash \exists x : s.\, B)$ where $A \equiv \exists x : s.\, B$ and $t : s \in \mathrm{LTerms}(\Sigma)$ (exists introduction);
($\exists$E) $(\Gamma \vdash A) = ((\Gamma \vdash \exists x : s.\, B)(\Gamma \vdash B \supset A) \vdash A)$ provided $x \in V_s$ is free in $\Gamma \vdash B \supset A$, that is, $x : s \not\in \bigcup_{C\in\Gamma} \mathrm{FV}(C)$ (exists elimination).

Usually, we will use the natural deduction presentation of the calculus, see, e.g., [TS00]: The rules in this form are presented in figure 1. Also, we omit the subscripts denoting left or right rules when it is clear from the context. Finally, the usual hypotheses on the language apply.

*Note* 2.1. Some inference rules, ($\vee$E), ($\bot$E), differ from the standard version shown in [TS00] because they do not use discharging of assumptions: it is immediate to prove that the chosen rules are equivalent to the standard ones. We prefer the non-standard version because it slightly simplify a few proofs and definitions in the following. Also, the rule ($\top$I) is usually omitted but it is needed to show the symmetry of introduction/elimination, a character of the calculus we will exploit.

$$A \text{ where } A \in \Gamma \qquad \frac{}{A}\,\text{Ax}$$

$$\frac{B \quad C}{B \wedge C}\,\wedge\text{I} \qquad \frac{A}{A \wedge B}\,\wedge\text{E}_{\text{l}} \qquad \frac{A}{B \wedge A}\,\wedge\text{E}_{\text{r}}$$

$$\frac{B \vee C}{B}\,\vee\text{I}_{\text{l}} \qquad \frac{B \vee C}{C}\,\vee\text{I}_{\text{r}} \qquad \frac{B \wedge C \quad B \supset A \quad C \supset A}{A}\,\vee\text{E}$$

$$\frac{\begin{array}{c}[C]\\ \vdots\\ B\end{array}}{C \supset B}\,\supset\text{I} \qquad \frac{B \supset A \quad B}{A}\,\supset\text{E}$$

$$\frac{}{\top}\,\top\text{I} \qquad \frac{}{\bot \supset B}\,\bot\text{E}$$

$$\frac{B}{\forall x\!:\!s.B}\,\forall\text{I}\,_{(*)} \qquad \frac{\forall x\!:\!s.B}{B[t/x]}\,\forall\text{E}$$

$$\frac{B[t/x]}{\exists x\!:\!s.B}\,\exists\text{I} \qquad \frac{\exists x\!:\!s.B \quad B \supset A}{A}\,\exists\text{E}\,_{(*)}$$

FIGURE 1. Inference rules in natural deduction: in $^{(*)}$, the bounded variable must be free in the corresponding proof.

## 3. First-order lambda-calculus

This section illustrates the syntax of the $\lambda$-calculi associated to the logical systems described in the previous section. For the details of the presentation, we follow the schema illustrated in [Joh02b].

**Definition 3.1** (Lambda signature)**.** A *$\lambda$-signature* $\Sigma = \langle S, F, R, \mathrm{Ax}\rangle$ is a structure where

(1) $\langle S, F, R\rangle$ is a logical signature;
(2) Ax is the set of *axiom symbols*, each one decorated as $a : A \to B$ where $A, B \in \lambda\mathrm{Types}(\Sigma)$ and $\mathrm{FV}(A \to B) = \emptyset$.

We assume the standard hypotheses on the language. The canonical logical signature $\langle S, F, R\rangle$ is usually denoted as $\Sigma$, which is non-ambiguous in most contexts.

**Definition 3.2** (Lambda type)**.** Fixed a $\lambda$-signature $\Sigma = \langle S, F, R, \mathrm{Ax}\rangle$, and with the standard hypotheses on the language, the *$\lambda$-types on* $\Sigma$ are inductively defined along with their *free variables* as follows:

(1) $0, 1 \in \lambda\mathrm{Types}(\Sigma)$ and $\mathrm{FV}(0) = \mathrm{FV}(1) = \emptyset$;
(2) if $p : s_1 \times \cdots \times s_n \in R$ and $t_1 : s_1, \ldots, t_n : s_n \in \mathrm{LTerms}(\Sigma)$ then $p(t_1, \ldots, t_n) \in \lambda\mathrm{Types}(\Sigma)$ and $\mathrm{FV}(p(t_1, \ldots, t_n)) = \bigcup_{i=1}^{n} \mathrm{FV}(t_i : s_i)$;
(3) if $A, B \in \lambda\mathrm{Types}(\Sigma)$ then $A \times B, A + B, A \to B \in \lambda\mathrm{Types}(\Sigma)$ and $\mathrm{FV}(A \times B) = \mathrm{FV}(A + B) = \mathrm{FV}(A \to B) = \mathrm{FV}(A) \cup \mathrm{FV}(B)$;
(4) if $x \in V_s$ and $A \in \lambda\mathrm{Types}(\Sigma)$ then $\forall x : s.\, A, \exists x : s.\, A \in \lambda\mathrm{Types}(\Sigma)$ and $\mathrm{FV}(\forall x : s.\, A) = \mathrm{FV}(\exists x : s.\, A) = \mathrm{FV}(A) \setminus \{x : s\}$.

**Fact 3.1.** *The definition of $\lambda$-signature is well-founded.*

*Proof.* Since the set Ax is not used in the definition of $\lambda$Types, circularity is avoided. □

There is an evident correspondence between $\lambda$-types and logical formulae, summarised in Table 1. This correspondence is referred to as the 'proposition as types' interpretation, or, equivalently, as the Curry-Howard isomorphism.

**Definition 3.3** (Lambda term)**.** Fixed a $\lambda$-signature $\Sigma = \langle S, F, R, \mathrm{Ax}\rangle$, for each type $t \in \lambda\mathrm{Types}(\Sigma)$, we assume there is a denumerable set $W_t$ of *(typed) variables*. The standard hypotheses on the language apply.

A *$\lambda$-term* is inductively defined together with its free variables as:

(1) if $x \in W_t$ then $x : t \in \lambda\mathrm{Terms}(\Sigma)$ and $\mathrm{FV}(x : t) = \{x : t\}$;
(2) if $f : A \to B \in \mathrm{Ax}$ and $t : A \in \lambda\mathrm{Terms}(\Sigma)$ then $f(t) : B \in \lambda\mathrm{Terms}(\Sigma)$ and $\mathrm{FV}(f(t) : B) = \mathrm{FV}(t : A)$;
(3) if $s : A, t : B \in \lambda\mathrm{Terms}(\Sigma)$ then $\langle s, t\rangle : A \times B \in \lambda\mathrm{Terms}(\Sigma)$ and $\mathrm{FV}(\langle s, t\rangle : A \times B) = \mathrm{FV}(s : A) \cup \mathrm{FV}(t : B)$;
(4) if $t : A \times B \in \lambda\mathrm{Types}(\Sigma)$ then $\mathrm{fst}(t) : A \in \lambda\mathrm{Terms}(\Sigma)$, $\mathrm{snd}(t) : B \in \lambda\mathrm{Terms}(\Sigma)$ and $\mathrm{FV}(\mathrm{fst}(t) : A) = \mathrm{FV}(\mathrm{snd}(t) : B) = \mathrm{FV}(t : A \times B)$;
(5) if $t : A \in \lambda\mathrm{Terms}(\Sigma)$ then $\mathrm{inl}_B(t) : A + B \in \lambda\mathrm{Terms}(\Sigma)$, $\mathrm{inr}_B(t) : B + A \in \lambda\mathrm{Terms}(\Sigma)$ and $\mathrm{FV}(\mathrm{inl}_B(t) : A + B) = \mathrm{FV}(\mathrm{inr}_B(t) : B + A) = \mathrm{FV}(t : A)$;

| | |
|---|---|
| $0$ | $\bot$ |
| $1$ | $\top$ |
| $\times$ | $\wedge$ |
| $+$ | $\vee$ |
| $\to$ | $\supset$ |
| $\forall$ | $\forall$ |
| $\exists$ | $\exists$ |

Table 1. The 'proposition as types' interpretation.

(6) if $s:A+B, t:A\to C, r:B\to C\in\lambda\text{Terms}(\Sigma)$ then $\text{when}(s,t,r):C\in\lambda\text{Terms}(\Sigma)$ and $\text{FV}(\text{when}(s,t,r):C)=\text{FV}(s:A+B)\cup\text{FV}(t:A\to C)\cup\text{FV}(r:B\to C)$;

(7) if $x\in W_A$ is a variable and $t:B\in\lambda\text{Terms}(\Sigma)$ then $(\lambda x:A.\,t):A\to B\in\lambda\text{Terms}(\Sigma)$ and $\text{FV}((\lambda x:A.\,t):A\to B)=\text{FV}(t:B)\setminus\{x:A\}$;

(8) if $s:A\to B, t:A\in\lambda\text{Terms}(\Sigma)$ then $s\cdot t:B\in\lambda\text{Terms}(\Sigma)$ and $\text{FV}(s\cdot t:B)=\text{FV}(s:A\to B)\cup\text{FV}(t:A)$;

(9) $*:1\in\lambda\text{Terms}(\Sigma)$ and $\text{FV}(*:1)=\emptyset$;

(10) $\text{F}_A:0\to A\in\lambda\text{Terms}(\Sigma)$ and $\text{FV}(\text{F}_A:0\to A)=\emptyset$;

(11) if $x\in V_s$ and $t:A\in\lambda\text{Terms}(\Sigma)$ where $x:s\notin\text{FV}^*(t:A)$, then it holds that the expression $\text{allI}(\lambda x:s.\,t):(\forall x:s.\,A)\in\lambda\text{Terms}(\Sigma)$ and $\text{FV}(\text{allI}(\lambda x:s.\,t):(\forall x:s.\,A))=\text{FV}(t:A)$;

(12) if $t:(\forall x:s.\,A)\in\lambda\text{Terms}(\Sigma)$ and $r:s\in\text{LTerms}(\Sigma)$ then it holds that the expression $\text{allE}(t,r):(A[r/x])\in\lambda\text{Terms}(\Sigma)$ and $\text{FV}(\text{allE}(t,r):(A[r/x]))=\text{FV}(t:(\forall x:s.\,A))$;

(13) if $x\in V_s$, $r:s\in\text{LTerms}(\Sigma)$ and $t:(A[r/x])\in\lambda\text{Terms}(\Sigma)$ then $\text{exI}_x(t):(\exists x:s.\,A)\in\lambda\text{Terms}(\Sigma)$ and $\text{FV}(\text{exI}_x(t):(\exists x:s.\,A))=\text{FV}(t:(A[r/x]))$;

(14) if $t:(\exists x:s.\,A),\in\lambda\text{Terms}(\Sigma)$ and $r:A\to B\in\lambda\text{Terms}(\Sigma)$ where $x:s\notin\text{FV}^*(r:A\to B)$, then it holds that $\text{exE}(t,(\lambda x:s.\,r)):B\in\lambda\text{Terms}(\Sigma)$ and $\text{FV}(\text{exE}(t,(\lambda x:s.\,r)):B)=\text{FV}(t:(\exists x:s.\,A))\cup\text{FV}(r:A\to B)$.

In the previous definition, $x:s\in\text{FV}^*(t:A)$ if and only if there is $r\in\lambda\text{Types}(\Sigma)$ and $y\in W_r$ such that $x:s\in\text{FV}(r)$ and $y:r\in\text{FV}(t:A)$.

There is an obvious correspondence between $\lambda$-terms and logical proofs: comparing definition 2.6 with the construction rules above, we obtain the explicit correspondence, summarised in Table 2. This correspondence is referred to as the ‘proofs as terms’ interpretation, or, equivalently, as the Curry-Howard isomorphism.

**Definition 3.4** (Context)**.** Fixed a $\lambda$-signature, *context* is any finite set of typed variables such that no variable occurs twice.

**Definition 3.5** (Term in context)**.** Fixed a $\lambda$-signature $\Sigma$, and given a context $\vec{x}:\vec{A}=x_1:A_1,\dots,x_n:A_n$, a *term-in-context* is a $\lambda$-term $t:B\in\lambda\text{Terms}(\Sigma)$, notation $\vec{x}:\vec{A}.\,t:B$, such that $\text{FV}(t:B)\subseteq\bigcup_{i=1}^n\text{FV}(x_i:A_i)$.

**Definition 3.6** (Equality in context)**.** Fixed a $\lambda$-signature $\Sigma$, and given a context $\vec{x}:\vec{A}=x_1:A_1,\dots,x_n:A_n$, an *equality-in-context*, notation $\vec{x}:\vec{A}.\,t=_B r$, is a pair $t:B\in\lambda\text{Terms}(\Sigma)$, $r:B\in\lambda\text{Terms}(\Sigma)$, such that $\vec{x}:\vec{A}.\,t:B$ and $\vec{x}:\vec{A}.\,r:B$ are terms-in-context.

**Definition 3.7** (Lambda calculus)**.** A *derivation* is inductively defined by the following inference rules, whose antecedents and consequents are equalities-in-context within a fixed $\lambda$-signature $\Sigma$ (the abbreviations $\vec{x}:\vec{A}=x_1:A_1,\dots,x_n:A_n$, and, similarly, $\vec{y}:\vec{B}=y_1:B_1,\dots,y_m:B_n$, as well as $[\vec{s}/\vec{y}]=[s_1/y_1,\dots,s_m/y_m]$ are used to lighten notation):

($\text{eq}_0$) $\vec{x}:\vec{A}.\,s=_C t\vdash\vec{y}:\vec{B}.\,s[r_1/x_1,\dots,r_n/x_n]=_C t[r_1/x_1,\dots,r_n/x_n]$ where, for any $1\le i\le n$, $\vec{y}:\vec{B}.\,r_i:A_i$ is a term-in-context;

($\text{eq}_1$) $$\left.\begin{array}{c}(\vec{x}:\vec{A}.\,s_1=_{B_1}t_1)\\ \vdots\\ (\vec{x}:\vec{A}.\,s_m=_{B_m}t_m)\end{array}\right\}\vdash\vec{x}:\vec{A}.\,r[\vec{s}/\vec{y}]=_C r[\vec{t}/\vec{y}];$$

($\text{eq}_2$) $\vdash x:A.\,x=_A x$;

($\text{eq}_3$) $x:A,y:A.\,x=_A y\vdash x:A,y:A.\,y=_A x$;

($\text{eq}_4$) $$\left.\begin{array}{c}(x:A,y:A,z:A.\,x=_A y)\\ (x:A,y:A,z:A.\,y=_A z)\end{array}\right\}\vdash x:A,y:A,z:A.\,x=_A z;$$

($\text{eq}_5$) $\vec{x}:\vec{A}.\,s=_C t\vdash\vec{x}:\vec{A}.\,(\lambda y:B.\,s)=_{B\to C}(\lambda y:B.\,t)$;

($\text{eq}_6$) $\vec{x}:\vec{A}.\,r=_C t\vdash\vec{x}:\vec{A}.\,\text{allI}(\lambda y:s.\,r)=_{(\forall y:s.\,C)}\text{allI}(\lambda y:s.\,t)$;

($\text{eq}_7$) $\vec{x}:\vec{A}.\,u=_C v\vdash\vec{x}:\vec{A}.\,\text{exE}(t,(\lambda y:s.\,u))=_C\text{exE}(t,(\lambda y:s.\,v))$;

($\times_0$) $\vdash x:1.\,x=_1 *$;

($\times_1$) $\vdash x\colon A, y\colon B.\, \mathrm{fst}(\langle x,y\rangle) =_A x$;
($\times_2$) $\vdash x\colon A, y\colon B.\, \mathrm{snd}(\langle x,y\rangle) =_B y$;
($\times_3$) $\vdash z\colon A\times B.\, \langle \mathrm{fst}(z), \mathrm{snd}(z)\rangle =_{A\times B} z$;
($+_0$) $\vdash \vec{x}\colon\vec{A}.\, \mathrm{when}(\mathrm{inl}_B(a), t, s) =_C t\cdot a$;
($+_1$) $\vdash \vec{x}\colon\vec{A}.\, \mathrm{when}(\mathrm{inr}_D(b), t, s) =_C s\cdot b$;
($+_2$) $\vdash x_0\colon A_1+A_2, x_1\colon A_1 \to (B_1+B_2), x_2\colon A_2 \to (B_1+B_2), x_3\colon B_1 \to C, x_4\colon B_2 \to C.$

$$\mathrm{when}(\mathrm{when}(x_0,x_1,x_2),x_3,x_4) =_C$$
$$=_C \mathrm{when}(x_0, (\lambda y\colon A_1.\, \mathrm{when}(x_1\cdot y, x_3, x_4)), (\lambda y\colon A_2.\, \mathrm{when}(x_2\cdot y, x_3, x_4)))$$

where $y\colon A_1 \not\in \mathrm{FV}(x_1\colon A_1+A_2) \cup \mathrm{FV}(x_3\colon B_1\to C) \cup \mathrm{FV}(x_4\colon B_2 \to C)$ and $y\colon A_2 \not\in \mathrm{FV}(x_2\colon A_1+A_2) \cup \mathrm{FV}(x_3\colon B_1\to C) \cup \mathrm{FV}(x_4\colon B_2\to C)$;
($+_3$) $\vdash x\colon A, y\colon 0.\, \mathrm{F}_A\cdot y =_A x$;
($\to_0$) $\vdash \vec{x}\colon\vec{A}.\, (\lambda y\colon C.\, s)\cdot t =_B s[t/y]$;
($\to_1$) $\vdash \vec{x}\colon\vec{A}.\, (\lambda y\colon C.\, t\cdot y) =_{C\to B} t$ where $y\colon C\not\in \mathrm{FV}(t\colon C\to B)$;
($\forall_0$) $\vdash \vec{x}\colon\vec{A}.\, \mathrm{allE}(\mathrm{allI}(\lambda z\colon s.\, t), r) =_{B[r/z]} t[r/z]$;
($\forall_1$) $\{\vec{x}\colon\vec{A}.\, \mathrm{allE}(u,r) =_B \mathrm{allE}(v,r)\}_{r\,:\,s\in\mathrm{LTerms}(\Sigma)} \vdash \vec{x}\colon\vec{A}.\, u =_{(\forall z\colon s.\, B)} v$;
($\exists_0$) $\vdash \vec{x}\colon\vec{A}.\, \mathrm{exE}(\mathrm{exI}_z(t), (\lambda z\colon s.\, v)) =_B (v[r/z])\cdot t$;
($\exists_1$) $\vec{x}\colon\vec{A}.\, \mathrm{exE}(u, (\lambda z\colon s.\, r)) =_B \mathrm{exE}(u,(\lambda z\colon s.\, t)) \vdash \vec{x}\colon\vec{A}.\, r =_{C\to B} t$ where $\mathrm{FV}(r\colon C\to B) = \mathrm{FV}(t\colon C\to B)$;
($\exists_2$) $\vdash v\colon(\exists y\colon s.\, A).\, w =_B \mathrm{exE}(v, (\lambda y\colon s.\, (\lambda z\colon A.\, w[\mathrm{exI}_y(z)/v])))$ with $z\colon A \not\in \mathrm{FV}(w\colon B)$;
($\exists_3$) $\vdash \vec{x}\colon\vec{A}.\, \mathrm{exE}(\mathrm{exE}(a, (\lambda y\colon s.\, (\lambda z\colon D.\, b))), (\lambda y\colon s.\, c)) =_C$
$=_C \mathrm{exE}(a, (\lambda y\colon s.\, (\lambda z\colon D.\, \mathrm{exE}(b, (\lambda y\colon s.\, c)))))$;
($\exists_4$) $\vdash \vec{x}\colon\vec{A}.\, \mathrm{exE}(a, (\lambda y\colon s.\, (\lambda z\colon C.\, b[\mathrm{exI}_y(z)/w]))) =_B b[a/w]$ with $z\colon C\not\in \mathrm{FV}(b\colon B)$.

*Note* 3.1. The complex rules of the $\lambda$-calculus require some remarks:

(1) rules ($\mathrm{eq}_0$) and ($\mathrm{eq}_1$) define the behaviour of equality with respect to substitutions. The format of these rules is taken from [Joh02b].
(2) rules ($\mathrm{eq}_2$), ($\mathrm{eq}_3$) and ($\mathrm{eq}_4$) say that equality is an equivalence relation;
(3) rules ($\mathrm{eq}_5$), ($\mathrm{eq}_6$) and ($\mathrm{eq}_7$) define equality when bounded variables are present in a term;
(4) rule ($\times_0$) is the reduction associated to the empty product;
(5) rules ($\times_1$), ($\times_2$) and ($\times_3$) are the standard reductions of binary products, see, e.g., [GLT89];
(6) rules ($+_0$) and ($+_1$) are the standard reductions of binary sums, see [GLT89];
(7) rule ($+_3$) defines the reduction of empty sums, modulo an abstraction;
(8) rule ($\to_0$) is the usual $\beta$-reduction, see [GLT89];
(9) rule ($\to_1$) is the usual $\eta$-reduction, , see [GLT89];
(10) rules ($+_2$), ($\forall_0$), ($\exists_0$), ($\exists_2$), ($\exists_3$) and ($\exists_4$) are the $\lambda$-coding via the Curry-Howard isomorphism, of some conversions used in the standard proof of the normalisation theorem, see [TS00, pp. 180–181];
(11) rule ($\exists_1$) encodes that two proofs, whose last step is ($\exists$E), which are equal, must have their parametric subproofs equal as well;
(12) the infinitary rule ($\forall_1$) says that two universal terms are equal when all their possible instances are equal.

It is worth noticing that equality of proofs is decidable since the only problematic rule, ($\forall_1$), can be simplified to:

$$\vec{x}\colon\vec{A}.\, \mathrm{allE}(u,y) =_B \mathrm{allE}(v,y) \vdash \vec{x}\colon\vec{A}.\, u =_{(\forall z\colon s.\, B)} v \enspace ,$$

provided that $y\in V_s$ is such that $y\colon s\not\in \mathrm{FV}(u\colon(\forall z\colon s.\, B)) \cup \mathrm{FV}(v\colon(\forall z\colon s.\, B))$. The proof of this fact, and its equivalence to ($\forall_1$) follows from the semantics, and it will appear in section 6.

| | |
|---|---|
| $x : A$ | the assumption $A$, whose occurrence is $x$ |
| $f(t) : B$ | $\dfrac{\overline{A \supset B}\ \text{Ax} \quad \overset{\vdots\, t}{A}}{B}$ |
| $\langle s, t\rangle : A \times B$ | $\dfrac{\overset{\vdots\, s}{A} \quad \overset{\vdots\, t}{B}}{A \wedge B}$ |
| $\mathrm{fst}(t) : A, \mathrm{snd}(t) : B$ | $\dfrac{\overset{\vdots\, t}{A \wedge B}}{A} \quad \dfrac{\overset{\vdots\, t}{A \wedge B}}{B}$ |
| $\mathrm{inl}_B(t) : A + B, \mathrm{inr}_B(t) : B + A$ | $\dfrac{\overset{\vdots\, t}{A}}{A \vee B} \quad \dfrac{\overset{\vdots\, t}{B}}{B \vee A}$ |
| $\mathrm{when}(s, t, r) : C$ | $\dfrac{\overset{\vdots\, s}{A \vee B} \quad \overset{\vdots\, t}{A \supset C} \quad \overset{\vdots\, r}{B \supset C}}{C}$ |
| $(\lambda x : s.\, t) : A \to B$ | $\dfrac{\overset{[A^x]}{\overset{\vdots\, t}{B}}}{A \supset B}$ |
| $s \cdot t : B$ | $\dfrac{\overset{\vdots\, s}{A \supset B} \quad \overset{\vdots\, t}{A}}{B}$ |
| $* : 1$ | the axiom $\top$ |
| $\mathrm{F}_A : 0 \to A$ | the axiom $\bot \supset A$ |
| $\mathrm{allI}(\lambda x : s.\, t) : (\forall x : s.\, A)$ | $\dfrac{\overset{\vdots\, t}{A}}{\forall x : s.\, A}$ |
| $\mathrm{allE}(t, r) : A[r/x]$ | $\dfrac{\overset{\vdots\, t}{\forall x : s.\, A}}{A[r/x]}$ |
| $\mathrm{exI}_x(t) : (\exists x : s.\, A)$ | $\dfrac{\overset{\vdots\, t}{A}}{\exists x : s.\, A}$ |
| $\mathrm{exE}(t, (\lambda x : s.\, r)) : B$ | $\dfrac{\overset{\vdots\, t}{\exists x : s.\, A} \quad \overset{\vdots\, r}{A \supset B}}{B}$ |

TABLE 2. The ‘proofs as terms’ interpretation.

## 4. Categorical semantics

This section defines the semantics of both the logic and the $\lambda$-calculus introduced in sections 2 and 3. The propositional part is based on [Tay99], while the first-order extension is novel.

**Definition 4.1** (Logically distributive category)**.** Fixed a $\lambda$-signature $\Sigma = \langle S, F, R, \mathrm{Ax}\rangle$, a category $\mathbb{C}$ together with a map $M\colon \lambda\mathrm{Types}(\Sigma) \to \mathrm{Obj}\,\mathbb{C}$ is said to be *logically distributive* if it satisfies the following seven conditions:

(1) $\mathbb{C}$ has finite products;
(2) $\mathbb{C}$ has finite co-products;
(3) $\mathbb{C}$ has exponentiation;
(4) $\mathbb{C}$ is *distributive*, i.e., for every $A, B, C \in \mathrm{Obj}\,\mathbb{C}$, the arrow $\Delta = [1_A \times \iota_1, 1_A \times \iota_2]\colon (A \times B) + (A \times C) \to A \times (B + C)$ has an inverse, where $[_,_]$ is the co-universal arrow of the $(A \times B) + (A \times C)$ co-product, $_ \times _$ is the product arrow, see [Gol06], $1_A$ is the identity arrow on $A$, and $\iota_1\colon B \to B + C$, $\iota_2\colon C \to B + C$ are the canonical injections of the $B + C$ co-product.

For every $s \in S$, $A \in \lambda\mathrm{Types}(\Sigma)$, and $x \in V_s$, let $\Sigma_A(x\colon s)\colon \mathrm{LTerms}(\Sigma)(s) \to \mathbb{C}$ be the functor from the discrete category $\mathrm{LTerms}(\Sigma)(s) = \{t\colon s \mid t\colon s \in \mathrm{LTerms}(\Sigma)\}$ to $\mathbb{C}$ defined by $t\colon s \mapsto M(A[t/x])$.

Also, for every $s \in S$, $A \in \lambda\mathrm{Types}(\Sigma)$, and $x \in V_s$, let $\mathbb{C}_{(\forall x\colon s.\, A)}$ be the subcategory of $\mathbb{C}$ whose objects are the vertexes of the cones on $\Sigma_A(x\colon s)$ such that they are of the form $MB$ for some $B \in \lambda\mathrm{Types}(\Sigma)$ and $x\colon s \not\in \mathrm{FV}(B)$. Moreover, the arrows of $\mathbb{C}_{(\forall x\colon s.\, A)}$, apart identities, are the arrows in the category of cones over $\Sigma_A(x\colon s)$ having the objects of $\mathbb{C}_{(\forall x\colon s.\, A)}$ as domain and $M(\forall x\colon s.\, A)$ as co-domain.

Finally, for every $s \in S$, $A \in \lambda\mathrm{Types}(\Sigma)$, and $x \in V_s$, let $\mathbb{C}_{(\exists x\colon s.\, A)}$ be the subcategory of $\mathbb{C}$ whose objects are the vertexes of the co-cones on $\Sigma_A(x\colon s)$ such that they are of the form $MB$ for some $B \in \lambda\mathrm{Types}(\Sigma)$ and $x\colon s \not\in \mathrm{FV}(B)$. Moreover, the arrows of $\mathbb{C}_{(\exists x\colon s.\, A)}$, apart identities, are the arrows in the category of co-cones over $\Sigma_A(x\colon s)$ having the objects of $\mathbb{C}_{(\exists x\colon s.\, A)}$ as co-domain and $M(\exists x\colon s.\, A)$ as domain.

(5) All the subcategories $\mathbb{C}_{(\forall x\colon s.\, A)}$ have terminal objects, and all the subcategories $\mathbb{C}_{(\exists x\colon s.\, A)}$ have initial objects;
(6) The $M$ map is such that
    (a) $M(0) = 0$, the initial object of $\mathbb{C}$;
    (b) $M(1) = 1$, the terminal object of $\mathbb{C}$;
    (c) $M(A \times B) = MA \times MB$, the binary product in $\mathbb{C}$;
    (d) $M(A + B) = MA + MB$, the binary co-product in $\mathbb{C}$;
    (e) $M(A \to B) = MB^{MA}$, the exponential object in $\mathbb{C}$;
    (f) $M(\forall x\colon s.\, A)$ is the terminal object in the subcategory $\mathbb{C}_{(\forall x\colon s.\, A)}$;
    (g) $M(\exists x\colon s.\, A)$ is the initial object in the subcategory $\mathbb{C}_{(\exists x\colon s.\, A)}$;
(7) For every $x \in V_s$, $A, B \in \lambda\mathrm{Types}(\Sigma)$ with $x\colon s \not\in \mathrm{FV}(A)$, $MA \times M(\exists x\colon s.\, B)$ is an object of $\mathbb{C}_{(\exists x\colon s.\, A \times B)}$ since, if $\big(M(\exists x\colon s.\, B), \{\delta_t\}_{t\colon s \in \mathrm{LTerms}(\Sigma)}\big)$ is a co-cone over $\Sigma_B(x\colon s)$, and there is one by condition (5), then $\big(MA \times M(\exists x\colon s.\, B), \{1_{MA} \times \delta_t\}_{t\colon s \in \mathrm{LTerms}(\Sigma)}\big)$ is a co-cone over $\Sigma_{A \times B}(x\colon s)$. Thus, there is a unique arrow $!\colon M(\exists x\colon s.\, A \times B) \to MA \times M(\exists x\colon s.\, B)$ in $\mathbb{C}_{(\exists x\colon s.\, A \times B)}$. Our last condition requires that the arrow $!$ has an inverse.

*Note* 4.1. The definition is well-founded as it does not **define** the map $M$, but it imposes a condition on its structure. This avoids the evident circularity. But, obviously, the definition is impredicative.

*Note* 4.2. The first four conditions are the ones used in [Tay99] to define the categorical semantics of propositional intuitionistic logic. The other three conditions are specific to first-order logic. It is worth noticing that every logically distributive category provides a model for a propositional system that considers quantified formulae as constants.

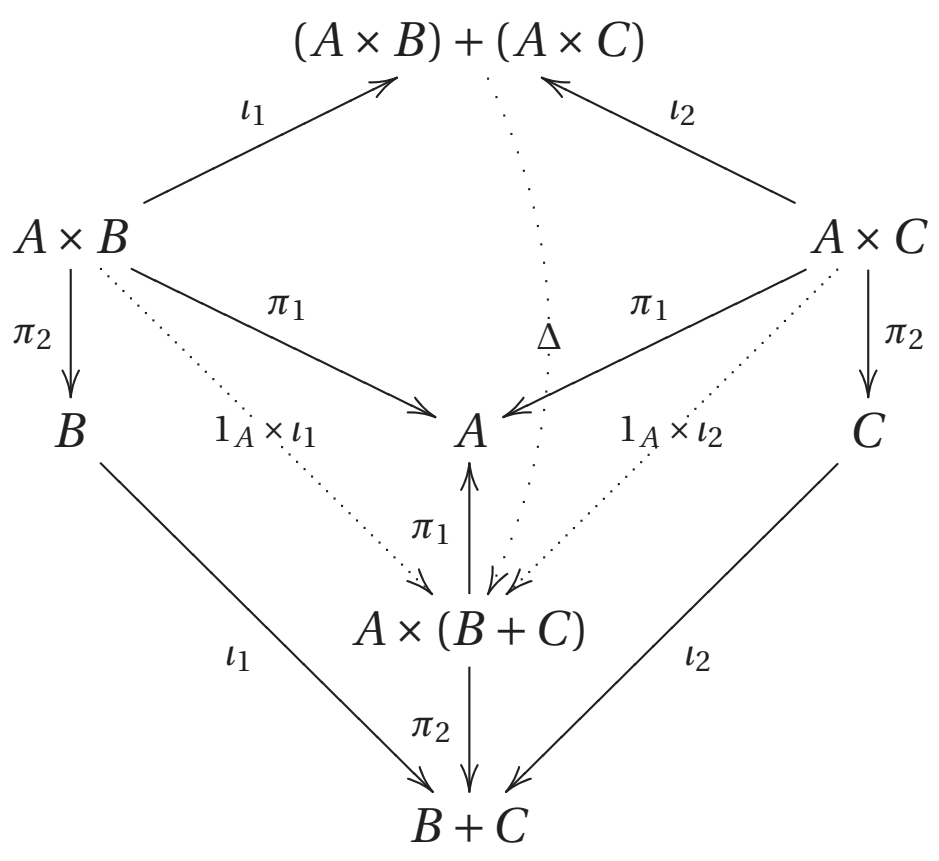

FIGURE 2. The construction of the $\Delta$ arrow.

*Note* 4.3. The arrow $\Delta$ is derived from the commutative diagram in figure 2. In the following, we will reserve the symbol $\Delta$ to denote this arrow.

*Note* 4.4. Condition (5) can be graphically visualised by noticing that the $\mathbb{C}_{(\forall x:s.\,A)}$ subcategories are 'star-shaped' with an unique arrow from the centre to each other object. Dually, the $\mathbb{C}_{(\exists x:s.\,A)}$ subcategories are 'co-star-shaped', see figures 3 and 4. The condition forces the fact that there are just enough arrows to build these 'stars'. It is worth noticing that, since $M0 = 0$ and $M1 = 1$ these categories are always inhabited.

*Note* 4.5. The last requirement is just distributivity of $\times$ over $\exists$, intuitively considered as an infinite disjunction.

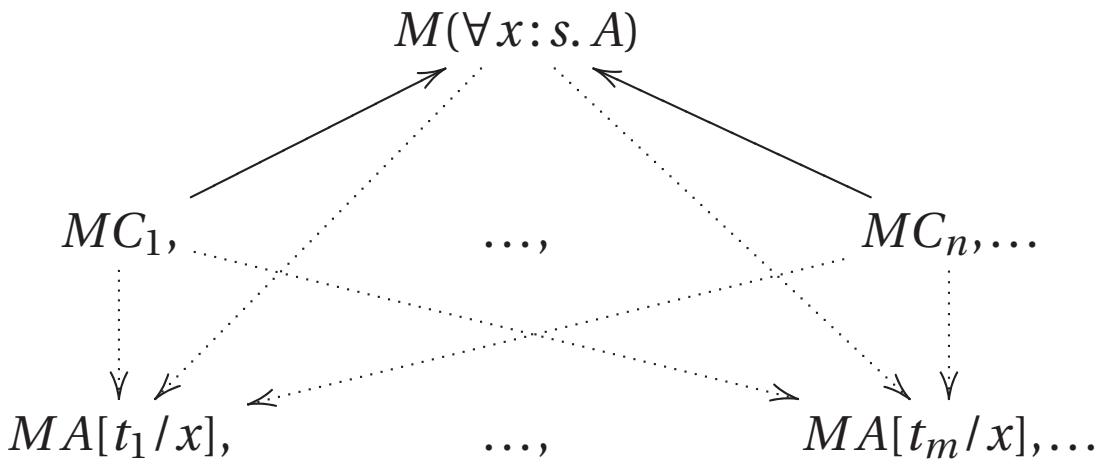

FIGURE 3. The graphical aspect of $\mathbb{C}_{\forall x:s.\,A}$.

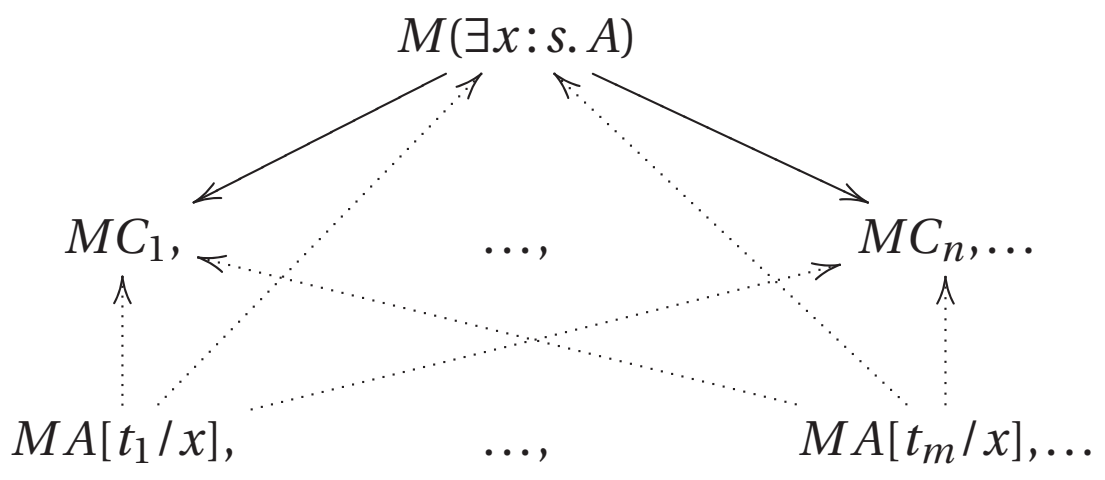

FIGURE 4. The graphical aspect of $\mathbb{C}_{\exists x:s.\,A}$.

The $\Sigma_A(x\colon s)$ functor is the 'glue' which links logical terms to the structure of the category. As we will show in the following sections, when soundness and completeness are proved, the only role terms have in the logical system is in the fact that they allow for substituting variables. Also, abstractly, the only role variables have is either being substituted, and this action is modelled by the $\Sigma_A$ family of functors, or being quantified, and this action is captured by the $\mathbb{C}_{\forall x\colon s.\,A}$ and $\mathbb{C}_{\exists x\colon s.\,A}$ subcategories.

**Proposition 4.1.** *When $z\colon s$ is a new variable,* $\Sigma_{\exists x\colon s.\,A}(y\colon r) = \Sigma_{\exists z\colon s.\,A[z/x]}(y\colon r)$.

*Proof.* For each $t\colon r \in \mathrm{LTerms}(\Sigma)$,

$$\Sigma_{\exists x\colon s.\,A}(y\colon r)(t\colon r) = M((\exists x\colon s.\,A)[t/y]) = M(\exists w\colon s.\,A[w/x][t/y])$$

with $w\colon s$ a new variable. Analogously,

$$\begin{aligned}\Sigma_{\exists z\colon s.\,A[z/x]}(y\colon r)(t\colon r) &= M((\exists z\colon s.\,A[z/x])[t/y]) = \\ &= M(\exists w\colon s.\,A[z/x][w/z][t/y]) = M(\exists w\colon s.\,A[w/x][t/y])\ .\end{aligned}$$ □

**Proposition 4.2.** *When $z\colon s$ is a new variable,* $\Sigma_{\forall x\colon s.\,A}(y\colon r) = \Sigma_{\forall z\colon s.\,A[z/x]}(y\colon r)$.

*Proof.* As proposition 4.1. □

*Note* 4.6. Propositions 4.1 and 4.2 allow the renaming of bounded variables in quantified types, ensuring that the denoted object is the same.

**Proposition 4.3.** $\Sigma_{A\times B}(x\colon s)(t\colon s) = \Sigma_A(x\colon s)(t\colon s) \times \Sigma_B(x\colon s)(t\colon s)$.

*Proof.* $\Sigma_{A\times B}(x\colon s)(t\colon s) = M((A\times B)[t/x]) = M(A[t/x]\times B[t/x]) =$

$$= M(A[t/x]) \times M(B[t/x]) = \Sigma_A(x\colon s)(t\colon s) \times \Sigma_B(x\colon s)(t\colon s)\ .$$ □

**Proposition 4.4.** $\Sigma_{A+B}(x\colon s)(t\colon s) = \Sigma_A(x\colon s)(t\colon s) + \Sigma_B(x\colon s)(t\colon s)$.

*Proof.* $\Sigma_{A+B}(x\colon s)(t\colon s) = M((A+B)[t/x]) = M(A[t/x]+B[t/x]) =$

$$= M(A[t/x]) + M(B[t/x]) = \Sigma_A(x\colon s)(t\colon s) + \Sigma_B(x\colon s)(t\colon s)\ .$$ □

**Proposition 4.5.** $\Sigma_0(x\colon s)(t\colon s) = 0$, $\Sigma_1(x\colon s)(t\colon s) = 1$.

*Proof.* $\Sigma_0(x\colon s)(t\colon s) = M(0[t/x]) = M0 = 0$ ;

$\Sigma_1(x\colon s)(t\colon s) = M(1[t/x]) = M1 = 1$ . □

**Proposition 4.6.** $\Sigma_{A\to B}(x\colon s)(t\colon s) = \Sigma_B(x\colon s)(t\colon s)^{\Sigma_A(x\colon s)(t\colon s)}$.

*Proof.* $\Sigma_{A\to B}(x\colon s)(t\colon s) = M((A\to B)[t/x]) = M(A[t/x]\to B[t/x]) =$

$$= M(B[t/x])^{M(A[t/x])} = \Sigma_B(x\colon s)(t\colon s)^{\Sigma_A(x\colon s)(t\colon s)}\ .$$ □

**Proposition 4.7.** $\Sigma_{\exists x\colon s.\,A}(y\colon r)(t\colon r)$ *is the initial object of* $\mathbb{C}_{\exists z\colon s.\,A[z/x][t/y]}$ *where $z\colon s$ is a new variable.*

*Proof.* $\Sigma_{\exists x\colon s.\,A}(y\colon r)(t\colon r) = M((\exists x\colon s.\,A)[t/y]) = M(\exists z\colon s.\,A[z/x][t/y])$ where $z\colon s$ a new variable. The conclusion follows by definition 4.1. □

**Proposition 4.8.** $\Sigma_{\forall x\colon s.\,A}(y\colon r)(t\colon r)$ *is the terminal object of* $\mathbb{C}_{\forall z\colon s.\,A[z/x][t/y]}$ *where $z\colon s$ is a new variable.*

*Proof.* $\Sigma_{\forall x\colon s.\,A}(y\colon r)(t\colon r) = M((\forall x\colon s.\,A)[t/y]) = M(\forall z\colon s.\,A[z/x][t/y])$ where $z\colon s$ a new variable. The conclusion follows by definition 4.1. □

*Note* 4.7. Propositions 4.3 to 4.8 show how the functors $\Sigma_A(x\colon s)$ depend only on the definition of $M$ on $\lambda$-terms of the form $p(t_1,\ldots,t_n)$ **and** on quantified terms—specifically, to identify the initial and terminal objects in the appropriate subcategories. Thus, these pieces of information are those characterising $M$.

This fact also explains in a different way why definition 4.1 is well-founded.

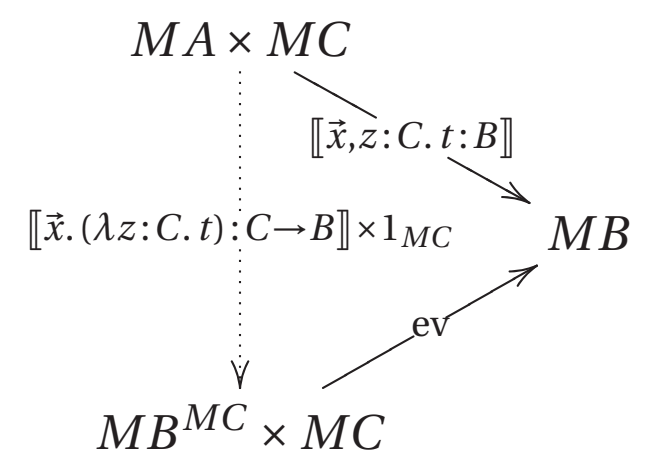

FIGURE 5. The semantics of $(\lambda z\!:\!C.\, t)\!:\!C \to B$.

**Definition 4.2** ($\Sigma$-structure)**.** Given a $\lambda$-signature $\Sigma = \langle S, F, R, \mathrm{Ax} \rangle$, a $\Sigma$-structure is a triple $\langle \mathbb{C}, M, M_{\mathrm{Ax}} \rangle$ such that $\mathbb{C}$ together with $M$ forms a logically distributive category and $M_{\mathrm{Ax}}$ is a map from Ax such that $M_{\mathrm{Ax}}(a\!:A \to B) \in \mathrm{Hom}_{\mathbb{C}}(MA, MB)$.

*Note* 4.8. Sometimes, we will refer to $M$ in a given $\Sigma$-structure $\langle \mathbb{C}, M, M_{\mathrm{Ax}} \rangle$ as the *semantics of* $\lambda\mathrm{Types}(\Sigma)$ or also as the semantics of the $\lambda$-types in the logically distributive category $\mathbb{C}$ with $M$.

Similarly, we will refer to $M_{\mathrm{Ax}}$ as the *semantics of* axioms in the logically distributive category $\mathbb{C}$ with $M$.

**Definition 4.3** ($\lambda$-terms semantics)**.** Fixed a $\Sigma$-structure $\langle \mathbb{C}, M, M_{\mathrm{Ax}} \rangle$, let $A \equiv A_1 \times \cdots A_n$, and let $\vec{x} \equiv x_1\!:\!A_1, \ldots, x_n\!:\!A_n$ be a context. The *semantics of a term-in-context* $\vec{x}.\, t\!:\!B$, notation $[\![\vec{x}.\, t\!:\!B]\!]$, is an arrow in $\mathrm{Hom}_{\mathbb{C}}(MA, MB)$ inductively defined as follows:

(1) $[\![\vec{x}.\, x_i : A_i]\!] = \pi_i$, the $i$-th projector of the product $MA = MA_1 \times \cdots \times MA_n$;
(2) if $a\!:\!C \to B \in \mathrm{Ax}$ then $[\![\vec{x}.\, a(t)\!:\!B]\!] = M_{\mathrm{Ax}} a \circ [\![\vec{x}.\, t\!:\!C]\!]$;
(3) $[\![\vec{x}.\, \langle s, t \rangle\!:\!B \times C]\!] = \big([\![\vec{x}.\, s\!:\!B]\!], [\![\vec{x}.\, t\!:\!C]\!]\big)$ where $(_,_)$ is the universal arrow of the product $MB \times MC$;
(4) $[\![\vec{x}.\, \mathrm{fst}(t)\!:\!B]\!] = \pi_1 \circ [\![\vec{x}.\, t\!:\!B \times C]\!]$ where $\pi_1$ is the first canonical projector of the product $MA \times MB$;
(5) $[\![\vec{x}.\, \mathrm{snd}(t)\!:\!C]\!] = \pi_2 \circ [\![\vec{x}.\, t\!:\!B \times C]\!]$ where $\pi_2$ is the second canonical projector of the product $MA \times MB$;
(6) $[\![\vec{x}.\, (\lambda z\!:\!C.\, t)\!:\!C \to B]\!]$ is the exponential transpose of $[\![\vec{x}, z\!:\!C.\, t\!:\!B]\!]\colon MA \times MC \to MB$, see figure 5;
(7) $[\![\vec{x}.\, s \cdot t\!:\!B]\!] = \mathrm{ev} \circ \big([\![\vec{x}.\, s\!:\!C \to B]\!], [\![\vec{x}.\, t\!:\!C]\!]\big)$ where ev is the exponential evaluation arrow;
(8) $[\![\vec{x}.\, \mathrm{inl}_B(t)\!:\!C + B]\!] = \iota_1 \circ [\![\vec{x}.\, t\!:\!C]\!]$ with $\iota_1$ the first canonical injection of the co-product $MC + MB$;
(9) $[\![\vec{x}.\, \mathrm{inr}_C(t)\!:\!C + B]\!] = \iota_2 \circ [\![\vec{x}.\, t\!:\!B]\!]$ with $\iota_2$ the second canonical injection of the co-product $MC + MB$;
(10)
$$\begin{aligned}[\![\vec{x}.\, \mathrm{when}(t, u, v)\!:\!B]\!] = \big[&\mathrm{ev} \circ \big([\![\vec{x}.\, u\!:\!C_1 \to B]\!] \times 1_{MC_1}\big),\\ &\mathrm{ev} \circ \big([\![\vec{x}.\, v\!:\!C_2 \to B]\!] \times 1_{MC_2}\big)\big] \circ\\ &\circ \Delta^{-1} \circ \big(1_{MA}, [\![\vec{x}.\, t\!:\!C_1 + C_2]\!]\big) \quad ,\end{aligned}$$
where $[_,_]$ is the co-universal arrow of the co-product $(MA \times MC_1) + (MA \times MC_2)$, $(_,_)$ is the universal arrow of the product $MA \times (MC_1 + MC_2)$, and the arrow $\Delta\colon (MA \times MC_1) + (MA \times MC_2) \to MA \times (MC_1 + MC_2)$ has an inverse because $\mathbb{C}$ with $M$ is logically distributive—figure 6 illustrates the commutative diagram behind this case;
(11) $[\![\vec{x}.\, *\!:\!1]\!] = !\colon MA \to 1$, the universal arrow of the terminal object;
(12) $[\![\vec{x}.\, \mathrm{F}_B\!:\!0 \to B]\!]$ is the exponential transpose of $(!\colon 0 \to MB) \circ (\pi_{n+1}\colon MA \times 0 \to 0)$;

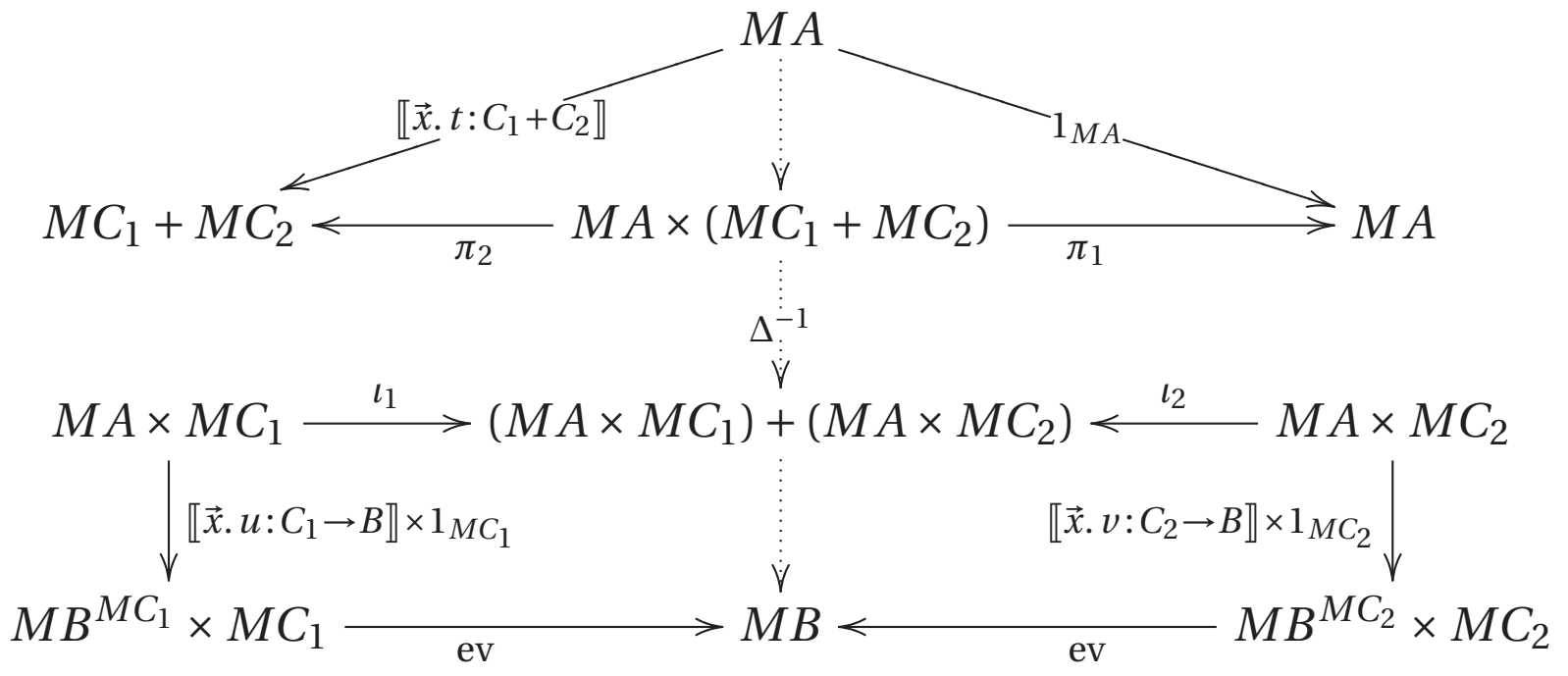

FIGURE 6. The semantics of when$(t,u,v)$.

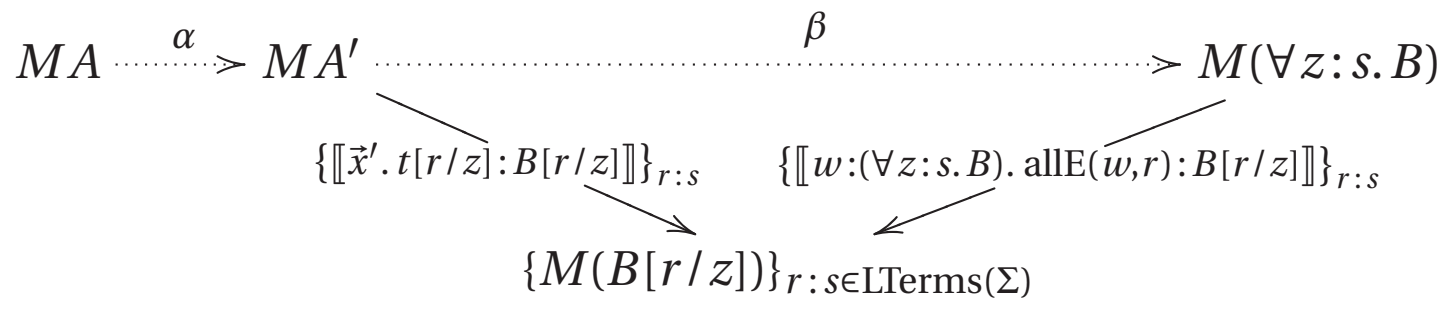

FIGURE 7. Semantics of allI$(\lambda z{:}\, s.\, t)$.

(13) $[\![\vec{x}.\ \mathrm{allI}(\lambda z{:}\, s.\, t){:}(\forall z{:}\, s.B)]\!] = \beta\circ\alpha$ where $\alpha \equiv 1_{MA_{i_1}} \times\cdots\times 1_{MA_{i_k}}\colon MA \to MA'$ with $A' \equiv A_{i_1}\times\cdots\times A_{i_k}$, where $\vec{x}' \equiv \{x_{i_1}{:}\,A_{i_1},\dots,x_{i_k}{:}\,A_{i_k}\} = \mathrm{FV}(t{:}\,B)$, and $\beta\colon MA' \to M(\forall z{:}\, s.B)$, the universal arrow from $MA'$ to the terminal object in $\mathbb{C}_{\forall z:s.B}$.

The definition makes sense, since $\left(MA', \left\{[\![\vec{x}'.\ t[r/z]{:}\,B[r/z]]\!]\right\}_{r:s\in\mathrm{LTerms}(\Sigma)}\right)$ is a cone on $\Sigma_B(z{:}\,s)$. In fact, $\vec{x}'.\ t{:}\,B$ is a term-in-context, so $\vec{x}'[r/z].\ t[r/z]{:}\,B[r/z]$ is term-in-context too, as it is immediate to prove by induction on the structure of terms. But $z{:}\,s \notin \mathrm{FV}^*(t{:}\,B) = \mathrm{FV}^*(\vec{x}') = \mathrm{FV}(A')$, by construction of the term allI$(\lambda z{:}\, s.\, t)$, see definition 3.3, so $\vec{x}'[r/z] \equiv \vec{x}'$. Thus, for every $r{:}\,s \in \mathrm{LTerms}(\Sigma)$, $[\![\vec{x}'.\ t[r/z]{:}\,B[r/z]]\!]\colon MA' \to M(B[r/z]) \equiv \Sigma_B(z{:}\,s)(r{:}\,s)$. Hence, $MA'$ is the vertex of a cone on $\Sigma_B(z{:}\,s)$ and $z{:}\,s\notin \mathrm{FV}(A')$, making $MA'$ an object of $\mathbb{C}_{\forall z:s.B}$. But $\mathbb{C}$ with $M$ is logically distributive, so $M(\forall z{:}\, s.B)$ is the terminal object of $\mathbb{C}_{\forall z:s.B}$, and $\beta$ is uniquely identified. Figure 7 illustrates the commutative diagram corresponding to this case.

(14) $[\![\vec{x}.\ \mathrm{allE}(t,r){:}\,B[r/z]]\!] = p_r\circ[\![\vec{x}.\ t{:}(\forall z{:}\, s.B)]\!]$ where $p_r\colon M(\forall z{:}\, s.B)\to M(B[r/z])$ is the $r$-th projector of the unique cone on $\Sigma_B(z{:}\,s)$ whose vertex is $M(\forall z{:}\, s.B)$.

The definition is sound: there is a cone $\left(M(\forall z{:}\, s.B), \{p_r\}_{r:s\in\mathrm{LTerms}(\Sigma)}\right)$ on $\Sigma_B(z{:}\,s)$ since $M(\forall z{:}\, s.B)$ is the terminal object of $\mathbb{C}_{M(\forall z:s.B)}$, as $\mathbb{C}$ with $M$ is logically distributive, see definition 4.1. For the same reason, there is a unique arrow $M(\forall z{:}\, s.B)\to M(\forall z{:}\, s.B)$ in $\mathbb{C}_{\forall z:s.B}$, which is $1_{M(\forall z:s.B)}$, by definition of $\mathbb{C}_{\forall z:s.B}$. Thus, if $\left(M(\forall z{:}\, s.B), \{q_r\}_{r:s\in\mathrm{LTerms}(\Sigma)}\right)$ is another cone on $\Sigma_B(z{:}\,s)$, then it induces an arrow $\phi\colon M(\forall z{:}\, s.B)\to M(\forall z{:}\, s.B)$ in $\mathbb{C}_{\forall z:s.B}$ such that $q_r = p_r\circ\phi$, because $M(\forall z{:}\, s.B)$ is an object of $\mathbb{C}_{\forall z:s.B}$ and, thus, there is an arrow from its cone to the terminal object's cone. Summing up, it follows that $\phi = 1_{M(\forall z:s.B)}$ and, thus, $q_r = p_r$.

It is worth noticing that $p_r = [\![w{:}(\forall z{:}\, s.B).\ \mathrm{allE}(w,r){:}\,B[r/z]]\!]$, as it immediately follows from the current case's definition. Thus, the $p_r$ arrows in the 'star-shaped' subcategory $\mathbb{C}_{\forall z:s.B}$ characterise the allE term constructor.

(15) $[\![\vec{x}.\, \mathrm{exI}_z(t)\colon(\exists z\colon s.\, B)]\!] = j_r \circ [\![\vec{x}.\, t\colon B[r/z]]\!]$ where $j_r\colon M(B[r/z]) \to M(\exists z\colon s.\, B)$ is the $r$-th injection of the unique co-cone on $\Sigma_B(z\colon s)$ whose vertex is $M(\exists z\colon s.\, B)$.

The definition is sound: there is a co-cone $\left(M(\exists z\colon s.\, B), \{j_r\}_{r\colon s\in\mathrm{LTerms}(\Sigma)}\right)$ on $\Sigma_B(z\colon s)$ since $M(\exists z\colon s.\, B)$ is the initial object of $\mathbb{C}_{M(\exists z\colon s.\, B)}$, as $\mathbb{C}$ with $M$ is logically distributive, see definition 4.1. For the same reason, there is a unique arrow $M(\exists z\colon s.\, B) \to M(\exists z\colon s.\, B)$ in $\mathbb{C}_{\exists z\colon s.\, B}$, which is $1_{M(\exists z\colon s.\, B)}$, by definition of $\mathbb{C}_{\exists z\colon s.\, B}$. Thus, if $\left(M(\exists z\colon s.\, B), \{q_r\}_{r\colon s\in\mathrm{LTerms}(\Sigma)}\right)$ is another co-cone on $\Sigma_B(z\colon s)$, then it induces an arrow $\phi\colon M(\exists z\colon s.\, B) \to M(\exists z\colon s.\, B)$ in $\mathbb{C}_{\exists z\colon s.\, B}$ such that $j_r = \phi \circ q_r$, because $M(\exists z\colon s.\, B)$ is an object of $\mathbb{C}_{\exists z\colon s.\, B}$ and, thus, there is an arrow from its co-cone to the initial object's co-cone. Summing up, it follows that $\phi = 1_{M(\exists z\colon s.\, B)}$ and, thus, $q_r = j_r$.

It is worth noticing that $j_r = [\![w\colon B[r/z].\, \mathrm{exI}_z(w)\colon(\exists z\colon s.\, B)]\!]$, as it immediately follows from the current case's definition. So, the $j_r$'s above are the arrows in the 'co-star-shaped' subcategory $\mathbb{C}_{\exists z\colon s.\, B}$ characterising the exI term constructor.

(16) $[\![\vec{x}.\, \mathrm{exE}(t, (\lambda z\colon s.\, r))\colon B]\!] = \gamma \circ \beta^{-1} \circ \left(\alpha, [\![\vec{x}.\, t\colon(\exists z\colon s.\, C)]\!]\right)$ where
   (a) $\alpha \equiv 1_{MA_{i_1}} \times \cdots \times 1_{MA_{i_k}}\colon MA \to MA'$ with $A' \equiv A_{i_1} \times \cdots \times A_{i_k}$, where $\vec{x}' \equiv \{x_{i_1}\colon A_{i_1}, \ldots, x_{i_k}\colon A_{i_k}\} = \mathrm{FV}(t\colon(\exists z\colon s.\, C)) \cup \mathrm{FV}(r\colon C \to B)$;
   (b) $\beta\colon M(\exists z\colon s.\, A' \times C) \to MA' \times M(\exists z\colon s.\, C)$ is the co-universal arrow in the subcategory $\mathbb{C}_{\exists z\colon s.\, A' \times C}$;
   (c) $\gamma\colon M(\exists z\colon s.\, A' \times C) \to MB$ is the co-universal arrow in $\mathbb{C}_{\exists z\colon s.\, A' \times C}$.

The definition is sound: in the first place, calling $j_u$ the $u$-th injection of the co-cone on $\Sigma_C(z\colon s)$ with vertex $M(\exists z\colon s.\, C)$, which is unique as seen in the previous case,

$$\left(MA' \times M(\exists z\colon s.\, C), \{1_{MA'} \times j_u\}_{u\colon s\in\mathrm{LTerms}(\Sigma)}\right)$$

is a co-cone on $\Sigma_{A'\times C}(z\colon s)$ since $\mathrm{dom}(1_{MA'} \times j_u) = MA' \times M(C[u/z]) = M((A' \times C)[u/z])$ as $z\colon s \notin \mathrm{FV}(A')$ by construction of the term $\mathrm{exE}(t, (\lambda z\colon s.\, r))$, see definition 3.3. So, $MA' \times M(\exists z\colon s.\, C)$ is an object of $\mathbb{C}_{\exists z\colon s.\, A'\times C}$ and, thus, there is an unique arrow $\beta$ in $\mathbb{C}_{\exists z\colon s.\, A'\times C}$ from it to the unique co-cone on $\Sigma_{A'\times C}(z\colon s)$ whose vertex is $M(\exists z\colon s.\, A' \times C)$. Moreover, $\beta$ has an inverse, since $\mathbb{C}$ with $M$ is logically distributive (last requirement).

In the second place, we notice that $z\colon s \notin \mathrm{FV}(B)$ and

$$\left(MB, \left\{\mathrm{ev} \circ \left([\![\vec{x}'.\, r[u/z]\colon C[u/z] \to B]\!] \times 1_{M(C[u/z])}\right)\right\}_{u\colon s\in\mathrm{LTerms}(\Sigma)}\right)$$

is a co-cone on $\Sigma_{A'\times C}(z\colon s)$, because $\vec{x}'.\, r\colon C \to B$ is a term-in-context, and so is $\vec{x}'[u/z].\, r[u/z]\colon C[u/z] \to B[u/z]$ for any $u\colon s \in \mathrm{LTerms}(\Sigma)$. But $z\colon s \notin \mathrm{FV}(\vec{x}') = \mathrm{FV}(A')$, so $\vec{x}'[u/z].\, r[u/z]\colon C[u/z] \to B[u/z] = \vec{x}'.\, r[u/z]\colon C[u/z] \to B$. Hence, the arrow $[\![\vec{x}'.\, r[u/z]\colon C[u/z] \to B]\!]\colon MA' \to MB^{M(C[u/z])}$ exists. It is then immediate to see that the exponential evaluation yields an arrow $MA' \times M(C[u/z]) \to MB$. So $MB$ is an object of $\mathbb{C}_{\exists z\colon s.\, A'\times C}$, thus $\gamma$ is uniquely defined.

Figure 8 shows the commutative diagram illustrating this last case.

*Note* 4.9. Semantics of $\lambda$-terms is well-founded, as we proved in the doubtful cases, and it uses all the conditions we imposed in definition 4.1 on logically distributive categories, showing that those conditions are needed.

**Definition 4.4** (Semantics of equalities)**.** With notation as above, the *semantics of an equality-in-context* associates to each $\vec{x}.\, s =_B t$ a value in $\{0, 1\}$ and, precisely,

$$[\![\vec{x}.\, s =_B t]\!] = \begin{cases} 1 & \text{if } [\![\vec{x}.\, s\colon B]\!] = [\![\vec{x}.\, t\colon B]\!] \\ 0 & \text{otherwise}\ . \end{cases}$$

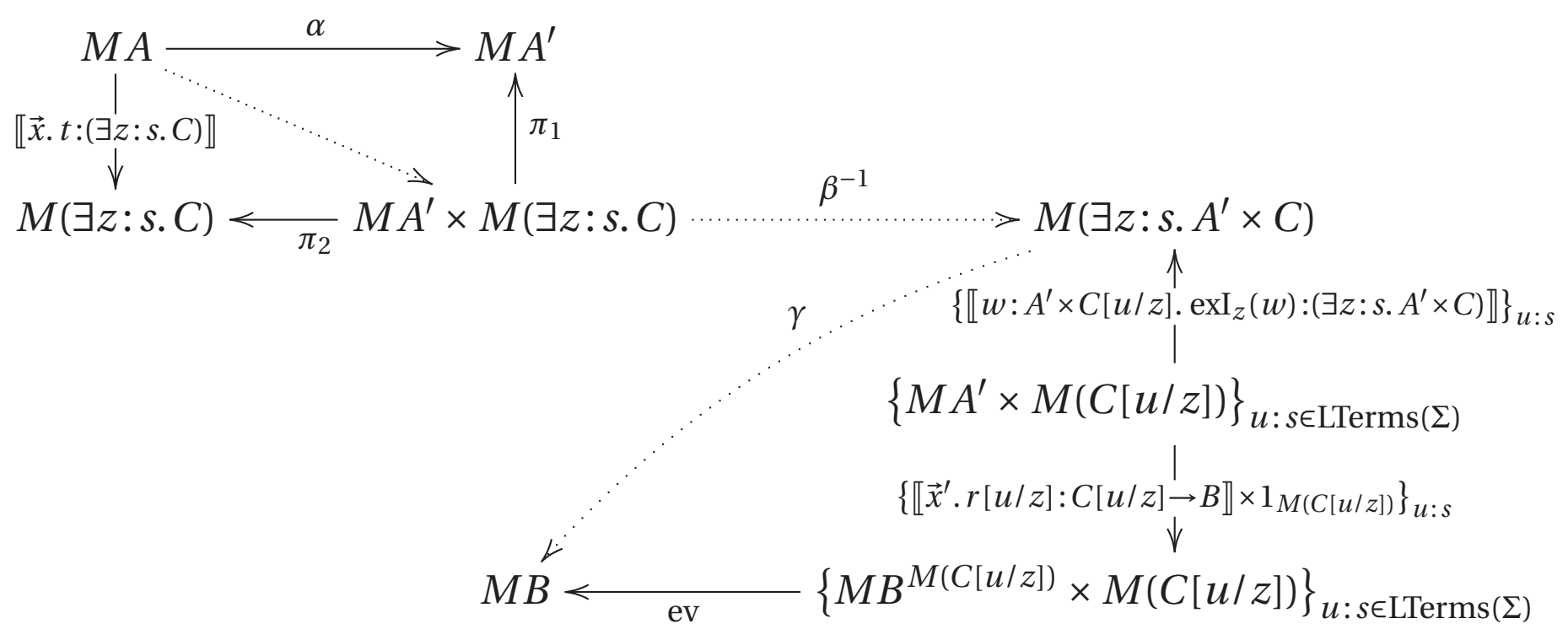

FIGURE 8. The semantics of $\mathrm{exE}(t, (\lambda z : s.\, r))$.

## 5. SOUNDNESS

In the following, when not otherwise stated, we assume to have a fixed $\Sigma$-structure whose elements are denoted as in the preceding sections.

First, we need to take care of $\alpha$-renaming in contexts.

**Proposition 5.1.** $[\![\vec{x}, z\colon B.\, u\colon C]\!] = [\![\vec{x}, w\colon B.\, u[w/z]\colon C]\!]$, *where* $w\colon A \notin \mathrm{FV}(u\colon C)$.

*Proof.* By induction on the structure of the term $t$:

(1) $t \equiv y \in W_B$: there are two sub-cases:
  (a) $y \equiv x_i$: $[\![\vec{x}, w\colon B.\, x_i[w/z]\colon A_i]\!] = [\![\vec{x}, w\colon B.\, x_i\colon A_i]\!] = \pi_i = [\![\vec{x}, z\colon B.\, x_i\colon A_i]\!]$;
  (b) $y \equiv z$: $[\![\vec{x}, w\colon B.\, z[w/z]\colon B]\!] = [\![\vec{x}, w\colon B.\, w\colon B]\!] = \pi_{n+1} = [\![\vec{x}, z\colon B.\, z\colon B]\!]$.

(2) $t \equiv a(u)$, $a\colon D \to C \in \mathrm{Ax}$:
$$\begin{aligned} & [\![\vec{x}, w\colon B.\, a(u)[w/z]\colon C]\!] \\ =& [\![\vec{x}, w\colon B.\, a(u[w/z])\colon C]\!] \\ =& M_{\mathrm{Ax}} a \circ [\![\vec{x}, w\colon B.\, u[w/z]\colon D]\!] \end{aligned}$$
and, by induction hypothesis,
$$\begin{aligned} =& M_{\mathrm{Ax}} a \circ [\![\vec{x}, z\colon B.\, u\colon D]\!] \\ =& [\![\vec{x}, z\colon B.\, a(u)\colon C]\!] \ ; \end{aligned}$$

(3) $t \equiv \langle u, v\rangle$:
$$\begin{aligned} & [\![\vec{x}, w\colon B.\, \langle u, v\rangle[w/z]\colon C_1 \times C_2]\!] \\ =& [\![\vec{x}, w\colon B.\, \langle u[w/z], v[w/z]\rangle\colon C_1 \times C_2]\!] \\ =& \left([\![\vec{x}, w\colon B.\, u[w/z]\colon C_1]\!], [\![\vec{x}, w\colon B.\, v[w/z]\colon C_2]\!]\right) \end{aligned}$$
and, by induction hypothesis,
$$\begin{aligned} =& \left([\![\vec{x}, z\colon B.\, u\colon C_1]\!], [\![\vec{x}, z\colon B.\, v\colon C_2]\!]\right) \\ =& [\![\vec{x}, z\colon B.\, \langle u, v\rangle\colon C_1 \times C_2]\!] \ ; \end{aligned}$$

(4) $t \equiv \mathrm{fst}(u)$:
$$\begin{aligned} & [\![\vec{x}, w\colon B.\, \mathrm{fst}(u)[w/z]\colon C_1]\!] \\ =& [\![\vec{x}, w\colon B.\, \mathrm{fst}(u[w/z])\colon C_1]\!] \\ =& \pi_1 \circ [\![\vec{x}, w\colon B.\, u[w/z]\colon C_1 \times C_2]\!] \end{aligned}$$
and, by induction hypothesis,
$$\begin{aligned} =& \pi_1 \circ [\![\vec{x}, z\colon B.\, u\colon C_1 \times C_2]\!] \\ =& [\![\vec{x}, z\colon B.\, \mathrm{fst}(u)\colon C_1]\!] \ ; \end{aligned}$$

(5) $t \equiv \mathrm{snd}(u)$:
$$\begin{aligned} & [\![\vec{x}, w\colon B.\, \mathrm{snd}(u)[w/z]\colon C_2]\!] \\ =& [\![\vec{x}, w\colon B.\, \mathrm{snd}(u[w/z])\colon C_2]\!] \\ =& \pi_2 \circ [\![\vec{x}, w\colon B.\, u[w/z]\colon C_1 \times C_2]\!] \end{aligned}$$
and, by induction hypothesis,
$$\begin{aligned} =& \pi_2 \circ [\![\vec{x}, z\colon B.\, u\colon C_1 \times C_2]\!] \\ =& [\![\vec{x}, z\colon B.\, \mathrm{snd}(u)\colon C_2]\!] \ ; \end{aligned}$$

(6) $t \equiv (\lambda y\colon C_1.\, u)$:
$$\begin{aligned} & [\![\vec{x}, w\colon B.\, (\lambda y\colon C_1.\, u)[w/z]\colon C_1 \to C_2]\!] \\ =& [\![\vec{x}, w\colon B.\, (\lambda v\colon C_1.\, u[v/y][w/z])\colon C_1 \to C_2]\!] \end{aligned}$$
where $v \notin \mathrm{FV}(u\colon C_2) \cup \{w\colon B\}$. This expression is the exponential transpose of
$$[\![\vec{x}, w\colon B, v\colon C_1.\, u[v/y][w/z]\colon C_2]\!]$$

by induction hypothesis,

$$= [\![\vec{x}, z : B, v : C_1 .\, u[v/y] : C_2]\!]$$

applying the induction hypothesis once again,

$$= [\![\vec{x}, z : B, y : C_1 .\, u : C_2]\!]$$

which is the exponential transpose of $[\![\vec{x}, z : B.\, (\lambda y : C_1 .\, u) : C_1 \to C_2]\!]$;

(7) $t \equiv u \cdot v$:
$$\begin{aligned} & [\![\vec{x}, w : B.\, (u \cdot v)[w/z] : C]\!] \\ =& [\![\vec{x}, w : B.\, u[w/z] \cdot v[w/z] : C]\!] \\ =& \mathrm{ev} \circ \big([\![\vec{x}, w : B.\, u[w/z] : D \to C]\!], [\![\vec{x}, w : B.\, v[w/z] : D]\!]\big) \end{aligned}$$

and, by induction hypothesis,

$$\begin{aligned} =& \mathrm{ev} \circ \big([\![\vec{x}, z : B.\, u : D \to C]\!], [\![\vec{x}, z : B.\, v : D]\!]\big) \\ =& [\![\vec{x}, z : B.\, u \cdot v : C]\!] \;; \end{aligned}$$

(8) $t \equiv \mathrm{inl}_D(u)$:
$$\begin{aligned} & [\![\vec{x}, w : B.\, \mathrm{inl}_D(u)[w/z] : C + D]\!] \\ =& [\![\vec{x}, w : B.\, \mathrm{inl}_D(u[w/z]) : C + D]\!] \\ =& \iota_1 \circ [\![\vec{x}, w : B.\, u[w/z] : C]\!] \end{aligned}$$

and, by induction hypothesis,

$$\begin{aligned} =& \iota_1 \circ [\![\vec{x}, z : B.\, u : C]\!] \\ =& [\![\vec{x}, z : B.\, \mathrm{inl}_D(u) : C + D]\!] \;; \end{aligned}$$

(9) $t \equiv \mathrm{inr}_C(u)$:
$$\begin{aligned} & [\![\vec{x}, w : B.\, \mathrm{inr}_C(u)[w/z] : C + D]\!] \\ =& [\![\vec{x}, w : B.\, \mathrm{inr}_C(u[w/z]) : C + D]\!] \\ =& \iota_2 \circ [\![\vec{x}, w : B.\, u[w/z] : C]\!] \end{aligned}$$

and, by induction hypothesis,

$$\begin{aligned} =& \iota_2 \circ [\![\vec{x}, z : B.\, u : C]\!] \\ =& [\![\vec{x}, z : B.\, \mathrm{inr}_C(u) : C + D]\!] \;; \end{aligned}$$

(10) $t \equiv \mathrm{when}(r, u, v)$:
$$\begin{aligned} & [\![\vec{x}, w : B.\, \mathrm{when}(r, u, v)[w/z] : C]\!] \\ =& [\![\vec{x}, w : B.\, \mathrm{when}(r[w/z], u[w/z], v[w/z]) : C]\!] \\ =& \big[\mathrm{ev} \circ \big([\![\vec{x}, w : B.\, u[w/z] : C_1 \to C]\!] \times 1_{MC_1}\big), \\ & \quad \mathrm{ev} \circ \big([\![\vec{x}, w : B.\, v[w/z] : C_2 \to C]\!] \times 1_{MC_2}\big)\big] \circ \\ & \circ \Delta^{-1} \circ \big(1_{MA \times MB}, [\![\vec{x}, w : B.\, r[w/z] : C_1 + C_2]\!]\big) \end{aligned}$$

and, by induction hypothesis,

$$\begin{aligned} =& \big[\mathrm{ev} \circ \big([\![\vec{x}, z : B.\, u : C_1 \to C]\!] \times 1_{MC_1}\big), \\ & \quad \mathrm{ev} \circ \big([\![\vec{x}, z : B.\, v : C_2 \to C]\!] \times 1_{MC_2}\big)\big] \circ \\ & \circ \Delta^{-1} \circ \big(1_{MA \times MB}, [\![\vec{x}, w : B.\, r : C_1 + C_2]\!]\big) \\ =& [\![\vec{x}, z : B.\, \mathrm{when}(r, u, v) : C]\!] \;; \end{aligned}$$

(11) $t \equiv *$: $[\![\vec{x}, w : B.\, *[w/z] : 1]\!] = [\![\vec{x}, w : B.\, * : 1]\!] = !\colon MA \times MB \to 1 = [\![\vec{x}, z : B.\, * : 1]\!]$;

(12) $t \equiv \mathrm{F}_D$: $[\![\vec{x}, w : B.\, \mathrm{F}[w/z] : 0 \to D]\!] = [\![\vec{x}, w : B.\, \mathrm{F} : 0 \to D]\!]$ is the exponential transpose of $(!\colon 0 \to MD) \circ (\pi_{n+2}\colon MA \times MD \times 0 \to 0)$, that is, $[\![\vec{x}, z : B.\, \mathrm{F} : 0 \to D]\!]$;

(13) $t : \mathrm{allI}(\lambda y : s.\, u)$ $\quad [\![\vec{x}, w : B.\, \mathrm{allI}(\lambda y : s.\, u)[w/z] : (\forall y : s.\, D)]\!]$

$$= [\![\vec{x}, w : B.\, \mathrm{allI}(\lambda y : s.\, u[w/z]) : (\forall y : s.\, D)]\!]$$
$$= \beta_w \circ \alpha_w \ ,$$

where $\alpha_w \equiv 1_{MA^w_{i_1}} \times \cdots \times 1_{MA^w_{i_k}} : MA \times MB \to MA'_w$ with $A'_w \equiv A^w_{i_1} \times \cdots \times A^w_{i_k}$ and $\vec{x}'_w \equiv \{x_{i_1} : A^w_{i_1}, \ldots, x_{i_k} : A^w_{i_k}\} = \mathrm{FV}(u[w/z] : D)$, and $\beta_w : MA'_w \to M(\forall y : s.\, D)$, the universal arrow from $MA'_w$ to the terminal object in $\mathbb{C}_{\forall y : s.\, D}$. Analogously

$$[\![\vec{x}, z : B.\, \mathrm{allI}(\lambda y : s.\, u) : (\forall y : s.\, D)]\!] = \beta_z \circ \alpha_z \ ,$$

where $\alpha_z \equiv 1_{MA^z_{i_1}} \times \cdots \times 1_{MA^z_{i_k}} : MA \times MB \to MA'_z$ with $A'_z \equiv A^z_{i_1} \times \cdots \times A^z_{i_k}$ and $\vec{x}'_z \equiv \{x_{i_1} : A^z_{i_1}, \ldots, x_{i_k} : A^z_{i_k}\} = \mathrm{FV}(u : D)$, and $\beta_z : MA'_z \to M(\forall y : s.\, D)$, the universal arrow from $MA'_z$ to the terminal object in $\mathbb{C}_{\forall y : s.\, D}$.

But $\mathrm{FV}(u : D)$ differs from $\mathrm{FV}(u[w/z] : D)$ only by the $z : B$ variable instead of $w : B$, thus $MA'_z = MA'_w$ and $\alpha_z = \alpha_w$. Moreover, since $MA'_z = MA'_w$, $\beta_z = \beta_w$ as universal arrows are unique;

(14) $t \equiv \mathrm{allE}(u, r)$: $\quad [\![\vec{x}, w : B.\, \mathrm{allE}(u, r)[w/z] : D[r/z]]\!]$

$$= [\![\vec{x}, w : B.\, \mathrm{allE}(u[w/z], r) : D[r/z]]\!]$$
$$= (p_r : M(\forall y : s.\, D) \to M(D[r/y])) \circ [\![\vec{x}, w : B.\, u[w/z] : (\forall y : s.\, D)]\!]$$

and, by induction hypothesis,

$$= (p_r : M(\forall y : s.\, D) \to M(D[r/y])) \circ [\![\vec{x}, z : B.\, u : (\forall y : s.\, D)]\!]$$
$$= [\![\vec{x}, z : B.\, \mathrm{allE}(u, r) : D[r/z]]\!] \ ;$$

(15) $t \equiv \mathrm{exI}_y(u)$: $\quad [\![\vec{x}, w : B.\, \mathrm{exI}_y(u)[w/z] : (\exists y : s.\, D)]\!]$

$$= [\![\vec{x}, w : B.\, \mathrm{exI}_y(u[w/z]) : (\exists y : s.\, D)]\!]$$
$$= (j_r : M(D[r/y]) \to M(\exists y : s.\, D)) \circ [\![\vec{x}, w : B.\, u[w/z] : D[r/y]]\!]$$

and, by induction hypothesis,

$$= (j_r : M(D[r/y]) \to M(\exists y : s.\, D)) \circ [\![\vec{x}, z : B.\, u : D[r/y]]\!]$$
$$= [\![\vec{x}, z : B.\, \mathrm{exI}_y(u) : (\exists y : s.\, D)]\!] \ ;$$

(16) $t \equiv \mathrm{exE}(u, (\lambda y : s.\, v))$: $\quad [\![\vec{x}, w : B.\, \mathrm{exE}(u, (\lambda y : s.\, v))[w/z] : C]\!]$

$$= [\![\vec{x}, w : B.\, \mathrm{exE}(u[w/z], (\lambda y : s.\, v[w/z])) : C]\!]$$
$$= \gamma_w \circ \beta_w^{-1} \circ \big(\alpha_w, [\![\vec{x}, w : B.\, u[w/z] : (\exists y : s.\, D)]\!]\big)$$

and, by induction hypothesis,

$$= \gamma_w \circ \beta_w^{-1} \circ \big(\alpha_w, [\![\vec{x}, z : B.\, u : (\exists y : s.\, D)]\!]\big) \ ,$$

where

(a) $\alpha_w \equiv 1_{MA^w_{i_1}} \times \cdots \times 1_{MA^w_{i_k}} : MA \times MB \to MA'_w$ with $A'_w \equiv A^w_{i_1} \times \cdots \times A^w_{i_k}$ and $\vec{x}'_w \equiv \{x_{i_1} : A^w_{i_1}, \ldots, x_{i_k} : A^w_{i_k}\} = \mathrm{FV}(u[w/z] : (\exists y : s.\, D)) \cup \mathrm{FV}(v[w/z] : D \to C)$;

(b) $\beta_w : M(\exists y : s.\, A'_w \times D) \to MA'_w \times M(\exists y : s.\, D)$ is the co-universal arrow in the subcategory $\mathbb{C}_{\exists y : s.\, A'_w \times D}$;

(c) $\gamma_w : M(\exists y : s.\, A'_w \times D) \to MC$ is the co-universal arrow in $\mathbb{C}_{\exists y : s.\, A'_w \times D}$.

Analogously, $\quad [\![\vec{x}, z : B.\, \mathrm{exE}(u, (\lambda y : s.\, v)) : C]\!]$

$$= \gamma_z \circ \beta_z^{-1} \circ \big(\alpha_z, [\![\vec{x}, z : B.\, u : (\exists y : s.\, D)]\!]\big) \ ,$$

where

(a) $\alpha_z \equiv 1_{MA^z_{i_1}} \times \cdots \times 1_{MA^z_{i_k}} : MA \times MB \to MA'_z$ with $A'_z \equiv A^z_{i_1} \times \cdots \times A^z_{i_k}$ and $\vec{x}'_z \equiv \{x_{i_1} : A^z_{i_1}, \ldots, x_{i_k} : A^z_{i_k}\} = \mathrm{FV}(u : (\exists y : s.\, D)) \cup \mathrm{FV}(v : D \to C)$;

(b) $\beta_z\colon M(\exists y\colon s.\, A'_z \times D) \to MA'_z \times M(\exists y\colon s.\, D)$ is the co-universal arrow in the subcategory $\mathbb{C}_{\exists y\colon s.\, A'_z \times D}$;

(c) $\gamma_z\colon M(\exists y\colon s.\, A'_z \times D) \to MC$ is the co-universal arrow in $\mathbb{C}_{\exists y\colon s.\, A'_z \times D}$.

But $\mathrm{FV}(u{:}(\exists y\colon s.\, D))$ differs from $\mathrm{FV}(u[w/z]{:}(\exists y\colon s.\, D))$ only by the $z{:}B$ variable instead of $w{:}B$, thus $MA'_z = MA'_w$ and $\alpha_z = \alpha_w$. Moreover, since $MA'_z = MA'_w$, $\beta_z = \beta_w$ and $\gamma_z = \gamma_w$ as co-universal arrows are unique. □

As a second step, we want to show that terms equivalent modulo renaming of bound variables, get interpreted in the same way

**Proposition 5.2.** *If $\vec{x}.\, t{:}A$ and $\vec{x}.\, s{:}A$ are two terms-in-context and $t$ is $\alpha$-equivalent to $s$ ($t =_\alpha s$), then $[\![\vec{x}.\, t{:}A]\!] = [\![\vec{x}.\, s{:}A]\!]$.*

*Proof.* By induction on the structure of the term $t$:

(1) $t \equiv y \in W_A$: by definition of $\alpha$-equivalence $s \equiv y$, so $s \equiv t$;

(2) $t \equiv f(u)$, $f\colon B \to C \in Ax$: by definition of $\alpha$-equivalence, $s \equiv f(v)$ and $u =_\alpha v$. By induction hypothesis, $[\![\vec{x}.\, u{:}B]\!] = [\![\vec{x}.\, v{:}B]\!]$, so $[\![\vec{x}.\, f(u){:}C]\!] = M_{\mathrm{Ax}} f \circ [\![\vec{x}.\, u{:}B]\!] = M_{\mathrm{Ax}} f \circ [\![\vec{x}.\, v{:}B]\!] = [\![\vec{x}.\, f(v){:}C]\!]$;

(3) $t \equiv \langle u, v\rangle$: by definition of $\alpha$-equivalence $s \equiv \langle u', v'\rangle$ and $u =_\alpha u'$, $v =_\alpha v'$. By induction hypothesis, $[\![\vec{x}.\, u{:}B]\!] = [\![\vec{x}.\, u'{:}B]\!]$, $[\![\vec{x}.\, v{:}C]\!] = [\![\vec{x}.\, v'{:}C]\!]$, so $[\![\vec{x}.\, \langle u, v\rangle{:}B \times C]\!] = \big([\![\vec{x}.\, u{:}B]\!], [\![\vec{x}.\, v{:}C]\!]\big) = \big([\![\vec{x}.\, u'{:}B]\!], [\![\vec{x}.\, v'{:}C]\!]\big) = [\![\vec{x}.\, \langle u', v'\rangle{:}B \times C]\!]$;

(4) $t \equiv \mathrm{fst}(u)$: by definition of $\alpha$-equivalence $s \equiv \mathrm{fst}(v)$ and $u =_\alpha v$. By induction hypothesis, $[\![\vec{x}.\, u{:}B \times C]\!] = [\![\vec{x}.\, v{:}B \times C]\!]$, so $[\![\vec{x}.\, \mathrm{fst}(u){:}B]\!] = \pi_1 \circ [\![\vec{x}.\, u{:}B \times C]\!] = \pi_1 \circ [\![\vec{x}.\, v{:}B \times C]\!] = [\![\vec{x}.\, \mathrm{fst}(v){:}B]\!]$;

(5) $t \equiv \mathrm{snd}(u)$: by definition of $\alpha$-equivalence $s \equiv \mathrm{snd}(v)$ and $u =_\alpha v$. By induction hypothesis, $[\![\vec{x}.\, u{:}B \times C]\!] = [\![\vec{x}.\, v{:}B \times C]\!]$, so $[\![\vec{x}.\, \mathrm{snd}(u){:}C]\!] = \pi_2 \circ [\![\vec{x}.\, u{:}B \times C]\!] = \pi_2 \circ [\![\vec{x}.\, v{:}B \times C]\!] = [\![\vec{x}.\, \mathrm{snd}(v){:}C]\!]$;

(6) $t \equiv (\lambda y{:}B.\, u)$: by definition of $\alpha$-equivalence, $s \equiv (\lambda y{:}B.\, v)$ and $u[z/y] = v[z/y]$, with $z{:}B$ a new variable. By proposition 5.1, $[\![\vec{x}, y{:}B.\, u{:}C]\!] = [\![\vec{x}, z{:}B.\, u[z/y]{:}C]\!]$; by induction hypothesis, $[\![\vec{x}, z{:}B.\, u[z/y]{:}C]\!] = [\![\vec{x}, z{:}B.\, v[z/y]{:}C]\!]$; by proposition 5.1 again, $[\![\vec{x}, z{:}B.\, v[z/y]{:}C]\!] = [\![\vec{x}, y{:}B.\, v{:}C]\!]$. Thus, the exponential transpose of $[\![\vec{x}, y{:}B.\, u{:}C]\!]$ equals the exponential transpose of $[\![\vec{x}, y{:}B.\, v{:}C]\!]$, that is $[\![\vec{x}.\, (\lambda y{:}B.\, u){:}B \to C]\!] = [\![\vec{x}.\, (\lambda y{:}B.\, v){:}B \to C]\!]$;

(7) $t \equiv u \cdot v$: by definition of $\alpha$-equivalence $s \equiv u' \cdot v'$ and $u =_\alpha u'$, $v =_\alpha v'$. By induction hypothesis, $[\![\vec{x}.\, u{:}B \to C]\!] = [\![\vec{x}.\, u'{:}B \to C]\!]$ and $[\![\vec{x}.\, v{:}B]\!] = [\![\vec{x}.\, v'{:}B]\!]$, so $[\![\vec{x}.\, u \cdot v{:}C]\!] = \mathrm{ev} \circ \big([\![\vec{x}.\, u{:}B \to C]\!], [\![\vec{x}.\, v{:}B]\!]\big) = \mathrm{ev} \circ \big([\![\vec{x}.\, u'{:}B \to C]\!], [\![\vec{x}.\, v'{:}B]\!]\big) = [\![\vec{x}.\, u' \cdot v'{:}C]\!]$;

(8) $t \equiv \mathrm{inl}_C(u)$: by definition of $\alpha$-equivalence $s \equiv \mathrm{inl}_C(v)$ and $u =_\alpha v$. By induction hypothesis, $[\![\vec{x}.\, u{:}B]\!] = [\![\vec{x}.\, v{:}B]\!]$, thus $[\![\vec{x}.\, \mathrm{inl}_C(u){:}B + C]\!] = \iota_1 \circ [\![\vec{x}.\, u{:}B]\!] = \iota_1 \circ [\![\vec{x}.\, v{:}B]\!] = [\![\vec{x}.\, \mathrm{inl}_C(v){:}B + C]\!]$;

(9) $t \equiv \mathrm{inr}_C(u)$: by definition of $\alpha$-equivalence $s \equiv \mathrm{inr}_C(v)$ and $u =_\alpha v$. By induction hypothesis, $[\![\vec{x}.\, u{:}B]\!] = [\![\vec{x}.\, v{:}B]\!]$, thus $[\![\vec{x}.\, \mathrm{inr}_C(u){:}C + B]\!] = \iota_2 \circ [\![\vec{x}.\, u{:}B]\!] = \iota_2 \circ [\![\vec{x}.\, v{:}B]\!] = [\![\vec{x}.\, \mathrm{inr}_C(v){:}C + B]\!]$;

(10) $t \equiv \mathrm{when}(r, u, v)$: by definition of $\alpha$-equivalence $s \equiv \mathrm{when}(r', u', v')$ and $r =_\alpha r'$, $u =_\alpha u'$ and $v =_\alpha v'$. By induction hypothesis, $[\![\vec{x}.\, r{:}B + C]\!] = [\![\vec{x}.\, r'{:}B + C]\!]$, $[\![\vec{x}.\, u{:}B \to D]\!] = [\![\vec{x}.\, u'{:}B \to D]\!]$ and $[\![\vec{x}.\, v{:}C \to D]\!] = [\![\vec{x}.\, v'{:}C \to D]\!]$, so:

$$\begin{aligned}
&[\![\vec{x}.\, \mathrm{when}(r, u, v){:}D]\!] \\
&= \big[\mathrm{ev} \circ \big([\![\vec{x}.\, u{:}B \to D]\!] \times 1_{MB}\big), \\
&\quad \mathrm{ev} \circ \big([\![\vec{x}.\, v{:}C \to D]\!] \times 1_{MC}\big)\big] \\
&\quad \circ \Delta^{-1} \circ \big(1_{MA}, [\![\vec{x}.\, r{:}B + C]\!]\big)
\end{aligned}$$

and, applying the induction hypothesis,

$$\begin{aligned}&= \big[\mathrm{ev}\circ\big([\![\vec{x}.\,u' : B \to D]\!] \times 1_{MB}\big),\\ &\quad \mathrm{ev}\circ\big([\![\vec{x}.\,v' : C \to D]\!] \times 1_{MC}\big)\big]\\ &\quad \circ \Delta^{-1} \circ \big(1_{MA}, [\![\vec{x}.\,r' : B + C]\!]\big)\\ &= [\![\vec{x}.\,\mathrm{when}(r', u', v') : D]\!]\ ;\end{aligned}$$

(11) $t \equiv *$: by definition of $\alpha$-equivalence $s \equiv *$, so $s \equiv t$;

(12) $t \equiv \mathrm{F}_B$: by definition of $\alpha$-equivalence $s \equiv \mathrm{F}_B$, so $s \equiv t$;

(13) $t \equiv \mathrm{allI}(\lambda x{:}\sigma.\,u)$: by definition of $\alpha$-equivalence $s \equiv \mathrm{allI}(\lambda y{:}\sigma.\,v)$ and $u[z/x] =_\alpha v[z/y]$ with $z{:}\sigma$ a new variable. Thus, for all $r{:}\sigma \in \mathrm{LTerms}(\Sigma)$, $u[z/x][r/z] =_\alpha v[z/y][r/z]$, that is, $u[r/x] =_\alpha v[r/y]$. So, by induction hypothesis,

$$[\![\vec{x}'.\,u[r/x] : B[r/x]]\!] = [\![\vec{x}'.\,v[r/x] : B'[r/x]]\!]\ ,$$

where $\vec{x}' = \mathrm{FV}(u{:}B) = \mathrm{FV}(v{:}B')$.

So, $[\![\vec{x}.\,\mathrm{allI}(\lambda x{:}\sigma.\,u) : (\forall x{:}\sigma.\,B)]\!] = \beta_x \circ \alpha$, with $\alpha \equiv 1 \times \cdots \times 1 \colon MA \to MA'$, the product corresponding to $\vec{x}'$, and $\beta_x \colon MA' \to M(\forall x{:}\sigma.\,B)$ universal in $\mathbb{C}_{\forall x:\sigma.\,B}$. In particular, $\beta_x$ is an arrow from the cone

$$\Big(MA', \big\{[\![\vec{x}'.\,u[r/x] : B[r/x]]\!]\big\}_{r:\sigma\in\mathrm{LTerms}(\Sigma)}\Big)$$

to the cone $\Big(M(\forall x{:}\sigma.\,B), \big\{p_r\big\}_{r:\sigma\in\mathrm{LTerms}(\Sigma)}\Big)$.

Also, $[\![\vec{x}.\,\mathrm{allI}(\lambda y{:}\sigma.\,v) : (\forall x{:}\sigma.\,B)]\!] = \beta_y \circ \alpha$, with $\alpha$ as before, and $\beta_y \colon MA' \to M(\forall y{:}\sigma.\,B')$, where $M(\forall y{:}\sigma.\,B') = M(\forall x{:}\sigma.\,B)$, and $\beta_y$ is universal in $\mathbb{C}_{\forall x:\sigma.\,B}$. In particular, $\beta_y$ is an arrow from the cone

$$\Big(MA', \big\{[\![\vec{x}'.\,v[r/x] : B'[r/x]]\!]\big\}_{r:\sigma\in\mathrm{LTerms}(\Sigma)}\Big)$$

to the cone $\Big(M(\forall x{:}\sigma.\,B), \big\{p_r\big\}_{r:\sigma\in\mathrm{LTerms}(\Sigma)}\Big)$.

We notice that the co-domains of $\beta_x$ and $\beta_y$ as arrows in the category of cones are the same. So, being both universal, $\beta_x = \beta_y$;

(14) $t \equiv \mathrm{allE}(u, r)$: by definition of $\alpha$-equivalence $s \equiv \mathrm{allE}(v, r)$ and $u =_\alpha v$. By induction hypothesis, $[\![\vec{x}.\,u : (\forall x{:}\sigma.\,B)]\!] = [\![\vec{x}.\,v : (\forall x{:}\sigma.\,B)]\!]$, so $[\![\vec{x}.\,\mathrm{allE}(u, r) : B[r/x]]\!] = p_r \circ [\![\vec{x}.\,u : (\forall x{:}\sigma.\,B)]\!] = p_r \circ [\![\vec{x}.\,v : (\forall x{:}\sigma.\,B)]\!] = [\![\vec{x}.\,\mathrm{allE}(v, r) : B[r/x]]\!]$;

(15) $t \equiv \mathrm{exI}_x(u)$: by definition of $\alpha$-equivalence $s \equiv \mathrm{exI}_x(v)$ and $u =_\alpha v$. By induction hypothesis, $[\![\vec{x}.\,u : B[r/x]]\!] = [\![\vec{x}.\,v : B[r/x]]\!]$, so $[\![\vec{x}.\,\mathrm{exI}_x(u) : (\exists x{:}\sigma.\,B)]\!] = j_r \circ [\![\vec{x}.\,u : B[r/x]]\!] = j_r \circ [\![\vec{x}.\,v : B[r/x]]\!] = [\![\vec{x}.\,\mathrm{exI}_x(v) : (\exists x{:}\sigma.\,B)]\!]$;

(16) $t \equiv \mathrm{exE}(u, (\lambda x{:}\sigma.\,v))$: by definition of $\alpha$-equivalence $s \equiv \mathrm{exE}(u', (\lambda y{:}\sigma.\,v'))$, $u =_\alpha u'$, and $v[z/x] =_\alpha v'[z/y]$ with $z{:}\sigma$ a new variable. We notice that

$$\mathrm{FV}(u : (\exists x{:}\sigma.\,B)) \cup \mathrm{FV}(v : B \to C) = \mathrm{FV}(u' : (\exists y{:}\sigma.\,B')) \cup \mathrm{FV}(v' : B' \to C)\ .$$

Then, $[\![\vec{x}.\,\mathrm{exE}(u, (\lambda x{:}\sigma.\,v)) : C]\!] = \gamma_x \circ \beta_x^{-1} \circ \big(\alpha, [\![\vec{x}.\,u : (\exists x{:}\sigma.\,B)]\!]\big)$ where

(a) $\alpha = 1 \times \cdots \times 1 \colon MA \to MA'$, the product related to $\vec{x}' = \mathrm{FV}(u : (\exists x{:}\sigma.\,B)) \cup \mathrm{FV}(v : B \to C)$ as for definition 4.3;

(b) $\beta_x \colon M(\exists x{:}\sigma.\,A' \times B) \to MA' \times M(\exists x{:}\sigma.\,B)$ co-universal in $\mathbb{C}_{\exists x:\sigma.\,A' \times B}$;

(c) $\gamma_x \colon M(\exists x{:}\sigma.\,A' \times B) \to MC$ co-universal in $\mathbb{C}_{\exists x:\sigma.\,A' \times B}$.

Analogously, $[\![\vec{x}.\,\mathrm{exE}(u', (\lambda y{:}\sigma.\,v')) : C]\!] = \gamma_y \circ \beta_y^{-1} \circ \big(\alpha_y, [\![\vec{x}.\,u' : (\exists y{:}\sigma.\,B')]\!]\big)$ where

(a) $\alpha = 1 \times \cdots \times 1 \colon MA \to MA'$, the product related to $\vec{x}' = \mathrm{FV}(u' : (\exists y{:}\sigma.\,B')) \cup \mathrm{FV}(v' : B' \to C)$ as for definition 4.3—notice how $\alpha_y$ equals $\alpha$.

(b) $\beta_y \colon M(\exists y{:}\sigma.\,A' \times B') \to MA' \times M(\exists y{:}\sigma.\,B')$ co-universal in $\mathbb{C}_{\exists y:\sigma.\,A' \times B'}$;

(c) $\gamma_y \colon M(\exists y{:}\sigma.\,A' \times B') \to MC$ co-universal in $\mathbb{C}_{\exists y:\sigma.\,A' \times B'}$.

Since $M(\exists x{:}\sigma.\,B) = M(\exists y{:}\sigma.\,B')$, $[\![\vec{x}.\,u : (\exists x{:}\sigma.\,B)]\!] = [\![\vec{x}.\,u' : (\exists y{:}\sigma.\,B')]\!]$ by induction hypothesis, so $\big(\alpha, [\![\vec{x}.\,u : (\exists x{:}\sigma.\,B)]\!]\big) = \big(\alpha, [\![\vec{x}.\,u' : (\exists y{:}\sigma.\,B')]\!]\big)$. Hence, it suffices to prove that $\gamma_x \circ \beta_x^{-1} = \gamma_y \circ \beta_y^{-1}$.

From $v[z/x] =_\alpha v'[z/y]$, we have that $v[z/x][r/z] =_\alpha v'[z/y][r/z]$, that is, $v[r/x] =_\alpha v'[r/y]$, for every $r:\sigma \in \mathrm{LTerms}(\Sigma)$. Thus, by induction hypothesis, $[\![\vec{x}'.\, v[r/x]:B[r/x] \to C]\!] = [\![\vec{x}'.\, v'[r/y]:B'[r/y] \to C]\!]$, from which $M(B[r/x]) = M(B'[r/y])$, an instance of which is $MB = M(B'[x/y])$, so $M(\exists x:\sigma.\, A' \times B) = M(\exists y:\sigma.\, A' \times B')$.

Considering the canonical co-cones:

$$\Big(MA' \times M(\exists x:\sigma.\, B), \{1_{MA'} \times [\![w:B[r/x].\, \mathrm{exI}_x(w):(\exists x:\sigma.\, B)]\!]\}_{r:\sigma\in\mathrm{LTerms}(\Sigma)}\Big) \enspace,$$
$$\Big(MA' \times M(\exists y:\sigma.\, B'), \{1_{MA'} \times [\![w:B'[r/y].\, \mathrm{exI}_y(w):(\exists y:\sigma.\, B')]\!]\}_{r:\sigma\in\mathrm{LTerms}(\Sigma)}\Big) \enspace,$$
$$\Big(M(\exists x:\sigma.\, A' \times B), \{[\![w:A' \times B[r/x].\, \mathrm{exI}_x(w):(\exists x:\sigma.\, A' \times B)]\!]\}_{r:\sigma\in\mathrm{LTerms}(\Sigma)}\Big) \enspace,$$
$$\Big(M(\exists y:\sigma.\, A' \times B'), \{[\![w:A' \times B'[r/y].\, \mathrm{exI}_y(w):(\exists y:\sigma.\, A' \times B')]\!]\}_{r:\sigma\in\mathrm{LTerms}(\Sigma)}\Big) \enspace,$$

from what we have shown till now, the first two co-cones are equal, and also the last two are equal. Thus, $\beta_x$ and $\beta_y$ share the same domain and the same co-domain; being co-universal arrows, it follows that $\beta_x = \beta_y$. Hence, it suffices to prove $\gamma_x = \gamma_y$.

But $\gamma_x$ is the unique arrow such that

$$\begin{aligned}
&\gamma_x \circ [\![w:A' \times B[r/x].\, \mathrm{exI}_x(w):(\exists x:\sigma.\, A' \times B)]\!] \\
&= \mathrm{ev} \circ \big([\![\vec{x}'.\, v[r/x]:B[r/x] \to C]\!] \times 1_{MA'}\big) \\
&= \mathrm{ev} \circ \big([\![\vec{x}'.\, v'[r/y]:B[r/y] \to C]\!] \times 1_{MA'}\big) \\
&= \gamma_y \circ [\![w:A' \times B[r/y].\, \mathrm{exI}_y(w):(\exists y:\sigma.\, A' \times B')]\!] \\
&= \gamma_y \circ [\![w:A' \times B[r/x].\, \mathrm{exI}_x(w):(\exists x:\sigma.\, A' \times B)]\!] \enspace.
\end{aligned}$$

Thus, $\gamma_x = \gamma_y$. □

The key to prove soundness is to notice how composing interpretations is the same as performing substitutions.

**Lemma 5.3.** *Let $\vec{x} \equiv x_1:A_1,\dots,x_n:A_n$ and $\vec{y} \equiv y_1:B_1,\dots,y_m:B_m$ be contexts, let $\vec{y}.\, r_i:A_i$ be terms-in-context for each $1 \le i \le n$, and let $\vec{x}.\, t:C$ be a term-in-context. Then,*

$$[\![\vec{y}.\, t[\vec{r}/\vec{x}]:C]\!] = [\![\vec{x}.\, t:C]\!] \circ [\![\vec{y}.\, \vec{r}:\vec{A}]\!] \enspace,$$

*where $t[\vec{r}/\vec{x}] \equiv t[r_1/x_1,\dots,r_n/x_n]$ and $[\![\vec{y}.\, \vec{r}:\vec{A}]\!] = \big([\![\vec{y}.\, r_1:A_1]\!],\dots,[\![\vec{y}.\, r_n:A_n]\!]\big)$.*

*Proof.* By induction on the structure of the $\lambda$-term $t$:

(1) $t:C \equiv x_i:A_i$: $[\![\vec{y}.\, x_i[\vec{r}/\vec{x}]:A_i]\!] = [\![\vec{y}.\, r_i:A_i]\!] = \pi_i \circ [\![\vec{y}.\, \vec{r}:\vec{A}]\!] = [\![\vec{x}.\, x_i:A_i]\!] \circ [\![\vec{y}.\, \vec{r}:\vec{A}]\!]$;

(2) $t:C \equiv (f:D \to C)(s:D)$:
$$\begin{aligned}
&[\![\vec{y}.\, f(s)[\vec{r}/\vec{x}]:C]\!] \\
&= [\![\vec{y}.\, f(s[\vec{r}/\vec{x}]):C]\!] \\
&= M_{\mathrm{Ax}} f \circ [\![\vec{y}.\, s[\vec{r}/\vec{x}]:D]\!]
\end{aligned}$$

and, by induction hypothesis,

$$\begin{aligned}
&= M_{\mathrm{Ax}} f \circ [\![\vec{x}.\, s:D]\!] \circ [\![\vec{y}.\, \vec{r}:\vec{A}]\!] \\
&= [\![\vec{x}.\, f(s):C]\!] \circ [\![\vec{y}.\, \vec{r}:\vec{A}]\!] \enspace;
\end{aligned}$$

(3) $t:C \equiv \langle u,v\rangle:C_1 \times C_2$:
$$\begin{aligned}
&[\![\vec{y}.\, \langle u,v\rangle[\vec{r}/\vec{x}]:C_1 \times C_2]\!] \\
&= [\![\vec{y}.\, \langle u[\vec{r}/\vec{x}], v[\vec{r}/\vec{x}]\rangle:C_1 \times C_2]\!] \\
&= \big([\![\vec{y}.\, u[\vec{r}/\vec{x}]:C_1]\!], [\![\vec{y}.\, v[\vec{r}/\vec{x}]:C_2]\!]\big)
\end{aligned}$$

and, by induction hypothesis,

$$
\begin{aligned}
&= \big( [\![\vec{x}.\, u : C_1]\!] \circ [\![\vec{y}.\, \vec{r} : \vec{A}]\!], [\![\vec{x}.\, v : C_2]\!] \circ [\![\vec{y}.\, \vec{r} : \vec{A}]\!] \big) \\
&= \big( [\![\vec{x}.\, u : C_1]\!], [\![\vec{x}.\, v : C_2]\!] \big) \circ [\![\vec{y}.\, \vec{r} : \vec{A}]\!] \\
&= [\![\vec{x}.\, \langle u, v\rangle : C_1 \times C_2]\!] \circ [\![\vec{y}.\, \vec{r} : \vec{A}]\!] \ ;
\end{aligned}
$$

(4) $t : C \equiv \mathrm{fst}(s : C_1 \times C_2) : C_1$:

$$
\begin{aligned}
&[\![\vec{y}.\, \mathrm{fst}(s)[\vec{r}/\vec{x}] : C_1]\!] \\
&= [\![\vec{y}.\, \mathrm{fst}(s[\vec{r}/\vec{x}]) : C_1]\!] \\
&= \pi_1 \circ [\![\vec{y}.\, s[\vec{r}/\vec{x}] : C_1 \times C_2]\!]
\end{aligned}
$$

and, by induction hypothesis,

$$
\begin{aligned}
&= \pi_1 \circ [\![\vec{x}.\, s : C_1 \times C_2]\!] \circ [\![\vec{y}.\, \vec{r} : \vec{A}]\!] \\
&= [\![\vec{x}.\, \mathrm{fst}(s) : C_1]\!] \circ [\![\vec{y}.\, \vec{r} : \vec{A}]\!] \ ;
\end{aligned}
$$

(5) $t : C \equiv \mathrm{snd}(s : C_1 \times C_2) : C_2$:

$$
\begin{aligned}
&[\![\vec{y}.\, \mathrm{snd}(s)[\vec{r}/\vec{x}] : C_2]\!] \\
&= [\![\vec{y}.\, \mathrm{snd}(s[\vec{r}/\vec{x}]) : C_2]\!] \\
&= \pi_2 \circ [\![\vec{y}.\, s[\vec{r}/\vec{x}] : C_1 \times C_2]\!]
\end{aligned}
$$

and, by induction hypothesis,

$$
\begin{aligned}
&= \pi_2 \circ [\![\vec{x}.\, s : C_1 \times C_2]\!] \circ [\![\vec{y}.\, \vec{r} : \vec{A}]\!] \\
&= [\![\vec{x}.\, \mathrm{snd}(s) : C_2]\!] \circ [\![\vec{y}.\, \vec{r} : \vec{A}]\!] \ ;
\end{aligned}
$$

(6) $t : C \equiv (\lambda z : D_1.\, s) : D_1 \to D_2$:

$$
\begin{aligned}
&[\![\vec{y}.\, (\lambda z : D_1.\, s)[\vec{r}/\vec{x}] : D_1 \to D_2]\!] \\
&= [\![\vec{y}.\, (\lambda w : D_1.\, s[w/z][\vec{r}/\vec{x}]) : D_1 \to D_2]\!] = h
\end{aligned}
$$

where $w : D_1$ is a new variable and $h$ is the exponential transpose of

$$g = [\![\vec{y}, w : D_1.\, s[w/z][\vec{r}/\vec{x}] : D_2]\!]$$

which, by induction hypothesis

$$
\begin{aligned}
&= [\![\vec{x}, w : D_1.\, s[w/z] : D_2]\!] \circ \\
&\quad \circ \big( [\![\vec{y}, w : D_1.\, r_1 : A_1]\!], \dots, [\![\vec{y}, w : D_1.\, r_n : A_n]\!], [\![\vec{y}, w : D_1.\, w : D_1]\!] \big) \\
&= [\![\vec{x}, w : D_1.\, s[w/z] : D_2]\!] \circ \big( [\![\vec{y}.\, \vec{r} : \vec{A}]\!] \times 1_{MD_1} \big)
\end{aligned}
$$

and, by proposition 5.1

$$= [\![\vec{x}, z : D_1.\, s : D_2]\!] \circ \big( [\![\vec{y}.\, \vec{r} : \vec{A}]\!] \times 1_{MD_1} \big) \ .$$

As illustrated in figure 9, $[\![\vec{x}.\, (\lambda z : D_1.\, s) : D_1 \to D_2]\!] = k$ is the exponential transpose of $[\![\vec{x}, z : D_1.\, s : D_2]\!]$. Thus,

$$
\begin{aligned}
g &= [\![\vec{x}, z : D_1.\, s : D_2]\!] \circ \big( [\![\vec{y}.\, \vec{r} : \vec{A}]\!] \times 1_{MD_1} \big) \\
&= \mathrm{ev} \circ (k \times 1_{MD_1}) \circ \big( [\![\vec{y}.\, \vec{r} : \vec{A}]\!] \times 1_{MD_1} \big) \\
&= \mathrm{ev} \circ \big( (k \circ [\![\vec{y}.\, \vec{r} : \vec{A}]\!]) \times 1_{MD_1} \big) \ .
\end{aligned}
$$

So, by exponentiation, $h = k \circ [\![\vec{y}.\, \vec{r} : \vec{A}]\!]$;

(7) $t : C \equiv (u : C_1 \to C_2) \cdot (v : C_1) : C_2$:

$$
\begin{aligned}
&[\![\vec{y}.\, (u \cdot v)[\vec{r}/\vec{x}] : C_2]\!] \\
&= [\![\vec{y}.\, u[\vec{r}/\vec{x}] \cdot v[\vec{r}/\vec{x}] : C_2]\!] \\
&= \mathrm{ev} \circ \big( [\![\vec{y}.\, u[\vec{r}/\vec{x}] : C_1 \to C_2]\!], [\![\vec{y}.\, v[\vec{r}/\vec{x}] : C_1]\!] \big)
\end{aligned}
$$

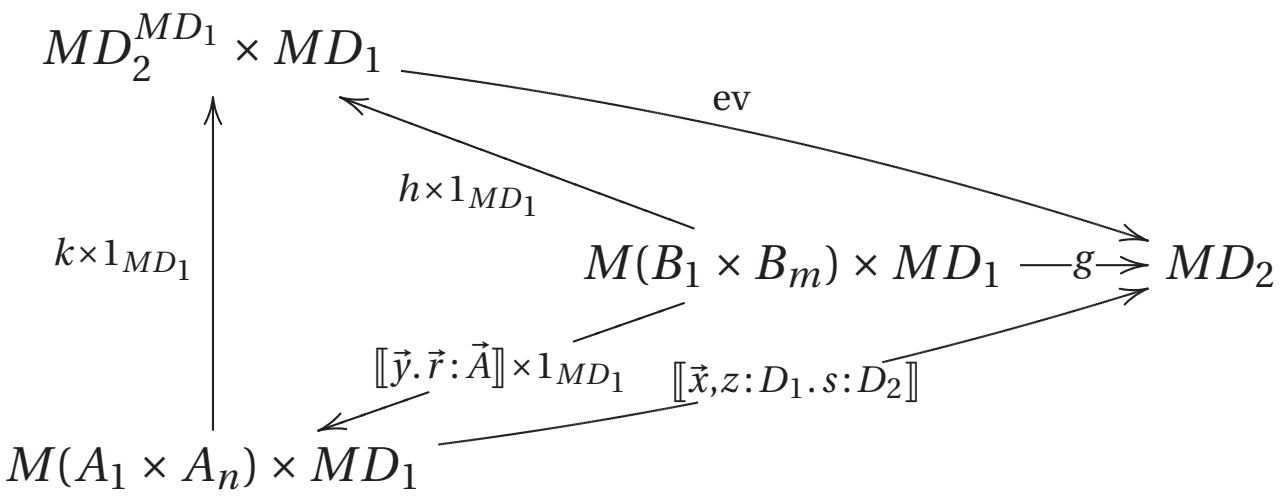

FIGURE 9. The case of $(\lambda z\!:\!D_1.\,s)$ in lemma 5.3.

and, by induction hypothesis,

$$
\begin{aligned}
&= \mathrm{ev} \circ \left([\![\vec{x}.\,u\!:\!C_1 \to C_2]\!] \circ [\![\vec{y}.\,\vec{r}:\vec{A}]\!], [\![\vec{x}.\,v\!:\!C_1]\!] \circ [\![\vec{y}.\,\vec{r}:\vec{A}]\!]\right) \\
&= \mathrm{ev} \circ \left([\![\vec{x}.\,u\!:\!C_1 \to C_2]\!], [\![\vec{x}.\,v\!:\!C_1]\!]\right) \circ [\![\vec{y}.\,\vec{r}:\vec{A}]\!] \\
&= [\![\vec{y}.\,u \cdot v\!:\!C_2]\!] \circ [\![\vec{y}.\,\vec{r}:\vec{A}]\!] \quad ;
\end{aligned}
$$

(8) $t\!:\!C \equiv \mathrm{inl}_{D_2}(s\!:\!D_1)\!:\!D_1 + D_2$:

$$
\begin{aligned}
&[\![\vec{y}.\,\mathrm{inl}_{D_2}(s)[\vec{r}/\vec{x}]\!:\!D_1 + D_2]\!] \\
=&[\![\vec{y}.\,\mathrm{inl}_{D_2}(s[\vec{r}/\vec{x}])\!:\!D_1 + D_2]\!] \\
=&\,\iota_1 \circ [\![\vec{y}.\,s[\vec{r}/\vec{x}]\!:\!D_1]\!]
\end{aligned}
$$

and, by induction hypothesis,

$$
\begin{aligned}
&= \iota_1 \circ [\![\vec{x}.\,s\!:\!D_1]\!] \circ [\![\vec{y}.\,\vec{r}:\vec{A}]\!] \\
&= [\![\vec{x}.\,\mathrm{inl}_{D_2}(s)\!:\!D_1 + D_2]\!] \circ [\![\vec{y}.\,\vec{r}:\vec{A}]\!] \quad ;
\end{aligned}
$$

(9) $t\!:\!C \equiv \mathrm{inr}_{D_1}(s\!:\!D_2)\!:\!D_1 + D_2$:
As in the preceding case,

$$
\begin{aligned}
&[\![\vec{y}.\,\mathrm{inr}_{D_1}(s)[\vec{r}/\vec{x}]\!:\!D_1 + D_2]\!] \\
=&[\![\vec{y}.\,\mathrm{inr}_{D_1}(s[\vec{r}/\vec{x}])\!:\!D_1 + D_2]\!] \\
=&\,\iota_2 \circ [\![\vec{y}.\,s[\vec{r}/\vec{x}]\!:\!D_2]\!]
\end{aligned}
$$

and, by induction hypothesis,

$$
\begin{aligned}
&= \iota_2 \circ [\![\vec{x}.\,s\!:\!D_2]\!] \circ [\![\vec{y}.\,\vec{r}:\vec{A}]\!] \\
&= [\![\vec{x}.\,\mathrm{inr}_{D_1}(s)\!:\!D_1 + D_2]\!] \circ [\![\vec{y}.\,\vec{r}:\vec{A}]\!] \quad ;
\end{aligned}
$$

(10) $t\!:\!C \equiv \mathrm{when}(s\!:\!D_1 + D_2, u\!:\!D_1 \to C, v\!:\!D_2 \to C)$:

$$
\begin{aligned}
&[\![\vec{y}.\,\mathrm{when}(s,u,v)[\vec{r}/\vec{x}]\!:\!C]\!] = [\![\vec{y}.\,\mathrm{when}(s[\vec{r}/\vec{x}], u[\vec{r}/\vec{y}], v[\vec{r}/\vec{x}])\!:\!C]\!] \\
=&\left[\mathrm{ev} \circ \left([\![\vec{y}.\,u[\vec{r}/\vec{x}]\!:\!D_1 \to C]\!] \times 1_{MD_1}\right), \mathrm{ev} \circ \left([\![\vec{x}.\,v[\vec{r}/\vec{x}]\!:\!D_2 \to C]\!] \times 1_{MD_2}\right)\right] \\
&\circ \Delta^{-1} \circ \left(1_{MB}, [\![\vec{y}.\,s[\vec{r}/\vec{x}]\!:\!D_1 + D_2]\!]\right)
\end{aligned}
$$

$\vec{x}' = \mathrm{FV}(u\colon D)$ and $\alpha'\colon MA \to MA'$ the corresponding projection. Thus, there is an arrow $\delta\colon MB' \to MA'$, and precisely $\delta = \alpha' \circ [\![\vec{y}'.\vec{r}:\vec{A}]\!]$ such that $\alpha' \circ [\![\vec{y}.\vec{r}:\vec{A}]\!] = \delta \circ \alpha$.

Thus,
$$\begin{aligned}&[\![\vec{x}.\,\mathrm{allI}(\lambda z\colon s.\,u)\colon(\forall z\colon s.\,D)]\!] \circ [\![\vec{y}.\vec{r}:\vec{A}]\!]\\ =&(\beta'\colon MA' \to M(\forall z\colon s.\,D)) \circ \alpha' \circ [\![\vec{y}.\vec{r}:\vec{A}]\!]\\ =&\beta' \circ \delta \circ \alpha\end{aligned}$$

where $\beta'$ is universal in $\mathbb{C}_{\forall z\colon s.\,D}$, but $\beta' \circ \delta = \beta$ as $\beta$ is universal in $\mathbb{C}_{\forall z\colon s.\,D}$

$$\begin{aligned}=&\beta \circ \alpha\\ =&[\![\vec{y}.\,\mathrm{allI}(\lambda z\colon s.\,u)[\vec{r}:\vec{A}]\colon(\forall z\colon s.\,D)]\!];\end{aligned}$$

(14) $t\colon C \equiv \mathrm{allE}(u,v)\colon D[v/z]$:
$$\begin{aligned}&[\![\vec{y}.\,\mathrm{allE}(u,v)[\vec{r}/\vec{x}]\colon D[v/z]]\!]\\ =&[\![\vec{y}.\,\mathrm{allE}(u[\vec{r}/\vec{x}],v)\colon D[v/z]]\!]\\ =&p_v \circ [\![\vec{y}.\,u[\vec{r}/\vec{x}]\colon(\forall z\colon s.\,D)]\!]\end{aligned}$$

and, by induction hypothesis,

$$\begin{aligned}=&p_v \circ [\![\vec{x}.\,u\colon(\forall z\colon s.\,D)]\!] \circ [\![\vec{y}.\vec{r}:\vec{A}]\!]\\ =&[\![\vec{x}.\,\mathrm{allE}(u,v)\colon D[v/z]]\!] \circ [\![\vec{y}.\vec{r}:\vec{A}]\!]\ ;\end{aligned}$$

(15) $t\colon C \equiv \mathrm{exI}_z(u)\colon(\exists z\colon s.\,D)$:
$$\begin{aligned}&[\![\vec{y}.\,\mathrm{exI}_z(u)[\vec{r}/\vec{x}]\colon(\exists z\colon s.\,D)]\!]\\ =&[\![\vec{y}.\,\mathrm{exI}_z(u[\vec{r}/\vec{x}])\colon(\exists z\colon s.\,D)]\!]\\ =&j_v \circ [\![\vec{y}.\,u[\vec{r}/\vec{x}]\colon D[v/z]]\!]\end{aligned}$$

and, by induction hypothesis,

$$\begin{aligned}=&j_v \circ [\![\vec{x}.\,u\colon D[v/z]]\!] \circ [\![\vec{y}.\vec{r}:\vec{A}]\!]\\ =&[\![\vec{x}.\,\mathrm{exI}_z(u)\colon(\exists z\colon s.\,D)]\!] \circ [\![\vec{y}.\vec{r}:\vec{A}]\!]\ ;\end{aligned}$$

(16) $t\colon C \equiv \mathrm{exE}(u,(\lambda z\colon s.\,v))\colon C$:
$$\begin{aligned}&[\![\vec{y}.\,\mathrm{exE}(u,(\lambda z\colon s.\,v))[\vec{r}/\vec{x}]\colon C]\!]\\ =&[\![\vec{y}.\,\mathrm{exE}(u[\vec{r}/\vec{x}],(\lambda w\colon s.\,v[w/z][\vec{r}/\vec{x}]))\colon C]\!]\end{aligned}$$

where $w\colon s$ is a new variable,

$$=\gamma \circ \beta^{-1} \circ \left(\alpha, [\![\vec{y}.\,u[\vec{r}/\vec{x}]\colon(\exists z\colon s.\,D)]\!]\right)\ ,$$

where
- (a) $\alpha\colon MB \to MB'$ is the projection corresponding to $MB'$, the product associated to $\vec{y}' = \mathrm{FV}(u[\vec{r}/\vec{x}]) \cup \mathrm{FV}(v[w/z][\vec{r}/\vec{x}])$;
- (b) $\beta\colon M(\exists z\colon s.\,B' \times D) \to MB' \times M(\exists z\colon s.\,D)$ co-universal in $\mathbb{C}_{\exists z\colon s.\,B'\times D}$;
- (c) $\gamma\colon M(\exists z\colon s.\,B' \times D) \to MC$ co-universal in $\mathbb{C}_{\exists z\colon s.\,B'\times D}$;

Let $\alpha'\colon MA \to MA'$ be the projection corresponding to the product $MA'$ associated with $\vec{x}' = \mathrm{FV}(u\colon(\exists z\colon s.\,D)) \cup \mathrm{FV}(v\colon D \to C)$.

But, $[\![\vec{y}.\vec{r}:\vec{A}]\!]\colon MB \to MA$ is a product of arrows and, in particular, by definition of $\alpha$, $[\![\vec{y}'.\vec{r}:\vec{A}]\!] \circ \alpha = [\![\vec{y}.\vec{r}:\vec{A}]\!]$. Thus, there is an arrow $\delta\colon MB' \to MA'$, and precisely $\delta = \alpha' \circ [\![\vec{y}'.\vec{r}:\vec{A}]\!]$ such that $\alpha' \circ [\![\vec{y}.\vec{r}:\vec{A}]\!] = \delta \circ \alpha$.

Now,
$$\begin{aligned}&[\![\vec{x}.\,\mathrm{exE}(u,(\lambda z\colon s.\,v))\colon C]\!] \circ [\![\vec{y}.\vec{r}:\vec{A}]\!]\\ =&\gamma' \circ \beta'^{-1} \circ \left(\alpha', [\![\vec{x}.\,u\colon(\exists z\colon s.\,D)]\!]\right) \circ [\![\vec{y}.\vec{r}:\vec{A}]\!]\ ,\end{aligned}$$

where
- (a) $\beta'\colon M(\exists z\colon s.\,A' \times D) \to MA' \times M(\exists z\colon s.\,D)$ co-universal in $\mathbb{C}_{\exists z\colon s.\,A'\times D}$;
- (b) $\gamma'\colon M(\exists z\colon s.\,A' \times D) \to MC$ co-universal in $\mathbb{C}_{\exists z\colon s.\,A'\times D}$;

Thus, $[\![\vec{x}.\, \mathrm{exE}(u,(\lambda x\colon s.\, v))\colon C]\!] \circ [\![\vec{y}.\vec{r}\colon\vec{A}]\!]$

$$=\gamma'\circ\beta'^{-1}\circ\big(\alpha',[\![\vec{x}.\,u\colon(\exists z\colon s.\,D)]\!]\big)\circ[\![\vec{y}.\vec{r}\colon\vec{A}]\!]$$

$$=\gamma'\circ\beta'^{-1}\circ\big(\alpha'\circ[\![\vec{y}.\vec{r}\colon\vec{A}]\!],[\![\vec{x}.\,u\colon(\exists z\colon s.\,D)]\!]\circ[\![\vec{y}.\vec{r}\colon\vec{A}]\!]\big)$$

$$=\gamma'\circ\beta'^{-1}\circ\big(\delta\circ\alpha,[\![\vec{x}.\,u\colon(\exists z\colon s.\,D)]\!]\circ[\![\vec{y}.\vec{r}\colon\vec{A}]\!]\big)$$

$$=\gamma'\circ\beta'^{-1}\circ\big(\delta\times 1_{M(\exists z\colon s.\,D)}\big)\circ\big(\alpha,[\![\vec{x}.\,u\colon(\exists z\colon s.\,D)]\!]\circ[\![\vec{y}.\vec{r}\colon\vec{A}]\!]\big)\ .$$

Consider the co-cone

$$\Big(M(\exists z\colon s.\,B'\times D),\big\{[\![w\colon B'\times D[h/z].\ \mathrm{exI}_z(w)\colon(\exists z\colon s.\,B'\times D)]\!]\big\}_{h\colon s\in\mathrm{LTerms}(\Sigma)}\Big)\ ,$$

and the co-cone

$$\Big(M(\exists z\colon s.\,A'\times D),\big\{[\![w\colon A'\times D[h/z].\ \mathrm{exI}_z(w)\colon(\exists z\colon s.\,A'\times D)]\!]\big\}_{h\colon s\in\mathrm{LTerms}(\Sigma)}\Big)\ :$$

evidently, $\delta\times 1_{M(D[h/z])}\colon MB'\times M(D[h/z])\to MA'\times M(D[h/z])$ for each $h\colon s\in\mathrm{LTerms}(\Sigma)$. So,

$$\big(M(\exists z\colon s.\,A'\times D),$$
$$\big\{[\![w\colon A'\times D[h/z].\ \mathrm{exI}_z(w)\colon(\exists z\colon s.\,A'\times D)]\!]\circ(\delta\times 1_{M(D[h/z])})\big\}_{h\colon s\in\mathrm{LTerms}(\Sigma)}\Big)$$

is a co-cone on $\Sigma_{B'\times D}(z\colon s)$ and $z\colon s\not\in\mathrm{FV}(\exists z\colon s.\,A'\times D)$. So, being $M(\exists z\colon s.\,B'\times D)$ initial in $\mathbb{C}_{\exists z\colon s.\,B'\times B}$, there is a unique arrow $\delta'\colon M(\exists z\colon s.\,B'\times D)\to M(\exists z\colon s.\,A'\times D)$ in $\mathbb{C}_{\exists z\colon s.\,B'\times B}$. Thus, $\delta'=\beta'^{-1}\circ(\delta\times 1_{M(\exists z\colon s.\,D)})\circ\beta$ as all these arrows are in the category of co-cones on $\Sigma_{B'\times D}(z\colon s)$.

Hence, $\gamma'\circ\beta'^{-1}\circ\big(\delta\times 1_{M(\exists z\colon s.\,D)}\big)=\gamma'\circ\delta'\circ\beta^{-1}$.

But, $\gamma$ is co-universal in $\mathbb{C}_{\exists z\colon s.\,B'\times B}$ and thus $\gamma'\circ\delta'=\gamma$ as all these arrows are in the category of co-cones on $\Sigma_{B'\times D}(z\colon s)$. □

Validity, i.e., truth, in a logically distributive category is naturally defined as one expects.

**Definition 5.1** (Validity)**.** An equality-in-context $\vec{x}.\,s=_A t$ is *valid* in the *$\lambda$-theory $T$*, a set of equalities-in-context, when, in every logically distributive category $\mathbb{C}$, each model $M$ of $T$ is also a model of $\vec{x}.\,s=_A t$.

A $\Sigma$-structure $M$ in $\mathbb{C}$ is a *model of a theory $T$* when it is a model of each $\phi$ in $T$.

Finally, $M$ is a *model of an equality-in-context $\phi$* if $[\![\phi]\!]=1$.

With the preceding propositions, the soundness theorem is, essentially, calculation.

**Theorem 5.4** (Soundness)**.** *If an equation-in-context $\vec{x}.\,s=_A t$ is derivable from a $\lambda$-theory $T$, then $\vec{x}.\,s=_A t$ is valid in each model of $T$ in every logically distributive category.*

*Proof.* By induction on the structure of the derivation of the equality-in-context (for clarity, we keep the notation of definition 3.7):

(eq$_0$) $$[\![\vec{y}\colon\vec{B}.\,s[\vec{r}/\vec{x}]\colon C]\!]$$

by lemma 5.3

$$=[\![\vec{x}\colon\vec{A}.\,s\colon C]\!]\circ[\![\vec{y}\colon\vec{B}.\vec{r}\colon\vec{A}]\!]$$

by induction hypothesis

$$=[\![\vec{x}\colon\vec{A}.\,t\colon C]\!]\circ[\![\vec{y}\colon\vec{B}.\vec{r}\colon\vec{A}]\!]$$

by lemma 5.3

$$=[\![\vec{y}\colon\vec{B}.\,t[\vec{r}/\vec{x}]\colon C]\!]\ ;$$

(eq$_1$) $$[\![\vec{x}\colon\vec{A}.\, r[\vec{s}/\vec{y}]\colon C]\!]$$

by lemma 5.3

$$=[\![\vec{y}\colon\vec{B}.\, r\colon C]\!]\circ[\![\vec{x}\colon\vec{A}.\,\vec{s}\colon\vec{B}]\!]$$

by induction hypothesis applied $m$ times,

$$=[\![\vec{y}\colon\vec{B}.\, r\colon C]\!]\circ[\![\vec{x}\colon\vec{A}.\,\vec{t}\colon\vec{B}]\!]$$

by lemma 5.3

$$=[\![\vec{x}\colon\vec{A}.\, r[\vec{t}/\vec{y}]\colon C]\!]\ ;$$

(eq$_2$) Immediate, as equality of arrows is reflexive;
(eq$_3$) Immediate, as equality of arrows is symmetric;
(eq$_4$) Immediate, as equality of arrows is transitive;
(eq$_5$) $[\![\vec{x}\colon\vec{A}.\,(\lambda y\colon B.\, s)\colon B\to C]\!]=[\![\vec{x}\colon\vec{A}.\,(\lambda y\colon B.\, t)\colon B\to C]\!]$ if and only if $[\![\vec{x}\colon\vec{A}, y\colon B.\, s\colon C]\!]=[\![\vec{x}\colon\vec{A}, y\colon B.\, t\colon C]\!]$. But

$$\begin{aligned}&[\![\vec{x}\colon\vec{A}, y\colon B.\, s\colon C]\!]\\ =&[\![\vec{x}\colon\vec{A}.\, s\colon C]\!]\circ\left(1_{MA_1}\times\cdots\times 1_{MA_n}\colon M\vec{A}\times MB\to M\vec{A}\right)\end{aligned}$$

by induction hypothesis

$$\begin{aligned}=&[\![\vec{x}\colon\vec{A}.\, t\colon C]\!]\circ\left(1_{MA_1}\times\cdots\times 1_{MA_n}\colon M\vec{A}\times MB\to M\vec{A}\right)\\ =&[\![\vec{x}\colon\vec{A}, y\colon B.\, t\colon C]\!]\ ;\end{aligned}$$

(eq$_6$) $[\![\vec{x}\colon\vec{A}.\,\mathrm{allI}(\lambda y\colon s.\, r)\colon(\forall y\colon s.\, C)]\!]=\beta_r\circ\alpha_r$ where $\alpha_r\colon M\vec{A}\to MA^r$ is the projection to the product corresponding to $\vec{x}^r=\mathrm{FV}(r\colon C)$, and $\beta_r\colon MA^r\to M(\forall y\colon s.\, C)$ is universal in $\mathbb{C}_{\forall y\colon s.\, C}$.

Analogously, $[\![\vec{x}\colon\vec{A}.\,\mathrm{allI}(\lambda y\colon s.\, t)\colon(\forall y\colon s.\, C)]\!]=\beta_t\circ\alpha_t$ where $\alpha_t\colon M\vec{A}\to MA^t$ is the projection to the product corresponding to $\vec{x}^t=\mathrm{FV}(t\colon C)$, and $\beta_r\colon MA^r\to M(\forall y\colon s.\, C)$ is universal in $\mathbb{C}_{\forall y\colon s.\, C}$.

Let $\alpha\colon M\vec{A}\to MA'$ be the projection to the product corresponding to $\vec{x}'=\mathrm{FV}(r\colon C)\cup\mathrm{FV}(t\colon C)$. Then, there are $\alpha'_r\colon MA'\to MA^r$ and $\alpha'_t\colon MA'\to MA^t$ projections such that $\alpha_r=\alpha'_r\circ\alpha$ and $\alpha_t=\alpha'_t\circ\alpha$.

Considering the cones

$$\left(MA',\{[\![\vec{x}'\colon A'.\, r[h/y]\colon C[h/y]]\!]\circ\alpha'_r\}_{h\colon s\in\mathrm{LTerms}(\Sigma)}\right)$$

and

$$\left(MA',\{[\![\vec{x}'\colon A'.\, t[h/y]\colon C[h/y]]\!]\circ\alpha'_t\}_{h\colon s\in\mathrm{LTerms}(\Sigma)}\right)$$

it follows that $\alpha'_r$ and $\alpha'_t$ are arrows in the category of cones on $\mathbb{C}$ and, moreover, $MA'$ is an object of $\mathbb{C}_{\forall y\colon s.\, C}$ as $y\colon s\not\in\mathrm{FV}(A')$ by construction.

Thus, there is $\beta\colon MA'\to M(\forall y\colon s.\, C)$ universal in $\mathbb{C}_{\forall y\colon s.\, C}$ and $\beta=\beta_r\circ\alpha'_r=\beta_t\circ\alpha'_t$. So, $\beta_r\circ\alpha_r=\beta_r\circ\alpha'_r\circ\alpha=\beta\circ\alpha=\beta_t\circ\alpha'_t\circ\alpha=\beta_t\circ\alpha_t$;

(eq$_7$) $[\![\vec{x}\colon\vec{A}.\,\mathrm{exE}(t,(\lambda y\colon s.\, u))\colon C]\!]=\gamma_u\circ\beta_u^{-1}\circ\left(\alpha_u,[\![\vec{x}\colon\vec{A}.\, t\colon(\exists y\colon s.\, D)]\!]\right)$ where

- $\alpha_u\colon M\vec{A}\to MA^u$ is the projection to the product corresponding to $\vec{x}_u=\mathrm{FV}(t\colon(\exists y\colon s.\, D))\cup\mathrm{FV}(u\colon D\to C)$;
- $\beta_u\colon M(\exists y\colon s.\, A^u\times D)\to MA^u\times M(\exists y\colon s.\, D)$ is co-universal in $\mathbb{C}_{\exists y\colon s.\, A^u\times D}$;
- $\gamma_u\colon M(\exists y\colon s.\, A^u\times D)\to MC$ is co-universal in $\mathbb{C}_{\exists y\colon s.\, A^u\times D}$.

Analogously, $[\![\vec{x}\colon\vec{A}.\,\mathrm{exE}(t,(\lambda y\colon s.\, v))\colon C]\!]=\gamma_v\circ\beta_v^{-1}\circ\left(\alpha_v,[\![\vec{x}\colon\vec{A}.\, t\colon(\exists y\colon s.\, D)]\!]\right)$ where

- $\alpha_v\colon M\vec{A}\to MA^v$ is the projection to the product corresponding to $\vec{x}_v=\mathrm{FV}(t\colon(\exists y\colon s.\, D))\cup\mathrm{FV}(v\colon D\to C)$;
- $\beta_v\colon M(\exists y\colon s.\, A^v\times D)\to MA^v\times M(\exists y\colon s.\, D)$ is co-universal in $\mathbb{C}_{\exists y\colon s.\, A^v\times D}$;
- $\gamma_v\colon M(\exists y\colon s.\, A^v\times D)\to MC$ is co-universal in $\mathbb{C}_{\exists y\colon s.\, A^v\times D}$.

Let $\alpha\colon M\vec{A} \to MA'$ be the projection to the product corresponding to

$$\vec{x}' = \mathrm{FV}(t\colon(\exists y\colon s.\,D)) \cup \mathrm{FV}(u\colon D \to C) \cup \mathrm{FV}(v\colon D \to C)\ .$$

Then, there are $\alpha'_u\colon MA' \to MA^u$ and $\alpha'_v\colon MA' \to MA^v$ projections such that $\alpha_u = \alpha'_u \circ \alpha$ and $\alpha_v = \alpha'_v \circ \alpha$.

Consider the co-cone

$$\left(MA' \times M(\exists y\colon s.\,D), \left\{1_{MA'} \times [\![w\colon D[h/y].\,\mathrm{exI}_y(w)\colon(\exists y\colon s.\,D)]\!]\right\}_{h\colon s\in \mathrm{LTerms}(\Sigma)}\right)\ ,$$

since $y\colon s \notin \mathrm{FV}(A' \times (\exists y\colon s.\,D))$, there is a co-universal arrows $\beta\colon M(\exists y\colon s.\,A' \times D) \to MA' \times M(\exists y\colon s.\,D)$ and, being $\mathbb{C}$ logically distributive, it admits inverse.

Considering the co-cones

$$\Big(M(\exists y\colon s.\,A^u \times D), \Big\{\ \beta_u^{-1} \circ (\alpha'_u \times 1_{M(\exists y\colon s.\,D)}) \circ \left(1_{MA'} \times [\![w\colon D[h/y].\,\mathrm{exI}_y(w)\colon(\exists y\colon s.\,D)]\!]\right)\Big\}_{h\colon s\in \mathrm{LTerms}(\Sigma)}\Big)$$

and

$$\Big(M(\exists y\colon s.\,A^v \times D), \Big\{\ \beta_v^{-1} \circ (\alpha'_v \times 1_{M(\exists y\colon s.\,D)}) \circ \left(1_{MA'} \times [\![w\colon D[h/y].\,\mathrm{exI}_y(w)\colon(\exists y\colon s.\,D)]\!]\right)\Big\}_{h\colon s\in \mathrm{LTerms}(\Sigma)}\Big)$$

there must be $\beta'_u\colon M(\exists y\colon s.\,A' \times D) \to M(\exists y\colon s.\,A^u \times D)$ and $\beta'_v\colon M(\exists y\colon s.\,A' \times D) \to M(\exists y\colon s.\,A^v \times D)$ both co-universal in $\mathbb{C}_{\exists y\colon s.\,A'\times D}$. Moreover, it holds that $\beta'_u = \beta_u^{-1} \circ (\alpha'_u \times 1_{M(\exists y\colon s.\,D)}) \circ \beta$ and $\beta'_v = \beta_v^{-1} \circ (\alpha'_v \times 1_{M(\exists y\colon s.\,D)}) \circ \beta$.

Considering the co-cone

$$\left(MC, \left\{\mathrm{ev} \circ \left([\![\vec{x}'\colon A'.\,u[h/y]\colon D[h/y] \to C]\!] \circ 1_{M(D[h/y])}\right)\right\}_{h\colon s\in \mathrm{LTerms}(\Sigma)}\right)$$

there must be $\gamma\colon M(\exists y\colon s.\,A' \times D) \to MC$ co-universal in $\mathbb{C}_{\exists y\colon s.\,A'\times D}$ and, moreover, $\gamma = \gamma_u \circ \beta'_u = \gamma_v \circ \beta'_v$.

So, applying all the equalities,

$$\begin{aligned}
&\gamma_u \circ \beta_u^{-1} \circ \left(\alpha_u, [\![\vec{x}\colon \vec{A}.\,t\colon(\exists y\colon s.\,D)]\!]\right)\\
=&\gamma_u \circ \beta_u^{-1} \circ \left(\alpha'_u \circ \alpha, [\![\vec{x}\colon \vec{A}.\,t\colon(\exists y\colon s.\,D)]\!]\right)\\
=&\gamma_u \circ \beta_u^{-1} \circ (\alpha'_u \times 1_{M(\exists y\colon s.\,D)}) \circ \left(\alpha, [\![\vec{x}\colon \vec{A}.\,t\colon(\exists y\colon s.\,D)]\!]\right)\\
=&\gamma_u \circ \beta'_u \circ \beta^{-1} \circ \left(\alpha, [\![\vec{x}\colon \vec{A}.\,t\colon(\exists y\colon s.\,D)]\!]\right)\\
=&\gamma \circ \beta^{-1} \circ \left(\alpha, [\![\vec{x}\colon \vec{A}.\,t\colon(\exists y\colon s.\,D)]\!]\right)\\
=&\gamma_v \circ \beta'_v \circ \beta^{-1} \circ \left(\alpha, [\![\vec{x}\colon \vec{A}.\,t\colon(\exists y\colon s.\,D)]\!]\right)\\
=&\gamma_v \circ \beta_v^{-1} \circ (\alpha'_v \times 1_{M(\exists y\colon s.\,D)}) \circ \left(\alpha, [\![\vec{x}\colon \vec{A}.\,t\colon(\exists y\colon s.\,D)]\!]\right)\\
=&\gamma_v \circ \beta_v^{-1} \circ \left(\alpha'_v \circ \alpha, [\![\vec{x}\colon \vec{A}.\,t\colon(\exists y\colon s.\,D)]\!]\right)\\
=&\gamma_v \circ \beta_v^{-1} \circ \left(\alpha_v, [\![\vec{x}\colon \vec{A}.\,t\colon(\exists y\colon s.\,D)]\!]\right)\ ;
\end{aligned}$$

$(\times_0)$ $[\![x\colon 1.\,x\colon]\!] = 1_1 = !\colon 1 \to 1 = [\![x\colon 1.\,*\colon 1]\!]$;

$(\times_1)$

$$\begin{aligned}
&[\![x\colon A, y\colon B.\,\mathrm{fst}(\langle x, y\rangle)\colon A]\!]\\
=&\pi_1 \circ [\![x\colon A, y\colon B.\,\langle x, y\rangle\colon A \times B]\!]\\
=&\pi_1 \circ \left([\![x\colon A, y\colon B.\,x\colon A]\!], [\![x\colon A, y\colon B.\,y\colon B]\!]\right)\\
=&\pi_1 \circ (\pi_1, \pi_2)\\
=&\pi_1 \circ 1_{M(A\times B)}\\
=&\pi_1\\
=&[\![x\colon A, y\colon B.\,x\colon A]\!]\ ;
\end{aligned}$$

$(\times_2)$
$$\begin{aligned}
&[\![x\colon A, y\colon B.\, \mathrm{snd}(\langle x,y\rangle) : B]\!]\\
&=\pi_2\circ[\![x\colon A, y\colon B.\, \langle x,y\rangle : A\times B]\!]\\
&=\pi_2\circ\big([\![x\colon A, y\colon B.\, x : A]\!],[\![x\colon A, y\colon B.\, y : B]\!]\big)\\
&=\pi_2\circ(\pi_1,\pi_2)\\
&=\pi_2\circ 1_{M(A\times B)}\\
&=\pi_2\\
&=[\![x\colon A, y\colon B.\, y : B]\!]\ ;
\end{aligned}$$

$(\times_3)$
$$\begin{aligned}
&[\![z\colon A\times B.\, \langle \mathrm{fst}(z),\mathrm{snd}(z)\rangle : A\times B]\!]\\
&=\big([\![z\colon A\times B.\, \mathrm{fst}(z) : A]\!],[\![z\colon A\times B.\, \mathrm{snd}(z) : B]\!]\big)\\
&=\big(\pi_1\circ[\![z\colon A\times B.\, z : A\times B]\!],\pi_2\circ[\![z\colon A\times B.\, z : A\times B]\!]\big)\\
&=\big(\pi_1\circ 1_{M(A\times B)},\pi_2\circ 1_{M(A\times B)}\big)\\
&=(\pi_1,\pi_2)\\
&=1_{M(A\times B)}\\
&=[\![z\colon A\times B.\, z : A\times B]\!]\ ;
\end{aligned}$$

$(+_0)$
$$\begin{aligned}
&[\![\vec{x}\colon\vec{A}.\, \mathrm{when}(\mathrm{inl}_B(a),t,s) : C]\!]\\
&=\big[\mathrm{ev}\circ\big([\![\vec{x}\colon\vec{A}.\, t : D\to C]\!],1_{MD}\big),\mathrm{ev}\circ\big([\![\vec{x}\colon\vec{A}.\, s : B\to C]\!],1_{MB}\big)\big]\\
&\quad\circ\Delta^{-1}\circ\big(1_{M\vec{A}},[\![\vec{x}\colon\vec{A}.\, \mathrm{inl}_B(a) : D+B]\!]\big)\\
&=\big[\mathrm{ev}\circ\big([\![\vec{x}\colon\vec{A}.\, t : D\to C]\!],1_{MD}\big),\mathrm{ev}\circ\big([\![\vec{x}\colon\vec{A}.\, s : B\to C]\!],1_{MB}\big)\big]\\
&\quad\circ\Delta^{-1}\circ\big(1_{M\vec{A}},\iota_1\circ[\![\vec{x}\colon\vec{A}.\, a : D]\!]\big)\\
&=\big[\mathrm{ev}\circ\big([\![\vec{x}\colon\vec{A}.\, t : D\to C]\!],1_{MD}\big),\mathrm{ev}\circ\big([\![\vec{x}\colon\vec{A}.\, s : B\to C]\!],1_{MB}\big)\big]\\
&\quad\circ\big(\iota_1\circ\big(1_{M\vec{A}},[\![\vec{x}\colon\vec{A}.\, a : D]\!]\big)\big)\\
&=\mathrm{ev}\circ\big([\![\vec{x}\colon\vec{A}.\, t : D\to C]\!],1_{MD}\big)\circ\big(1_{M\vec{A}},[\![\vec{x}\colon\vec{A}.\, a : D]\!]\big)\\
&=\mathrm{ev}\circ\big([\![\vec{x}\colon\vec{A}.\, t : D\to C]\!],[\![\vec{x}\colon\vec{A}.\, a : D]\!]\big)\\
&=[\![\vec{x}\colon\vec{A}.\, t\cdot a : C]\!]\ ;
\end{aligned}$$

$(+_1)$
$$\begin{aligned}
&[\![\vec{x}\colon\vec{A}.\, \mathrm{when}(\mathrm{inr}_D(b),t,s) : C]\!]\\
&=\big[\mathrm{ev}\circ\big([\![\vec{x}\colon\vec{A}.\, t : D\to C]\!],1_{MD}\big),\mathrm{ev}\circ\big([\![\vec{x}\colon\vec{A}.\, s : B\to C]\!],1_{MB}\big)\big]\\
&\quad\circ\Delta^{-1}\circ\big(1_{M\vec{A}},[\![\vec{x}\colon\vec{A}.\, \mathrm{inr}_D(b) : D+B]\!]\big)\\
&=\big[\mathrm{ev}\circ\big([\![\vec{x}\colon\vec{A}.\, t : D\to C]\!],1_{MD}\big),\mathrm{ev}\circ\big([\![\vec{x}\colon\vec{A}.\, s : B\to C]\!],1_{MB}\big)\big]\\
&\quad\circ\Delta^{-1}\circ\big(1_{M\vec{A}},\iota_2\circ[\![\vec{x}\colon\vec{A}.\, b : B]\!]\big)\\
&=\big[\mathrm{ev}\circ\big([\![\vec{x}\colon\vec{A}.\, t : D\to C]\!],1_{MD}\big),\mathrm{ev}\circ\big([\![\vec{x}\colon\vec{A}.\, s : B\to C]\!],1_{MB}\big)\big]\\
&\quad\circ\big(\iota_2\circ\big(1_{M\vec{A}},[\![\vec{x}\colon\vec{A}.\, b : B]\!]\big)\big)\\
&=\mathrm{ev}\circ\big([\![\vec{x}\colon\vec{A}.\, s : B\to C]\!],1_{MB}\big)\circ\big(1_{M\vec{A}},[\![\vec{x}\colon\vec{A}.\, b : B]\!]\big)\\
&=\mathrm{ev}\circ\big([\![\vec{x}\colon\vec{A}.\, s : B\to C]\!],[\![\vec{x}\colon\vec{A}.\, b : B]\!]\big)\\
&=[\![\vec{x}\colon\vec{A}.\, s\cdot b : C]\!]\ ;
\end{aligned}$$

$(+_2)$ Let

$$\vec{x}\equiv x_0 : A_1+A_2, x_1 : A_1\to B_1+B_2, x_2 : A_2\to B_1+B_2, x_3 : B_1\to C, x_4 : B_2\to C$$

and

$$E\equiv (A_1+A_2)\times(A_1\to B_1+B_2)\times(A_2\to B_1+B_2)\times(B_1\to C)\times(B_2\to C)\ .$$

Then,

$$
\begin{aligned}
&[\![\vec{x}.\ \mathrm{when}(x_0, (\lambda y\colon A_1.\ \mathrm{when}(x_1\cdot y, x_3, x_4)),\\
&\qquad (\lambda y\colon A_2.\ \mathrm{when}(x_2\cdot y, x_3, x_4)))\colon C]\!]\\
&= \big[\mathrm{ev}\circ\big([\![\vec{x}.(\lambda y\colon A_1.\ \mathrm{when}(x_1\cdot y, x_3, x_4))\colon C]\!]\times 1_{MA_1}\big),\\
&\quad \mathrm{ev}\circ\big([\![\vec{x}.(\lambda y\colon A_2.\ \mathrm{when}(x_2\cdot y, x_3, x_4))\colon C]\!]\times 1_{MA_2}\big)\big]\\
&\circ\Delta^{-1}\circ\big(1_{ME}, [\![\vec{x}.x_0\colon A_1+A_2]\!]\big)
\end{aligned}
$$

by definition of exponential transpose

$$
\begin{aligned}
&= \big[[\![\vec{x}, y\colon A_1.\ \mathrm{when}(x_1\cdot y, x_3, x_4)\colon C]\!],\\
&\quad [\![\vec{x}, y\colon A_2.\ \mathrm{when}(x_2\cdot y, x_3, x_4)\colon C]\!]\big]\\
&\circ\Delta^{-1}\circ(1_{ME}, \pi_1\colon M(E\to A_1+A_2))\\
&= \big[\big[\mathrm{ev}\circ\big([\![\vec{x}, y\colon A_1.x_3\colon B_1\to C]\!]\times 1_{MB_1}\big), \mathrm{ev}\circ\big([\![\vec{x}, y\colon A_1.x_4\colon B_2\to C]\!]\times 1_{MB_2}\big)\big]\\
&\qquad\circ\Delta^{-1}\circ\big(1_{ME\times MA_1}, [\![\vec{x}, y\colon A_1.x_1\cdot y\colon B_1+B_2]\!]\big),\\
&\quad \big[\mathrm{ev}\circ\big([\![\vec{x}, y\colon A_2.x_3\colon B_1\to C]\!]\times 1_{MB_1}\big), \mathrm{ev}\circ\big([\![\vec{x}, y\colon A_2.x_4\colon B_2\to C]\!]\times 1_{MB_2}\big)\big]\\
&\qquad\circ\Delta^{-1}\circ\big(1_{ME\times MA_2}, [\![\vec{x}, y\colon A_2.x_2\cdot y\colon B_1+B_2]\!]\big)\big]\\
&\circ\Delta^{-1}\circ(1_{ME}, \pi_1)\\
&= \big[\big[\mathrm{ev}\circ\big(\pi_4\colon M(E\times A_1\to(B_1\to C))\times 1_{MB_1}\big),\\
&\qquad\quad \mathrm{ev}\circ\big(\pi_5\colon M(E\times A_1\to(B_2\to C))\times 1_{MB_2}\big)\big]\\
&\qquad\circ\Delta^{-1}\circ\big(1_{ME\times MA_1}, \mathrm{ev}\circ\big([\![\vec{x}, y\colon A_1.\ x_1\colon A_1\to B_1+B_2]\!], [\![\vec{x}, y\colon A_1.\ y\colon A_1]\!]\big)\big),\\
&\quad \big[\mathrm{ev}\circ\big(\pi_4\colon M(E\times A_2\to(B_1\to C))\times 1_{MB_1}\big),\\
&\qquad\quad \mathrm{ev}\circ\big(\pi_5\colon M(E\times A_2\to(B_2\to C))\times 1_{MB_2}\big)\big]\\
&\qquad\circ\Delta^{-1}\circ\big(1_{ME\times MA_2}, \mathrm{ev}\circ\big([\![\vec{x}, y\colon A_2.\ x_2\colon A_2\to B_1+B_2]\!], [\![\vec{x}, y\colon A_2.\ y\colon A_2]\!]\big)\big)\big]\\
&\circ\Delta^{-1}\circ(1_{ME}, \pi_1)\\
&= \big[\big[\mathrm{ev}\circ\big(\pi_4\times 1_{MB_1}\big), \mathrm{ev}\circ\big(\pi_5\times 1_{MB_2}\big)\big]\circ\Delta^{-1}\circ\big(1_{ME\times MA_1}, \mathrm{ev}\circ(\pi_2,\pi_6)\big),\\
&\quad \big[\mathrm{ev}\circ\big(\pi_4\times 1_{MB_1}\big), \mathrm{ev}\circ\big(\pi_5\times 1_{MB_2}\big)\big]\circ\Delta^{-1}\circ\big(1_{ME\times MA_2}, \mathrm{ev}\circ(\pi_3,\pi_6)\big)\big]\\
&\circ\Delta^{-1}\circ(1_{ME}, \pi_1)
\end{aligned}
$$

by composition of co-universal arrows,

$$
\begin{aligned}
&= \big[\mathrm{ev}\circ\big(\pi_4\times 1_{MB_1}\big), \mathrm{ev}\circ\big(\pi_5\times 1_{MB_2}\big)\big]\\
&\circ(\Delta^{-1}\colon M(E\times(B_1+B_2))\to M((E\times B_1)+(E\times B_2)))\\
&\circ[(\pi\colon M(E\times A_1)\to ME, \mathrm{ev}\circ(\pi_2,\pi_6)), (\pi\colon M(E\times MA_2)\to E, \mathrm{ev}\circ(\pi_3,\pi_6))]\\
&\circ\Delta^{-1}\circ(1_{ME}, \pi_1)
\end{aligned}
$$

by composition of universal arrows,

$$
\begin{aligned}
&= \big[\mathrm{ev}\circ\big(\pi_4\times 1_{MB_1}\big), \mathrm{ev}\circ\big(\pi_5\times 1_{MB_2}\big)\big]\circ\Delta^{-1}\\
&\circ\big(1_{ME}, [\mathrm{ev}\circ(\pi_2,\pi_6), \mathrm{ev}\circ(\pi_3,\pi_6)]\circ\Delta^{-1}\circ(1_{ME},\pi_1)\big)
\end{aligned}
$$

and, expliciting,

$$
\begin{aligned}
&= \big[\mathrm{ev}\circ\big(\pi_4\times 1_{MB_1}\big), \mathrm{ev}\circ\big(\pi_5\times 1_{MB_2}\big)\big]\circ\Delta^{-1}\\
&\circ\big(1_{ME}, \big[\mathrm{ev}\circ\big([\![\vec{x}.x_1\colon A_1\to B_1+B_2]\!]\times 1_{MA_1}\big),\\
&\qquad \mathrm{ev}\circ\big([\![\vec{x}.x_2\colon A_2\to B_1+B_2]\!]\times 1_{MA_2}\big)\big]\circ\Delta^{-1}\circ\big(1_{ME}, [\![\vec{x}.x_0\colon A_1+A_2]\!]\big)\big)
\end{aligned}
$$

expliciting again, and recognising the when,

$$
\begin{aligned}
&= \left[\mathrm{ev}\circ\left(\llbracket\vec{x}.\,x_3 : B_1 \to C\rrbracket \times 1_{MB_1}\right), \mathrm{ev}\circ\left(\llbracket\vec{x}.\,x_4 : B_2 \to C\rrbracket \times 1_{MB_2}\right)\right] \\
&\quad \circ \Delta^{-1} \circ \left(1_{ME}, \llbracket\vec{x}.\,\mathrm{when}(x_0, x_1, x_2) : B_1 + B_2\rrbracket\right) \\
&= \llbracket\vec{x}.\,\mathrm{when}(\mathrm{when}(x_0, x_1, x_2), x_3, x_4) : C\rrbracket \ ;
\end{aligned}
$$

$(+_3)$

$$
\begin{aligned}
&\llbracket x : A, y : 0.\, \mathrm{F}_A \cdot y : A\rrbracket \\
&= \mathrm{ev}\circ\left(\llbracket x : A, y : 0.\, \mathrm{F}_A : 0 \to A\rrbracket, \llbracket x : A, y : 0.\, y : 0\rrbracket\right) \\
&= \mathrm{ev}\circ(h, \pi_2 : MA \times 0 \to 0)
\end{aligned}
$$

where $h$ is the exponential transpose of $(! : 0 \to MA) \circ (\pi_2 : MA \times 0 \to 0)$,

$$
\begin{aligned}
&= \mathrm{ev}\circ(h \times 1_0) \circ (\pi_1 : MA \times 0 \to MA, \pi_2 : MA \times 0 \to 0) \\
&= \mathrm{ev}\circ(h \times 1_0) \circ 1_{MA\times 0} \\
&= \mathrm{ev}\circ(h \times 1_0) \\
&= (! : 0 \to MA) \circ (\pi_2 : MA \times 0 \to 0) \\
&= (\pi_1 : MA \times 0 \to 0) \circ (1_0, ! : 0 \to MA) \circ (\pi_2 : MA \times 0 \to 0) \\
&= (\pi_1 : MA \times 0 \to 0) \circ 1_{MA\times 0} \\
&= \pi_1 : MA \times 0 \to 0 \\
&= \llbracket x : A, y : 0.\, x : A\rrbracket \ ;
\end{aligned}
$$

$(\to_0)$

$$
\begin{aligned}
&\llbracket\vec{x} : \vec{A}.\,(\lambda y : C.\, s) \cdot t : B\rrbracket \\
&= \mathrm{ev}\circ\left(\llbracket\vec{x} : \vec{A}.\,(\lambda y : C.\, s) : C \to B\rrbracket, \llbracket\vec{x} : \vec{A}.\, t : C\rrbracket\right) \\
&= \mathrm{ev}\circ\left(h, \llbracket\vec{x} : \vec{A}.\, t : C\rrbracket\right)
\end{aligned}
$$

where $h$ is the exponential transpose of $\llbracket\vec{x} : \vec{A}, y : C.\, s : B\rrbracket$

$$
\begin{aligned}
&= \mathrm{ev}\circ(h \times 1_{MC}) \circ \left(1_{M\vec{A}}, \llbracket\vec{x} : \vec{A}.\, t : C\rrbracket\right) \\
&= \llbracket\vec{x} : \vec{A}, y : C.\, s : B\rrbracket \circ \left(1_{M\vec{A}}, \llbracket\vec{x} : \vec{A}.\, t : C\rrbracket\right)
\end{aligned}
$$

and, by lemma 5.3,

$$
= \llbracket\vec{x} : \vec{A}.\, s[t/y] : B\rrbracket \ ;
$$

$(\to_1)$ $\llbracket\vec{x} : \vec{A}.\,(\lambda y : C.\, t \cdot y) : C \to B\rrbracket = h$ where $h$ is the exponential transpose of

$$
\begin{aligned}
&\llbracket\vec{x} : \vec{A}, y : C.\, t \cdot y : B\rrbracket \\
&= \mathrm{ev}\circ\left(\llbracket\vec{x} : \vec{A}, y : C.\, t : C \to B\rrbracket, \llbracket\vec{x} : \vec{A}, y : C.\, y : C\rrbracket\right) \\
&= \mathrm{ev}\circ\left(\llbracket\vec{x} : \vec{A}, y : C.\, t : C \to B\rrbracket, \pi_{n+1}\right) \\
&= \mathrm{ev}\circ\left(\llbracket\vec{x} : \vec{A}.\, t : C \to B\rrbracket \circ \pi_{\vec{A}}, \pi_{n+1}\right)
\end{aligned}
$$

where $\pi_{\vec{A}} : \vec{A} \times C \to \vec{A}$ is the obvious projector—notice how the term-in-context is well-formed thanks to the proviso in the $(\to_1)$ rule

$$
\begin{aligned}
&= \mathrm{ev}\circ\left(\llbracket\vec{x} : \vec{A}.\, t : C \to B\rrbracket \times 1_{MC}\right) \circ \left(\pi_{\vec{A}}, \pi_{n+1}\right) \\
&= \mathrm{ev}\circ\left(\llbracket\vec{x} : \vec{A}.\, t : C \to B\rrbracket \times 1_{MC}\right) \circ 1_{M\vec{A}\times MC} \\
&= \mathrm{ev}\circ\left(\llbracket\vec{x} : \vec{A}.\, t : C \to B\rrbracket \times 1_{MC}\right) \ .
\end{aligned}
$$

Thus, by definition of exponential transpose, $h = \llbracket\vec{x} : \vec{A}.\, t : C \to B\rrbracket$;

$(\forall_0)$
$$\begin{aligned}&[\![\vec{x}:\vec{A}.\ \mathrm{allE}(\mathrm{allI}(\lambda z\colon s.\ t),r):B[r/z]]\!]\\ =&p_r\circ[\![\vec{x}:\vec{A}.\ \mathrm{allI}(\lambda z\colon s.\ t):(\forall z\colon s.B)]\!]\\ =&p_r\circ\beta\circ\alpha\end{aligned}$$

where $p_r$ is the projection of the (unique) cone with vertex $M(\forall z\colon s.B)$, the arrow $\alpha\colon M\vec{A}\to MA'$ is the projection to the product corresponding to $\vec{x}'=\mathrm{FV}(t\colon B)$ and $\beta\colon MA'\to M(\forall z\colon s.B)$ is universal in $\mathbb{C}_{\forall z\colon s.C}$

$$=[\![\vec{x}'.\ t[r/z]:B[r/z]]\!]\circ\alpha$$

since $p_r\circ\beta$ equals the $r$-th projection in the cone whose vertex is $MA'$, as $\beta$ is an arrow between cones,

$$=[\![\vec{x}:\vec{A}.\ t[r/z]:B[r/z]]\!]\ ;$$

$(\forall_1)$ By induction hypothesis,

$$\begin{aligned}&[\![\vec{x}:\vec{A}.\ \mathrm{allE}(u,r):B]\!]\\ =&p_r\circ[\![\vec{x}:\vec{A}.\ u:(\forall z\colon s.B)]\!]\\ =&p_r\circ[\![\vec{x}:\vec{A}.\ v:(\forall z\colon s.B)]\!]\\ =&[\![\vec{x}:\vec{A}.\ \mathrm{allE}(v,r):B]\!]\end{aligned}$$

where $p_r$ is the $r$-th projection of the unique cone with vertex $M(\forall z\colon s.B)$. So,

$$\begin{aligned}&p_r\circ[\![\vec{x}'.\ u:(\forall z\colon s.B)]\!]\circ\alpha\\ =&p_r\circ[\![\vec{x}'.\ v:(\forall z\colon s.B)]\!]\circ\alpha\ ,\end{aligned}$$

where $\alpha\colon M\vec{A}\to MA'$ is the projection to the product corresponding to $\vec{x}'=\mathrm{FV}(u:(\forall z\colon s.B))\cup\mathrm{FV}(v:(\forall z\colon s.B))$.

Evidently,

$$\left(MA',\{p_r\circ[\![\vec{x}'.\ u:(\forall z\colon s.B)]\!]\}_{r\colon s\in\mathrm{LTerms}(\Sigma)}\right)$$

and

$$\left(MA',\{p_r\circ[\![\vec{x}'.\ v:(\forall z\colon s.B)]\!]\}_{r\colon s\in\mathrm{LTerms}(\Sigma)}\right)$$

are both cones on $\Sigma_B(z\colon s)$ and $z\colon s\not\in\mathrm{FV}^*(A')$.

Thus, there is a unique arrow $MA'\to M(\forall z\colon s.B)$ in $\mathbb{C}_{\forall z\colon s.B}$ and it is

$$[\![\vec{x}'.\ u:(\forall z\colon s.B)]\!]=[\![\vec{x}'.\ v:(\forall z\colon s.B)]\!]\ ;$$

$(\exists_0)$
$$\begin{aligned}&[\![\vec{x}:\vec{A}.\ \mathrm{exE}(\mathrm{exI}_z(t),(\lambda z\colon s.\ v)):B]\!]\\ =&\gamma\circ\beta^{-1}\circ\left(\alpha,[\![\vec{x}:\vec{A}.\ \mathrm{exI}_z(t):(\exists z\colon s.\ C)]\!]\right)\\ =&\gamma\circ\beta^{-1}\circ\left(\alpha,j_r\circ[\![\vec{x}:\vec{A}.\ t:C[r/z]]\!]\right)\\ =&\gamma\circ\beta^{-1}\circ\left(1_{MA'},j_r\circ[\![\vec{x}'.\ t:C[r/z]]\!]\right)\circ\alpha\end{aligned}$$

where $j_r\colon M(C[r/z])\to M(\exists z\colon s.C)$ is the $r$-th injection in the unique co-cone with vertex $M(\exists z\colon s.C)$, $\alpha\colon M\vec{A}\to MA'$ is the projection to the product corresponding to $\vec{x}'=\mathrm{FV}(t\colon C[r/z])\cup\mathrm{FV}(v\colon C\to B)$, $\beta\colon M(\exists z\colon s.A'\times C)\to MA'\times M(\exists z\colon s.C)$ is co-universal in $\mathbb{C}_{\exists z\colon s.\ A'\times C}$, and $\gamma\colon M(\exists z\colon s.A'\times C)\to MB$ is co-universal in $\mathbb{C}_{\exists z\colon s.\ A'\times C}$.

Let $\phi\colon M(\exists z\colon s.A'\times C)\to MB^{M(C[r/z])}\times M(C[r/z])$ be the co-universal arrow in $\mathbb{C}_{\exists z\colon s.\ A'\times C}$ to the co-cone

$$\Big(MB^{M(C[r/z])}\times M(C[r/z]),\Big\{\ \left([\![\vec{x}'.\ v[r/z]:C[r/z]\to B]\!],[\![\vec{x}'.\ t:C[r/z]]\!]\right)\circ\pi_1\circ\beta\circ\delta_u\Big\}_{u\colon s\in\mathrm{LTerms}(\Sigma)}\Big)$$

where $\delta_u\colon MA' \times M(C[u/z]) \to M(\exists z\colon s.\, A' \times C)$ is the $u$-th injection of the unique co-cone of vertex $M(\exists z\colon s.\, A' \times C)$.

The arrow $\phi$ exists because $\mathbb{C}$ is logically distributive, and, being co-universal,

$$\phi = \left(\llbracket \vec{x}'.\, v[r/z] : C[r/z] \to B \rrbracket, \llbracket \vec{x}'.\, t : C[r/z] \rrbracket\right) \circ \pi_1 \circ \beta \ .$$

Moreover, by exponentiation, $\gamma \circ \delta_r = \mathrm{ev} \circ \phi \circ \delta_r$, thus, being $\gamma$ and $\phi$ both co-universal, $\gamma = \mathrm{ev} \circ \phi$. So,

$$\begin{aligned}
&\llbracket \vec{x} : \vec{A}.\, \mathrm{exE}(\mathrm{exI}_z(t), (\lambda z\colon s.\, v)) : B \rrbracket \\
=&\gamma \circ \beta^{-1} \circ \left(1_{MA'}, j_r \circ \llbracket \vec{x}'.\, t : C[r/z] \rrbracket\right) \circ \alpha \\
=&\mathrm{ev} \circ \phi \circ \beta^{-1} \circ \left(1_{MA'}, j_r \circ \llbracket \vec{x}'.\, t : C[r/z] \rrbracket\right) \circ \alpha \\
=&\mathrm{ev} \circ \left(\llbracket \vec{x}'.\, v[r/z] : C[r/z] \to B \rrbracket, \llbracket \vec{x}'.\, t : C[r/z] \rrbracket\right) \\
&\circ \pi_1 \circ \left(1_{MA'}, j_r \circ \llbracket \vec{x}'.\, t : C[r/z] \rrbracket\right) \circ \alpha \\
=&\mathrm{ev} \circ \left(\llbracket \vec{x}'.\, v[r/z] : C[r/z] \to B \rrbracket, \llbracket \vec{x}'.\, t : C[r/z] \rrbracket\right) \circ 1_{MA'} \circ \alpha \\
=&\mathrm{ev} \circ \left(\llbracket \vec{x}'.\, v[r/z] : C[r/z] \to B \rrbracket, \llbracket \vec{x}'.\, t : C[r/z] \rrbracket\right) \circ \alpha \\
=&\llbracket \vec{x}'.\, v[r/z] \cdot t : B \rrbracket \circ \alpha \\
=&\llbracket \vec{x} : \vec{A}.\, v[r/z] \cdot t : B \rrbracket \ ;
\end{aligned}$$

($\exists_1$) $$\llbracket \vec{x} : \vec{A}.\, \mathrm{exE}(u, (\lambda z\colon s.\, r)) : B \rrbracket = \gamma_r \circ \beta_r^{-1} \circ \left(\alpha_r, \llbracket \vec{x} : \vec{A}.\, u : (\exists z\colon s.\, C) \rrbracket\right)$$

and

$$\llbracket \vec{x} : \vec{A}.\, \mathrm{exE}(u, (\lambda z\colon s.\, t)) : B \rrbracket = \gamma_t \circ \beta_t^{-1} \circ \left(\alpha_t, \llbracket \vec{x} : \vec{A}.\, u : (\exists z\colon s.\, C) \rrbracket\right) \ ,$$

where $\alpha_r\colon M\vec{A} \to MA^r$, $\alpha_t\colon M\vec{A} \to MA^t$ are the projections to the products corresponding, respectively, to $\vec{x}^r = \mathrm{FV}(u : (\exists z\colon s.\, C)) \cup \mathrm{FV}(r : C \to B)$ and $\vec{x}^t = \mathrm{FV}(u : (\exists z\colon s.\, C)) \cup \mathrm{FV}(t : C \to B)$; moreover, the arrows

$$\begin{aligned}
&\beta_r\colon M(\exists z\colon s.\, A^r \times C) \to MA^r \times M(\exists z\colon s.\, C) \ , \\
&\beta_t\colon M(\exists z\colon s.\, A^t \times C) \to MA^t \times M(\exists z\colon s.\, C) \ , \\
&\gamma_r\colon M(\exists z\colon s.\, A^r \times C) \to MB \ , \\
&\gamma_t\colon M(\exists z\colon s.\, A^t \times C) \to MB
\end{aligned}$$

are co-universal in $\mathbb{C}_{\exists z\colon s.\, A^r \times C}$ and $\mathbb{C}_{\exists z\colon s.\, A^t \times C}$, respectively.

By hypothesis, $\mathrm{FV}(r\colon C \to B) = \mathrm{FV}(t : C \to B)$, so $\alpha = \alpha_r = \alpha_t$, corresponding to $\vec{x}'$. In turn, it immediately follows that $\beta = \beta_r = \beta_t$ and $\gamma = \gamma_r = \gamma_t$.

Considering the product of arrows

$$\llbracket \vec{x}'.\, r : C \to B \rrbracket \times 1_{MC}\colon MA' \times MC \to MB^{MC} \times MC$$

by exponentiation

$$=\gamma \circ j_x$$

where $j_x\colon MA' \times MC \to M(\exists z\colon s.\, A' \times C)$ is the $x$-th injection of the unique co-cone of vertex $M(\exists z\colon s.\, A' \times C)$.

Similarly, it holds that $\llbracket \vec{x}'.\, t : C \to B \rrbracket \times 1_{MC} = \gamma \circ j_x$.

So, by exponentiation, $\llbracket \vec{x}'.\, r : C \to B \rrbracket = \llbracket \vec{x}'.\, t : C \to B \rrbracket$, from which it immediately follows that $\llbracket \vec{x} : \vec{A}.\, r : C \to B \rrbracket = \llbracket \vec{x} : \vec{A}.\, t : C \to B \rrbracket$;

($\exists_2$) $$\begin{aligned}
&\llbracket v : (\exists x\colon s.\, A).\, \mathrm{exE}(v, (\lambda x\colon s.\, (\lambda z\colon A.\, w[\mathrm{exI}_x(z)/v]))) : B \rrbracket \\
=&\gamma \circ \beta^{-1} \circ \left(\alpha, \llbracket v : (\exists x\colon s.\, A).\, v : (\exists x\colon s.\, A) \rrbracket\right) \\
=&\gamma \circ \beta^{-1} \circ \left(\alpha, 1_{M(\exists x\colon s.\, A)}\right)
\end{aligned}$$

and, since $\alpha\colon M(\exists x\colon s.\,A) \to M(\exists x\colon s.\,A)$ is a projection,

$$=\gamma\circ\beta^{-1}\circ\left(1_{M(\exists x\colon s.\,A)}, 1_{M(\exists x\colon s.\,A)}\right) \quad ,$$

where both arrows

$$\beta\colon M(\exists x\colon s.\,(\exists x\colon s.\,A)\times A) \to M(\exists x\colon s.\,A)\times M(\exists x\colon s.\,A)$$

and

$$\gamma\colon M(\exists x\colon s.\,(\exists x\colon s.\,A)\times A) \to MB$$

are co-universal in $\mathbb{C}_{\exists x\colon s.\,(\exists x\colon s.\,A)\times A}$.

But

$$\left(M(\exists x\colon s.\,A), \{\pi_1\colon M(\exists x\colon s.\,A)\times M(A[r/x]) \to M(\exists x\colon s.\,A)\}_{r\colon s\in \mathrm{LTerms}(\Sigma)}\right)$$

is a co-cone on $\Sigma_{(\exists x\colon s.\,A)\times A}(x\colon s)$ so there exists

$$\delta\colon M(\exists x\colon s.\,(\exists x\colon s.\,A)\times A) \to M(\exists x\colon s.\,A)$$

co-universal in $\mathbb{C}_{\exists x\colon s.\,(\exists x\colon s.\,A)\times A}$.

By exponentiation, there exists a unique $\theta_r\colon M(\exists x\colon s.\,A) \to MB^{M(A[r/x])}$) such that $[\![v\colon(\exists x\colon s.\,A).\,w\colon B]\!]\circ\pi_1 = \mathrm{ev}\circ(\theta_r\times 1_{M(A[r/x])})$ for each $r\colon s\in\mathrm{LTerms}(\Sigma)$. Thus,

$$\left(MB, \left\{\mathrm{ev}\circ(\theta_r\times 1_{M(A[r/x])})\right\}_{r\colon s\in\mathrm{LTerms}(\Sigma)}\right)$$

is a co-cone on $\Sigma_{(\exists x\colon s.\,A)\times A}(x\colon s)$ and so $[\![v\colon(\exists x\colon s.\,A).\,w\colon B]\!]\circ\delta$ is in $\mathbb{C}_{\exists x\colon s.\,(\exists x\colon s.\,A)\times A}$.

But $[\![v\colon(\exists x\colon s.\,A).\,w\colon B]\!]\circ\pi_1\circ\beta$ is in $\mathbb{C}_{\exists x\colon s.\,(\exists x\colon s.\,A)\times A}$ too, so by co-universality, $\gamma = [\![v\colon(\exists x\colon s.\,A).\,w\colon B]\!]\circ\delta = [\![v\colon(\exists x\colon s.\,A).\,w\colon B]\!]\circ\pi_1\circ\beta$.

Thus, $[\![v\colon(\exists x\colon s.\,A).\,\mathrm{exE}(v,(\lambda x\colon s.\,(\lambda z\colon A.\,w[\mathrm{exI}_x(z)/v]))) \colon B]\!]$

$$\begin{aligned}
&=\gamma\circ\beta^{-1}\circ\left(1_{M(\exists x\colon s.\,A)}, 1_{M(\exists x\colon s.\,A)}\right)\\
&=[\![v\colon(\exists x\colon s.\,A).\,w\colon B]\!]\circ\pi_1\circ\beta\circ\beta^{-1}\circ\left(1_{M(\exists x\colon s.\,A)}, 1_{M(\exists x\colon s.\,A)}\right)\\
&=[\![v\colon(\exists x\colon s.\,A).\,w\colon B]\!]\circ\pi_1\circ\left(1_{M(\exists x\colon s.\,A)}, 1_{M(\exists x\colon s.\,A)}\right)\\
&=[\![v\colon(\exists x\colon s.\,A).\,w\colon B]\!]\circ 1_{M(\exists x\colon s.\,A)}\\
&=[\![v\colon(\exists x\colon s.\,A).\,w\colon B]\!] \quad ;
\end{aligned}$$

($\exists_3$)

$$\begin{aligned}
&[\![\vec{x}\colon\vec{A}.\,\mathrm{exE}(\mathrm{exE}(a,(\lambda y\colon s.\,(\lambda z\colon D.\,b))),(\lambda y\colon s.\,c))\colon C]\!]\\
&=\gamma_1\circ\beta_1^{-1}\circ\left(\alpha_1, [\![\vec{x}\colon\vec{A}.\,\mathrm{exE}(a,(\lambda y\colon s.\,(\lambda z\colon D.\,b)))\colon(\exists y\colon s.\,B)]\!]\right)\\
&=\gamma_1\circ\beta_1^{-1}\circ\left(\alpha_1, \gamma_2\circ\beta_2^{-1}\circ\left(\alpha_2, [\![\vec{x}\colon\vec{A}.\,a\colon(\exists y\colon s.\,D)]\!]\right)\right)
\end{aligned}$$

where

- $\alpha_1\colon M\vec{A}\to MA^1$ and $\alpha_2\colon M\vec{A}\to MA^2$ are projections to the products corresponding to $\vec{x}^1 = \mathrm{FV}(\mathrm{exE}(a,(\lambda y\colon s.\,(\lambda z\colon D.\,b)))\colon(\exists y\colon s.\,B))\cup\mathrm{FV}(c\colon B\to C) = \mathrm{FV}(a\colon(\exists y\colon s.\,D))\cup(\mathrm{FV}(b\colon(\exists y\colon s.\,B))\setminus\{z\colon D\})\cup\mathrm{FV}(c\colon B\to C)$ (for easiness of reference, $\mathrm{FV}(a)\cup\mathrm{FV}(b)\cup\mathrm{FV}(c)$) and $\vec{x}^2 = \mathrm{FV}(a\colon(\exists y\colon s.\,D))\cup\mathrm{FV}((\lambda z\colon D.\,b)\colon D\to(\exists y\colon s.\,B)) = \mathrm{FV}(a\colon(\exists y\colon s.\,D)\cup(\mathrm{FV}(b\colon(\exists y\colon s.\,B))\setminus\{z\colon D\})$ (in the abbreviated form, $\mathrm{FV}(a)\cup\mathrm{FV}(b)$), respectively;
- $\beta^1\colon M(\exists y\colon s.\,A^1\times B)\to MA^1\times M(\exists y\colon s.\,B)$ and $\gamma^1\colon M(\exists y\colon s.\,A^1\times B)\to MC$ are co-universal in $\mathbb{C}_{\exists y\colon s.\,A^1\times B}$;
- $\beta^2\colon M(\exists y\colon s.\,A^2\times D)\to MA^2\times M(\exists y\colon s.\,D)$ and $\gamma^2\colon M(\exists y\colon s.\,A^2\times D)\to MC$ are co-universal in $\mathbb{C}_{\exists y\colon s.\,A^2\times D}$.

Similarly,

$$\begin{aligned}
&[\![\vec{x}\colon\vec{A}.\,\mathrm{exE}(a,(\lambda y\colon s.\,(\lambda z\colon D.\,\mathrm{exE}(b,(\lambda y\colon s.\,c)))))\colon C]\!]\\
&=\gamma_3\circ\beta_3^{-1}\circ\left(\alpha_3, [\![\vec{x}\colon\vec{A}.\,a\colon(\exists y\colon s.\,D)]\!]\right)
\end{aligned}$$

where

- $\alpha_3\colon M\vec{A} \to MA^3$ is the projection to the product corresponding to $\vec{x}^3 = \mathrm{FV}(a{:}(\exists y{:}\, s.\, D)) \cup \mathrm{FV}(\lambda z{:}\, D.\ \mathrm{exE}(b, (\lambda y{:}\, s.\, c)))$, i.e., $\mathrm{FV}(a) \cup \mathrm{FV}(b) \cup \mathrm{FV}(c)$ in abbreviated form;
- $\beta_3\colon M(\exists y{:}\, s.\, A^3 \times D) \to MA^3 \times M(\exists y{:}\, s.\, D)$ is co-universal in $\mathbb{C}_{\exists y:\, s.\, A^3 \times D}$;
- $\gamma_3\colon M(\exists y{:}\, s.\, A^3 \times D) \to MC$ is co-universal in $\mathbb{C}_{\exists y:\, s.\, A^3 \times D}$;

We notice that $\alpha_3 = \alpha_1$ and $A^3 = A^1$ as they correspond to the same set of free variables. Also, there is $\alpha_4\colon MA^1 \to MA^2$ such that $\alpha_3 = \alpha_4 \circ \alpha_1$ because $\vec{x}^3 \subseteq \vec{x}^1$.

Being $MA^2 \times M(\exists y{:}\, s.D)$ a product, $\alpha_4 \times 1_{M(\exists y:\, s.D)}\colon MA^1 \times M(\exists y{:}\, s.D) \to MA^2 \times M(\exists y{:}\, s.B)$ is universal for such a product. Thus,

$$\left(\pi_1, \gamma_2 \circ \beta_2^{-1} \circ (\alpha_4 \times 1_{M(\exists y:\, s.\, D)})\right)\colon MA^1 \times M(\exists y{:}\, s.\, D) \to MA^1 \times M(\exists y{:}\, s.\, B)$$

is universal for the product $MA^1 \times M(\exists y{:}\, s.\, B)$.

Moreover, $\gamma_2 \circ \beta_2^{-1} \circ (\alpha_4 \times 1_{M(\exists y:\, s.\, D)})$ preserves co-cones, as it is immediate to prove, being $\gamma_2$ and $\beta_2$ co-universal in $\mathbb{C}_{\exists y:\, s.\, A^2 \times D}$, so it is co-universal in $\mathbb{C}_{\exists y:\, s.\, B}$. Thus, $\left(\pi_1, \gamma_2 \circ \beta_2^{-1} \circ (\alpha_4 \times 1_{M(\exists y:\, s.\, D)})\right)$ preserves co-cones, too.

As a consequence, we have that

$$\gamma_3 \circ \beta_3^{-1} = \gamma_1 \circ \beta_1^{-1} \circ \left(\pi_1, \gamma_2 \circ \beta_2^{-1} \circ (\alpha_4 \times 1_{M(\exists y:\, s.\, D)})\right) \ .$$

So, applying these equalities, we obtain

$$\begin{aligned}
&[\![\vec{x}{:}\, \vec{A}.\ \mathrm{exE}(\mathrm{exE}(a, (\lambda y{:}\, s.\, (\lambda z{:}\, D.\, b))), (\lambda y{:}\, s.\, c)) : C]\!]\\
=&\gamma_1 \circ \beta_1^{-1} \circ \left(\alpha_1, \gamma_2 \circ \beta_2^{-1} \circ \left(\alpha_2, [\![\vec{x}{:}\, \vec{A}.\, a{:}(\exists y{:}\, s.\, D)]\!]\right)\right)\\
=&\gamma_1 \circ \beta_1^{-1} \circ \left(\alpha_1, \gamma_2 \circ \beta_2^{-1} \circ \left(\alpha_4 \circ \alpha_1, [\![\vec{x}{:}\, \vec{A}.\, a{:}(\exists y{:}\, s.\, D)]\!]\right)\right)\\
=&\gamma_1 \circ \beta_1^{-1} \circ \left(\alpha_1, \gamma_2 \circ \beta_2^{-1} \circ \left(\alpha_4 \times 1_{M(\exists y:\, s.\, D)}\right) \circ \left(\alpha_1, [\![\vec{x}{:}\, \vec{A}.\, a{:}(\exists y{:}\, s.\, D)]\!]\right)\right)\\
=&\gamma_1 \circ \beta_1^{-1} \circ \left(\pi_1, \gamma_2 \circ \beta_2^{-1} \circ \left(\alpha_4 \times 1_{M(\exists y:\, s.\, D)}\right)\right) \circ \left(\alpha_1, [\![\vec{x}{:}\, \vec{A}.\, a{:}(\exists y{:}\, s.\, D)]\!]\right)\\
=&\gamma_3 \circ \beta_3^{-1} \circ \left(\alpha_1, [\![\vec{x}{:}\, \vec{A}.\, a{:}(\exists y{:}\, s.\, D)]\!]\right)\\
=&\gamma_3 \circ \beta_3^{-1} \circ \left(\alpha_3, [\![\vec{x}{:}\, \vec{A}.\, a{:}(\exists y{:}\, s.\, D)]\!]\right)\\
=&[\![\vec{x}{:}\, \vec{A}.\ \mathrm{exE}(a, (\lambda y{:}\, s.\, (\lambda z{:}\, D.\ \mathrm{exE}(b, (\lambda y{:}\, s.\, c))))) : C]\!] \ ;
\end{aligned}$$

$(\exists_4)$
$$\begin{aligned}
&[\![\vec{x}{:}\, \vec{A}.\ \mathrm{exE}(a, (\lambda y{:}\, s.\, (\lambda z{:}\, C.\, b[\mathrm{exI}_y(z)/w]))) : B]\!]\\
=&\gamma \circ \beta^{-1} \circ \left(\alpha, [\![\vec{x}{:}\, \vec{A}.\, a{:}(\exists y{:}\, s.\, B)]\!]\right)
\end{aligned}$$

where

- $\alpha\colon M\vec{A} \to MA'$ is the projection to the product corresponding to

$$\begin{aligned}
\vec{x}' &= \mathrm{FV}(a{:}(\exists y{:}\, s.\, C)) \cup \mathrm{FV}((\lambda z{:}\, C.\, b[\mathrm{exI}_y(z)/w]))\\
&= \mathrm{FV}(a{:}(\exists y{:}\, s.\, C)) \cup (\mathrm{FV}(b{:}\, B) \setminus \{z{:}\, C\}) \ ;
\end{aligned}$$

- $\beta\colon M(\exists y{:}\, s.\, A' \times C) \to MA' \times M(\exists y{:}\, s.\, C)$ is co-universal in $\mathbb{C}_{\exists y:\, s.\, A' \times C}$;
- $\gamma\colon M(\exists y{:}\, s.\, A' \times C) \to MB$ is co-universal in $\mathbb{C}_{\exists y:\, s.\, A' \times C}$.

By definition, $\beta$ is the arrow which maps the co-cone

$$\left(M(\exists y{:}\, s.\, A' \times C), \{[\![v{:}\, A' \times C[r/z].\ \mathrm{exI}_y(v) : (\exists y{:}\, s.\, A' \times C)]\!]\}_{r:\, s \in \mathrm{LTerms}(\Sigma)}\right)$$

in the co-cone

$$\left(MA' \times M(\exists y{:}\, s.\, C), \{1_{MA'} \times [\![u{:}\, C[r/z].\ \mathrm{exI}_y(u) : (\exists y{:}\, s.\, C)]\!]\}_{r:\, s \in \mathrm{LTerms}(\Sigma)}\right) \ .$$

In turn, $[\![\vec{x}', w{:}(\exists y{:}\, s.\, C).\, b{:}\, B]\!]$ maps the co-cone above with vertex $MA' \times M(\exists y{:}\, s.\, C)$ in

$$\begin{aligned}
&\Big(MB, \Big\{\ [\![\vec{x}', w{:}(\exists y{:}\, s.\, C).\, b{:}\, B]\!]\\
&\qquad \circ \left(1_{MA'} \times [\![u{:}\, C[r/z].\ \mathrm{exI}_y(u) : (\exists y{:}\, s.\, C)]\!]\right)\Big\}_{r:\, s \in \mathrm{LTerms}(\Sigma)}\Big) \ .
\end{aligned}$$

Thus, $[\![\vec{x}', w{:}(\exists y{:}s.\,C).\,b{:}B]\!] \circ \beta$ is an arrow of $\mathbb{C}_{\exists y:s.\,A'\times C}$, and it equals $\gamma$ by co-universality. So,

$$
\begin{aligned}
&[\![\vec{x}{:}\vec{A}.\,\mathrm{exE}(a,(\lambda y{:}s.\,(\lambda z{:}C.\,b[\mathrm{exI}_y(z)/w]))){:}B]\!]\\
=&\gamma \circ \beta^{-1} \circ \big(\alpha, [\![\vec{x}{:}\vec{A}.\,a{:}(\exists y{:}s.\,B)]\!]\big)\\
=&[\![\vec{x}', w{:}(\exists y{:}s.\,C).\,b{:}B]\!] \circ \beta \circ \beta^{-1} \circ \big(\alpha, [\![\vec{x}{:}\vec{A}.\,a{:}(\exists y{:}s.\,B)]\!]\big)\\
=&[\![\vec{x}', w{:}(\exists y{:}s.\,C).\,b{:}B]\!] \circ \big(\alpha, [\![\vec{x}{:}\vec{A}.\,a{:}(\exists y{:}s.\,B)]\!]\big)
\end{aligned}
$$

and, by lemma 5.3,

$$
=[\![\vec{x}{:}\vec{A}.\,b[a/w]{:}B]\!] \quad . \qquad \square
$$

## 6. Completeness

In this section, we will prove that the $\lambda$-calculi of section 5 are complete with respect to models in logically distributive categories. Moreover, from this result and the soundness theorem, we will prove that the logical systems of section 4 are sound and complete with respect to the same class of models.

The first part of the completeness proof aims at constructing a canonical model based on the language of the $\lambda$-calculus.

**Definition 6.1** (Syntactical equivalence)**.** Given a $\lambda$-theory $T$, the *syntactical equivalence* of two terms-in-context is defined by fixing the generated equivalence classes. Precisely, the equivalence class $[x\colon A.\, t\colon B]$ is defined as the minimal set, composed by terms-in-context, such that

(1) $x\colon A.\, t\colon B \in [x\colon A.\, t\colon B]$—reflexivity;
(2) if $T \vdash \vec{y}\colon\vec{D}.\, s =_C r$, where $\vec{y}\colon\vec{D}.\, s =_C r$ is an equality-in-context, and $\vec{y}\colon\vec{D}.\, s\colon C \in [x\colon A.\, t\colon B]$, then $\vec{y}\colon\vec{D}.\, r\colon C \in [x\colon A.\, t\colon B]$—closure under provable equivalence;
(3) if $\vec{y}\colon\vec{D}.\, s\colon C$ is a term-in-context and, for some $1 \leq i < m$ and $z\colon D_i \times D_{i+1} \not\in \mathrm{FV}(s\colon C) \cup \{y_1\colon D_1, \dots, y_m\colon D_m\}$, it happens that
$$y_1\colon D_1, \dots, y_{i-1}\colon D_{i-1}, z\colon D_i \times D_{i+1}, y_{i+1}\colon D_{i+2}, \\ \dots, y_m\colon D_m.\, s[\mathrm{fst}(z)/y_i][\mathrm{snd}(z)/y_{i+1}]\colon C \in [x\colon A.\, t\colon B] \enspace ,$$
then $\vec{y}\colon\vec{D}.\, s\colon C \in [x\colon A.\, t\colon B]$—closure under associativity in contexts;
(4) if $\vec{y}\colon\vec{D}.\, s\colon C$ is a term-in-context and, for some $1 \leq i < m$ and $z\colon D_{i+1} \times D_i \not\in \mathrm{FV}(s\colon C) \cup \{y_1\colon D_1, \dots, y_m\colon D_m\}$, it happens that
$$y_1\colon D_1, \dots, y_{i-1}\colon D_{i-1}, z\colon D_{i+1} \times D_i, y_{i+1}\colon D_{i+2}, \\ \dots, y_m\colon D_m.\, s[\mathrm{snd}(z)/y_i][\mathrm{fst}(z)/y_{i+1}]\colon C \in [x\colon A.\, t\colon B] \enspace ,$$
then $\vec{y}\colon\vec{D}.\, s\colon C \in [x\colon A.\, t\colon B]$—closure under commutativity in contexts;
(5) if $\vec{y}\colon\vec{D}.\, s\colon C \in [x\colon A.\, y\colon B]$ and $z\colon D_i \not\in \mathrm{FV}(s\colon C) \cup \{y_1\colon D_1, \dots, y_m\colon D_m\}$ for some $1 \leq i \leq m$, then
$$y_1\colon D_1, \dots, y_{i-1}\colon D_{i-1}, z\colon D_i, y_{i+1}\colon D_{i+1}, \dots, y_m\colon D_m.\, s[z/y_i]\colon C \in [x\colon A.\, t\colon B]$$
—closure under $\alpha$-renaming in contexts.

Since representatives of equivalence classes have a contest composed by only one variable, we have to prove that it suffices.

**Proposition 6.1.** *Each term-in-context $x_1\colon A_1, \dots, x_n\colon A_n.\, t\colon B$ appears in some equivalence class $[y\colon C.\, s\colon B]$.*

*Proof.* Let $C = A_1 \times (A_2 \times (\cdots (A_{n-1} \times A_n) \cdots))$ and
$$s \equiv t[\mathrm{fst}(y)/x_1][\mathrm{fst}(\mathrm{snd}(y))/x_2] \cdots [\mathrm{snd}(\mathrm{snd}(\cdots(\mathrm{snd}(y))\cdots))/x_n] \enspace ,$$
where $y\colon C$ is a new variable.

By closure under $\alpha$-renaming in contexts, the class $[y\colon C.\, s\colon B]$ is uniquely defined.

By iteration of closure under associativity of contexts, $\vec{x}\colon\vec{A}.\, t\colon B \in [y\colon C.\, s\colon B]$. □

Also, we want to show that equivalent terms in the same context are in the same equivalence class. These facts are needed to show that equivalence classes form a partition of the terms-in-context.

**Proposition 6.2.** *If $\vec{x}\colon\vec{A}.\, t\colon B$ and $\vec{x}\colon\vec{A}.\, s\colon B$ are two terms-in-context in the same equivalence class $[y\colon C.\, r\colon D]$, then $T \vdash \vec{x}\colon\vec{A}.\, t =_B s$ and $B = D$.*

*Proof.* An immediate induction on the definition of equivalence classes shows that $B = D$, $T \vdash \vec{x}\colon\vec{A}, y\colon C.\, t =_B r$, and $T \vdash \vec{x}\colon\vec{A}, y\colon C.\, s =_B r$. The conclusion follows by $(\mathrm{eq}_4)$ and $(\mathrm{eq}_0)$. □

**Proposition 6.3.** *Syntactical equivalence is an equivalence relation.*

*Proof.* By proposition 6.1, we know that each term-in-context belongs to some equivalence class. Suppose $\vec{x}:\vec{A}.\,t:B$ belongs to $[y_1:C_1.\,r_1:B]$ and $[y_2:C_2.\,r_2:B]$. Then, by proposition 6.2, $T\vdash \vec{x}:\vec{A}, y_1:C_1.\,t=_B r_1$ and $T\vdash \vec{x}:\vec{A}, y_2:C_2.\,t=_B r_2$. By ($eq_4$) and ($eq_0$), $T\vdash \vec{x}:\vec{A}, y_1:C_1, y_2:C_2.\,r_1=_B r_2$. Thus, both $\vec{x}:\vec{A}, y_1:C_1, y_2:C_2.\,r_1:B\in[y_2:C_2.\,r_2:B]$ and $\vec{x}:\vec{A}, y_1:C_1, y_2:C_2.\,r_2:B\in[y_1:C_1.\,r_1:B]$.

Now, $\vec{x}:\vec{A}, y_1:C_1, y_2:C_2.\,r_1:B\in[y_1:C_1.\,r_1:B]$, because of ($eq_3$) and ($eq_0$).

Also, $\vec{x}:\vec{A}, y_1:C_1, y_2:C_2.\,r_2:B\in[y_2:C_2.\,r_2:B]$, for the same reason.

Thus, by proposition 6.2, every term-in-context in $[y_1:C_1.\,r_1:B]$ is in $[y_2:C_2.\,r_2:B]$ and every term-in-context in $[y_2:C_2.\,r_2:B]$ belongs to $[y_1:C_1.\,r_1:B]$, i.e., the two equivalence classes are the same.

Since syntactical equivalence induces a partition, it is an equivalence relation. □

The canonical logically distributive category is now immediate to define. Of course, we have to prove that it is logically distributive.

**Definition 6.2** (Syntactical category)**.** Given a $\lambda$-theory $T$. the *syntactical category* $\mathbb{C}_T$ has $\lambda$Types($\Sigma$) as objects, where $\Sigma$ is the $\lambda$-signature of $T$, and the equivalence classes $[x:A.\,t:B]\colon A\to B$ as arrows. Identities are given by the classes $[x:A.\,x:A]\colon A\to A$ for each $\lambda$-type $A$, and composition is given by substitution:

$$[y:B.\,s:C]\circ[x:A.\,t:B]=[x:A.\,s[t/y]:C]\ \ .$$

Moreover, the map $M_T\colon \lambda\mathrm{Types}(\Sigma)\to \mathrm{Obj}\,\mathbb{C}_T$ is defined as $M_T A=A$. We will omit subscripts when they are clear from the context.

**Proposition 6.4.** *The* $\mathbb{C}_T$ *category has finite products.*

*Proof.* It suffices to prove that $\mathbb{C}_T$ has terminal object and binary products.

The terminal object of $\mathbb{C}_T$ is the $\lambda$-type 1: for any object $A$, there is a morphism $A\to 1$, i.e., $[x:A.\,*:1]$. Also, if $[x:A.\,t:1]\colon A\to 1$, then $T\vdash x:A.\,t=_1 *$ by ($\times_0$) and ($eq_0$).

The product of $A$ and $B$ is $A\times B$ with projections $[z:A\times B.\,\mathrm{fst}(z):A]\colon A\times B\to A$ and $[z:A\times B.\,\mathrm{snd}(z):B]\colon A\times B\to B$. The universal arrow induced by $[w:C.\,s:A]\colon C\to A$ and $[w:C.\,t:B]\colon C\to B$ is $[w:C.\,\langle s,t\rangle:A\times B]\colon C\to A\times B$, because

$$\begin{aligned}&[z:A\times B.\,\mathrm{fst}(z):A]\circ[w:C.\,\langle s,t\rangle:A\times B]\\=&[w:C.\,\mathrm{fst}(\langle s,t\rangle):A]\\=&[w:C.\,s:A]\end{aligned}$$

by ($\times_1$) and

$$\begin{aligned}&[z:A\times B.\,\mathrm{snd}(z):B]\circ[w:C.\,\langle s,t\rangle:A\times B]\\=&[w:C.\,\mathrm{snd}(\langle s,t\rangle):B]\\=&[w:C.\,t:B]\end{aligned}$$

by ($\times_2$). If $[w:C.\,r:A\times B]$ is such that $[z:A\times B.\,\mathrm{fst}(z):A]\circ[w:C.\,r:A\times B]=[w:C.\,s:A]$ and $[z:A\times B.\,\mathrm{snd}(z):B]\circ[w:C.\,r:A\times B]=[w:C.\,t:B]$, then $[w:C.\,\mathrm{fst}(r):A]=[w:C.\,s:A]$ and $[w:C.\,\mathrm{snd}(r):B]=[w:C.\,t:B]$, that is, $T\vdash w:C.\,\mathrm{fst}(r)=_A s$ and $T\vdash w:C.\,\mathrm{snd}(r)=_B t$. Then $T\vdash w:C.\,\langle s,t\rangle=_{A\times B} r$ by ($\times_3$), ($eq_0$) and ($eq_1$). □

**Proposition 6.5.** *The* $\mathbb{C}_T$ *category has finite co-products.*

*Proof.* It suffices to prove that $\mathbb{C}_T$ has initial object and binary co-products.

The initial object is the $\lambda$-type 0: for any object $A$, $[x:0.\,\mathrm{F}_A\cdot x:A]\colon 0\to A$. Also, if $[x:0.\,t:A]$ then $[x:0.\,t:A]=[x:0.\,\mathrm{F}_A\cdot x:A]$ by ($+_3$) and ($eq_0$).

The co-product of $A$ and $B$ is $A+B$ with injections $[z\colon A.\, \mathrm{inl}_B(z)\colon A+B]\colon A\to A+B$ and $[z\colon B.\, \mathrm{inr}_A(z)\colon A+B]\colon B\to A+B$. The co-universal arrow induced by $[u\colon A.\, t\colon C]\colon A\to C$ and $[v\colon B.\, s\colon C]\colon B\to C$ is $[z\colon A+B.\, \mathrm{when}(z,(\lambda u\colon A.\, t),(\lambda v\colon B.\, s))\colon C]\colon A+B\to C$ since

$$\begin{aligned}
&[z\colon A+B.\, \mathrm{when}(z,(\lambda u\colon A.\, t),(\lambda v\colon B.\, s))\colon C]\circ[u\colon A.\, \mathrm{inl}_B(u)\colon A+B]\\
=&[u\colon A.\, \mathrm{when}(\mathrm{inl}_B(u),(\lambda u\colon A.\, t),(\lambda v\colon B.\, s))\colon C]\\
=&[u\colon A.\, (\lambda u\colon A.\, t)\cdot u\colon C]\\
=&[u\colon A.\, t\colon C]\quad ,
\end{aligned}$$

by $(+_0)$ and $(\to_0)$, and

$$\begin{aligned}
&[z\colon A+B.\, \mathrm{when}(z,(\lambda u\colon A.\, t),(\lambda v\colon B.\, s))\colon C]\circ[v\colon B.\, \mathrm{inr}_A(v)\colon A+B]\\
=&[v\colon B.\, \mathrm{when}(\mathrm{inr}_A(v),(\lambda u\colon A.\, t),(\lambda v\colon B.\, s))\colon C]\\
=&[v\colon B.\, (\lambda v\colon B.\, s)\cdot v\colon C]\\
=&[v\colon B.\, s\colon C]
\end{aligned}$$

by $(+_1)$ and $(\to_0)$. Also, if $[z\colon A+B.\, r\colon C]$ is such that $[z\colon A+B.\, r\colon C]\circ[u\colon A.\, \mathrm{inl}_B(u)\colon A+B]=[u\colon A.\, t\colon C]$ and $[z\colon A+B.\, r\colon C]\circ[v\colon B.\, \mathrm{inr}_A(v)\colon A+B]\circ[v\colon B.\, s\colon C]$, then

$$\begin{aligned}
&[u\colon A.\, r[\mathrm{inl}_B(u)/z]\colon C]\\
=&[z\colon A+B.\, r\colon C]\circ[u\colon A.\, \mathrm{inl}_B(u)\colon A+B]\\
=&[u\colon A.\, t\colon C]
\end{aligned}$$

by $(\to_0)$

$$=[u\colon A.\, (\lambda u\colon A.\, t)\cdot u\colon C]\quad ,$$

and

$$\begin{aligned}
&[v\colon B.\, r[\mathrm{inr}_A(v)/z]\colon C]\\
=&[z\colon A+B.\, r\colon C]\circ[v\colon B.\, \mathrm{inr}_A(v)\colon A+B]\\
=&[v\colon B.\, s\colon C]
\end{aligned}$$

by $(\to_0)$

$$=[v\colon B.\, (\lambda v\colon B.\, s)\cdot v\colon C]\quad .$$

Thus, $T\vdash z\colon A+B.\, r=_C \mathrm{when}(z,(\lambda u\colon A.\, t),(\lambda v\colon B.\, s))$ by $(+_3)$. □

**Proposition 6.6.** *The $\mathbb{C}_T$ category has exponentiation.*

*Proof.* The exponential object $B^A$ is the $\lambda$-type $A\to B$ and the corresponding evaluation arrow is $[w\colon(A\to B)\times A.\, \mathrm{fst}(w)\cdot\mathrm{snd}(w)\colon B]\colon (A\to B)\times A\to B$. Given any morphism $[z\colon C\times A.\, t\colon B]\colon C\times A\to B$, its exponential transpose is $[w\colon C.\, (\lambda x\colon A.\, t[\langle w,x\rangle/z])\colon A\to B]$. In fact,

$$\begin{aligned}
&\mathrm{ev}\circ([w\colon C.\, (\lambda x\colon A.\, t[\langle w,x\rangle/z])\colon A\to B]\times 1_A)\\
=&[w\colon(A\to B)\times A.\, \mathrm{fst}(w)\cdot\mathrm{snd}(w)\colon B]\\
&\circ([w\colon C.\, (\lambda x\colon A.\, t[\langle w,x\rangle/z])\colon A\to B]\times[x\colon A.\, x\colon A])\\
=&[w\colon(A\to B)\times A.\, \mathrm{fst}(w)\cdot\mathrm{snd}(w)\colon B]\\
&\circ[w\colon C\times A.\, \langle(\lambda x\colon A.\, t[\langle\mathrm{fst}(w),x\rangle/z]),\mathrm{snd}(w)\rangle\colon(A\to B)\times A]\\
=&[w\colon C\times A.\, (\lambda x\colon A.\, t[\langle\mathrm{fst}(w),x\rangle/z])\cdot\mathrm{snd}(w)\colon B]\\
=&[w\colon C\times A.\, t[\langle\mathrm{fst}(w),\mathrm{snd}(w)\rangle/z]\colon B]\\
=&[w\colon C\times A.\, t[w/z]\colon B]\\
=&[z\colon C\times A.\, t\colon B]\quad .
\end{aligned}$$

Also, let $[z\colon C.\, r\colon A\to B]\colon C\to(A\to B)$ be such that $\mathrm{ev}\circ([z\colon C.\, r\colon A\to B]\times 1_A)=[z\colon C\times A.\, t\colon B]$. Then,

$$
\begin{aligned}
&\mathrm{ev}\circ([z\colon C.\, r\colon A\to B]\times 1_A)\\
&=[w\colon (A\to B)\times A.\, \mathrm{fst}(w)\cdot \mathrm{snd}(w)\colon B]\\
&\quad\circ([z\colon C.\, r\colon A\to B]\times[x\colon A.\, x\colon A])\\
&=[w\colon (A\to B)\times A.\, \mathrm{fst}(w)\cdot \mathrm{snd}(w)\colon B]\\
&\quad\circ[w\colon C\times A.\, \langle r[\mathrm{fst}(w)/z], \mathrm{snd}(w)\rangle\colon (A\to B)\times B]\\
&=[w\colon C\times A.\, r[\mathrm{fst}(w)/z]\cdot \mathrm{snd}(w)\colon B]\\
&=[w\colon C\times A.\, (\lambda z\colon C.\, r)\cdot \mathrm{fst}(w))\cdot \mathrm{snd}(w)\colon B]
\end{aligned}
$$

and, by hypothesis,

$$
\begin{aligned}
&=[w\colon C\times A.\, t[w/z]\colon B]\\
&=[w\colon C\times A.\, t[\langle \mathrm{fst}(w), \mathrm{snd}(w)\rangle/z]\colon B]\\
&=[w\colon C\times A.\, (\lambda x\colon A.\, t[\langle \mathrm{fst}(w), x\rangle/z])\cdot \mathrm{snd}(w)\colon B]\\
&=[w\colon C\times A.\, ((\lambda v\colon C.\, (\lambda x\colon A.\, t[\langle v, x\rangle/z]))\cdot \mathrm{fst}(w))\cdot \mathrm{snd}(w)\colon B]
\end{aligned}
$$

Thus, $r$ equals the exponential transpose modulo $\alpha$-renaming. □

**Proposition 6.7.** *The $\mathbb{C}_T$ category is distributive.*

*Proof.* Define

$$
\begin{aligned}
\Delta^{-1} = [w\colon A\times(B+C).\, &\mathrm{when}(\mathrm{snd}(w), (\lambda x\colon B.\, \mathrm{inl}_{A\times C}(\langle \mathrm{fst}(w), x\rangle)),\\
&(\lambda y\colon C.\, \mathrm{inr}_{A\times B}(\langle \mathrm{fst}(w), y\rangle)))\colon (A\times B)+(A\times C)]\ .
\end{aligned}
$$

It holds that

$$
\begin{aligned}
\Delta &= [1_{MA}\times\iota_1, 1_{MA}\times\iota_2]\\
&= [[w\colon A\times B.\, \langle \mathrm{fst}(w), \mathrm{inl}_C(\mathrm{snd}(w))\rangle\colon A\times(B+C)],\\
&\quad [w\colon A\times C.\, \langle \mathrm{fst}(w), \mathrm{inr}_B(\mathrm{snd}(w))\rangle\colon A\times(B+C)]]\\
&= [z\colon (A\times B)+(A\times C).\, \mathrm{when}(z,\\
&\qquad (\lambda w\colon A+B.\, \langle \mathrm{fst}(w), \mathrm{inl}_C(\mathrm{snd}(w))\rangle),\\
&\qquad (\lambda w\colon A+C.\, \langle \mathrm{fst}(w), \mathrm{inr}_B(\mathrm{snd}(w))\rangle))\colon A\times(B+C)]\ .
\end{aligned}
$$

Considering

$$
\begin{aligned}
&\Delta^{-1}\circ\Delta\circ[u\colon A\times B.\, \mathrm{inl}_{A\times C}(u)\colon (A\times B)+(A\times C)]\\
&=\Delta^{-1}\circ[z\colon (A\times B)+(A\times C).\, \mathrm{when}(z,\\
&\qquad (\lambda w\colon A+B.\, \langle \mathrm{fst}(w), \mathrm{inl}_C(\mathrm{snd}(w))\rangle),\\
&\qquad (\lambda w\colon A+C.\, \langle \mathrm{fst}(w), \mathrm{inr}_B(\mathrm{snd}(w))\rangle))\colon A\times(B+C)]\\
&\quad\circ[u\colon A\times B.\, \mathrm{inl}_{A\times C}(u)\colon (A\times B)+(A\times C)]\\
&=\Delta^{-1}\circ[u\colon A\times B.\, \mathrm{when}(\mathrm{inl}_{A\times C}(u),\\
&\qquad (\lambda w\colon A+B.\, \langle \mathrm{fst}(w), \mathrm{inl}_C(\mathrm{snd}(w))\rangle),\\
&\qquad (\lambda w\colon A+C.\, \langle \mathrm{fst}(w), \mathrm{inr}_B(\mathrm{snd}(w))\rangle))\colon A\times(B+C)]
\end{aligned}
$$

by $(\exists_0)$

$$
\begin{aligned}
&=\Delta^{-1}\circ[u\colon A\times B.\, (\lambda w\colon A+B.\, \langle \mathrm{fst}(w), \mathrm{inl}_C(\mathrm{snd}(w))\rangle)\cdot u\colon A\times(B+C)]\\
&=\Delta^{-1}\circ[u\colon A\times B.\, \langle \mathrm{fst}(u), \mathrm{inl}_C(\mathrm{snd}(u))\rangle\colon A\times(B+C)]
\end{aligned}
$$

expanding $\Delta^{-1}$

$$\begin{aligned}
&=[w\colon A\times(B+C).\,\text{when}(\text{snd}(w),(\lambda x\colon B.\,\text{inl}_{A\times C}(\langle\text{fst}(w),x\rangle)),\\
&\qquad(\lambda y\colon C.\,\text{inr}_{A\times B}(\langle\text{fst}(w),y\rangle))):(A\times B)+(A\times C)]\\
&\quad\circ[u\colon A\times B.\,\langle\text{fst}(u),\text{inl}_C(\text{snd}(u))\rangle : A\times(B+C)]\\
&=[u\colon A\times B.\,\text{when}(\text{snd}(\langle\text{fst}(u),\text{inl}_C(\text{snd}(u))\rangle),\\
&\qquad(\lambda x\colon B.\,\text{inl}_{A\times C}(\langle\text{fst}(\langle\text{fst}(u),\text{inl}_C(\text{snd}(u))\rangle),x\rangle)),\\
&\qquad(\lambda y\colon C.\,\text{inr}_{A\times B}(\langle\text{fst}(\langle\text{fst}(u),\text{inl}_C(\text{snd}(u))\rangle),y\rangle))):(A\times B)+(A\times C)]\\
&=[u\colon A\times B.\,\text{when}(\text{inl}_C(\text{snd}(u)),(\lambda x\colon B.\,\text{inl}_{A\times C}(\langle\text{fst}(u),x\rangle)),\\
&\qquad(\lambda y\colon C.\,\text{inr}_{A\times B}(\langle\text{fst}(u),y\rangle))):(A\times B)+(A\times C)]
\end{aligned}$$

by $(\exists_0)$

$$\begin{aligned}
&=[u\colon A\times B.\,(\lambda x\colon B.\,\text{inl}_{A\times C}(\langle\text{fst}(u),x\rangle))\cdot\text{snd}(u):(A\times B)+(A\times C)]\\
&=[u\colon A\times B.\,\text{inl}_{A\times C}(\langle\text{fst}(u),\text{snd}(u)\rangle):(A\times B)+(A\times C)]\\
&=[u\colon A\times B.\,\text{inl}_{A\times C}(u):(A\times B)+(A\times C)]\quad.
\end{aligned}$$

Similarly, $\Delta^{-1}\circ\Delta\circ[v\colon A\times C.\,\text{inr}_{A\times B}(v):(A\times B)+(A\times C)]=[v\colon A\times C.\,\text{inr}_{A\times B}(v):(A\times B)+(A\times C)]$. So, by co-universality of $1_{(A\times B)+(A\times C)}$, $\Delta^{-1}\circ\Delta=1_{(A\times B)+(A\times C)}$.

It holds that,

$$\begin{aligned}
&\Delta^{-1}\circ[w\colon A\times B.\,\langle\text{fst}(w),\text{inl}_C(\text{snd}(w))\rangle : A\times(B+C)]\\
&=[w\colon A\times(B+C).\,\text{when}(\text{snd}(w),\\
&\qquad(\lambda x\colon B.\,\text{inl}_{A\times C}(\langle\text{fst}(w),x\rangle)),\\
&\qquad(\lambda y\colon C.\,\text{inr}_{A\times B}(\langle\text{fst}(w),y\rangle))):(A\times B)+(A\times C)]\\
&\quad\circ[w\colon A\times B.\,\langle\text{fst}(w),\text{inl}_C(\text{snd}(w))\rangle : A\times(B+C)]\\
&=[w\colon A\times B.\,\text{when}(\text{snd}(\langle\text{fst}(w),\text{inl}_C(\text{snd}(w))),\\
&\qquad(\lambda x\colon B.\,\text{inl}_{A\times C}(\langle\text{fst}(w),x\rangle)),\\
&\qquad(\lambda y\colon C.\,\text{inr}_{A\times B}(\langle\text{fst}(w),y\rangle))):(A\times B)+(A\times C)]\\
&=[w\colon A\times B.\,\text{when}(\text{inl}_C(\text{snd}(w)),\\
&\qquad(\lambda x\colon B.\,\text{inl}_{A\times C}(\langle\text{fst}(w),x\rangle)),\\
&\qquad(\lambda y\colon C.\,\text{inr}_{A\times B}(\langle\text{fst}(w),y\rangle))):(A\times B)+(A\times C)]\\
&=[w\colon A\times B.\,(\lambda x\colon B.\,\text{inl}_{A\times C}(\langle\text{fst}(w),x\rangle))\cdot\text{snd}(w):(A\times B)+(A\times C)]\\
&=[w\colon A\times B.\,\text{inl}_{A\times C}(\langle\text{fst}(w),\text{snd}(w)\rangle):(A\times B)+(A\times C)]\\
&=[w\colon A\times B.\,\text{inl}_{A\times C}(w):(A\times B)+(A\times C)]=\iota_1\quad.
\end{aligned}$$

Since $\Delta$ is co-universal, it is the unique arrow such that

$$\Delta\circ\iota_1=[w\colon A\times B.\,\text{inl}_{A\times C}(w):(A\times B)+(A\times C)]\quad,$$

so, calling $\alpha=[w\colon A\times B.\,\langle\text{fst}(w),\text{inl}_C(\text{snd}(w))\rangle : A\times(B+C)]$, $\Delta\circ\Delta^{-1}\circ\alpha=\alpha$. But $\alpha$ is universal for the product $A\times(B+C)$, so $\Delta\circ\Delta^{-1}=1_{A\times(B+C)}$. □

**Proposition 6.8.** *In the $\mathbb{C}_T$ category, $M_T(\forall x\colon s.\,A)$ is the terminal object of $\mathbb{C}_{\forall x\colon s.\,A}$.*

*Proof.* Since $x\colon s\notin\text{FV}(\forall x\colon s.\,A)$ and

$$\left(M(\forall x\colon s.\,A),\{[u\colon(\forall x\colon s.\,A).\,\text{allE}(u,r):A[r/x]]\}_{r\colon s\in\text{LTerms}(\Sigma)}\right)$$

is a cone on $\Sigma_A(x\colon s)$, $M(\forall x\colon s.\,A)\in\text{Obj}\,\mathbb{C}_{\forall x\colon s.\,A}$.

If $MB\in\text{Obj}\,\mathbb{C}_{\forall x\colon s.\,A}$, then there is a cone on $\Sigma_A(x\colon s)$:

$$\left(MB,\{[v\colon B.\,t_r : A[r/x]]\}_{r\colon s\in\text{LTerms}(\Sigma)}\right)\quad.$$

But also

$$\big(MB, \{[v:B.\, t_x[r/x] : A[r/x]]\}_{r:s \in \mathrm{LTerms}(\Sigma)}\big)$$

is a cone on $\Sigma_A(x:s)$ since $x:s \notin \mathrm{FV}(B)$ and so

$$[(u:B)[r/x].\, t_x[r/x] : A[r/x]] = [u:B.\, t_x[r/x] : A[r/x]] \ .$$

Consider $[v:B.\, \mathrm{allI}(\lambda x:s.\, t_x) :(\forall x:s.\, A)]$: it holds that

$$\begin{aligned} &[u:(\forall x:s.\ A).\, \mathrm{allE}(u,r) : A[r/x]] \circ [v:B.\, \mathrm{allI}(\lambda x:s.\, t_x) :(\forall x:s.\, A)] \\ =&[v:B.\, \mathrm{allE}(\mathrm{allI}(\lambda x:s.\, t_x), r) : A[r/x]] \\ =&[v:B.\, t_x[r/x] : A[r/x]] \ . \end{aligned}$$

Thus, $[v:B.\, \mathrm{allI}(\lambda x:s.\, t_x) :(\forall x:s.\, A)]$ is an arrow in $\mathbb{C}_{\forall x:s.\,A}$.

Let $[v:B.\, w :(\forall x:s.\, A)]$ be an arrow in $\mathbb{C}_{\forall x:s.\,A}$. For each $r:s \in \mathrm{LTerms}(\Sigma)$,

$$\begin{aligned} &[v:B.\, \mathrm{allE}(w,r) : A[r/x]] \\ =&[u:(\forall x:s.\, A).\, \mathrm{allE}(u,r) : A[r/x]] \circ [v:B.\, w :(\forall x:s.\, A)] \end{aligned}$$

and, being an arrow between cones,

$$=[v:B.\, t_x[r/x] : A[r/x]] \ .$$

Thus $[v:B.\, \mathrm{allE}(\mathrm{allI}(\lambda x:s.\, t_x), r) : A[r/x]] = [v:B.\, \mathrm{allE}(w,r) : A[r/x]]$.

Applying $(\forall_1)$, $[v:B.\, w :(\forall x:s.\, A)] = [v:B.\, \mathrm{allI}(\lambda x:s.\, t_x) :(\forall x:s.\, A)]$. □

**Proposition 6.9.** *In the $\mathbb{C}_T$ category, $M_T(\exists x:s.\, A)$ is the initial object of $\mathbb{C}_{\exists x:s.\,A}$.*

*Proof.* Since $x:s \notin \mathrm{FV}(\exists x:s.\, A)$ and

$$\big(M(\exists x:s.\, A), \{[u:a[r/x].\, \mathrm{exI}_x(u) :(\exists x:s.\ A)]\}_{r:s \in \mathrm{LTerms}(\Sigma)}\big)$$

is a co-cone on $\Sigma_A(x:s)$, $M(\exists x:s.\, A) \in \mathrm{Obj}\,\mathbb{C}_{\exists x:s.\,A}$.

If $MB \in \mathrm{Obj}\,\mathbb{C}_{\exists x:s.\,A}$, then $x:s \notin \mathrm{FV}(B)$ and there is a co-cone on $\Sigma_A(x:s)$:

$$\big(MB, \{[u:A[r/x].\, t_r \cdot u : B]\}_{r:s \in \mathrm{LTerms}(\Sigma)}\big) \ .$$

The peculiar form $t_r \cdot u$ is general: we can always substitute a generic term $t'_r$ with $t_r \cdot u$, posing $t_r \equiv (\lambda w:A[r/x].\, t'_r)$ where $w:A[r/x] \notin \mathrm{FV}(t'_r)$.

But also

$$\big(MB, \{[u:A[r/x].\, (t_x[r/x]) \cdot u : B]\}_{r:s \in \mathrm{LTerms}(\Sigma)}\big)$$

is a co-cone on $\Sigma_A(x:s)$ since $x:s \notin \mathrm{FV}(B)$ and so

$$[(u:A)[r/x].\, (t_x \cdot u)[r/x] : B[r/x]] = [u:A[r/x].\, (t_x[r/x]) \cdot u : B] \ .$$

Consider $[v:(\exists x:s.\, A).\, \mathrm{exE}(v, (\lambda x:s.\, t_x)) : B]$: it holds that

$$\begin{aligned} &[v:(\exists x:s.\, A).\, \mathrm{exE}(v, (\lambda x:s.\, t_x)) : B] \circ [u:A[r/x].\, \mathrm{exI}_x(u) :(\exists x:s.\, A)] \\ =&[u:A[r/x].\, \mathrm{exE}(\mathrm{exI}_x(u), (\lambda x:s.\, t_x)) : B] \\ =&[u:A[r/x].\, (t_x[r/x]) \cdot u : B] \ . \end{aligned}$$

Thus, $[v:(\exists x:s.\, A).\, \mathrm{exE}(v, (\lambda x:s.\, t_x)) : B]$ is an arrow in $\mathbb{C}_{\exists x:s.\,A}$.

Let $[v:(\exists x:s.\, A).\, w : B]$ be an arrow in $\mathbb{C}_{\exists x:s.\,A}$. Clearly, $x:s \notin \mathrm{FV}^*(w:B)$ as $x:s$ is not in $\mathrm{FV}^*(v:(\exists x:s.\, A))$. For each $r:s \in \mathrm{LTerms}(\Sigma)$,

$$[v:(\exists x:s.\, A).\, w : B] \circ [u:A[r/x].\, \mathrm{exI}_x(u) :(\exists x:s.\, A)]$$

by $(\exists_2)$, with $z:A$ a new variable,

$$\begin{aligned} =&[v:(\exists x:s.\, A).\, \mathrm{exE}(u, (\lambda x:s.\, (\lambda z:A.\, w[\mathrm{exI}_x(z)/v]))) : B] \\ &\circ [u:A[r/x].\, \mathrm{exI}_x(u) :(\exists x:s.\, A)] \\ =&[u:A[r/x].\, \mathrm{exE}(\mathrm{exI}_x(u), (\lambda x:s.\, (\lambda z:A.\, w[\mathrm{exI}_x(z)/v]))) : B] \end{aligned}$$

by $(\exists_0)$

$$\begin{aligned}
&= [u\colon A[r/x].\,(\lambda z\colon A.\,w[\mathrm{exI}_x(z)/v])[r/z]\cdot u\colon B] \\
&= [u\colon A[r/x].\,(\lambda z\colon A.\,w[\mathrm{exI}_x(z)/v][r/z])\cdot u\colon B] \\
&= [u\colon A[r/x].\,w[\mathrm{exI}_x(u)/v]\colon B]
\end{aligned}$$

and, since $[v\colon(\exists x\colon s.\,A).\,w\colon B]$ is an arrow of $\mathbb{C}_{\exists x\colon s.\,A}$,

$$\begin{aligned}
&= [u\colon A[r/x].\,(t_x[r/x])\cdot u\colon B] \\
&= [u\colon A[r/x].\,\mathrm{exE}(\mathrm{exI}_x(u),(\lambda x\colon s.\,t_x))\colon B]\ .
\end{aligned}$$

Thus, in particular

$$\begin{aligned}
&[u\colon A[r/x].\,\mathrm{exE}(\mathrm{exI}_x(u),(\lambda x\colon s.\,(\lambda z\colon A.\,w[\mathrm{exI}_x(z)/v])))\colon B] \\
=&[u\colon A[r/x].\,\mathrm{exE}(\mathrm{exI}_x(u),(\lambda x\colon s.\,t_x))\colon B]\ .
\end{aligned}$$

Then $[u\colon A.\,(\lambda z\colon A.\,w[\mathrm{exI}_x(z)/v])\colon A\to B] = [u\colon A.\,t_x\colon A\to B]$ by $(\exists_1)$. So,

$$\begin{aligned}
&[v\colon(\exists x\colon s.\,A).\,\mathrm{exE}(v,(\lambda x\colon s.\,t_x))\colon B] \\
=&[v\colon(\exists x\colon s.\,A).\,\mathrm{exE}(v,(\lambda x\colon s.\,(\lambda z\colon A.\,w[\mathrm{exI}_x(z)/v])))\colon B] \\
=&[v\colon(\exists x\colon s.\,A).\,w\colon B]\ ,
\end{aligned}$$

proving the uniqueness of the co-universal arrow. □

**Proposition 6.10.** *In the $\mathbb{C}_T$ category, the arrow $\delta\colon M(\exists x\colon s.\,A\times B)\to MA\times M(\exists x\colon s.\,B)$, co-universal in $\mathbb{C}_{\exists x\colon s.\,A}$, has an inverse.*

*Proof.* Define

$$\theta = [v\colon(\exists x\colon s.\,A\times B).\,\mathrm{exE}(v,(\lambda x\colon s.\,(\lambda z\colon A\times B.\,\langle \mathrm{fst}(z),\mathrm{exI}_x(\mathrm{snd}(z))\rangle)))\colon A\times(\exists x\colon s.\,B)]\ .$$

It holds that

$$\begin{aligned}
&\theta\circ[u\colon A\times B[r/x].\,\mathrm{exI}_x(u)\colon(\exists x\colon s.\,A\times B)] \\
=&[u\colon A\times B[r/x].\,\mathrm{exE}(\mathrm{exI}_x(u), \\
&\qquad(\lambda x\colon s.\,(\lambda z\colon A\times B.\,\langle \mathrm{fst}(z),\mathrm{exI}_x(\mathrm{snd}(z))\rangle)))\colon A\times(\exists x\colon s.\,B)]
\end{aligned}$$

by $(\exists_0)$

$$\begin{aligned}
&= [u\colon A\times B[r/x].\,(\lambda z\colon A\times B.\,\langle \mathrm{fst}(z),\mathrm{exI}_x(\mathrm{snd}(z))\rangle)[r/x]\cdot u\colon A\times(\exists x\colon s.\,B)] \\
&= [u\colon A\times B[r/x].\,(\lambda z\colon A\times B[r/x].\,\langle \mathrm{fst}(z),\mathrm{exI}_x(\mathrm{snd}(z))\rangle)\cdot u\colon A\times(\exists x\colon s.\,B)] \\
&= [u\colon A\times B[r/x].\,\langle \mathrm{fst}(u),\mathrm{exI}_x(\mathrm{snd}(u))\rangle\colon A\times(\exists x\colon s.\,B)] \\
&= [u\colon A.\,u\colon A]\times[u\colon B[r/x].\,\mathrm{exI}_x(u)\colon(\exists x\colon s.\,B)] \\
&= 1_{MA}\times[u\colon B[r/x].\,\mathrm{exI}_x(u)\colon(\exists x\colon s.\,B)]\ .
\end{aligned}$$

So, $\theta$ is in $\mathbb{C}_{\exists x\colon s.\,A\times B}$. Thus, by co-universality of $\delta$, $\theta=\delta$.

Define

$$\phi = [v\colon A\times(\exists x\colon s.\,B).\,\mathrm{exE}(\mathrm{snd}(v),(\lambda x\colon s.\,(\lambda z\colon B.\,\mathrm{exI}_x(\langle \mathrm{fst}(v),z\rangle))))\colon(\exists x\colon s.\,A\times B)]\ .$$

It holds that

$$\begin{aligned}
&\phi\circ\theta\circ[u\colon A\times B[r/x].\,\mathrm{exI}_x(u)\colon(\exists x\colon s.\,A\times B)] \\
=&\phi\circ(1_{MA}\times[u\colon B[r/x].\,\mathrm{exI}_x(u)\colon(\exists x\colon s.\,B)]) \\
=&\phi\circ[u\colon A\times B[r/x].\,\langle \mathrm{fst}(u),\mathrm{exI}_x(\mathrm{snd}(u))\rangle\colon A\times(\exists x\colon s.\,B)] \\
=&[u\colon A\times B[r/x].\,\mathrm{exE}(\mathrm{snd}(u),(\lambda x\colon s.\,(\lambda z\colon B.\,\mathrm{exI}_x(\langle \mathrm{fst}(u),z\rangle))))\colon(\exists x\colon s.\,A\times B)]
\end{aligned}$$

by $(\exists_0)$

$$\begin{aligned}
&=[u\colon A\times B[r/x].\,(\lambda z\colon B.\ \mathrm{exI}_x(\langle \mathrm{fst}(u), z\rangle))[r/z]\cdot \mathrm{snd}(u)\colon(\exists x\colon s.\,A\times B)]\\
&=[u\colon A\times B[r/x].\,(\lambda z\colon B[r/z].\ \mathrm{exI}_x(\langle \mathrm{fst}(u), z\rangle))\cdot \mathrm{snd}(u)\colon(\exists x\colon s.\,A\times B)]\\
&=[u\colon A\times B[r/x].\ \mathrm{exI}_x(\langle \mathrm{fst}(u), \mathrm{snd}(u)\rangle)\colon(\exists x\colon s.\,A\times B)]\\
&=[u\colon A\times B[r/x].\ \mathrm{exI}_x(u)\colon(\exists x\colon s.\,A\times B)]\quad .
\end{aligned}$$

So, $\phi$ is an arrow between co-cones.

It follows that $\phi\circ\delta = 1_{M(\exists x\colon s.\,A\times B)}$ by co-universality of $1_{M(\exists x\colon s.\,A\times B)}$ in $\mathbb{C}_{\exists x\colon s.\,A\times B}$.

On the other way,

$$\begin{aligned}
&\delta\circ\phi\\
&=[u\colon(\exists x\colon s.\,A\times B).\ \mathrm{exE}(u,(\lambda x\colon s.\,(\lambda z\colon A\times B.\\
&\qquad \langle \mathrm{fst}(z), \mathrm{exI}_x(\mathrm{snd}(z))\rangle))))\colon A\times(\exists x\colon s.\,B)]\\
&\quad\circ[v\colon A\times(\exists x\colon s.\,B).\ \mathrm{exE}(\mathrm{snd}(v),\\
&\qquad (\lambda x\colon s.\,(\lambda z\colon B.\ \mathrm{exI}_x(\langle \mathrm{fst}(v), z\rangle))))\colon(\exists x\colon s.\,A\times B)]\\
&=[v\colon A\times(\exists x\colon s.\,B).\ \mathrm{exE}(\mathrm{exE}(\mathrm{snd}(v),(\lambda x\colon s.\,(\lambda z\colon B.\ \mathrm{exI}_x(\langle \mathrm{fst}(v), z\rangle)))),\\
&\qquad (\lambda x\colon s.\,(\lambda z\colon A\times B.\,\langle \mathrm{fst}(z), \mathrm{exI}_x(\mathrm{snd}(z))\rangle)))\colon A\times(\exists x\colon s.\,B)]
\end{aligned}$$

by $(\exists_3)$

$$\begin{aligned}
&=[v\colon A\times(\exists x\colon s.\,B).\ \mathrm{exE}(\mathrm{snd}(v),(\lambda x\colon s.\,(\lambda w\colon B.\ \mathrm{exE}(\mathrm{exI}_x(\langle \mathrm{fst}(v), w\rangle),\\
&\qquad (\lambda x\colon s.\,(\lambda z\colon A\times B.\,\langle \mathrm{fst}(z), \mathrm{exI}_x(\mathrm{snd}(z))\rangle)))))))\colon A\times(\exists x\colon s.\,B)]
\end{aligned}$$

by $(\exists_0)$

$$\begin{aligned}
&=[v\colon A\times(\exists x\colon s.\,B).\ \mathrm{exE}(\mathrm{snd}(v),(\lambda x\colon s.\,(\lambda w\colon B.\\
&\qquad (\lambda z\colon A\times B.\,\langle \mathrm{fst}(z), \mathrm{exI}_x(\mathrm{snd}(z))\rangle)\cdot\langle \mathrm{fst}(v), w\rangle)))\colon A\times(\exists x\colon s.\,B)]\\
&=[v\colon A\times(\exists x\colon s.\,B).\ \mathrm{exE}(\mathrm{snd}(v),(\lambda x\colon s.\,(\lambda w\colon B.\\
&\qquad \langle \mathrm{fst}(v), \mathrm{exI}_x(w)\rangle)))\colon A\times(\exists x\colon s.\,B)]
\end{aligned}$$

by $(\exists_4)$

$$\begin{aligned}
&=[v\colon A\times(\exists x\colon s.\,B).\,\langle \mathrm{fst}(v), \mathrm{snd}(v)\rangle\colon A\times(\exists x\colon s.\,B)]\\
&=[v\colon A\times(\exists x\colon s.\,B).\ v\colon A\times(\exists x\colon s.\,B)]\\
&=1_{MA\times M(\exists x\colon s.\,B)}\quad .
\end{aligned}$$

Thus, $\phi=\delta^{-1}$. □

**Proposition 6.11.** *The $\mathbb{C}_T$ category is logically distributive.*

*Proof.* Immediate consequence of propositions 6.4 to 6.10. □

Having a logically distributive category $\mathbb{C}_T$, we can immediately identify a model for $T$ in it.

**Proposition 6.12.** *Given a $\lambda$-theory $T$ on the $\Sigma$ signature, the $\Sigma$-structure $\langle\mathbb{C}_T, M_T, M_{\mathrm{Ax}}\rangle$ on the corresponding syntactical category is defined by $M_{\mathrm{Ax}}$ which maps $f\colon A\to B\in\mathrm{Ax}$ to $[x\colon A.\,f(x)\colon B]$. This $\Sigma$-structure is a model for $T$ and, moreover, it satisfies exactly those equalities-in-context which are provable in $T$.*

*Proof.* By induction on the structure of $t$, we prove that $[\![x\colon A.\,t\colon B]\!]=[x\colon A.\,t\colon B]$:

(1) $t\in W_B$, i.e., $t\colon B$ is a variable: then, necessarily, $t\equiv x$ and $B\equiv A$, thus

$$[\![x\colon A.\,t\colon B]\!]=[\![x\colon A.\,x\colon A]\!]=1_{MA}=[x\colon A.\,x\colon A]=[x\colon A.\,t\colon B]\quad ;$$

(2) $t : B \equiv f(s)$ with $f : C \to B \in Ax$:

$$\begin{aligned}&[\![x : A.\, f(s) : B]\!]\\ =&M_{\mathrm{Ax}} f \circ [\![x : A.\, s : C]\!]\end{aligned}$$

by induction hypothesis

$$\begin{aligned}=&[x : A.\, f(x) : B] \circ [x : A.\, s : C]\\ =&[x : A.\, f(s) : B] \quad ;\end{aligned}$$

(3) $t : B \equiv \langle s, r\rangle : B_1 \times B_2$:

$$\begin{aligned}&[\![x : A.\, \langle s, r\rangle : B_1 \times B_2]\!]\\ =&\left([\![x : A.\, s : B_1]\!], [\![x : A.\, r : B_2]\!]\right)\end{aligned}$$

by induction hypothesis

$$= ([x : A.\, s : B_1], [x : A.\, r : B_2])$$

by proposition 6.4

$$= [x : A.\, \langle s, r\rangle : B_1 \times B_2] \quad ;$$

(4) $t : B \equiv \mathrm{fst}(s) : B$:

$$\begin{aligned}&[\![x : A.\, \mathrm{fst}(s) : B]\!]\\ =&\pi_1 \circ [\![x : A.\, s : B \times C]\!]\end{aligned}$$

by induction hypothesis

$$= \pi_1 \circ [x : A.\, s : B \times C]$$

by proposition 6.4

$$\begin{aligned}=&[z : B \times C.\, \mathrm{fst}(z) : B] \circ [x : A.\, s : B \times C]\\ =&[x : A.\, \mathrm{fst}(s) : B] \quad ;\end{aligned}$$

(5) $t : B \equiv \mathrm{snd}(s) : B$:

$$\begin{aligned}&[\![x : A.\, \mathrm{snd}(s) : B]\!]\\ =&\pi_2 \circ [\![x : A.\, s : C \times B]\!]\end{aligned}$$

by induction hypothesis

$$= \pi_2 \circ [x : A.\, s : C \times B]$$

by proposition 6.4

$$\begin{aligned}=&[z : C \times B.\, \mathrm{snd}(z) : B] \circ [x : A.\, s : C \times B]\\ =&[x : A.\, \mathrm{snd}(s) : B] \quad ;\end{aligned}$$

(6) $t : B \equiv (\lambda z : C.\, s) : C \to D$: $[\![x : A.\, (\lambda z : C.\, s) : C \to D]\!]$ is the exponential transpose of $[\![x : A, z : C.\, s : D]\!] = [x : A, z : C.\, s : D]$ by induction hypothesis. So,

$$[\![x : A.\, (\lambda z : C.\, s) : C \to D]\!]$$

by proposition 6.6

$$\begin{aligned}=&[y : A.\, (\lambda u : C.\, s[\mathrm{fst}(w)/x][\mathrm{snd}(w)/z])[\langle y, u\rangle / w] : C \to D]\\ =&[y : A.\, (\lambda u : C.\, s[y/x][u/z]) : C \to D]\\ =&[x : A.\, (\lambda z : C.\, s) : C \to D] \quad ;\end{aligned}$$

(7) $t : B \equiv s \cdot r : B$:

$$\begin{aligned}&[\![x : A.\, s \cdot r : B]\!]\\ &= \mathrm{ev} \circ \big([\![x : A.\, s : C \to B]\!], [\![x : A.\, r : C]\!]\big)\end{aligned}$$

by induction hypothesis

$$= \mathrm{ev} \circ ([x : A.\, s : C \to B], [x : A.\, r : C])$$

by proposition 6.4

$$= \mathrm{ev} \circ [x : A.\, \langle s, r\rangle : (C \to B) \times C]$$

by proposition 6.6

$$\begin{aligned}&= [w : (C \to B) \times C.\, \mathrm{fst}(w) \cdot \mathrm{snd}(w) : B] \circ [x : A.\, \langle s, r\rangle : (C \to B) \times C]\\ &= [x : A.\, \mathrm{fst}(\langle s, r\rangle) \cdot \mathrm{snd}(\langle s, r\rangle) : B]\\ &= [x : A.\, s \cdot r : B] \quad ;\end{aligned}$$

(8) $t : B \equiv \mathrm{inl}_D(s) : C + D$:

$$\begin{aligned}&[\![x : A.\, \mathrm{inl}_D(s) : C + D]\!]\\ &= \iota_1 \circ [\![x : A.\, s : C]\!]\end{aligned}$$

by induction hypothesis

$$= \iota_1 \circ [x : A.\, s : C]$$

by proposition 6.5

$$\begin{aligned}&= [z : C.\, \mathrm{inl}_D(z) : C + D] \circ [x : A.\, s : C]\\ &= [x : A.\, \mathrm{inl}_D(s) : C + D] \quad ;\end{aligned}$$

(9) $t : B \equiv \mathrm{inr}_D(s) : D + C$:

$$\begin{aligned}&[\![x : A.\, \mathrm{inr}_D(s) : D + C]\!]\\ &= \iota_2 \circ [\![x : A.\, s : C]\!]\end{aligned}$$

by induction hypothesis

$$= \iota_2 \circ [x : A.\, s : C]$$

by proposition 6.5

$$\begin{aligned}&= [z : C.\, \mathrm{inr}_D(z) : D + C] \circ [x : A.\, s : C]\\ &= [x : A.\, \mathrm{inr}_D(s) : D + C] \quad ;\end{aligned}$$

(10) $t : B \equiv \mathrm{when}(s, u, v) : B$:

$$\begin{aligned}&[\![x : A.\, \mathrm{when}(s, u, v) : B]\!]\\ &= \big[\mathrm{ev} \circ \big([\![x : A.\, u : B_1 \to B]\!] \times 1_{MB_1}\big),\\ &\quad \mathrm{ev} \circ \big([\![x : A.\, v : B_2 \to B]\!] \times 1_{MB_2}\big)\big]\\ &\quad \circ \Delta^{-1} \circ \big(1_{MA}, [\![x : A.\, s : B_1 + B_2]\!]\big)\end{aligned}$$

by induction hypothesis

$$\begin{aligned}&= \big[\mathrm{ev} \circ \big([x : A.\, u : B_1 \to B] \times [y : B_1.\, y : B_1]\big),\\ &\quad \mathrm{ev} \circ \big([x : A.\, v : B_2 \to B] \times [y : B_2.\, y : B_2]\big)\big]\\ &\quad \circ \Delta^{-1} \circ ([x : A.\, x : A], [x : A.\, s : B_1 + B_2])\end{aligned}$$

by propositions 6.4, 6.5, 6.6, and 6.7

$$
\begin{aligned}
&= \bigl[\mathrm{ev}\circ\bigl([x\colon A.\, u\colon B_1 \to B] \times [y\colon B_1.\, y\colon B_1]\bigr),\\
&\quad \mathrm{ev}\circ\bigl([x\colon A.\, v\colon B_2 \to B] \times [y\colon B_2.\, y\colon B_2]\bigr)\bigr]\\
&\quad \circ [w\colon A \times (B_1 + B_2).\, \mathrm{when}(\mathrm{snd}(w), (\lambda x\colon B_1.\, \mathrm{inl}_{A\times B_2}(\langle \mathrm{fst}(w), x\rangle)),\\
&\qquad (\lambda y\colon B_2.\, \mathrm{inr}_{A\times B_1}(\langle \mathrm{fst}(w), y\rangle)))\colon (A\times B_1) + (A \times B_2)]\\
&\quad \circ [x\colon A.\, \langle x, s\rangle \colon A \times (B_1 + B_2)]\\
&= [\mathrm{ev}\circ [w\colon A\times B_1.\, \langle u[\mathrm{fst}(w)/x], \mathrm{snd}(w)\rangle \colon (B_1 \to B)\times B_1],\\
&\quad \mathrm{ev}\circ [w\colon A\times B_2.\, \langle v[\mathrm{fst}(w)/x], \mathrm{snd}(w)\rangle \colon (B_2 \to B)\times B_2]]\\
&\quad \circ [z\colon A.\, \mathrm{when}(s[z/x], (\lambda x\colon B_1.\, \mathrm{inl}_{A\times B_2}(\langle z, x\rangle)),\\
&\qquad (\lambda y\colon B_2.\, \mathrm{inr}_{A\times B_1}(\langle z, y\rangle)))\colon (A\times B_1) + (A\times B_2)]\\
&= [[w\colon A\times B_1.\, u[\mathrm{fst}(w)/x]\cdot \mathrm{snd}(w)\colon B],\\
&\quad [w\colon A\times B_2.\, v[\mathrm{fst}(w)/x]\cdot \mathrm{snd}(w)\colon B]]\\
&\quad \circ [z\colon A.\, \mathrm{when}(s[z/x], (\lambda x\colon B_1.\, \mathrm{inl}_{A\times B_2}(\langle z, x\rangle)),\\
&\qquad (\lambda y\colon B_2.\, \mathrm{inr}_{A\times B_1}(\langle z, y\rangle)))\colon (A\times B_1) + (A\times B_2)]\\
&= [z\colon (A\times B_1) + (A\times B_2).\, \mathrm{when}(z,\\
&\quad (\lambda w\colon A\times B_1.\, u[\mathrm{fst}(w)/x]\cdot \mathrm{snd}(w)),\\
&\quad (\lambda w\colon A\times B_2.\, v[\mathrm{fst}(w)/x]\cdot \mathrm{snd}(w)))\colon B]\\
&\quad \circ [z\colon A.\, \mathrm{when}(s[z/x], (\lambda x\colon B_1.\, \mathrm{inl}_{A\times B_2}(\langle z, x\rangle)),\\
&\qquad (\lambda y\colon B_2.\, \mathrm{inr}_{A\times B_1}(\langle z, y\rangle)))\colon (A\times B_1) + (A\times B_2)]\\
&= [z\colon A.\, \mathrm{when}(\mathrm{when}(s[z/x], (\lambda x\colon B_1.\, \mathrm{inl}_{A\times B_2}(\langle z, x\rangle)),\\
&\qquad (\lambda y\colon B_2.\, \mathrm{inr}_{A\times B_1}(\langle z, y\rangle))),\\
&\quad (\lambda w\colon A\times B_1.\, u[\mathrm{fst}(w)/x]\cdot \mathrm{snd}(w)),\\
&\quad (\lambda w\colon A\times B_2.\, v[\mathrm{fst}(w)/x]\cdot \mathrm{snd}(w)))\colon B]
\end{aligned}
$$

by $(+_2)$

$$
\begin{aligned}
&= [z\colon A.\, \mathrm{when}(s[z/x],\\
&\quad (\lambda y\colon B_1.\, \mathrm{when}(\mathrm{inl}_{A\times B_2}(\langle z, y\rangle),\\
&\qquad (\lambda w\colon A\times B_1.\, u[\mathrm{fst}(w)/x]\cdot \mathrm{snd}(w)),\\
&\qquad (\lambda w\colon A\times B_2.\, v[\mathrm{fst}(w)/x]\cdot \mathrm{snd}(w))))\\
&\quad (\lambda y\colon B_2.\, \mathrm{when}(\mathrm{inr}_{A\times B_1}(\langle z, y\rangle),\\
&\qquad (\lambda w\colon A\times B_1.\, u[\mathrm{fst}(w)/x]\cdot \mathrm{snd}(w)),\\
&\qquad (\lambda w\colon A\times B_2.\, v[\mathrm{fst}(w)/x]\cdot \mathrm{snd}(w)))))\colon B]\\
&= [z\colon A.\, \mathrm{when}(s[z/x],\\
&\quad (\lambda y\colon B_1.\, (\lambda w\colon A\times B_1.\, u[\mathrm{fst}(w)/x]\cdot \mathrm{snd}(w))\cdot \langle z, y\rangle),\\
&\quad (\lambda y\colon B_2.\, (\lambda w\colon A\times B_2.\, v[\mathrm{fst}(w)/x]\cdot \mathrm{snd}(w))\cdot \langle z, y\rangle))\colon B]\\
&= [z\colon A.\, \mathrm{when}(s[z/x], (\lambda y\colon B_1.\, u[z/x]\cdot y), (\lambda y\colon B_2.\, v[z/x]\cdot y))\colon B]\\
&= [z\colon A.\, \mathrm{when}(s[z/x], u[z/x], v[z/x])\colon B]\\
&= [x\colon A.\, \mathrm{when}(s, u, v)\colon B]\ \ ;
\end{aligned}
$$

(11) $t\colon B \equiv *\colon 1$: $[\![x\colon A.\, *\colon 1]\!] = !\colon\ A \to 1 = [x\colon A.\, *\colon 1]$ by proposition 6.4;
(12) $t\colon B \equiv \mathrm{F}_C\colon 0 \to C$: $[\![x\colon A.\, \mathrm{F}_C\colon 0 \to C]\!]$ is the exponential transpose of

$$(!\colon\ 0 \to C)\circ \pi_2$$

by proposition 6.4 and 6.5

$$\begin{aligned}&=[x\!:\!0.\ \mathrm{F}_C\cdot x\!:\!C]\circ[w\!:\!A\times 0.\ \mathrm{snd}(w)\!:\!0]\\&=[w\!:\!A\times 0.\ \mathrm{F}_C\cdot\mathrm{snd}(w)\!:\!C]\ .\end{aligned}$$

By proposition 6.6, the exponential transpose of $[w\!:\!A\times 0.\ \mathrm{F}_C\cdot\mathrm{snd}(w)\!:\!C]$ is

$$\begin{aligned}&[x\!:\!A.\,(\lambda z\!:\!0.\,(\mathrm{F}_C\cdot\mathrm{snd}(w))[\langle x,z\rangle/w])\!:\!0\to C]\\=&[x\!:\!A.\,(\lambda z\!:\!0.\ \mathrm{F}_C\cdot\mathrm{snd}(\langle x,z\rangle))\!:\!0\to C]\\=&[x\!:\!A.\,(\lambda z\!:\!0.\ \mathrm{F}_C\cdot z)\!:\!0\to C]\end{aligned}$$

by $(\to_1)$

$$=[x\!:\!A.\ \mathrm{F}_C\!:\!0\to C]\ ;$$

(13) $t\!:\!B\equiv\mathrm{allI}(\lambda z\!:\!s.\,u)\!:\!(\forall z\!:\!s.\,C)$: $[\![x\!:\!A.\ \mathrm{allI}(\lambda z\!:\!s.\,u)\!:\!(\forall z\!:\!s.\,C)]\!]=\beta\circ\alpha$, where $\alpha$ is the projection from $MA$ to $MA'$, the product corresponding to $x'=\mathrm{FV}(u\!:\!C)$, and $\beta\colon MA'\to M(\forall z\!:\!s.\,C)$ is universal in $\mathbb{C}_{\forall z:s.\,C}$.

So, either $\alpha=1_{MA}$ or $\alpha\colon MA\to 1$ as $\mathrm{FV}(u\!:\!C)\subseteq\{x\!:\!A\}$. In the former case, $\alpha=[x\!:\!A.,\ x\!:\!A]$, in the latter $\alpha=[x\!:\!A.\ *\!:\!1]$, see proposition 6.4.

Also, $\beta$ is such that $[v\!:\!(\forall z\!:\!s.\,C).\ \mathrm{allE}(v,r)\!:\!C[r/z]]\circ\beta=[\![x\!:\!A'.\,u[r/z]\!:\!C[r/z]$ for every $r\!:\!s\in\mathrm{LTerms}(\Sigma)$, since $\beta$ is an arrow between cones. But,

$$\begin{aligned}&[v\!:\!(\forall z\!:\!s.\,C).\ \mathrm{allE}(v,r)\!:\!C[r/z]]\circ[x\!:\!A'.\ \mathrm{allI}(\lambda z\!:\!s.\,u)\!:\!(\forall z\!:\!s.\,C)]\\=&[x\!:\!A'.\ \mathrm{allE}(\mathrm{allI}(\lambda z\!:\!s.\,u),r)\!:\!C[r/z]]\end{aligned}$$

by $(\forall_0)$

$$=[x\!:\!A'.\,u[r/z]\!:\!C[r/z]]$$

by induction hypothesis

$$=[\![x\!:\!A'.\,u[r/z]\!:\!C[r/z]]\!]\ .$$

So, by universality of $\beta$, $\beta=[x\!:\!A'.\ \mathrm{allI}(\lambda z\!:\!s.\,u)\!:\!(\forall z\!:\!s.\,C)]$.

It is immediate to verify that $\beta\circ\alpha=[x\!:\!A.\ \mathrm{allI}(\lambda z\!:\!s.\,u)\!:\!(\forall z\!:\!s.\,C)]$ for both possible choices of $\alpha$;

(14) $t\!:\!B\equiv\mathrm{allE}(u,r)\!:\!C[r/z]$:

$$\begin{aligned}&[\![x\!:\!A.\ \mathrm{allE}(u,r)\!:\!C[r/z]]\!]\\=&p_r\circ[\![x\!:\!A.\,u\!:\!(\forall z\!:\!s.\,C)]\!]\end{aligned}$$

where $p_r$ is the $r$-th projector of the canonical cone with vertex $M(\forall z\!:\!s.\,C)$, so, by proposition 6.8,

$$=[v\!:\!(\forall z\!:\!s.\,C).\ \mathrm{allE}(v,t)\!:\!C[r/z]]\circ[\![x\!:\!A.\,u\!:\!(\forall z\!:\!s.\,C)]\!]$$

by induction hypothesis

$$\begin{aligned}&=[v\!:\!(\forall z\!:\!s.\,C).\ \mathrm{allE}(v,t)\!:\!C[r/z]]\circ[x\!:\!A.\,u\!:\!(\forall z\!:\!s.\,C)]\\&=[x\!:\!A.\ \mathrm{allE}(u,t)\!:\!C[r/z]]\ ;\end{aligned}$$

(15) $t\!:\!B\equiv\mathrm{exI}_z(u)\!:\!(\exists z\!:\!s.\,C)$:

$$\begin{aligned}&[\![x\!:\!A.\ \mathrm{exI}_z(u)\!:\!(\exists z\!:\!s.\,C)]\!]\\=&j_r\circ[\![x\!:\!A.\,u\!:\!C[r/z]]\!]\end{aligned}$$

where $j_r$ is the $r$-th injection of the canonical co-cone with vertex $M(\exists z\colon s.\, C)$, so, by proposition 6.9,

$$=[v\colon C[r/z].\, \mathrm{exI}_z(v)\colon(\exists z\colon s.\, C)]\circ [\![x\colon A.\, u\colon C[r/z]]\!]$$

by induction hypothesis

$$\begin{aligned}&=[v\colon C[r/z].\, \mathrm{exI}_z(v)\colon(\exists z\colon s.\, C)]\circ [x\colon A.\, u\colon C[r/z]]\\ &=[x\colon A.\, \mathrm{exI}_z(u)\colon C[r/z]]\ ;\end{aligned}$$

(16) $t\colon B\equiv \mathrm{exE}(u,(\lambda z\colon s.\, v))\colon B$:

$$\begin{aligned}&[\![x\colon A.\, \mathrm{exE}(u,(\lambda z\colon s.\, v))\colon B]\!]\\ &=\gamma\circ\beta^{-1}\circ\big(\alpha,[\![x\colon A.\, u\colon(\exists z\colon s.\, C)]\!]\big)\end{aligned}$$

by induction hypothesis

$$=\gamma\circ\beta^{-1}\circ(\alpha,[x\colon A.\, u\colon(\exists z\colon s.\, C)])\ ,$$

where $\alpha$ is the projection to $MA'$, the product corresponding to

$$x'=\mathrm{FV}(u\colon(\exists z\colon s.\, C))\cup\mathrm{FV}(v\colon C\to B)\subseteq\{x\colon A\}\ ,$$

$\beta\colon M(\exists z\colon s.\, A'\times C)\to MA'\times M(\exists z\colon s.\, C)$ and $\gamma\colon M(\exists z\colon s.\, A'\times C)\to MB$ are co-universal in $\mathbb{C}_{\exists z\colon s.\, A'\times C}$.

So, either $\alpha=1_{MA}$ or $\alpha\colon MA\to 1$ as $\mathrm{FV}(u\colon C)\subseteq\{x\colon A\}$. In the former case, $\alpha=[x\colon A.,\ x\colon A]$, in the latter $\alpha=[x\colon A.\ *\colon 1]$, see proposition 6.4. Thus,

$$\begin{aligned}&\beta^{-1}\circ(\alpha,[x\colon A.\, u\colon(\exists z\colon s.\, C)])\\ &=\beta^{-1}\circ(1_{MA'}\times[x\colon A.\, u\colon(\exists z\colon s.\, C)])\circ(\alpha,1_{MA})\end{aligned}$$

by proposition 6.10,

$$\begin{aligned}&=[v\colon A'\times(\exists z\colon s.\, B).\, \mathrm{exE}(\mathrm{snd}(v),\\ &\quad(\lambda z\colon s.\,(\lambda y\colon B.\, \mathrm{exI}_z(\langle\mathrm{fst}(v),y\rangle))))\colon(\exists z\colon s.\, A'\times B)]\\ &\quad\circ(1_{MA'}\times[x\colon A.\, u\colon(\exists z\colon s.\, C)])\circ(\alpha,1_{MA})\\ &=[v\colon A'\times(\exists z\colon s.\, B).\, \mathrm{exE}(\mathrm{snd}(v),\\ &\quad(\lambda z\colon s.\,(\lambda y\colon B.\, \mathrm{exI}_z(\langle\mathrm{fst}(v),y\rangle))))\colon(\exists z\colon s.\, A'\times B)]\\ &\quad\circ[w\colon A'\times A.\, \langle\mathrm{fst}(w),u\rangle\colon A'\times(\exists z\colon s.\, C)]\circ(\alpha,1_{MA})\\ &=[w\colon A'\times A.\, \mathrm{exE}(u,\\ &\quad(\lambda z\colon s.\,(\lambda y\colon B.\, \mathrm{exI}_z(\langle\mathrm{fst}(w),y\rangle))))\colon(\exists z\colon s.\, A'\times B)]\circ(\alpha,1_{MA})\ .\end{aligned}$$

Consider the co-cone on $\Sigma_{A'\times C}(z\colon s)$

$$\Big(MB,\big\{\mathrm{ev}\circ\big([\![x\colon A'.\, v[r/z]\colon C[r/z]\to B]\!]\times 1_{M(C[r/z])}\big)\big\}_{r\colon s\in\mathrm{LTerms}(\Sigma)}\Big)\ .$$

It holds that

$$\mathrm{ev}\circ\big([\![x\colon A'.\, v[r/z]\colon C[r/z]\to B]\!]\times 1_{M(C[r/z])}\big)$$

by induction hypothesis

$$=\mathrm{ev}\circ\big([x\colon A'.\, v[r/z]\colon C[r/z]\to B]\times[w\colon C[r/z].\, w\colon C[r/z]]\big)$$

by proposition 6.6

$$\begin{aligned}&=[w\colon(C[r/z]\to B)\times C[r/z].\, \mathrm{fst}(w)\cdot\mathrm{snd}(w)\colon B]\\ &\quad\circ[w\colon A'\times C[r/z].\, \langle v[r/z][\mathrm{fst}(w)/x],\mathrm{snd}(w)\rangle\colon(C[r/z]\to B)\times C[r/z]]\\ &=[w\colon A'\times C[r/z].\, v[r/z][\mathrm{fst}(w)/x]\cdot\mathrm{snd}(w)\colon B]\ .\end{aligned}$$

So, from proposition 6.9, it follows that

$$\gamma = [w\colon(\exists z\colon s.\, A' \times C).\, \mathrm{exE}(w, (\lambda z\colon s.\, (\lambda y\colon A' \times C.\, v[\mathrm{fst}(y)/x] \cdot \mathrm{snd}(y)))) \colon B] \ .$$

Summarising, we can calculate

$$\begin{aligned}
&[\![x\colon A.\, \mathrm{exE}(u, (\lambda z\colon s.\, v)) \colon B]\!]\\
=&\gamma \circ \beta^{-1} \circ (\alpha, [x\colon A.\, u\colon(\exists z\colon s.\, C)])\\
=&\gamma \circ [w\colon A' \times A.\, \mathrm{exE}(u,\\
&(\lambda z\colon s.\, (\lambda y\colon B.\, \mathrm{exI}_z(\langle \mathrm{fst}(w), y\rangle)))) \colon (\exists z\colon s.\, A' \times B)] \circ (\alpha, 1_{MA})\\
=&[w\colon(\exists z\colon s.\, A' \times C).\, \mathrm{exE}(w,\\
&(\lambda z\colon s.\, (\lambda y\colon A' \times C.\, v[\mathrm{fst}(y)/x] \cdot \mathrm{snd}(y)))) \colon B]\\
&\circ [w\colon A' \times A.\, \mathrm{exE}(u,\\
&(\lambda z\colon s.\, (\lambda y\colon B.\, \mathrm{exI}_z(\langle \mathrm{fst}(w), y\rangle)))) \colon (\exists z\colon s.\, A' \times B)] \circ (\alpha, 1_{MA})\\
=&[w\colon A' \times A.\, \mathrm{exE}(\mathrm{exE}(u, (\lambda z\colon s.\, (\lambda y\colon B.\, \mathrm{exI}_z(\langle \mathrm{fst}(w), y\rangle)))),\\
&(\lambda z\colon s.\, (\lambda y\colon A' \times C.\, v[\mathrm{fst}(y)/x] \cdot \mathrm{snd}(y)))) \colon B] \circ (\alpha, 1_{MA})
\end{aligned}$$

by $(\exists_3)$

$$\begin{aligned}
=&[w\colon A' \times A.\, \mathrm{exE}(u, (\lambda z\colon s.\, (\lambda y\colon B.\, \mathrm{exE}(\mathrm{exI}_z(\langle \mathrm{fst}(w), y\rangle),\\
&(\lambda z\colon s.\, (\lambda y\colon A' \times C.\, v[\mathrm{fst}(y)/x] \cdot \mathrm{snd}(y))))))) \colon B] \circ (\alpha, 1_{MA})
\end{aligned}$$

by $(\exists_0)$

$$\begin{aligned}
=&[w\colon A' \times A.\, \mathrm{exE}(u, (\lambda z\colon s.\, (\lambda y\colon B.\, (\lambda y'\colon A' \times C.\\
&v[\mathrm{fst}(y')/x] \cdot \mathrm{snd}(y')) \cdot \langle \mathrm{fst}(w), y\rangle))) \colon B] \circ (\alpha, 1_{MA})\\
=&[w\colon A' \times A.\, \mathrm{exE}(u, (\lambda z\colon s.\, (\lambda y\colon B.\, v[\mathrm{fst}(w)/x] \cdot y))) \colon B] \circ (\alpha, 1_{MA})\\
=&[w\colon A' \times A.\, \mathrm{exE}(u, (\lambda z\colon s.\, v[\mathrm{fst}(w)/x])) \colon B] \circ (\alpha, 1_{MA})\\
=&[x\colon A.\, \mathrm{exE}(u, (\lambda z\colon s.\, v)) \colon B] \ .
\end{aligned}$$

So, in the given $\Sigma$-structure, $[\![x\colon A.\, t\colon B]\!] = [\![x\colon A.\, s\colon B]\!]$ if and only if $[x\colon A.\, t\colon B] = [x\colon A.\, s\colon B]$ and, by proposition 6.2, this happens if and only if $T \vdash s =_B t$. □

The identified model $M_T$ in $\mathbb{C}_T$ is canonical because it allows to uniformly generate all other models in any other logically distributive category.

**Proposition 6.13.** *For every logically distributive category $\mathbb{C}$, there is a bijection between equivalence classes, modulo natural equivalences, of structure-preserving functors $\mathbb{C}_T \to \mathbb{C}$ and equivalence classes, modulo isomorphisms, of models of $T$ in $\mathbb{C}$, induced by the map $F \mapsto F(M_T)$.*

*Proof.* A functor preserves logical distributivity when it preserves finite products, finite co-products and exponentiation also in all the subcategories $\mathbb{C}_{\forall x\colon s.\, A}$ and $\mathbb{C}_{\exists x\colon s.\, A}$. As any functor preserves isomorphisms, it automatically preserves the existence of inverses as required in the definition of logically distributive category.

So, given a model $\langle \mathbb{N}, N_T, N_{\mathrm{Ax}}\rangle$ for $T$, the functor $F_{\mathbb{N}}\colon \mathbb{C}_T \to \mathbb{N}$ mapping $A \mapsto NA$ for each $A \in \mathrm{Obj}\,\mathbb{C}_T$ and $[x\colon A.\, t\colon B] \mapsto [\![x\colon A.\, t\colon B]\!]$ is clearly structure-preserving. Moreover, $F_N(M_T) = N_T$ and $F_N(M_{\mathrm{Ax}}) = N_{\mathrm{Ax}}$.

On the other side, every structure-preserving functor $F\colon \mathbb{N} \to \mathbb{C}_T$ must also preserve the interpretation of terms-in-context, so it is naturally isomorphic to $F_N$, defined on the $\Sigma$-structure $\langle \mathbb{N}, F(M_T), F(M_{\mathrm{Ax}})\rangle$. □

*Note* 6.1. Proposition 6.13 is an extension of proposition 4.2.5 in [Joh02b, page 955].

It can be regarded as a *classifying* result as it says, essentially, that each model $M$ of $T$ in a category $\mathbb{N}$ can be uniquely (up to natural isomorphisms) identified with the structure-preserving functor $F_M$.

Now, the completeness result follows immediately.

**Theorem 6.14** (Completeness)**.** *If $\vec{x}.\, s =_A t$ is an equality-in-context valid in every model for $T$ in each logically distributive category, then $T \vdash \vec{x}.\, s =_A t$.*

*Proof.* It is safe to assume, thanks to the notion of syntactical equivalence, that $\vec{x} \equiv x\colon B$.

Now, if $x\colon B.\, s =_A t$ is valid in each model of every logically distributive category, it is valid in particular in $M_T$ of $\mathbb{C}_T$, the canonical model. Thus,

$$[\![x\colon B.\, s\colon A]\!]_{M_T} = [\![x\colon B.\, t\colon A]\!]_{M_T} \ ,$$

that is, $[x\colon B.\, s\colon A] = [x\colon B.\, t\colon A]$. So, by proposition 6.12, $T \vdash \vec{x}.\, s =_A t$. □

**Corollary 6.15.** *An equality-in-context is valid in every model of every logically distributive category if and only if it is valid in $M_T$ of $\mathbb{C}_T$, which, in turn, happens if and only if $T$ proves the equality.*

*Proof.* Immediate consequence of theorem 6.14, of theorem 5.4 and of the classifying property expressed in proposition 6.13. □

6.1. **Soundness and completeness in logic.** From the soundness and completeness of the $\lambda$-theory, it follows that the corresponding logical theory via the Curry-Howard isomorphism is sound and complete as well, with an appropriate notion of truth inside logically distributive categories.

**Definition 6.3** (Valid type)**.** A $\lambda$-type $A$ is *valid in the model* $\mathcal{N} = \langle \mathbb{N}, N, N_{\mathrm{Ax}} \rangle$ when there exists $1 \to NA$ in $\mathbb{N}$.

A $\lambda$-type $A$ is a *logical consequence in the model* $\mathcal{N}$ of the $\lambda$-types $B_1, \ldots, B_n$ when there exists $N(B_1 \times \cdots \times B_n) \to NA$ in $\mathbb{N}$.

A $\lambda$-type $A$ is a *logical consequence* of $B_1, \ldots, B_n$ when it is a logical consequence of $B_1, \ldots, B_n$ in every model in every logically distributive category.

**Proposition 6.16.** *A $\lambda$-type $A$ is a logical consequence of $B_1, \ldots, B_n$ if and only if there exists a term-in-context $x\colon B_1 \times \cdots \times B_n.\, t\colon A$.*

*Proof.* Suppose $x\colon B_1 \times \cdots B_n.\, t\colon A$ is a term-in-context.

Then, in every model $\mathcal{N} = \langle \mathbb{N}, N, N_{\mathrm{Ax}} \rangle$,

$$[\![x\colon B_1 \times \cdots B_n.\, t\colon A]\!]_{\mathcal{N}}\colon N(B_1 \times \cdots \times B_n) \to NA$$

is an arrow of $\mathbb{N}$, so $A$ is a logical consequence of $B_1, \ldots, B_n$ in $\mathcal{N}$.

Vice versa, if $A$ is a logical consequence of $B_1, \ldots, B_n$, then it is so in the canonical model. Thus, there is an arrow

$$[x\colon B_1 \times \cdots \times B_n.\, t\colon A]\colon :\ M_T(B_1 \times \cdots \times B_n) \to M_T A \ .$$

But, for each logically distributive category $\mathbb{N}$, every model $\mathcal{N}$ in $\mathbb{N}$ is representable, up to isomorphisms, as a structure-preserving functor $F\colon \mathbb{C}_T \to \mathbb{N}$ by proposition 6.13. So,

$$F([x\colon B_1 \times \cdots \times B_n.\, t\colon A])\colon N(B_1 \times \cdots \times B_n) \to NA \ ,$$

and, moreover, $F([x\colon B_1 \times \cdots \times B_n.\, t\colon A])$ is isomorphic to $[\![x\colon B_1 \times \cdots \times B_n.\, t\colon A]\!]_{\mathcal{N}}$. □

**Corollary 6.17.** *A $\lambda$-type $A$ is a logical consequence of $B_1, \ldots, B_n$ if and only if there is a proof of $A$ from the hypotheses $B_1, \ldots, B_n$, when $\lambda$-types are interpreted as logical formulae and $\lambda$-terms as logical proofs, according to the Curry-Howard isomorphism.*

*Proof.* Obvious from proposition 6.16, as $t$ is a logical proof of $A$. □

*Note* 6.2. Corollary 6.17 is a soundness and completeness proof for the first-order intuitionistic logic with respect to the semantics given by the notion of logical consequence on the models in logically distributive categories.

It is worth noticing that the soundness and completeness theorem holds for any logical theory $T$. In particular, a contradictory theory $T$ has a model, differently from other semantics. But, in this case, the model identifies $\top$ and $\bot$, making it trivial: this fact is evident, since it happens in the canonical model in $\mathbb{C}_T$ thus, by proposition 6.13, in any other model.

We conclude this section showing that the $(\forall_1)$ rule of definition 3.7, as remarked at the end of section 5, can be simplified to

$$\vec{x}.\, \mathrm{allE}(u,y) =_B \mathrm{allE}(v,y) \vdash \vec{x}.\, u =_{(\forall z:s.\,B)} v \quad ,$$

where $y:s \not\in \mathrm{FV}(u:(\forall z:s.\,B)) \cup \mathrm{FV}(v:(\forall z:s.\,B))$.

Equivalence in one direction is evident: if this rule is assumed, and we know that $\vec{x}.\, \mathrm{allE}(u,r) =_B \mathrm{allE}(v,r)$ for all $r:s \in \mathrm{LTerms}(\Sigma)$, then, for $r \equiv y$, the usual $(\forall_1)$ rule follows. Conversely, if we know $\vec{x}.\, \mathrm{allE}(u,y) =_B \mathrm{allE}(v,y)$, then, by soundness, it holds that $[\![\vec{x}.\, \mathrm{allE}(u,y):B]\!] = [\![\vec{x}.\, \mathrm{allE}(v,y):B]\!]$. By the semantics of allE, we can safely assume that $y:s$ does not appear in $\vec{x}$—if this is the case, it must be absorbed by the $\alpha$ arrow, as this piece of information is not needed to correctly derive $u$ and $v$—so, substituting $r:s$ for $y:s$ yields $[\![\vec{x}.\, \mathrm{allE}(u,y):B]\!] = [\![\vec{x}.\, \mathrm{allE}(v,y):B]\!]$, i.e., $\vec{x}.\, \mathrm{allE}(u,r) =_B \mathrm{allE}(v,r)$ for any $r:s \in \mathrm{LTerms}(\Sigma)$. Thus, applying $(\forall_1)$, the simplified rule is derived.

## 7. Heyting categories

This section describes the semantics of the logical systems of section 2 in the Heyting categories, which are also defined in the following. The results here are a summary of well-known facts: appropriate references are given, even if the we mainly follow the presentation in [Joh02a, Joh02b].

First, we remind the definition of Heyting algebra.

**Definition 7.1** (Poset)**.** Let $P$ be a set and $\leq \subseteq P \times P$ a relation. $P$ is a *poset* if and only if

(1) for all $a \in P$, $a \leq a$;
(2) for all $a, b \in P$, $a \leq b$ and $b \leq a$ implies $a = b$;
(3) for all $a, b, c \in P$, $a \leq b$ and $b \leq c$ implies $a \leq c$.

If the second clause is omitted, we say that $P$ is a pre-order.

As usual, a poset $P$ can be regarded as a category $\mathbb{P}$ where $\mathrm{Obj}\,\mathbb{P} = P$ and $\mathrm{Hom}(a, b)$ contains one element exactly when $a \leq b$ in P.

**Definition 7.2** (Lattice)**.** A poset $\langle P; \leq\rangle$ is a *lattice* if, for every $a, b \in P$,

(1) there is the greatest lower bound of $a$ and $b$, called the *meet* of $a$ and $b$, and denoted as $a \wedge b$; $c = a \wedge b$ if and only if $c \leq a$ and $c \leq b$, and, for all $d \in P$ such that $d \leq a$ and $d \leq b$, $d \leq c$;
(2) there is the least upper bound of $a$ and $b$, called the *join* of $a$ and $b$, and denoted as $a \vee b$; $c = a \vee b$ if and only if $a \leq c$ and $b \leq c$, and, for all $d \in P$ such that $a \leq d$ and $b \leq d$, $c \leq d$.

It is routine to show that meets and joins are unique.

**Definition 7.3** (Bounded lattice)**.** A lattice $\langle P; \leq\rangle$ is *bounded* when

(1) there is $\top \in P$ such that, for all $a \in P$, $a \leq \top$;
(2) there is $\bot \in P$ such that, for all $a \in P$, $\bot \leq a$.

It is routine to show that the top and bottom elements are unique. Moreover, it is immediate to show that a bounded lattice has all the finite meets and joins.

**Definition 7.4** (Heyting algebra)**.** A bounded lattice $\langle P; \leq, \top, \bot\rangle$ is a *Heyting algebra* if and only if, for all $a, b \in P$, there is $x \in P$ such that $a \wedge x \leq b$ and, for all $y \in P$ such that $a \wedge y \leq b$, $y \leq x$. This element $x$ is called the *relative pseudo-complement* of $a$ with respect to $b$, and is denoted as $a \Rightarrow b$. Sometimes, the corresponding operation is called the *lattice implication*.

**Proposition 7.1.** *A Heyting algebra $\mathbb{H}$ is a bounded lattice (expressed as a category) in which, for each $b \in \mathrm{Obj}\,\mathbb{H}$, the functor*

$$\begin{array}{rcl} - \wedge b \colon \mathbb{H} & \to & \mathbb{H} \\ a & \mapsto & a \wedge b \end{array}$$

*has a right adjoint*

$$\begin{array}{rcl} b \Rightarrow - \colon \mathbb{H} & \to & \mathbb{H} \\ c & \mapsto & b \Rightarrow c \end{array} .$$

*Proof.* Standard. □

**Proposition 7.2.** *A Heyting algebra is a distributive lattice.*

*Proof.* See [Bor94b, p. 6]. □

To define Heyting categories, we need a number of auxiliary notions which have an interest by themselves.

**Definition 7.5** (Cartesian category)**.** A category is *Cartesian* when its has all the finite limits.

**Definition 7.6** (Having images)**.** A Cartesian category $\mathbb{C}$ *has images* if, for each $f \in \text{Hom}_{\mathbb{C}}(a,b)$, there exists $\text{im}(f)\colon c \rightarrowtail b$ such that $f = \text{im}(f) \circ r$ for some $r \in \text{Hom}_{\mathbb{C}}(a,c)$ and, for each $g\colon d \rightarrowtail b$ such that $f = g \circ h$ for an appropriate $h \in \text{Hom}_{\mathbb{C}}(a,d)$, there is $k \in \text{Hom}_{\mathbb{C}}(c,d)$ such that $\text{im}(f) = g \circ k$.

**Definition 7.7** (Ordering sub-objects)**.** Given $m_1\colon a_1 \rightarrowtail c$ and $m_2\colon a_2 \rightarrowtail c$, $m_1 \leq m_2$ if and only if there is $r\colon a_1 \to a_2$, necessarily unique, such that $m_1 = m_2 \circ r$.

**Definition 7.8** (Stability under pullback of images)**.** In a Cartesian category $\mathbb{C}$ with images, the images are *stable under pullback* whenever the pullback $f'\colon P \to C$ of $f\colon A \to B$ along $g\colon C \to B$ has as image $\text{im}(f')\colon P' \rightarrowtail C$ which is the pullback of $\text{im}(f)\colon D \rightarrowtail B$ along $g$. Briefly, the image of a pullback is the pullback of the image.

**Definition 7.9** (Regular category)**.** A Cartesian category is *regular* when it has images and they are stable under pullbacks.

**Definition 7.10** (Sub-object category)**.** Given a regular category $\mathbb{C}$, $\text{Sub}_{\mathbb{C}}(a)$, for any $a \in \text{Obj}\,\mathbb{C}$, denotes the full subcategory of $\mathbb{C}/a$ whose objects are the sub-objects of $a$. We remind that $\mathbb{C}/a$ is the category whose objects are the arrows having $a$ as co-domain and whose arrows are the obvious arrows of $\mathbb{C}$ making the triangles to commute. Also, a sub-object of $a$ is a mono-morphism in $\mathbb{C}$ having $a$ as co-domain.

A sub-object category forms a Heyting algebra when the base category is Heyting. Moreover, a Heyting category can be intuitively visualised as a collection of Heyting algebras linked by appropriate arrows so to have enough pullbacks.

**Proposition 7.3.** $\text{Sub}_{\mathbb{C}}(a)$ *is a pre-order.*

*Proof.* As before, given $m_1, m_2 \in \text{Obj}\,\text{Sub}_{\mathbb{C}}(a)$, i.e., $m_1\colon b_1 \rightarrowtail a$ and $m_2\colon b_2 \rightarrowtail a$, $m_1 \leq m_2$ if and only if there is $r \in \text{Hom}\,\text{Sub}_{\mathbb{C}}(m_1, m_2)$, i.e., $r\colon b_1 \to b_2$ in $\mathbb{C}$, such that $m_1 = m_2 \circ r$.

- $m \leq m$ since $(m\colon b \rightarrowtail a) = m \circ 1_b$;
- if $m_1 \leq m_2$ and $m_2 \leq m_1$, then $(m_1\colon b_1 \rightarrowtail a) = (m_2\colon b_2 \rightarrowtail a) \circ (h_1\colon b_1 \to b_2)$ and $m_2 = m_1 \circ (h_2\colon b_2 \to b_1)$, so $m_1 \circ 1_{b_1} = m_1 = m_2 \circ h_1 = m_1 \circ h_1 \circ h_1$. Being $m_1$ mono, $h_2 \circ h_1 = 1_{b_1}$. Conversely, $m_2 \circ 1_{b_2} = m_2 = m_1 \circ h_2 = m_2 \circ h_1 \circ h_2$, so $h_1 \circ h_2 = 1_{b_2}$. Thus, $h_2 = h_1^{-1}$.
- if $m_1 \leq m_2$ and $m_2 \leq m_3$, then $(m_1\colon b_1 \rightarrowtail a) = (m_2\colon b_2 \rightarrowtail a) \circ (h_1\colon b_1 \to b_2)$ and $m_2 = (m_3\colon b_3 \rightarrowtail a) \circ (h_2\colon b_2 \to b_3)$, so $m_1 = m_2 \circ h_1 = m_3 \circ h_2 \circ h_1$, i.e., $m_1 \leq m_3$. □

**Corollary 7.4.** $\text{Sub}_{\mathbb{C}}(a)$ *is a poset, modulo isomorphisms.*

**Proposition 7.5.** *If* $\mathbb{C}$ *is Cartesian, so is* $\mathbb{C}/a$.

*Proof.* Immediate, noticing that $\mathbb{C}$ is Cartesian when it has a terminal object and all the pullbacks. □

**Proposition 7.6.** *If* $\mathbb{C}$ *is Cartesian, so is* $\text{Sub}_{\mathbb{C}}(a)$.

*Proof.* See [Joh02a, A.1.3]. □

**Proposition 7.7.** *If* $\mathbb{C}$ *is Cartesian,* $\text{Sub}_{\mathbb{C}}(a)$ *has a top element* $\top$ *and finite meets.*

*Proof.* Obvious, since $\text{Sub}_{\mathbb{C}}(a)$ is Cartesian and thus has a terminal object $\top$ and binary products, the meets. □

**Definition 7.11** (Pullback functor)**.** If $\mathbb{C}$ is regular and $f \in \text{Hom}_{\mathbb{C}}(A,B)$, then the *pullback functor* $f^*\colon \text{Sub}_{\mathbb{C}}(B) \to \text{Sub}_{\mathbb{C}}(A)$ is defined, up to isomorphisms, as:

(1) let $g \in \mathrm{Obj}\,\mathrm{Sub}_{\mathbb{C}}(B)$, i.e., $g\colon X \rightarrowtail B$, then $f^*(g)$ is the pullback of $g$ along $f$

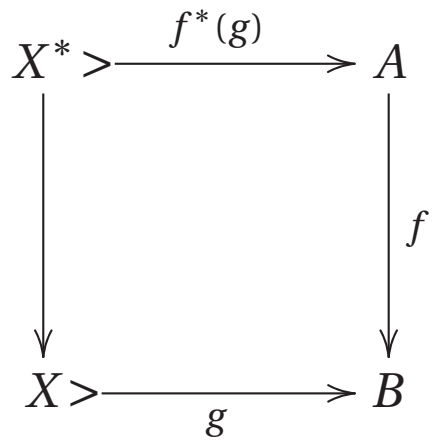

(2) let $h \in \mathrm{Hom}_{\mathrm{Sub}_{\mathbb{C}}(B)}(g_1, g_2)$, i.e., $g_1\colon X_1 \rightarrowtail B$, $g_2\colon X_2 \rightarrowtail B$, $h\colon X_1 \to X_2$ and $g_1 = g_2 \circ h$, then $f^*(h)$ is the universal arrow of the pullback as in the diagram:

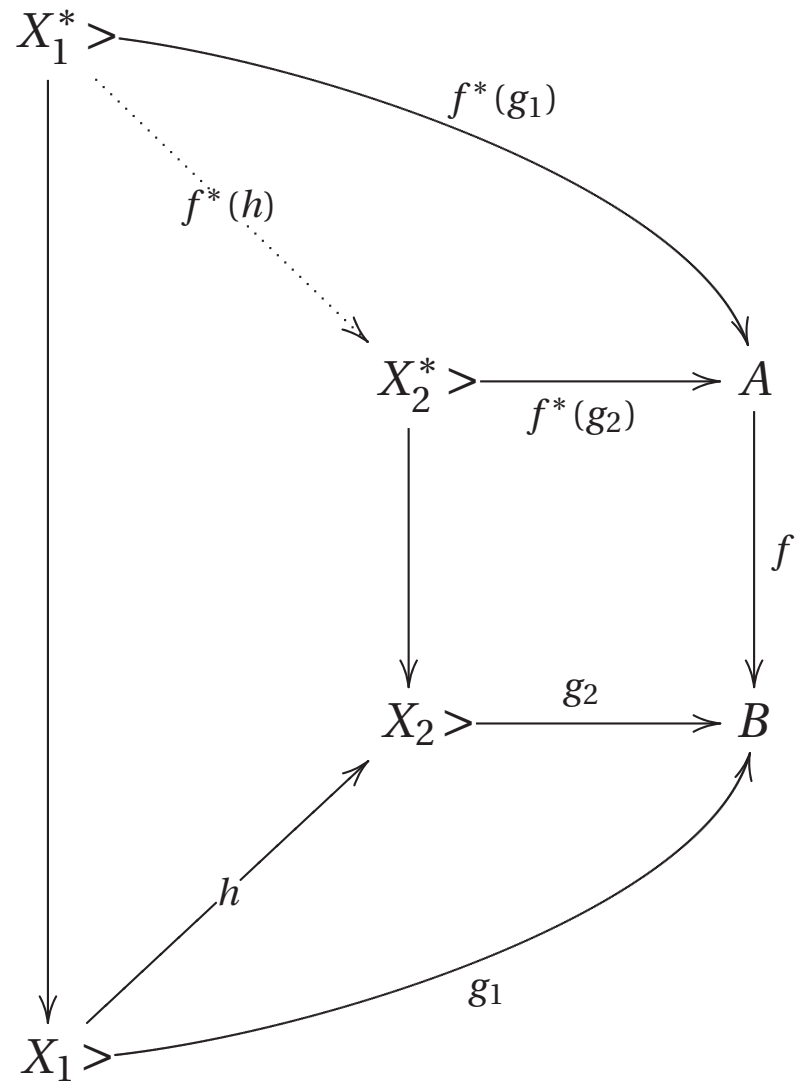

**Definition 7.12** (Coherent category)**.** A regular category $\mathbb{C}$ is *coherent* when, for each $c \in \mathrm{Obj}\,\mathbb{C}$, $\mathrm{Sub}_{\mathbb{C}}(c)$ seen as a pre-order category, has finite joins and each pullback functor preserves them.

**Proposition 7.8.** *In a coherent category* $\mathbb{C}$*, each* $\mathrm{Sub}_{\mathbb{C}}(a)$ *is a bounded lattice up to isomorphisms.*

*Proof.* Identifying isomorphic sub-objects, $\mathrm{Sub}_{\mathbb{C}}(a)$ is a poset, as seen in Corollary 7.4. Having all the finite meets, see Proposition 7.7, and joins, the statement follows. □

**Proposition 7.9.** *If the category* $\mathbb{C}$ *is regular, for each* $f\colon A \to B$ *in* $\mathbb{C}$*, the pullback functor* $f^*\colon \mathrm{Sub}_{\mathbb{C}}(B) \to \mathrm{Sub}_{\mathbb{C}}(A)$ *has a left adjoint* $\exists_f$*.*

*Proof.* See [Joh02a, Lemma A.1.3.1]. □

It is worth to explicitly construct the left adjoint $\exists_f$ from $f^*$: let $f\colon B \to A$ in $\mathbb{C}$ and let $b \in \mathrm{Obj}\,\mathrm{Sub}_{\mathbb{C}}(B)$, i.e., $b\colon C \rightarrowtail B$ in $\mathbb{C}$. Then, a reflection $(R_b, \eta_b)$ of $b$ along $f^*$ is such that

(1) $R_b\colon D \rightarrowtail A$ in $\mathbb{C}$;
(2) $\eta_b\colon C \to P$ in $\mathbb{C}$, where $f^*(R_b)\colon P \rightarrowtail B$ and $b = f^*(R_b) \circ \eta_b$;
(3) for each $a\colon X \rightarrowtail A$ in $\mathbb{C}$ and $\beta\colon C \to X^*$ where $f^*(a)\colon X^* \rightarrowtail B$ and $b = f^*(a) \circ \beta$, there is a unique $\alpha\colon D \to X$ such that $R_b = a \circ \alpha$ and $\beta = f^*(\alpha) \circ \eta_b$.

In a diagram:

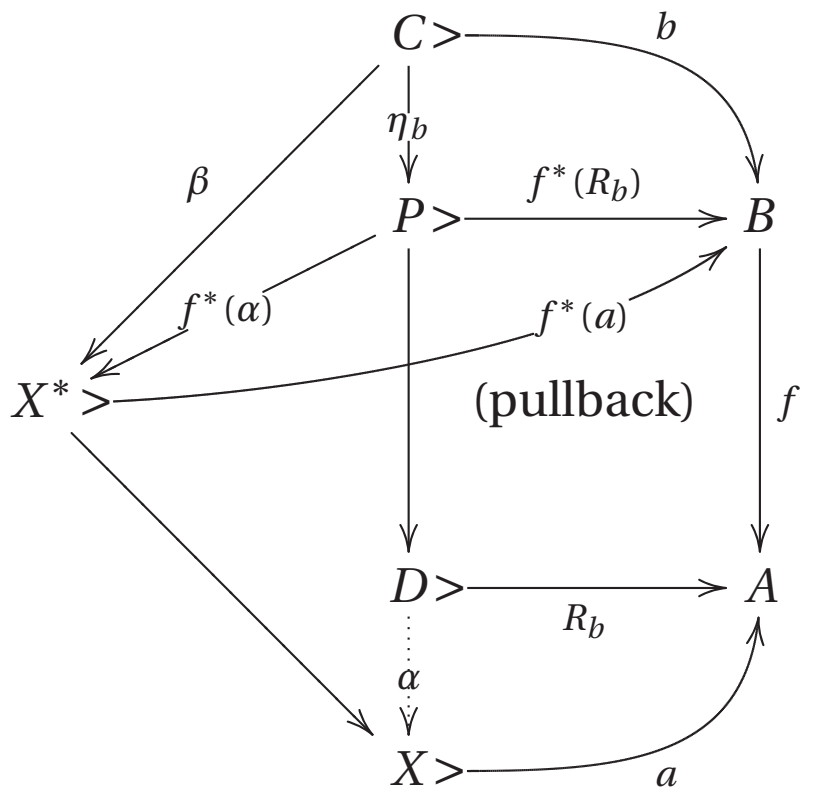

Thus, a reflection $(R_b, \eta_b)$ is the least sub-object of $A$ whose pullback along $f$ is greater than $b$. Then, the functor $\exists_f \colon \mathrm{Sub}_{\mathbb{C}}(B) \to \mathrm{Sub}_{\mathbb{C}}(A)$ sends each sub-object $b$ of $B$ to $R_b$, the least sub-object of $A$ whose pullback along $f$ is greater than $b$.

**Proposition 7.10.** *In a regular category* $\mathbb{C}$, *for each* $f\colon A \to B$, $\mathrm{im}(f) \cong \exists_f(1_A)$ *in* $\mathrm{Sub}_{\mathbb{C}}(A)$.

*Proof.* See [Joh02a, Lemma A.1.3.1]. □

**Proposition 7.11.** *In a coherent category* $\mathbb{C}$*, the initial object* $0$ *of* $\mathbb{C}$ *exists and it is the initial object of* $\mathrm{Sub}_{\mathbb{C}}(1)$.

*Proof.* See [Joh02a, Lemma A.1.4.1]. □

**Proposition 7.12.** *Let* $A_1, A_2 \in \mathrm{Obj}\,\mathrm{Sub}_{\mathbb{C}}(A)$ *in a coherent category* $\mathbb{C}$. *When* $A_1 \wedge A_2 \cong 0$, $A_1 \vee A_2$ *is a co-product of* $A_1$ *and* $A_2$.

*Proof.* See [Joh02a, Corollary A.1.4.4]. □

**Definition 7.13** (Heyting category)**.** A coherent category $\mathbb{C}$ is a *Heyting category* when, for each $f \in \mathrm{Hom}_{\mathbb{C}}(a, b)$, the pullback functor $f^*\colon \mathrm{Sub}_{\mathbb{C}}(b) \to \mathrm{Sub}_{\mathbb{C}}(a)$ has a right adjoint $\forall_f\colon \mathrm{Sub}_{\mathbb{C}}(a) \to \mathrm{Sub}_{\mathbb{C}}(b)$.

It is worth to explicitly construct the right adjoint $\forall_f$ from $f^*$: let $f\colon B \to A$ in $\mathbb{C}$ and let $b \in \mathrm{Obj}\,\mathrm{Sub}_{\mathbb{C}}(B)$, i.e., $b\colon C \rightarrowtail B$ in $\mathbb{C}$. Then, a co-reflection $(R_b, \varepsilon_b)$ of $b$ along $f^*$ is such that

(1) $R_b\colon D \rightarrowtail A$ in $\mathbb{C}$;
(2) $\varepsilon_b\colon P \to C$ in $\mathbb{C}$, where $f^*(R_b)\colon P \rightarrowtail B$ and $f^*(R_b) = b \circ \varepsilon_b$;
(3) for each $a\colon X \rightarrowtail A$ in $\mathbb{C}$ and $\beta\colon X^* \to C$ where $f^*(a)\colon X^* \rightarrowtail B$ and $f^*(a) = b \circ \beta$, there is a unique $\alpha\colon X \to D$ such that $a = R_b \circ \alpha$ and $\beta = \varepsilon_b \circ f^*(\alpha)$.

In a diagram:

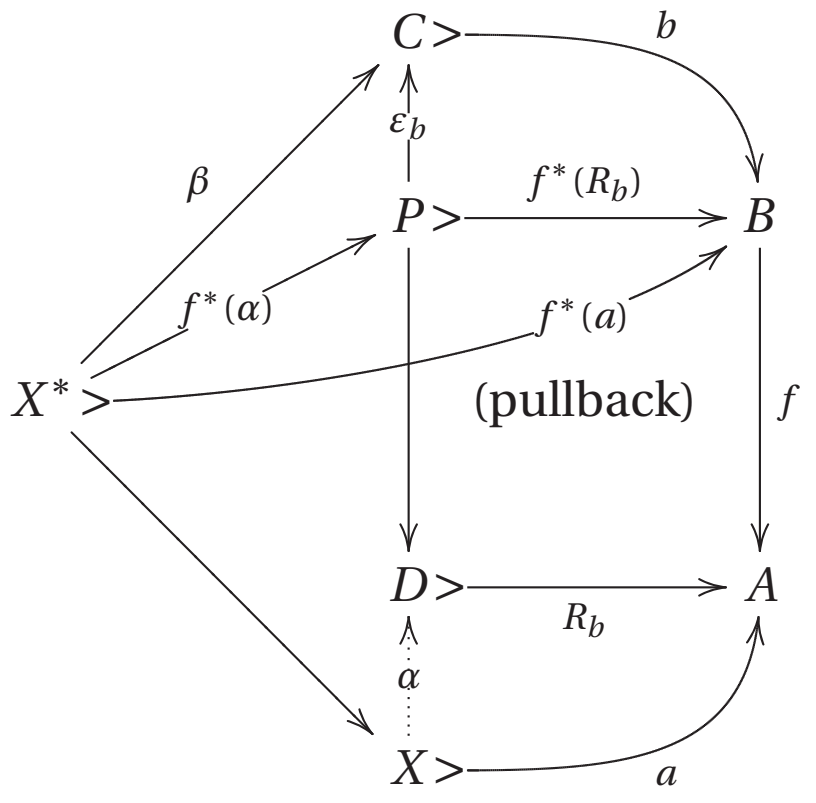

Thus, a co-reflection $(R_b, \varepsilon_b)$ is the greatest sub-object of $A$ whose pullback along $f$ is less than $b$. Then, the functor $\forall_f \colon \mathrm{Sub}_{\mathbb{C}}(B) \to \mathrm{Sub}_{\mathbb{C}}(A)$ sends each sub-object $b$ of $B$ to $R_b$, the greatest sub-object of $A$ whose pullback along $f$ is less than $b$.

**Proposition 7.13.** *Let $A_1 \rightarrowtail A$ and $A_2 \rightarrowtail A$ be sub-objects in a Heyting category. Then there exists a largest sub-object $A_3 \rightarrowtail A$ such that $A_3 \wedge A_1 \leq A_2$. Denoting this sub-object by $(A_1 \Rightarrow A_2) \rightarrowtail A$, the binary operation thus defined is stable under pullback.*

*Proof.* See [Joh02a, Lemma A.1.4.13]. □

**Corollary 7.14.** *In a Heyting category $\mathbb{C}$, for each $a \in \mathrm{Obj}\,\mathbb{C}$, $\mathrm{Sub}_{\mathbb{C}}(a)$ is a Heyting algebra.*

The purpose of Heyting categories is to provide a framework where it is possible to give meaning to logical formulae in such a way that a soundness and completeness theorem holds.

**Definition 7.14** ($\Sigma$-H-structure)**.** Given a Heyting category $\mathbb{C}$ and a logical signature $\Sigma$, a $\Sigma$*-H-structure* $M$ in $\mathbb{C}$ is specified by:

(1) a function assigning to each sort symbol $A$ in $\Sigma$ an object $MA$ in $\mathbb{C}$;
(2) a function assigning to each function symbol $f \colon A_1 \times \cdots \times A_n \to B$ in $\Sigma$ an arrow $Mf \colon MA_1 \times \cdots \times MA_n \to MB$ in $\mathbb{C}$;
(3) a function assigning to each relation symbol $R \colon A_1 \times \cdots \times A_n$ in $\Sigma$ a sub-object $MR \rightarrowtail MA_1 \times \cdots \times MA_n$ in $\mathbb{C}$.

For each finite string of sorts, we abbreviate $M(A_1 \times \cdots \times A_n) \equiv MA_1 \times \cdots \times MA_n$ and $M([]) \equiv 1$. Also, we speak of $\Sigma$-structure when it is clear from the context that we are dealing with a Heyting category instead of a logically distributive category.

**Definition 7.15** (Semantics of terms)**.** Let $M$ be a $\Sigma$-structure in a Heyting category $\mathbb{C}$. If $\vec{x}.\,t{:}B$ is a logical term-in-context over $\Sigma$ with $\vec{x} \equiv x_1{:}A_1,\ldots,x_n{:}A_n$, then $[\![\vec{x}.\,t{:}B]\!] \colon MA_1 \times \cdots \times MA_n \to MB$ is recursively defined as:

(1) if $t$ is a variable, it is necessarily $x_i$ for some $1 \leq i \leq n$, and $[\![\vec{x}.\,t{:}B]\!] = \pi_i$, the $i$-th projector of $MA_1 \times \cdots \times MA_n$;
(2) if $t \equiv f(t_1,\ldots,t_m)$ where $t_i{:}C_i$ is a term for all $1 \leq i \leq n$, then $[\![\vec{x}.\,t{:}B]\!] = Mf \circ \left([\![\vec{x}.\,t_1{:}C_1]\!],\ldots,[\![\vec{x}.\,t_m{:}C_m]\!]\right)$.

**Definition 7.16** (Semantics of logical formulas)**.** Let $M$ be a $\Sigma$-structure in a Heyting category $\mathbb{C}$. A formula-in-context $\vec{x}.\phi$ over $\Sigma$, where $\vec{x} \equiv x_1{:}A_1,\ldots,x_n{:}A_n$, is interpreted as a sub-object $[\![\vec{x}.\phi]\!] \rightarrowtail MA_1 \times \cdots \times MA_n$ according to the following recursive definition:

(1) $\phi \equiv R(t_1,\ldots,t_m)$ with $t_i{:}B_i$ a term for each $1 \leq i \leq m$ and $R{:}B_1 \times \cdots \times B_m$ a relation symbol in $\Sigma$: $[\![\vec{x}.\phi]\!]$ is the pullback of $MR \rightarrowtail MB_1 \times \cdots \times MB_m$ along $\left([\![\vec{x}.\,t_1{:}B_1]\!],\ldots,[\![\vec{x}.\,t_m{:}B_m]\!]\right)$

$$\begin{array}{ccc}
[\![\vec{x}.\phi]\!] & \rightarrowtail & MA_1 \times \cdots \times MA_n \\
\downarrow & & \downarrow{\scriptstyle ([\![\vec{x}.\,t_1:B_1]\!],\ldots,[\![\vec{x}.\,t_m:B_m]\!])} \\
MR & \rightarrowtail & MB_1 \times \cdots \times MB_m
\end{array}$$

(2) $\phi \equiv \top$: $[\![\vec{x}.\phi]\!] = 1_{\mathrm{Sub}_{\mathbb{C}}(MA_1 \times \cdots \times MA_n)}$;
(3) $\phi \equiv \bot$: $[\![\vec{x}.\phi]\!] = 0_{\mathrm{Sub}_{\mathbb{C}}(MA_1 \times \cdots \times MA_n)}$;
(4) $\phi \equiv \psi_1 \wedge \psi_2$: $[\![\vec{x}.\phi]\!] = [\![\vec{x}.\psi_1]\!] \wedge [\![\vec{x}.\psi_2]\!]$, the meet of $[\![\vec{x}.\psi_1]\!]$ and $[\![\vec{x}.\psi_2]\!]$ in the Heyting algebra $\mathrm{Sub}_{\mathbb{C}}(MA_1 \times \cdots \times MA_n)$;
(5) $\phi \equiv \psi_1 \vee \psi_2$: $[\![\vec{x}.\phi]\!] = [\![\vec{x}.\psi_1]\!] \vee [\![\vec{x}.\psi_2]\!]$, the join of $[\![\vec{x}.\psi_1]\!]$ and $[\![\vec{x}.\psi_2]\!]$ in the Heyting algebra $\mathrm{Sub}_{\mathbb{C}}(MA_1 \times \cdots \times MA_n)$;
(6) $\phi \equiv \psi_1 \supset \psi_2$: $[\![\vec{x}.\phi]\!] = [\![\vec{x}.\psi_1]\!] \Rightarrow [\![\vec{x}.\psi_2]\!]$, the implication from $[\![\vec{x}.\psi_1]\!]$ to $[\![\vec{x}.\psi_2]\!]$ in the Heyting algebra $\mathrm{Sub}_{\mathbb{C}}(MA_1 \times \cdots \times MA_n)$;

(7) $\phi \equiv \exists y\colon s.\,\psi$: $[\![\vec{x}.\,\phi]\!] = \exists_\pi\left([\![\vec{x}, y\colon s.\,\psi]\!]\right)$, where $\pi\colon MA_1 \times \cdots \times MA_n \times MS \to MA_1 \times \cdots \times MA_n$ is the canonical projection of the first $n$ factors;

(8) $\phi \equiv \forall y\colon s.\,\psi$: $[\![\vec{x}.\,\phi]\!] = \forall_\pi\left([\![\vec{x}, y\colon s.\,\psi]\!]\right)$, where $\pi\colon MA_1 \times \cdots \times MA_n \times MS \to MA_1 \times \cdots \times MA_n$ is the canonical projection of the first $n$ factors.

**Definition 7.17** (Heyting validity)**.** Let $M$ be a $\Sigma$-structure in a Heyting category $\mathbb{C}$. If $\sigma \equiv \phi \vdash_{\vec{x}} \psi$ is a sequent over $\Sigma$ in the context $\vec{x}$, we say that $\sigma$ is *valid* in $M$ when $[\![\vec{x}.\,\phi]\!] \leq [\![\vec{x}.\,\psi]\!]$ in $\mathrm{Sub}_{\mathbb{C}}(MA_1 \times \cdots \times MA_n)$.

**Definition 7.18** (Heyting model)**.** Let $M$ be a $\Sigma$-structure in a Heyting category $\mathbb{C}$. A theory $T$ over $\Sigma$ has $M$ as a *model* if and only if all the axioms of $T$ are valid in $M$.

**Theorem 7.15** (Soundness)**.** *Let $T$ be a first-order theory over a logical signature $\Sigma$ and let $M$ be a model for $T$ in a Heyting category $\mathbb{C}$. If $\sigma$ is a sequent which is provable in $T$, then $\sigma$ is valid in $M$.*

*Proof.* See [Joh02b, Proposition D.1.3.2]. □

**Theorem 7.16** (Completeness)**.** *Let $T$ be a first-order theory over a logical signature $\Sigma$. If a sequent $\sigma$ is valid in all models of $T$ in any Heyting category, then $\sigma$ is provable in $T$.*

*Proof.* See [Joh02b, Proposition D.1.4.11]. □

*Note* 7.1. It is worth remarking that the generic model $M_T$, see [Joh02b, p. 845], has **not** the universal property that all $T$-models are images of it under Heyting functors, as remarked after [Joh02b, Proposition D.1.4.11]. Thus, the generic model is not a classifying object for $T$-models.

## 8. Heyting versus logically distributive categories

In this section, we relate the semantics based on Heyting categories with the one based on logically distributive categories. After presenting a couple of useful properties of the interpretation in Heyting categories, we will use it to build a functor which maps models in Heyting categories into models in logically distributive categories preserving truth.

The first property says that, given a formula $\phi$, its interpretation in a Heyting category can be restricted to the minimal context $\mathrm{FV}(\phi)$, as the interpretations in richer contexts can be recovered by an appropriate pullback.

**Proposition 8.1.** *Let $\phi$ be a logical formula, and $M$ a $\Sigma$-structure in a Heyting category $\mathbb{C}$, where $\phi$ is written in the language of $\Sigma$. Then, for every context $\vec{x} \equiv x_1 : A_1, \dots, x_n : A_n$ such that $\mathrm{FV}(\phi) \subseteq \vec{x}$, $[\![\vec{x}.\phi]\!] \colon X \rightarrowtail MA_1 \times \cdots \times MA_n$. It holds that $[\![\vec{x}.\phi]\!]$ is the pullback of $[\![\mathrm{FV}(\phi).\phi]\!]$ along the canonical projection $\pi \colon MA_1 \times \cdots \times MA_n \to M(\mathrm{FV}(\phi))$.*

*Proof.* By induction on the structure of $\phi$:

(1) $\phi \equiv R(t_1, \dots, t_m)$ with $t_1 : B_i$ a term for each $1 \le i \le m$ and $R : B_1 \times \cdots \times B_m$ a relation symbol in $\Sigma$: $[\![\vec{x}.\phi]\!]$ is the pullback of $MR \rightarrowtail MB_1 \times \cdots \times MB_m$, the interpretation of $R$ in $\mathbb{C}$, along $\big([\![\vec{x}.\, t_1 : B_1]\!], \dots, [\![\vec{x}.\, t_m : B_m]\!]\big)$. Similarly, $[\![\mathrm{FV}(\phi).\phi]\!]$ is the pullback of $MR \rightarrowtail MB_1 \times \cdots \times MB_m$ along $\big([\![\mathrm{FV}(\phi).\, t_1 : B_1]\!], \dots, [\![\mathrm{FV}(\phi).\, t_m : B_m]\!]\big)$.

Universality ensures that

$$\big([\![\vec{x}.\, t_1 : B_1]\!], \dots, [\![\vec{x}.\, t_m : B_m]\!]\big) = \big([\![\mathrm{FV}(\phi).\, t_1 : B_1]\!], \dots, [\![\mathrm{FV}(\phi).\, t_m : B_m]\!]\big) \circ \pi \ .$$

Thus, $a$ is the universal arrow of the pullback as in the following diagram:

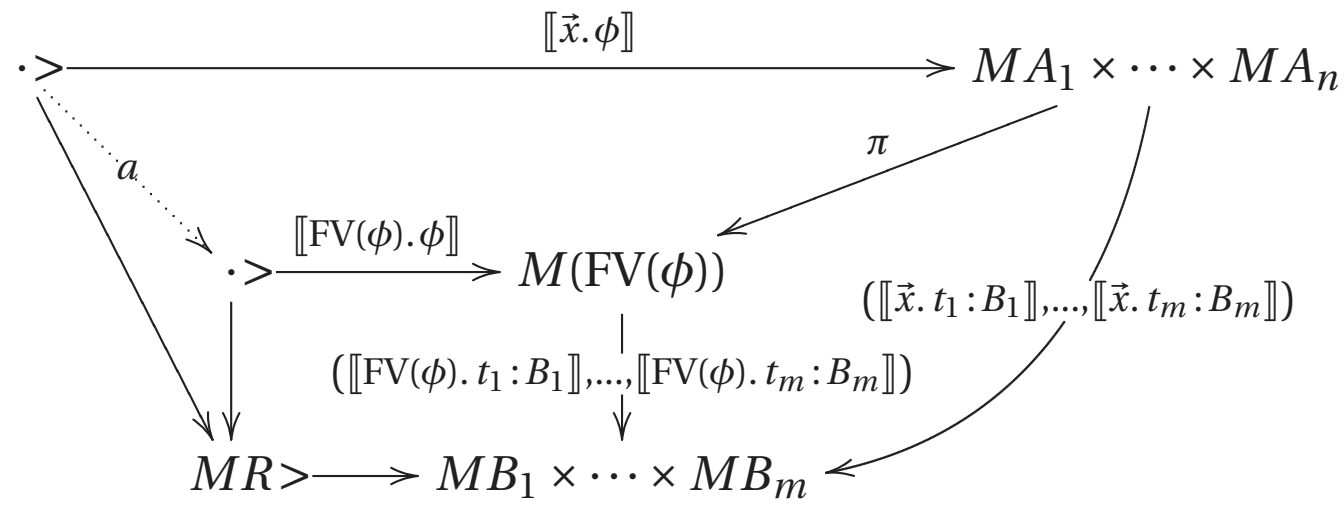

So, by the pullback lemma [Gol06], the conclusion follows.

(2) $\phi \equiv \top$: since $\mathrm{FV}(\phi) = \emptyset$, $M(\mathrm{FV}(\phi)) = 1$, the terminal object of $\mathbb{C}$. This fact makes evident that the following is a pullback diagram:

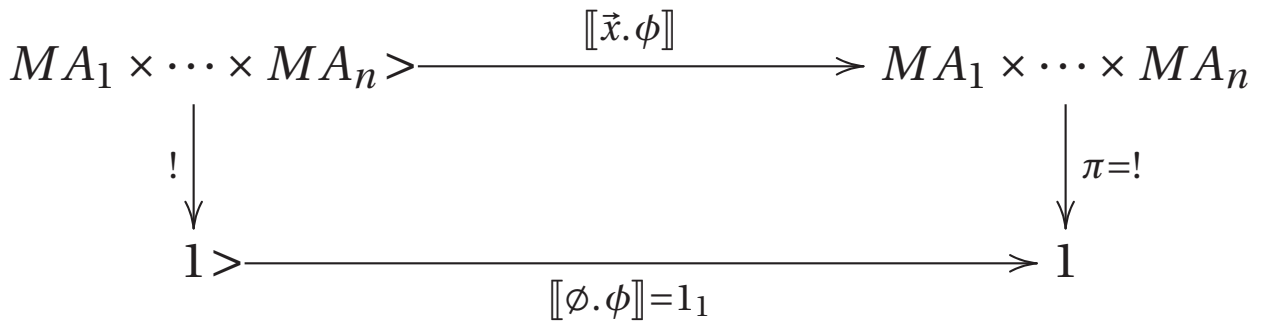

(3) $\phi \equiv \bot$: since $\mathrm{FV}(\phi) = \emptyset$, $M(\mathrm{FV}(\phi)) = 1$. This fact makes evident that the following is a pullback diagram:

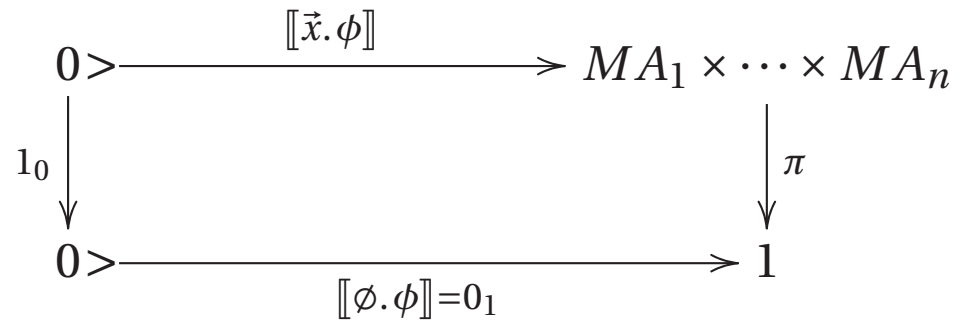

(4) $\phi \equiv \psi_1 \wedge \psi_2$: by induction hypothesis, we know that:

- $[\![\mathrm{FV}(\phi).\psi_1]\!]$ is the pullback of $[\![\mathrm{FV}(\psi_1).\psi_1]\!]$ along the canonical projection $M(\mathrm{FV}(\phi)) \to M(\mathrm{FV}(\psi_1))$;

- $[\![\mathrm{FV}(\phi).\psi_2]\!]$ is the pullback of $[\![\mathrm{FV}(\psi_2).\psi_2]\!]$ along $M(\mathrm{FV}(\phi)) \to M(\mathrm{FV}(\psi_2))$;
- $[\![\vec{x}.\psi_1]\!]$ is the pullback of $[\![\mathrm{FV}(\psi_1).\psi_1]\!]$ along the evident projection $MA_1 \times \cdots \times MA_n \to M(\mathrm{FV}(\psi_1))$;
- $[\![\vec{x}.\psi_2]\!]$ is the pullback of $[\![\mathrm{FV}(\psi_2).\psi_2]\!]$ along $MA_1 \times \cdots \times MA_n \to M(\mathrm{FV}(\psi_2))$.

Thus, by the pullback lemma, $[\![\vec{x}.\psi_1]\!]$ is the pullback of $[\![\mathrm{FV}(\phi).\psi_1]\!]$ along $\pi$, and $[\![\vec{x}.\psi_2]\!]$ is the pullback of $[\![\mathrm{FV}(\phi).\psi_2]\!]$ along $\pi$.

Consider $\pi_1 \colon [\![\mathrm{FV}(\phi).\phi]\!] \to [\![\mathrm{FV}(\phi).\psi_1]\!]$ and $\pi_2 \colon [\![\mathrm{FV}(\phi).\phi]\!] \to [\![\mathrm{FV}(\phi).\psi_2]\!]$ in $\mathrm{Sub}_{\mathbb{C}}(M(\mathrm{FV}(\phi)))$, and the corresponding $\pi_1^* \colon [\![\vec{x}.\phi]\!] \to [\![\vec{x}.\psi_1]\!]$ and $\pi_2^* \colon [\![\vec{x}.\phi]\!] \to [\![\vec{x}.\psi_2]\!]$ in $\mathrm{Sub}_{\mathbb{C}}(MA_1 \times \cdots \times MA_n)$.

Composing the pullback $g_1$ of $\pi$ along $[\![\mathrm{FV}(\phi).\psi_1]\!]$ with $\pi_1^*$, and the pullback $g_2$ of $\pi$ along $[\![\mathrm{FV}(\phi).\psi_2]\!]$ with $\pi_2^*$, there is an arrow $f \colon X^* \to X$, where $X^*$ is the codomain of $[\![\vec{x}.\phi]\!]$ and $X$ is the codomain of $[\![\mathrm{FV}(\phi).\phi]\!]$, since $X$ is a product in $\mathbb{C}$, such that $\pi_1 \circ f = g_1 \circ \pi_1^*$ and $\pi_2 \circ f = g_2 \circ \pi_2^*$. By diagram chasing and the properties of products, the following diagram commutes:

$$\begin{array}{ccc} X^* & \xrightarrow{[\![\vec{x}.\phi]\!]} & MA_1 \times \cdots \times MA_n \\ {\scriptstyle f}\big\downarrow & & \big\downarrow{\scriptstyle \pi} \\ X & \xrightarrow[{[\![\mathrm{FV}(\phi).\phi]\!]}]{} & M(\mathrm{FV}(\phi)) \end{array}$$

Suppose $h \colon P \to X$ and $k \colon P \to MA_1 \times \cdots \times MA_n$ are such that $[\![\mathrm{FV}(\phi).\phi]\!] \circ h = \pi \circ k$. Then $\pi \circ k = [\![\mathrm{FV}(\phi).\psi_1]\!] \circ \pi_1 \circ h$ and $\pi \circ k = [\![\mathrm{FV}(\phi).\psi_2]\!] \circ \pi_2 \circ h$, thus there are universal pullback arrows from $P$ to the domains of $[\![\vec{x}.\psi_1]\!]$ and $[\![\vec{x}.\psi_2]\!]$.

But $X^*$ is a product, so there is a universal arrow from $P$ to $X^*$, showing that the diagram above is a pullback.

(5) $\phi \equiv \psi_1 \vee \psi_2$: by induction hypothesis, we know that:

- $[\![\mathrm{FV}(\phi).\psi_1]\!]$ is the pullback of $[\![\mathrm{FV}(\psi_1).\psi_1]\!]$ along the canonical projection $M(\mathrm{FV}(\phi)) \to M(\mathrm{FV}(\psi_1))$;
- $[\![\mathrm{FV}(\phi).\psi_2]\!]$ is the pullback of $[\![\mathrm{FV}(\psi_2).\psi_2]\!]$ along $M(\mathrm{FV}(\phi)) \to M(\mathrm{FV}(\psi_2))$;
- $[\![\vec{x}.\psi_1]\!]$ is the pullback of $[\![\mathrm{FV}(\psi_1).\psi_1]\!]$ along the evident projection $MA_1 \times \cdots \times MA_n \to M(\mathrm{FV}(\psi_1))$;
- $[\![\vec{x}.\psi_2]\!]$ is the pullback of $[\![\mathrm{FV}(\psi_2).\psi_2]\!]$ along $MA_1 \times \cdots \times MA_n \to M(\mathrm{FV}(\psi_2))$.

Thus, by the pullback lemma, $[\![\vec{x}.\psi_1]\!]$ is the pullback of $[\![\mathrm{FV}(\phi).\psi_1]\!]$ along $\pi$, and $[\![\vec{x}.\psi_2]\!]$ is the pullback of $[\![\mathrm{FV}(\phi).\psi_2]\!]$ along $\pi$.

Consider the injections $\iota_1 \colon [\![\mathrm{FV}(\phi).\psi_1]\!] \to [\![\mathrm{FV}(\phi).\phi]\!]$ and $\iota_2 \colon [\![\mathrm{FV}(\phi).\psi_2]\!] \to [\![\mathrm{FV}(\phi).\phi]\!]$ in $\mathrm{Sub}_{\mathbb{C}}(M(\mathrm{FV}(\phi)))$, and the corresponding $\iota_1^* \colon [\![\vec{x}.\psi_1]\!] \to [\![\vec{x}.\phi]\!]$ and $\iota_2^* \colon [\![\vec{x}.\psi_2]\!] \to [\![\vec{x}.\phi]\!]$ in $\mathrm{Sub}_{\mathbb{C}}(MA_1 \times \cdots \times MA_n)$.

Composing $\iota_1$ with the pullback $g_1$ of $\pi$ along $[\![\mathrm{FV}(\phi).\psi_1]\!]$, and $\iota_2$ with the pullback $g_2$ of $\pi$ along $[\![\mathrm{FV}(\phi).\psi_2]\!]$, there is the co-universal arrow $f \colon X^* \to X$ of the $X^*$ co-product, where $X^*$ is the domain of $[\![\vec{x}.\phi]\!]$ and $X$ is the codomain of $[\![\mathrm{FV}(\phi).\phi]\!]$ in $\mathbb{C}$. By diagram chasing and the properties of co-products, the following diagram commutes:

$$\begin{array}{ccc} X^* & \xrightarrow{[\![\vec{x}.\phi]\!]} & MA_1 \times \cdots \times MA_n \\ {\scriptstyle f}\big\downarrow & & \big\downarrow{\scriptstyle \pi} \\ X & \xrightarrow[{[\![\mathrm{FV}(\phi).\phi]\!]}]{} & M(\mathrm{FV}(\phi)) \end{array}$$

But $\mathbb{C}$ is a Heyting category so, in particular, it is a coherent category, thus the pullback functor $\pi^*$ transforms $X$ in $X^*$, $\iota_1$ in $\iota_1^*$ and $\iota_2$ in $\iota_2^*$, showing that the above diagram is a pullback, up to isomorphisms.

(6) $\phi \equiv \psi_1 \supset \psi_2$: by induction hypothesis, we know that:

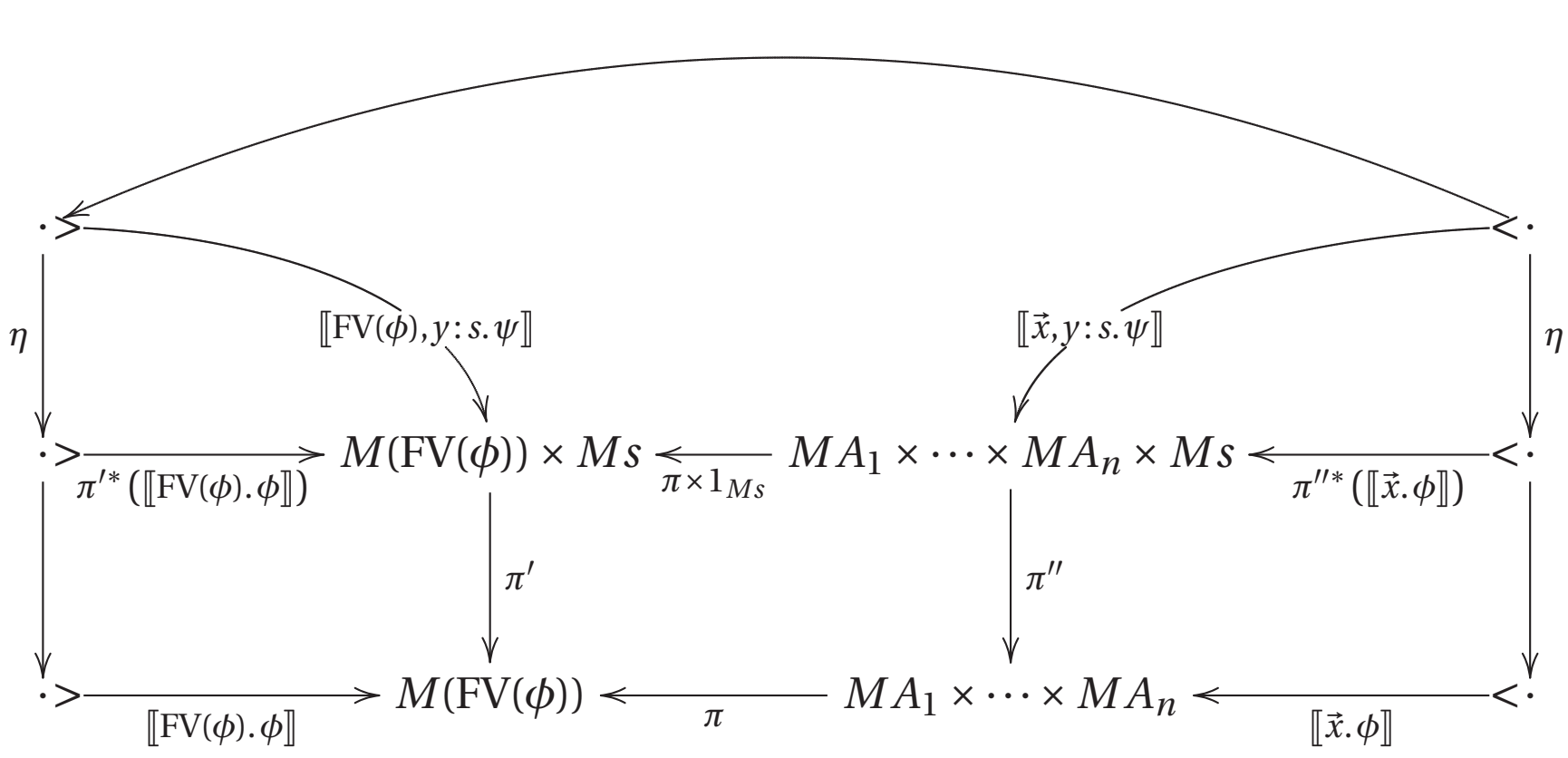

FIGURE 10. The existential case in proposition 8.1.

- $[\![\mathrm{FV}(\phi).\psi_1]\!]$ is the pullback of $[\![\mathrm{FV}(\psi_1).\psi_1]\!]$ along the canonical projection $M(\mathrm{FV}(\phi)) \to M(\mathrm{FV}(\psi_1))$;
- $[\![\mathrm{FV}(\phi).\psi_2]\!]$ is the pullback of $[\![\mathrm{FV}(\psi_2).\psi_2]\!]$ along $M(\mathrm{FV}(\phi)) \to M(\mathrm{FV}(\psi_2))$;
- $[\![\vec{x}.\psi_1]\!]$ is the pullback of $[\![\mathrm{FV}(\psi_1).\psi_1]\!]$ along the evident projection $MA_1 \times \cdots \times MA_n \to M(\mathrm{FV}(\psi_1))$;
- $[\![\vec{x}.\psi_2]\!]$ is the pullback of $[\![\mathrm{FV}(\psi_2).\psi_2]\!]$ along $MA_1 \times \cdots \times MA_n \to M(\mathrm{FV}(\psi_2))$.

Thus, by the pullback lemma, $[\![\vec{x}.\psi_1]\!]$ is the pullback of $[\![\mathrm{FV}(\phi).\psi_1]\!]$ along $\pi$, and $[\![\vec{x}.\psi_2]\!]$ is the pullback of $[\![\mathrm{FV}(\phi).\psi_2]\!]$ along $\pi$.

So, by proposition 7.13 the conclusion follows.

(7) $\phi \equiv \exists y\colon s.\psi$: By induction hypothesis, the arrow $[\![\vec{x}, y\colon s.\psi]\!]$ is the pullback of $[\![\mathrm{FV}(\phi), y\colon s.\psi]\!]$ along the projection $\pi \times 1_{Ms}$.

Denoting by $\pi'\colon M(\mathrm{FV}(\phi)) \times Ms \to M(\mathrm{FV}(\phi))$ and $\pi''\colon MA_1 \times \cdots \times MA_n \times Ms \to MA_1 \times \cdots \times MA_n$ the canonical projections, by definition of the $\exists_{\pi'}$ and $\exists_{\pi''}$ functors, the diagram in figure 10 commutes: Now, $[\![\mathrm{FV}(\phi).\phi]\!]$ is the image of $\pi' \circ [\![\mathrm{FV}(\phi), y\colon s.\psi]\!]$, and $[\![\vec{x}.\phi]\!]$ is the image of $\pi'' \circ [\![\vec{x}, y\colon s.\psi]\!]$. Moreover, the central square in the diagram is a pullback as it is immediate to verify from the properties of products and projections. Thus, by the pullback lemma, the composition of the upper pullback with the central square yields a pullback. Its left and right sides have $[\![\mathrm{FV}(\phi).\phi]\!]$ and $[\![\vec{x}.\phi]\!]$ as images, respectively. But $\mathbb{C}$ is a Heyting category, hence regular, so images are preserved by pullbacks, thus the conclusion follows.

(8) $\phi \equiv \forall y\colon s.\psi$: Consider the diagram in figure 11:

Let $P$ be the pullback of $[\![\mathrm{FV}(\phi).\phi]\!]$ and $\pi \circ [\![\vec{x}.\phi]\!]$, and let $P^*$ be the pullback of $\pi'^*\left([\![\mathrm{FV}(\phi).\phi]\!]\right)$ and $(\pi \times 1_{Ms}) \circ \pi''^*\left([\![\vec{x}.\phi]\!]\right)$.

It is immediate to see that the $(*)$ square is a pullback by the properties of products and projections. Also, it is evident that $\pi'^*\left([\![\mathrm{FV}(\phi).\phi]\!]\right) \circ a^*$ is the pullback of $[\![\mathrm{FV}(\phi).\phi]\!] \circ a$ along $\pi'$, and, symmetrically, that $\pi''^*\left([\![\vec{x}.\phi]\!]\right) \circ b^*$ is the pullback of $[\![\vec{x}.\phi]\!] \circ b$ along $\pi''$.

Being $\mathbb{C}$ a Heyting category, it is also regular, so we can construct the images of $\pi \circ [\![\vec{x}.\phi]\!]$ and $(\pi \times 1_{Ms}) \circ \pi''^*\left([\![\vec{x}.\phi]\!]\right))$, and, by regularity, it holds that $\mathrm{im}\left((\pi \times 1_{Ms}) \circ \pi''^*\left([\![\vec{x}.\phi]\!]\right))\right)$ is the pullback of $\mathrm{im}\left(\pi \circ [\![\vec{x}.\phi]\!]\right)$ along $\pi'$.

By induction hypothesis, $[\![\vec{x}, y\colon s.\psi]\!]$ is the pullback of $[\![\mathrm{FV}(\phi), y\colon s.\psi]\!]$ along $\pi \times 1_{Ms}$, thus $h$ as in figure exists and makes the top-centre square a pullback.

But $(\pi \times 1_{Ms}) \circ \pi''^*\left([\![\vec{x}.\phi]\!]\right) = [\![\mathrm{FV}(\phi), y\colon s.\psi]\!] \circ h \circ \varepsilon$ so, by minimality of the image, $\mathrm{im}\left((\pi \times 1_{Ms}) \circ \pi''^*\left([\![\vec{x}.\phi]\!]\right)\right) \leq [\![\mathrm{FV}(\phi), y\colon s.\psi]\!]$.

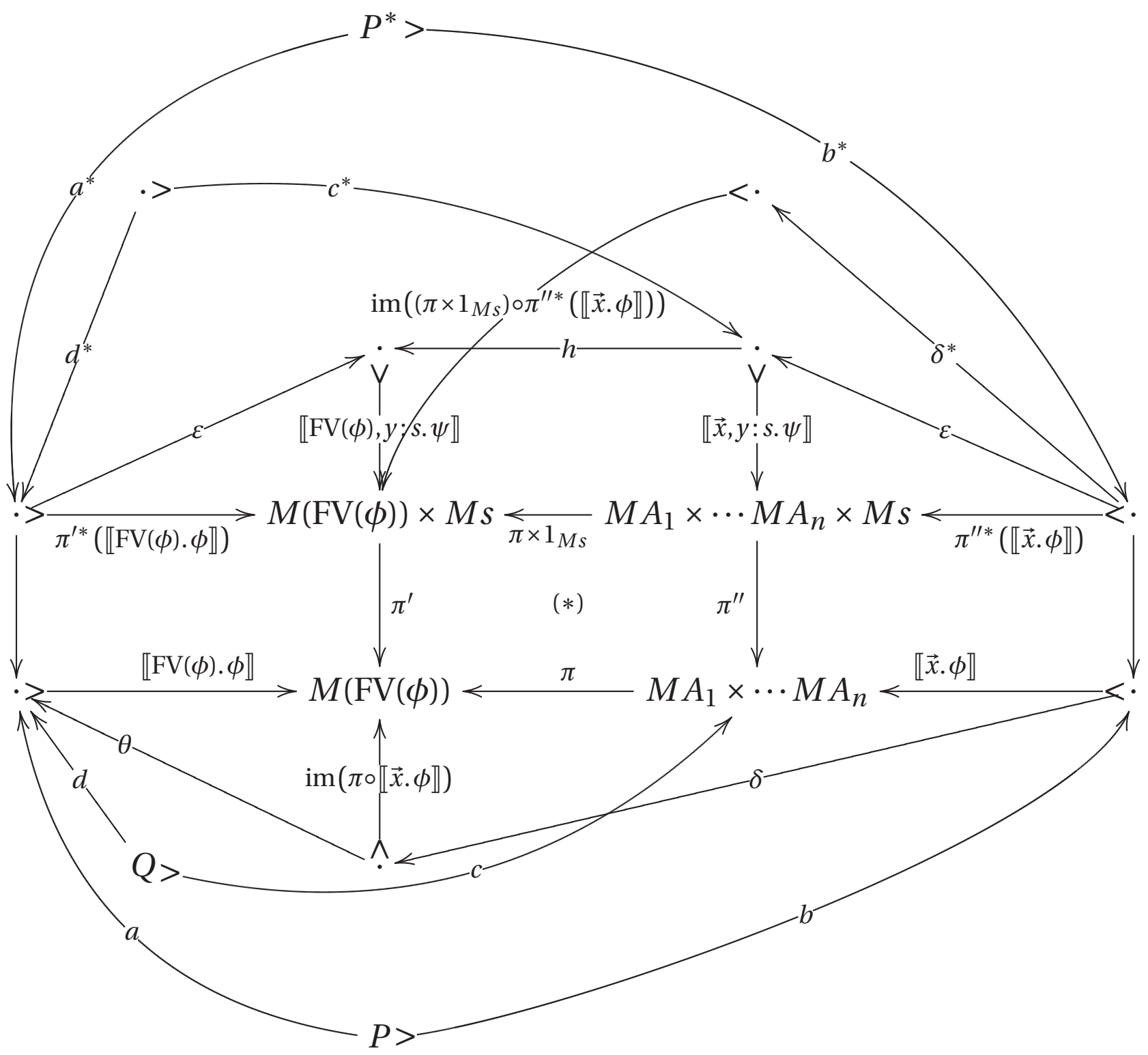

FIGURE 11. The universal case in proposition 8.1.

Since $\mathrm{im}\left(\pi\circ[\![\vec{x}.\phi]\!]\right)$ is a sub-object of $M(\mathrm{FV}(\phi))$ and its pullback along $\pi'$, $\mathrm{im}\left((\pi\times 1_{Ms})\circ\pi''^*\left([\![\vec{x}.\phi]\!]\right)\right)\le[\![\mathrm{FV}(\phi), y{:}s.\psi]\!]$, by right-adjointness of the functor $\forall_{\pi'}$, $\mathrm{im}\left(\pi\circ[\![\vec{x}.\phi]\!]\right)\le[\![\mathrm{FV}(\phi).\phi]\!]$, i.e., there is an arrow $\theta$ as in the diagram. Let $\alpha=\theta\circ\delta$.

Let $Q$ be the pullback of $[\![\mathrm{FV}(\phi).\phi]\!]$ and $\pi$, and let $c^*$ be the pullback of $c$ along $\pi''\circ[\![\vec{x}, y{:}s.\psi]\!]$. Since $Q$ is a pullback and the pair $\left(\alpha,[\![\vec{x}.\phi]\!]\right)$ makes $\left([\![\mathrm{FV}(\phi).\phi]\!],\pi\right)$ to commute, the universal arrow of $Q$ ensures that $[\![\vec{x}.\phi]\!]\le c$.

But $[\![\vec{x}, y{:}s.\psi]\!]\circ c^*\le[\![\vec{x}, y{:}s.\psi]\!]$ so, by right-adjointness of $\forall_{\pi''}$, $c\le[\![\vec{x}.\phi]\!]$. Thus, $c=[\![\vec{x}.\phi]\!]$ up to isomorphisms, that is, $[\![\vec{x}.\phi]\!]$ is the pullback of $[\![\mathrm{FV}(\phi).\phi]\!]$ along $\pi$, proving the statement. □

*Note* 8.1. Although proposition 8.1 as well as proposition 8.2, have not been explicitly stated in [Joh02b], they have been implicitly proved there. Nevertheless, having explicit, concrete proofs provides some insight on how Heyting categories are structured.

Next, we want to show that arrows between interpretations of formulae in some context are preserved by context restrictions.

**Proposition 8.2.** *Let $M$ be a $\Sigma$-structure in a Heyting category $\mathbb{C}$ and let $[\![\vec{x}.\phi]\!]\le[\![\vec{x}.\psi]\!]$ in $\mathrm{Sub}_{\mathbb{C}}(MA_1\times\cdots\times MA_n)$, where $\vec{x}.\phi$ and $\vec{x}.\psi$ are logical formulae in the context $\vec{x}\equiv x_1{:}A_1,\ldots,x_n{:}A_n$. Then, $[\![\mathrm{FV}(\phi)\cup\mathrm{FV}(\psi).\phi]\!]\le[\![\mathrm{FV}(\phi)\cup\mathrm{FV}(\psi).\psi]\!]$ in $\mathrm{Sub}_{\mathbb{C}}(\mathrm{FV}(\phi)\cup\mathrm{FV}(\psi))$.*

*Proof.* Let $\vec{z}=\mathrm{FV}(\phi)\cup\mathrm{FV}(\psi)$ and $\vec{y}=\vec{x}\setminus\vec{z}$. Calling $\pi\colon MA_1\times\cdots\times MA_n\to M\vec{z}$, $\pi_\phi\colon M\vec{z}\to M(\mathrm{FV}(\phi))$ and $\pi_\psi\colon M\vec{z}\to M(\mathrm{FV}(\psi))$ the evident projections, by proposition 8.1, $[\![\vec{z}.\phi]\!]$ is the pullback of $[\![\mathrm{FV}(\phi).\phi]\!]$ along $\pi_\phi$, $[\![\vec{z}.\psi]\!]$ is the pullback of $[\![\mathrm{FV}(\psi).\psi]\!]$ along $\pi_\psi$, $[\![\vec{x}.\phi]\!]$ is the pullback of $[\![\mathrm{FV}(\phi).\phi]\!]$ along $\pi_\phi\circ\pi$, and $[\![\vec{x}.\psi]\!]$ is the pullback of $[\![\mathrm{FV}(\psi).\psi]\!]$ along $\pi_\psi \circ \pi$. By the pullback lemma, $[\![\vec{x}.\phi]\!]$ is the pullback of $[\![\vec{z}.\phi]\!]$ along $\pi$, and $[\![\vec{x}.\psi]\!]$ is the pullback of $[\![\vec{z}.\psi]\!]$ along $\pi$

Since, when $\vec{y} = \emptyset$, the conclusion is trivial, we may safely assume that $\vec{y} \neq \emptyset$. We write $\forall \vec{y}.\theta$ for $\forall y_1 : s_1.(\forall y_2 : s_2.(\ldots(\forall y_m : s_m.\theta)\ldots))$ where $\vec{y} = y_1 : s_1, \ldots, y_m : s_m$.

Noticing that $\forall_\pi([\![\vec{x}.\psi]\!]) = [\![\vec{z}.\forall \vec{y}.\psi]\!]$ which exists, being $\mathbb{C}$ a Heyting category, we can draw the following commutative diagram:

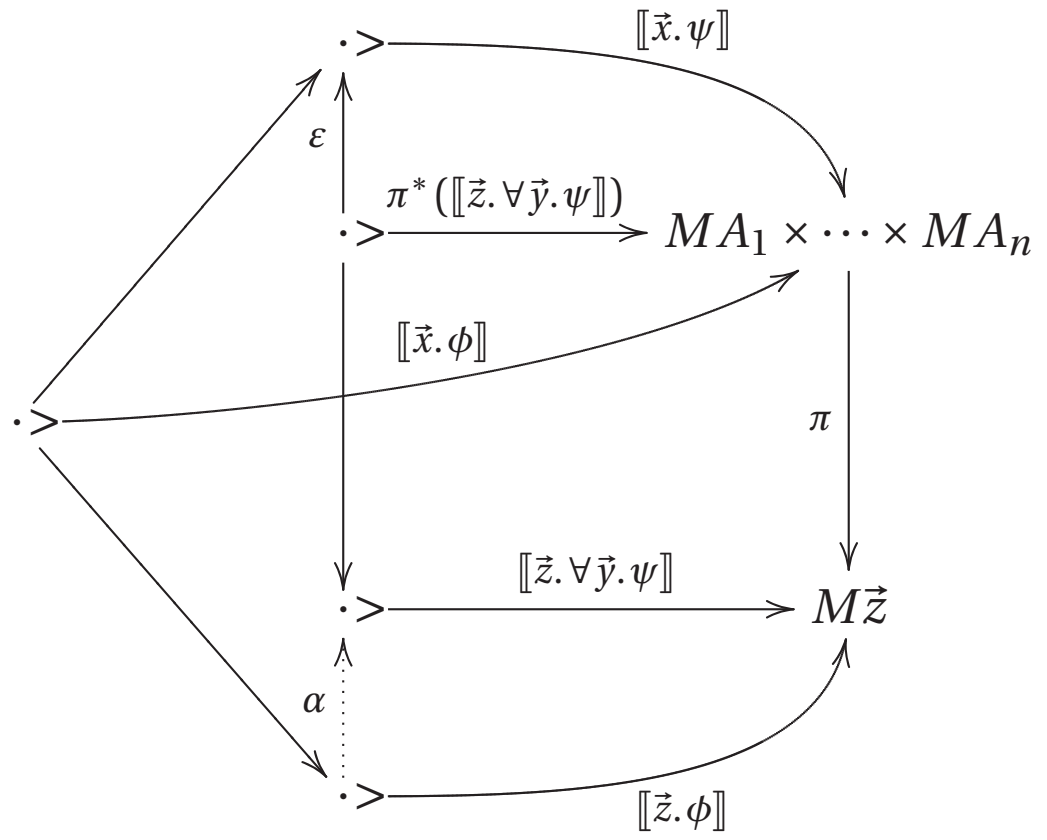

Being $\forall_\pi$ the right adjoint of the pullback functor $\pi^*$, there is a unique $\alpha$ making the above diagram to commute, thus $[\![\vec{z}.\phi]\!] \leq [\![\vec{z}.\forall\vec{y}.\psi]\!]$.

But the sequent $\forall\vec{y}.\psi \vdash_{\vec{z}} \psi$ holds in the empty theory, thus, by theorem 7.15 (soundness), it is valid in $M$, that is $[\![\vec{z}.\forall\vec{y}.\psi]\!] \leq [\![\vec{z}.\psi]\!]$. □

To every $\Sigma$-structure in a Heyting category, one can associate an essentially equivalent logical distributive category.

**Definition 8.1** (Associated logically distributive category)**.** Given $M$, a $\Sigma$-structure in a Heyting category $\mathbb{C}$, $\mathrm{LD}(M) = \langle \mathbb{D}, N\rangle$ is the *logically distributive category associated with* $M$, defined as follows (to lighten notation, we will identify logical formulae with the corresponding $\lambda$-types from now on):

(1) the objects of $\mathbb{D}$ are $[\![\mathrm{FV}(\phi).\phi]\!]_\mathbb{C}$ for every logical formula $\phi$ in the language over $\Sigma$;
(2) $\mathrm{Hom}_\mathbb{D}([\![\mathrm{FV}(\phi).\phi]\!]_\mathbb{C}, [\![\mathrm{FV}(\psi).\psi]\!]_\mathbb{C})$ contains at most one arrow, and this happens exactly when $[\![\mathrm{FV}(\phi)\cup\mathrm{FV}(\psi).\phi]\!]_\mathbb{C} \leq [\![\mathrm{FV}(\phi)\cup\mathrm{FV}(\psi).\psi]\!]_\mathbb{C}$ in the sub-object category $\mathrm{Sub}_\mathbb{C}(M(\mathrm{FV}(\phi)\cup\mathrm{FV}(\psi)))$;
(3) identities in $\mathbb{D}$ are the obvious arrows generated by reflexivity of $\leq$;
(4) composition in $\mathbb{D}$ reduces to transitivity of $\leq$;
(5) $N$ maps every $\lambda$-type $\phi$ to $[\![\mathrm{FV}(\phi).\phi]\!]_\mathbb{C}$.

*Note* 8.2. It is worth remarking that the map $N$ in the above definition, is surjective.

**Proposition 8.3.** *If $\phi \vdash_{\vec{x}} \psi$ is a sequent in the language over $\Sigma$ which is provable from the empty theory, i.e., in the pure logic, then, for every $M$ $\Sigma$-structure in a Heyting category $\mathbb{C}$, there is an arrow $N\phi \to N\psi$ in* $\mathrm{LD}(M) = \langle\mathbb{D}, N\rangle$.

*Proof.* By theorem 7.15, since $\phi \vdash_{\vec{x}} \psi$ is provable in the empty theory, $[\![\vec{x}.\phi]\!]_\mathbb{C} \leq [\![\vec{x}.\psi]\!]_\mathbb{C}$ in $\mathrm{Sub}_\mathbb{C}(MA_1 \times \cdots \times MA_n)$, where $\vec{x} \equiv x_1 : A_1, \ldots, x_n : A_n$. By proposition 8.2, it holds that $[\![\mathrm{FV}(\phi)\cup\mathrm{FV}(\psi).\phi]\!]_\mathbb{C} \leq [\![\mathrm{FV}(\phi)\cup\mathrm{FV}(\psi).\psi]\!]_\mathbb{C}$ in $\mathrm{Sub}_\mathbb{C}(M(\mathrm{FV}(\phi)) \cup M(\mathrm{FV}(\psi)))$. Thus, by definition 8.1, there is an arrow $N\phi \to N\psi$ in $\mathbb{D}$. □

**Proposition 8.4.** *For any $M$ $\Sigma$-structure in some Heyting category,* $\mathrm{LD}(M) = \langle\mathbb{D}, N\rangle$ *and* $\mathbb{D}$ *is a category.*

*Proof.* Evident from definition 8.1. □

**Proposition 8.5.** *For any $M$ $\Sigma$-structure in a Heyting category $\mathbb{C}$, $\mathrm{LD}(M) = \langle \mathbb{D}, N\rangle$ is such that $\mathbb{D}$ has finite products.*

*Proof.* It suffices to prove that $\mathbb{D}$ has a terminal object and binary products.

Since $\top$ is an axiom, the sequent $\phi \vdash_{\mathrm{FV}(\phi)} \top$ is trivially provable in the empty theory for any logical formula $\phi$. By proposition 8.3, there is an arrow $N\phi \to N\top$ in $\mathbb{D}$ which is necessarily unique, making $N\top$ the terminal object of $\mathbb{D}$.

Let $N\phi$ and $N\psi$ be objects of $\mathbb{D}$, and consider $\theta = \phi\wedge\psi$. Since the sequents $\theta \vdash_{\mathrm{FV}(\theta)} \phi$ and $\theta \vdash_{\mathrm{FV}(\theta)} \psi$ are immediately provable in the empty theory, by proposition 8.3, there are $N\theta \to N\phi$ and $N\theta \to N\psi$ in $\mathbb{D}$. Suppose $N\delta$ is an object of $\mathbb{D}$ such that $N\delta \to N\phi$ and $N\delta \to N\psi$ are arrows in $\mathbb{D}$. Then $[\![\mathrm{FV}(\delta)\cup\mathrm{FV}(\phi).\delta]\!]_{\mathbb{C}} \le [\![\mathrm{FV}(\delta)\cup\mathrm{FV}(\phi).\phi]\!]_{\mathbb{C}}$ in $\mathrm{Sub}_{\mathbb{C}}\big(M\big(\mathrm{FV}(\delta)\cup\mathrm{FV}(\phi)\big)\big)$, and similarly $[\![\mathrm{FV}(\delta)\cup\mathrm{FV}(\psi).\delta]\!]_{\mathbb{C}} \le [\![\mathrm{FV}(\delta)\cup\mathrm{FV}(\psi).\psi]\!]_{\mathbb{C}}$ in $\mathrm{Sub}_{\mathbb{C}}\big(M\big(\mathrm{FV}(\delta)\cup\mathrm{FV}(\psi)\big)\big)$, thus $[\![\mathrm{FV}(\delta)\cup\mathrm{FV}(\theta).\delta]\!]_{\mathbb{C}} \le [\![\mathrm{FV}(\delta)\cup\mathrm{FV}(\theta).\phi]\!]_{\mathbb{C}}$ and $[\![\mathrm{FV}(\delta)\cup\mathrm{FV}(\theta).\delta]\!]_{\mathbb{C}} \le [\![\mathrm{FV}(\delta)\cup\mathrm{FV}(\theta).\psi]\!]_{\mathbb{C}}$ in $\mathrm{Sub}_{\mathbb{C}}(M(\mathrm{FV}(\delta)\cup\mathrm{FV}(\theta)))$. By definition of meet, it holds that $[\![\mathrm{FV}(\delta)\cup\mathrm{FV}(\theta).\delta]\!]_{\mathbb{C}} \le [\![\mathrm{FV}(\delta)\cup\mathrm{FV}(\theta).\phi\wedge\psi]\!]_{\mathbb{C}}$ in $\mathrm{Sub}_{\mathbb{C}}(M(\mathrm{FV}(\delta)\cup\mathrm{FV}(\theta)))$, so, by definition of $\mathrm{LD}(M)$, $N\delta \to N\theta$ uniquely, making $N\theta$ the product $N\phi\times N\psi$ in $\mathbb{D}$. □

**Proposition 8.6.** *For any $M$ $\Sigma$-structure in a Heyting category $\mathbb{C}$, $\mathrm{LD}(M) = \langle \mathbb{D}, N\rangle$ is such that $\mathbb{D}$ ha finite co-products.*

*Proof.* It suffices to prove that $\mathbb{D}$ has an initial object and binary co-products.

The sequent $\bot \vdash_{\mathrm{FV}(\phi)} \phi$ is immediately provable in the empty theory for any logical formula $\phi$. By proposition 8.3, there is an arrow $N\bot \to N\phi$ in $\mathbb{D}$ which is necessarily unique, making $N\bot$ the initial object of $\mathbb{D}$.

Let $N\phi$ and $N\psi$ be objects of $\mathbb{D}$, and consider $\theta = \phi\vee\psi$. Since the sequents $\phi \vdash_{\mathrm{FV}(\theta)} \theta$ and $\psi \vdash_{\mathrm{FV}(\theta)} \theta$ are immediately provable in the empty theory, by proposition 8.3, there are $N\phi \to N\theta$ and $N\psi \to N\theta$ in $\mathbb{D}$. Suppose $N\delta$ is an object of $\mathbb{D}$ such that $N\phi \to N\delta$ and $N\psi \to N\delta$ are arrows in $\mathbb{D}$. Then $[\![\mathrm{FV}(\delta)\cup\mathrm{FV}(\phi).\phi]\!]_{\mathbb{C}} \le [\![\mathrm{FV}(\delta)\cup\mathrm{FV}(\phi).\delta]\!]_{\mathbb{C}}$ in $\mathrm{Sub}_{\mathbb{C}}\big(M\big(\mathrm{FV}(\delta)\cup\mathrm{FV}(\phi)\big)\big)$, and similarly $[\![\mathrm{FV}(\delta)\cup\mathrm{FV}(\psi).\psi]\!]_{\mathbb{C}} \le [\![\mathrm{FV}(\delta)\cup\mathrm{FV}(\psi).\delta]\!]_{\mathbb{C}}$ in $\mathrm{Sub}_{\mathbb{C}}\big(M\big(\mathrm{FV}(\delta)\cup\mathrm{FV}(\psi)\big)\big)$, thus $[\![\mathrm{FV}(\delta)\cup\mathrm{FV}(\theta).\phi]\!]_{\mathbb{C}} \le [\![\mathrm{FV}(\delta)\cup\mathrm{FV}(\theta).\delta]\!]_{\mathbb{C}}$ and $[\![\mathrm{FV}(\delta)\cup\mathrm{FV}(\theta).\psi]\!]_{\mathbb{C}} \le [\![\mathrm{FV}(\delta)\cup\mathrm{FV}(\theta).\delta]\!]_{\mathbb{C}}$ in $\mathrm{Sub}_{\mathbb{C}}(M(\mathrm{FV}(\delta)\cup\mathrm{FV}(\theta)))$. By definition of join, it holds that $[\![\mathrm{FV}(\delta)\cup\mathrm{FV}(\theta).\phi\vee\psi]\!]_{\mathbb{C}} \le [\![\mathrm{FV}(\delta)\cup\mathrm{FV}(\theta).\delta]\!]_{\mathbb{C}}$ in $\mathrm{Sub}_{\mathbb{C}}(M(\mathrm{FV}(\delta)\cup\mathrm{FV}(\theta)))$, so, by definition of $\mathrm{LD}(M)$, $N\theta \to N\delta$ uniquely, making $N\theta$ the co-product $N\phi+N\psi$ in $\mathbb{D}$. □

**Proposition 8.7.** *For any $M$ $\Sigma$-structure in a Heyting category $\mathbb{C}$, $\mathrm{LD}(M) = \langle \mathbb{D}, N\rangle$ is such that $\mathbb{D}$ has exponentiation.*

*Proof.* Let $N\phi$ and $N\psi$ be objects in $\mathbb{D}$ and consider $\theta = \phi\supset\psi$. Since $\phi\wedge\theta \vdash_{\mathrm{FV}(\theta)} \psi$ is immediately provable in the empty theory, by proposition 8.3 there is $N\big(\phi\wedge(\phi\supset\psi)\big) \to N\psi$ in $\mathbb{D}$. But $N\big(\phi\wedge(\phi\supset\psi)\big) = N\phi\times N(\phi\supset\psi)$ by proposition 8.5, so there is $N\phi\times N(\phi\supset\psi) \to N\psi$ in $\mathbb{D}$.

Let $N\delta \in \mathrm{Obj}\,\mathbb{D}$ such that $N\phi\times N\delta \to N\psi$ is an arrow of $\mathbb{D}$, that is, $N(\phi\wedge\delta) \to N\psi$. Then, $[\![\mathrm{FV}(\delta)\cup\mathrm{FV}(\theta).\phi\wedge\delta]\!]_{\mathbb{C}} \le [\![\mathrm{FV}(\delta)\cup\mathrm{FV}(\theta).\psi]\!]_{\mathbb{C}}$ in $\mathrm{Sub}_{\mathbb{C}}(M(\mathrm{FV}(\delta)\cup\mathrm{FV}(\theta)))$, so, by definition of implication in a Heyting algebra, $[\![\mathrm{FV}(\delta)\cup\mathrm{FV}(\theta).\delta]\!]_{\mathbb{C}} \le [\![\mathrm{FV}(\delta)\cup\mathrm{FV}(\theta).\phi\supset\psi]\!]_{\mathbb{C}}$ in $\mathrm{Sub}_{\mathbb{C}}(M(\mathrm{FV}(\delta)\cup\mathrm{FV}(\theta)))$. Thus, by definition of $\mathrm{LD}(M)$, there is a unique $N\delta \to N\theta$ in $\mathbb{D}$, proving that $N\theta = N\psi^{N\phi}$. □

**Proposition 8.8.** *For any $M$ $\Sigma$-structure in a Heyting category $\mathbb{C}$, $\mathrm{LD}(M) = \langle \mathbb{D}, N\rangle$ is such that $\mathbb{D}$ is a distributive category.*

*Proof.* Using the notation of definition 4.1, let $NA, NB, NC \in \mathrm{Obj}\,\mathbb{D}$ and let $\vec{x} \equiv \mathrm{FV}(A)\cup\mathrm{FV}(B)\cup\mathrm{FV}(C)$. By an easy proof, it holds that $A\wedge(B\vee C) \vdash_{\vec{x}} (A\wedge B)\vee(A\wedge C)$ is provable in the empty theory. By propositions 8.5, 8.6 and 8.3, $NA\times(NB+BC) \to (NA\times NB)+(NA\times NC)$ is an arrow in $\mathbb{D}$, which is the inverse of $\Delta$ by uniqueness. □

**Proposition 8.9.** *For any $M$ $\Sigma$-structure in a Heyting category $\mathbb{C}$, $\mathrm{LD}(M) = \langle \mathbb{D}, N\rangle$ is such that all the sub-categories $\mathbb{D}_{(\forall x:s.\,A)}$ have terminal objects.*

*Proof.* Because $\forall x\colon s.\, A \vdash_{\vec{x}} A[t/x]$, where $\vec{x} \equiv \mathrm{FV}(\forall x\colon s.\, A) \cup \mathrm{FV}(t\colon s)$, is provable in the empty theory for any logical formula $A$ over $\Sigma$ and any $t\colon s \in \mathrm{LTerms}(\Sigma)$, from proposition 8.3, there is $N(\forall x\colon s.\, A) \to NA[t/x] = \Sigma_A(x\colon s)(t\colon s)$ in $\mathbb{D}$, showing that $N(\forall x\colon s.\, A)$ is the vertex of a cone on $\Sigma_A(x\colon s)$, i.e., $N(\forall x\colon s.\, A) \in \mathrm{Obj}\,\mathbb{D}_{\forall x\colon s.\, A}$.

Let $NB \in \mathrm{Obj}\,\mathbb{D}_{\forall x\colon s.\, A}$: by definition of the sub-category, $B$ is a logical formula over $\Sigma$ such that $x\colon s \notin \mathrm{FV}(B)$ and, for each $t\colon s \in \mathrm{LTerms}(\Sigma)$, $NB \to NA[t/x]$ is an arrow in $\mathbb{D}$. Then, calling $\vec{y} \equiv \mathrm{FV}(B) \cup \mathrm{FV}(\forall x\colon s.\, A) \cup \mathrm{FV}(t\colon s)$, it holds that $[\![\vec{y}.\, B]\!]_{\mathbb{C}} \le [\![\vec{y}.\, A[t/x]]\!]_{\mathbb{C}}$ in $\mathrm{Sub}_{\mathbb{C}}\,(M\,(\mathrm{FV}(B) \cup \mathrm{FV}(\forall x\colon s.\, A) \cup \mathrm{FV}(t\colon s)))$. Thus, by definition, $B \vdash_{\vec{y}} A[t/x]$ for any $t\colon s \in \mathrm{LTerms}(\Sigma)$. In particular, for $x\colon s$, $B \vdash_{\vec{y}} A$ and, since $x\colon s \notin \mathrm{FV}(B)$, $B \vdash_{\vec{y}} \forall x\colon s.\, A$ holds in the empty theory. Hence, by proposition 8.3, there is a unique $NB \to N(\forall x\colon s.\, A)$ in $\mathbb{D}$, that is, $N(\forall x\colon s.\, A)$ is the terminal object of $\mathbb{D}_{\forall x\colon s.\, A}$. □

**Proposition 8.10.** *For any $M$ $\Sigma$-structure in a Heyting category $\mathbb{C}$, $\mathrm{LD}(M) = \langle \mathbb{D}, N\rangle$ is such that all the sub-categories $\mathbb{D}_{(\exists x\colon s.\, A)}$ have initial objects.*

*Proof.* Because $A[t/x] \vdash_{\vec{x}} \exists x\colon s.\, A$, where $\vec{x} \equiv \mathrm{FV}(\exists x\colon s.\, A) \cup \mathrm{FV}(t\colon s)$, is provable in the empty theory for any logical formula $A$ over $\Sigma$ and any $t\colon s \in \mathrm{LTerms}(\Sigma)$, by proposition 8.3, there is $\Sigma_A(x\colon s)(t\colon s) = NA[t/x] \to N(\exists x\colon s.\, A)$ in $\mathbb{D}$, showing that $N(\exists x\colon s.\, A)$ is the vertex of a co-cone on $\Sigma_A(x\colon s)$, i.e., $N(\exists x\colon s.\, A) \in \mathrm{Obj}\,\mathbb{D}_{\exists x\colon s.\, A}$.

Let $NB \in \mathrm{Obj}\,\mathbb{D}_{\exists x\colon s.\, A}$: by definition of the sub-category, $B$ is a logical formula over $\Sigma$ such that $x\colon s \notin \mathrm{FV}(B)$ and, for each $t\colon s \in \mathrm{LTerms}(\Sigma)$, $NA[t/x] \to NB$ is an arrow in $\mathbb{D}$. Then, calling $\vec{y} \equiv \mathrm{FV}(B) \cup \mathrm{FV}(\exists x\colon s.\, A) \cup \mathrm{FV}(t\colon s)$, it holds that $[\![\vec{y}.\, A[t/x]]\!]_{\mathbb{C}} \le [\![\vec{y}.\, B]\!]_{\mathbb{C}}$ in $\mathrm{Sub}_{\mathbb{C}}\,(M\,(\mathrm{FV}(B) \cup \mathrm{FV}(\exists x\colon s.\, A) \cup \mathrm{FV}(t\colon s)))$. Thus, by definition, $A[t/x] \vdash_{\vec{y}} B$ for any $t\colon s \in \mathrm{LTerms}(\Sigma)$. In particular, for $x\colon s$, $A \vdash_{\vec{y}} B$ and, since $x\colon s \notin \mathrm{FV}(B)$, $\exists x\colon s.\, A \vdash_{\vec{y}} B$ holds in the empty theory. Hence, by proposition 8.3, there is a unique $N(\exists x\colon s.\, A) \to NB$ in $\mathbb{D}$, that is, $N(\exists x\colon s.\, A)$ is the initial object of $\mathbb{D}_{\exists x\colon s.\, A}$. □

**Proposition 8.11.** *For any $M$ $\Sigma$-structure in a Heyting category $\mathbb{C}$, $\mathrm{LD}(M) = \langle \mathbb{D}, N\rangle$ is such that $N$ satisfies the condition expressed in point* (6) *of definition* 4.1.

*Proof.* By propositions 8.5 to 8.10, the claim is evident. □

**Proposition 8.12.** *For any $M$ $\Sigma$-structure in a Heyting category $\mathbb{C}$, $\mathrm{LD}(M) = \langle \mathbb{D}, N\rangle$ is such that the unique map $N(\exists x\colon s.\, A \times B) \to NA \times N(\exists x\colon s.\, B)$ in $\mathbb{D}$, where $x\colon s \notin \mathrm{FV}(A)$, has an inverse.*

*Proof.* Since the sequent $A \wedge \exists x\colon s.\, B \vdash_{\vec{x}} \exists x\colon s.\, A \wedge B$, where $\vec{x} \equiv \mathrm{FV}(A) \cup \mathrm{FV}(\exists x\colon s.\, B)$ and $x\colon s \notin \mathrm{FV}(A)$, holds in the empty theory, by propositions 8.3 and 8.5, there is an arrow $NA \times N(\exists x\colon s.\, B) \to N(\exists x\colon s.\, A \times B)$ in $\mathbb{D}$. By uniqueness of this arrow, it is immediate to see that it is the inverse of $N(\exists x\colon s.\, A \times B) \to NA \times N(\exists x\colon s.\, B)$. □

**Proposition 8.13.** *For any $M$ $\Sigma$-structure in a Heyting category $\mathbb{C}$, $\mathrm{LD}(M) = \langle \mathbb{D}, N\rangle$ is a logically distributive category.*

*Proof.* Obvious by propositions 8.4 to 8.12. □

When we say that $\mathrm{LD}(M)$ is essentially equivalent to $M$, we mean that these structures validates the same sequents which can be expressed in the language of logic.

**Proposition 8.14.** *Let $M$ be a $\Sigma$-structure in a Heyting category $\mathbb{C}$. If $\phi \vdash_{\vec{x}} \psi$ is a sequent in the language over $\Sigma$, and $\phi \vdash_{\vec{x}} \psi$ is valid in $M$, then $\phi \vdash_{\vec{x}} \psi$ is valid is* $\mathrm{LD}(M)$.

*Proof.* By definition, $\phi \vdash_{\vec{x}} \psi$ is valid in $M$ if and only if $[\![\vec{x}.\, \phi]\!]_M \le [\![\vec{x}.\, \psi]\!]_M$ in $\mathrm{Sub}_{\mathbb{C}}(MA_1 \times \cdots \times MA_n)$, where $\vec{x} \equiv x_1\colon A_1, \ldots, x_n\colon A_n$. So, by proposition 8.2, $[\![\mathrm{FV}(\phi) \cup \mathrm{FV}(\psi).\, \phi]\!]_M \le [\![\mathrm{FV}(\phi) \cup \mathrm{FV}(\psi).\, \psi]\!]_M$ in $\mathrm{Sub}_{\mathbb{C}}\left(M\left(\mathrm{FV}(\phi) \cup \mathrm{FV}(\psi)\right)\right)$. Calling $\mathrm{LD}(M) = \langle \mathbb{D}, N\rangle$, by definition of $\mathrm{LD}(M)$, there is an arrow $N\phi \to N\psi$, that is, the sequent $\phi \vdash_{\vec{x}} \psi$ is valid in $\mathbb{D}$. □

**Corollary 8.15.** *Let $M$ be a $\Sigma$-structure in a Heyting category $\mathbb{C}$. If $\phi \vdash_{\vec{x}} \psi$ is a sequent in the language over $\Sigma$, then $\phi \vdash_{\vec{x}} \psi$ is valid in $M$ if and only if $\phi \vdash_{\vec{x}} \psi$ is valid is* $\mathrm{LD}(M)$.

*Proof.* On one side this is just proposition 8.14. The other side is immediate noticing that the proof of proposition 8.14 is reversible. □

Since we want to consider LD as a functor, so to relate models in Heyting categories with models in logically distributive categories, we need to define appropriate categories grouping models, and the morphisms between models as well.

**Definition 8.2** (Σ-H-homomorphism)**.** Let $\Sigma = \langle S, F, R\rangle$ be a logical signature. Given a Σ-H-structure $M$ in a Heyting category $\mathbb{C}_M$, and a Σ-H-structure $N$ in a Heyting category $\mathbb{C}_N$, a *Σ-H-homomorphism* $h\colon M \to N$ is a Heyting functor $h\colon \mathbb{C}_M \to \mathbb{C}_N$ such that $h(\llbracket \vec{x}.\phi \rrbracket_M) = \llbracket \vec{x}.\phi \rrbracket_N$ for every formula-in-context $\vec{x}.\phi$.

*Note* 8.3. This definition differs from [Joh02b] where homomorphisms between Σ-H-structures are considered only inside the same Heyting category. Also, the definition does not reduce to the one in [Joh02b] when restricted. We believe that the above definition is more natural than the standard one in literature, and much well-behaved.

**Definition 8.3.** Fixed a Heyting category $\mathbb{C}$ and a logical signature Σ, Σ-H-Str($\mathbb{C}$) is the category whose objects are the Σ-H-structures in $\mathbb{C}$, and whose arrows are the Σ-H-homomorphisms between objects.

**Proposition 8.16.** *If* $\mathbb{C}$ *is a small Heyting category, so is* Σ-H-Str($\mathbb{C}$).

*Proof.* Elementary cardinal arithmetic provides a bound to the size of Σ-H-Str($\mathbb{C}$). □

**Definition 8.4.** Given a $\lambda$-signature Σ, Σ-Str is the category whose objects are the logically distributive categories $\langle \mathbb{D}, M\rangle$ such that $\mathrm{Obj}\,\mathbb{D} \subseteq \lambda\mathrm{Types}(\Sigma)$ and $M$ is surjective, and whose arrows are the Σ-homomorphisms between objects.

**Proposition 8.17.** *The* Σ-Str *category is small.*

*Proof.* Elementary cardinal arithmetic provides a bound to the size of it. □

**Proposition 8.18.** *For each* Σ*-H-structure* $M$ *in a Heyting category,* $\mathrm{LD}(M) \in \Sigma$-Str.

*Proof.* Evident by definition of LD($M$). □

**Proposition 8.19.** *If* $T$ *is a* $\lambda$*-theory over the signature* Σ*, and* $\mathbb{C}_T$ *is the syntactical category, see definition* 6.2*, then* $\mathbb{C}_t \in \Sigma$-Str.

*Proof.* Evident by definition of $\mathbb{C}_T$. □

**Definition 8.5** (LD functor)**.** Given a Heyting category $\mathbb{C}$ and a $\lambda$-signature Σ, the functor LD extends the map LD, see definition 8.1, as follows: LD: Σ-H-Str($\mathbb{C}$) → Σ-Str and, for each Σ-H-structure $M$ in $\mathbb{C}$, LD($M$) is defined as previously seen, while, for any Σ-H-homomorphism $h\colon M \to N$ between the Σ-H-structures $M$ and $N$ on $\mathbb{C}$, the Σ-homomorphism $\mathrm{LD}(h)\colon \mathrm{LD}(M) = \langle \mathbb{D}, K\rangle \to \mathrm{LD}(N) = \langle \mathbb{E}, H\rangle$ is defined as

$$\mathrm{LD}(h)\left(\llbracket \mathrm{FV}(\phi).\phi \rrbracket_M\right) = h\left(\llbracket \mathrm{FV}(\phi).\phi \rrbracket_M\right)$$

on the objects of LD($M$), while

$$\begin{aligned}\mathrm{LD}(h)\left(\llbracket \mathrm{FV}(\phi)\cup\mathrm{FV}(\psi).\phi \rrbracket_M \le \llbracket \mathrm{FV}(\phi)\cup\mathrm{FV}(\psi).\psi \rrbracket_M\right) &= \\ \left(h\left(\llbracket \mathrm{FV}(\phi)\cup\mathrm{FV}(\psi).\phi \rrbracket_M\right) \le h\left(\llbracket \mathrm{FV}(\phi)\cup\mathrm{FV}(\psi).\psi \rrbracket_M\right)\right)&\end{aligned}$$

on the arrows of LD($M$).

**Proposition 8.20.** LD *is a functor.*

*Proof.* By proposition 8.13, $\mathrm{LD}(M)$ is a logically distributive category and, by proposition 8.18, $\mathrm{LD}(M) \in \Sigma\text{-Str}$, making sound the definition of LD on objects.

Since every object of $\mathrm{LD}(M) = \langle \mathbb{D}, K \rangle$ has the form $\llbracket \mathrm{FV}(\phi).\,\phi \rrbracket_M$, and because $h$ is a $\Sigma$-H-homomorphism, $\mathrm{LD}\,(h)\left(\llbracket \mathrm{FV}(\phi).\,\phi \rrbracket_M\right) = h\left(\llbracket \mathrm{FV}(\phi).\,\phi \rrbracket_M\right) = \llbracket \mathrm{FV}(\phi).\,\phi \rrbracket_N$, satisfying the condition for $\mathrm{LD}(h)$ being a $\Sigma$-homomorphism.

Let $\vec{x} = \mathrm{FV}(\phi) \cup \mathrm{FV}(\psi)$. Since all the arrows of $\mathbb{D}$ have the form $f \colon \llbracket \vec{x}.\phi \rrbracket_M \leq \llbracket \vec{X}.\,\psi \rrbracket_M$, and $\mathrm{LD}(h)(f) = \left(h\left(\llbracket \vec{x}.\,\phi \rrbracket_M\right) \leq h\left(\llbracket \vec{x}.\,\psi \rrbracket_M\right)\right) = \left(\llbracket \vec{x}.\,\phi \rrbracket_N \leq \llbracket \vec{x}.\,\psi \rrbracket_N\right)$, LD correctly maps arrows into arrows, preserving domains and co-domains.

An immediate expansion of definitions shows that identities and composition are preserved by $\mathrm{LD}(h)$. Thus, $\mathrm{LD}(h)$ is a $\Sigma$-homomorphism, as claimed in definition 8.5.

Finally, since $1_M$ is a $\Sigma$-H-homomorphism, $\mathrm{LD}(1_M)$ immediately equals $1_{\mathrm{LD}(M)}$, the identity of $\mathrm{LD}(M)$, which is a $\Sigma$-homomorphism. Also, $\mathrm{LD}(f \circ g) = \mathrm{LD}(f) \circ \mathrm{LD}(g)$ as this reduces to composition of functors preserving interpretations. □

Combining proposition 8.20 with corollary 8.15, one obtains that the LD functors generate validity preserving maps from generic $\Sigma$-H-structures in any Heyting category, to some $\Sigma$-structure in a logically distributive category. Moreover, these maps also reflects validity. These facts show that logically distributive categories are at least as expressive as Heyting categories when interpreting first-order logic. Of course, being both semantics sound and complete, their expressive power 'inside' the logical language is the same; but it seems that the ability to model pieces of the mathematical world, is wider when they are expressed as logically distributive categories, rather than as Heyting categories—as the results in this section show that any model in a Heyting category can be converted into a model in a logically distributed category, while the converse is still an open question.

*Note* 8.4. One may want to 'invert' the LD functor, e.g., by trying to find a left adjoint to it. This seems to be impossible as the only natural candidate, the category whose arrows are the interpretation of terms and formulae in contexts, is not Cartesian. For example, the pullback of $\llbracket x : A.\,\phi \rrbracket$ along $\llbracket \vec{y}.\, t : A \rrbracket$, where $\phi$ is any formula in the context $x : A$ and $t : B$ is any closed term,

$$\begin{array}{ccc} \cdot & \xrightarrow{\llbracket x:A.\,\phi \rrbracket} & MA \\ \uparrow & & \uparrow{\scriptstyle \llbracket \vec{y}.\,t:A \rrbracket} \\ P & \xrightarrow[?]{} & M(\mathrm{FV}(\vec{y})) \end{array}$$

has no logical meaning, thus it does not appear as an object of $\mathrm{LD}(K)$ for any $K$.

## 9. Models in a topos

In this section, we introduce the topos-theoretic interpretation of first-order intuitionistic logic. Although the topic is well-known, the basic definitions and results are stated, following [Gol06]. This material is in preparation to the following section, where this semantics is compared with the one based on logically distributive categories.

**Definition 9.1** (Model in a topos)**.** Let $\Sigma = \langle S, F, R\rangle$ be a logical signature and let $\mathscr{E}$ be an elementary topos. A *$\mathscr{E}$-model* on $\Sigma$ is a structure

$$\begin{aligned}\mathscr{U} = \langle &\{A_s \in \mathrm{Obj}\,\mathscr{E} \colon \mathrm{Hom}_{\mathscr{E}}(1, A_s) \neq \emptyset\}_{s\in S},\\ &\{f_{\mathscr{E}} \in \mathrm{Hom}_{\mathscr{E}}(MA_1 \times \cdots \times MA_n, MA_0)\}_{f\colon A_1\times\cdots\times A_n \to A_0 \in F},\\ &\{r_{\mathscr{E}} \in \mathrm{Hom}_{\mathscr{E}}(MA_1 \times \cdots \times MA_n, \Omega)\}_{r\colon A_1\times\cdots\times A_n \in R}\rangle\end{aligned}$$

together with the interpretation $M$ such that

(1) $Ms = A_s$ for each $s \in S$;
(2) $Mf = f_{\mathscr{E}}$ for each $f \in F$;
(3) $Mr = r_{\mathscr{E}}$ for each $r \in R$.

**Definition 9.2** (Interpretation of terms)**.** Given a $\mathscr{E}$-model $\mathscr{U}$ with interpretation $M$ over $\Sigma = \langle S, F, R\rangle$, a signature, a logical term-in-context $\vec{x}.t\colon s$, with $\vec{x} \equiv x_1\colon A_1, \ldots, x_n\colon A_n$, is interpreted as an arrow $M(t\colon s)\colon MA_1 \times \cdots \times MA_n \to Ms$ in $\mathscr{E}$ according to the following inductive definition:

(1) $t\colon s \equiv x_i\colon s_i$, a variable: $M(t\colon s) = \pi_i\colon MA_1 \times \cdots \times MA_n \to MA_i$, the canonical $i$-th projection;
(2) $t\colon s \equiv f(t_1\colon s_1, \ldots, t_m\colon s_m)\colon s_0$, where $f\colon s_1 \times \cdots \times s_m \to s_0 \in F$ and $t_1\colon s_1, \ldots, t_m\colon s_m \in \mathrm{LTerms}(\Sigma)$: $M(t\colon s) = Mf \circ (M(t_1\colon s_1), \ldots, M(t_m\colon s_m))\colon MA_1 \times \cdots \times MA_n \to Ms_0$.

**Definition 9.3.** Fixed a topos $\mathscr{E}$, and $a \in \mathrm{Obj}\,\mathscr{E}$,

(1) $\top_a = \top \circ !\colon a \to \Omega$ where $!\colon a \to 1$ and $\top\colon 1 \rightarrowtail \Omega$ is the truth arrow;
(2) $\bot_a = \bot \circ !\colon a \to \Omega$ where $!\colon a \to 1$ and $\bot\colon 1 \to \Omega$ is the character of $0 \rightarrowtail 1$:

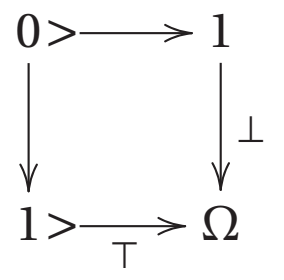

(3) $\cap\colon \Omega^2 \to \Omega$ is the character of $(\top, \top)$:

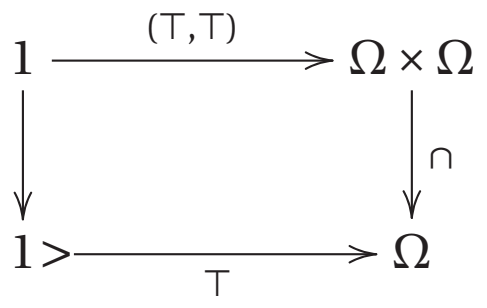

(4) $\cup\colon \Omega^2 \to \Omega$ is the character of the image of $[(\top_\Omega, 1_\Omega), (1_\Omega, \top_\Omega)]\colon \Omega + \Omega \to \Omega^2$:

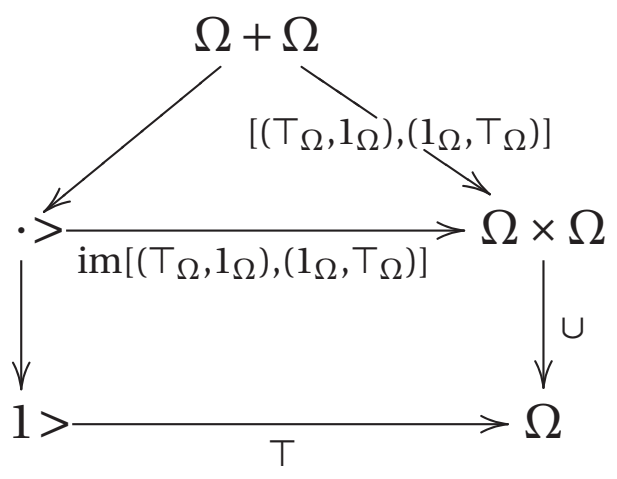

(5) $\Rightarrow\colon \Omega^2 \to \Omega$ is the character of the equaliser of $\cap$ and $\pi_1\colon \Omega^2 \to \Omega$:

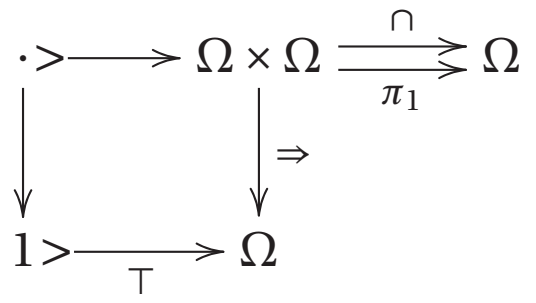

(6) $\forall_a\colon \Omega^a \to \Omega$ is the character of the exponential transpose of $\top_a \circ \pi_2\colon 1 \times a \to \Omega$:

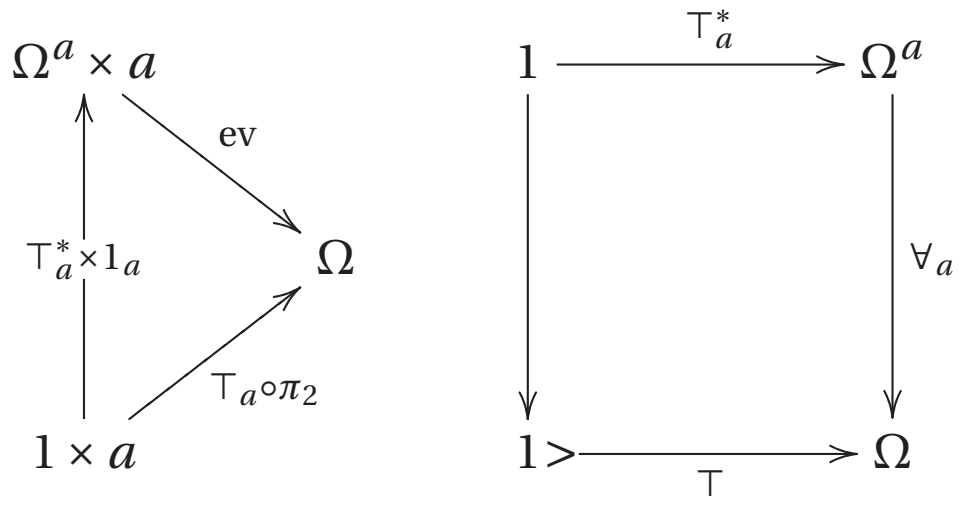

(7) $\exists_a\colon \Omega^a \to \Omega$ is the character of $\mathrm{im}(\pi_1 \circ \mathrm{ev}_a)$, where $\mathrm{ev}_a$ is the sub-object of $\Omega^a \times a$ whose character is $\mathrm{ev}\colon \Omega^a \times a \to \Omega$:

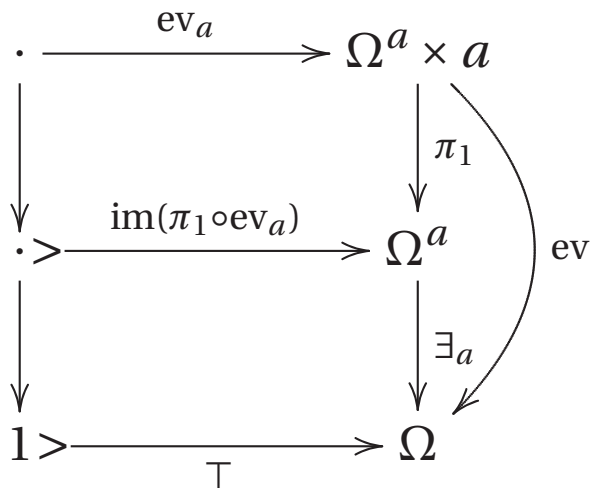

**Definition 9.4** (Interpretation of formulae)**.** Given a $\mathscr{E}$-model $\mathscr{U}$ on $\Sigma$, a logical formula-in-context $\vec{x}.\phi$, with $\vec{x} \equiv x_1 : A_1, \ldots, x_n : A_n$, is interpreted according to the following inductive definition:

(1) $\phi \equiv \top$: $[\![\vec{x}.\phi]\!] = \top_{MA_1 \times \cdots \times MA_n}$;
(2) $\phi \equiv \bot$: $[\![\vec{x}.\phi]\!] = \bot_{MA_1 \times \cdots \times MA_n}$;
(3) $\phi \equiv r(t_1, \ldots, t_m)$, where $r : B_1 \times \cdots \times B_m \in R$ and $t_1 : B_1, \ldots, t_m : B_m \in \mathrm{LTerms}(\Sigma)$: $[\![\vec{x}.\phi]\!] = Mr \circ (M(t_1 : B_1), \ldots, M(t_m : B_m))$;
(4) $\phi \equiv \psi_1 \wedge \psi_2$: $[\![\vec{x}.\phi]\!] = \cap \circ \big([\![\vec{x}.\psi_1]\!], [\![\vec{x}.\psi_2]\!]\big)$;
(5) $\phi \equiv \psi_1 \vee \psi_2$: $[\![\vec{x}.\phi]\!] = \cup \circ \big([\![\vec{x}.\psi_1]\!], [\![\vec{x}.\psi_2]\!]\big)$;
(6) $\phi \equiv \psi_1 \supset \psi_2$: $[\![\vec{x}.\phi]\!] = \Rightarrow \circ \big([\![\vec{x}.\psi_1]\!], [\![\vec{x}.\psi_2]\!]\big)$;
(7) $\phi \equiv \exists y : s.\psi$: $[\![\vec{x}.\phi]\!] = \exists_s \circ \theta$ with $\theta$ the exponential transpose of $[\![\vec{x}, y : s.\psi]\!]\colon MA_1 \times \cdots MA_n \times Ms \to \Omega$;
(8) $\phi \equiv \forall y : s.\psi$: $[\![\vec{x}.\phi]\!] = \forall_s \circ \theta$ with $\theta$ the exponential transpose of $[\![\vec{x}, y : s.\psi]\!]$.

We say that $\mathscr{U}$ *validates* the logical formula $\phi$, notation $\mathscr{U} \models_{\mathscr{E}} \phi$, when $[\![\mathrm{FV}(\phi).\phi]\!] = \top_{M(\mathrm{FV}(\phi))}$. A formula $\phi$ is $\mathscr{E}$*-valid*, notation $\mathscr{E} \models \phi$, if $\mathscr{U} \models_{\mathscr{E}} \phi$ for every $\mathscr{E}$-model $\mathscr{U}$.

**Proposition 9.1.** *Let $\phi$ be a logical formula over the signature $\Sigma$, and let $\mathscr{U}$ be a $\mathscr{E}$-model on $\Sigma$ in the topos $\mathscr{E}$. Then, for every context $\vec{x} \equiv x_1 : A_1, \ldots, x_n : A_n$ such that $\mathrm{FV}(\phi) \subseteq \vec{x}$, $[\![\vec{x}.\phi]\!]\colon MA_1 \times \cdots \times MA_n \to \Omega$. It holds that $[\![\vec{x}.\phi]\!] = [\![\mathrm{FV}(\phi).\phi]\!] \circ \pi$ where $\pi\colon MA_1 \times \cdots \times MA_n \to M(\mathrm{FV}(\phi))$ is the canonical projection.*

*Proof.* See [Gol06, p. 247]. □

**Theorem 9.2** (Soundness and completeness)**.** *If $\phi$ is a logical formula provable in first-order intuitionistic logic, then $\mathscr{E} \models \phi$ for any elementary topos $\mathscr{E}$, and vice versa.*

*Proof.* See [Gol06, p. 249] and [Gol06, p. 265] □

## 10. TOPOI VERSUS LOGICALLY DISTRIBUTIVE CATEGORIES

Analogously to section 8, here we present a way to embed topos-based model in suitable logically distributive categories in such a way that validity is preserved and reflected. As a rule of thumb, the development of these embedding results can be pursued along the lines already followed in section 8 and hereafter remarked.

**Definition 10.1** (Associated logically distributive category)**.** Given $\mathcal{U}$, a $\mathcal{E}$-model on $\Sigma$ in the topos $\mathcal{E}$, $\mathrm{LD}(\mathcal{U}) = \langle \mathbb{D}, N\rangle$ is the *logically distributive category associated with* $\mathcal{U}$, defined as:

(1) the objects of $\mathbb{D}$ are $[\![\mathrm{FV}(\phi).\phi]\!]_{\mathcal{U}}$ for every formula $\phi$ on $\Sigma$;
(2) $\mathrm{Hom}_{\mathbb{D}}\left([\![\mathrm{FV}(\phi).\phi]\!]_{\mathcal{U}}, [\![\mathrm{FV}(\psi).\psi]\!]_{\mathcal{U}}\right)$ contains at most one arrow, and this happens precisely when

$$\begin{aligned}&[\![\mathrm{FV}(\phi)\cup\mathrm{FV}(\psi).\phi\supset\psi]\!]_{\mathcal{U}} =\\ &= \Rightarrow\circ\left([\![\mathrm{FV}(\phi)\cup\mathrm{FV}(\psi).\phi]\!]_{\mathcal{U}}, [\![\mathrm{FV}(\phi)\cup\mathrm{FV}(\psi).\psi]\!]_{\mathcal{U}}\right) =\\ &= \top_{M(\mathrm{FV}(\phi)\cup\mathrm{FV}(\psi))} \quad ,\end{aligned}$$

or, equivalently, when $\mathcal{U}\models_{\mathcal{E}}\phi\supset\psi$, or, in other terms, when $\mathcal{U}$ validates the sequent $\phi\vdash_{\mathrm{FV}(\phi)\cup\mathrm{FV}(\psi)}\psi$;
(3) identities are the evident arrows given by $\mathcal{U}\models_{\mathcal{E}}\phi\supset\phi$;
(4) given $f\colon [\![\mathrm{FV}(\phi).\phi]\!]_{\mathcal{U}}\to[\![\mathrm{FV}(\psi).\psi]\!]_{\mathcal{U}}$ and $g\colon[\![\mathrm{FV}(\psi).\psi]\!]_{\mathcal{U}}\to[\![\mathrm{FV}(\theta).\theta]\!]_{\mathcal{U}}$ in $\mathbb{D}$, $g\circ f$ is defined by $\mathcal{U}\models_{\mathcal{E}}\phi\supset\theta$;
(5) $N$ maps every $\lambda$-type $\phi$ to $[\![\mathrm{FV}(\phi).\phi]\!]_{\mathcal{U}}$.

*Note* 10.1. The equivalent formulations of clause (2) are immediate consequences of definition 9.4 and theorem 9.2.

**Proposition 10.1.** *The sequent $\phi\vdash_{\vec{x}}\psi$ on the logical signature $\Sigma$ is provable in the empty theory if and only if, for each $\mathcal{E}$-model $\mathcal{U}$ on $\Sigma$ in a topos $\mathcal{E}$, there is an arrow $N\phi\to N\psi$ in* $\mathrm{LD}(\mathcal{U}) = \langle\mathbb{D},N\rangle$.

*Proof.* When $\phi\vdash_{\vec{x}}\psi$ is provable, $\emptyset\vdash_{\vec{x}}\phi\supset\psi$ is provable as well. By theorem 9.2, this is equivalent to say that $\mathcal{U}\models_{\mathcal{E}}\phi\supset\psi$, i.e., $N\phi\to N\psi$ in $\mathbb{D}$. □

**Proposition 10.2.** *For any $\mathcal{U}$, $\mathcal{E}$-model on $\Sigma$ in some topos $\mathcal{E}$,* $\mathrm{LD}(\mathcal{U}) = \langle\mathbb{D},N\rangle$ *where $\mathbb{D}$ is a category.*

*Proof.* Immediate, by unfolding definitions. □

**Proposition 10.3.** *For any $\mathcal{U}$, $\mathcal{E}$-model on $\Sigma$ in some topos $\mathcal{E}$,* $\mathrm{LD}(\mathcal{U}) = \langle\mathbb{D},N\rangle$ *is such that $\mathbb{D}$ has finite products.*

*Proof.* Since $\top$ is an axiom, $\phi\vdash_{\mathrm{FV}(\phi)}\top$ is trivially provable for any logical formula $\phi$. By proposition 10.1, there is $N\phi\to N\top$ in $\mathbb{D}$, which is necessarily unique, making $N\top$ the terminal object of $\mathbb{D}$.

Let $N\phi, N\psi\in\mathrm{Obj}\,\mathbb{D}$ and consider $\theta=\phi\wedge\psi$. Since the sequents $\theta\vdash_{\mathrm{FV}(\theta)}\phi$ and $\theta\vdash_{\mathrm{FV}(\theta)}\psi$ are immediately provable, by proposition 10.1, there are $N\theta\to N\phi$ and $N\theta\to N\psi$ in $\mathbb{D}$. Suppose $N\delta\in\mathrm{Obj}\,\mathbb{D}$ is such that $N\delta\to N\phi$ and $N\delta\to N\psi$ are arrows in $\mathbb{D}$. Then $\mathcal{U}\models_{\mathcal{E}}\delta\supset\phi$ and $\mathcal{U}\models_{\mathcal{E}}\delta\supset\psi$, i.e., the sequents $\delta\vdash_{\mathrm{FV}(\delta)\cup\mathrm{FV}(\theta)}\phi$ and $\delta\vdash_{\mathrm{FV}(\delta)\cup\mathrm{FV}(\theta)}\psi$ are valid in $\mathcal{U}$. Since $(\delta\supset\phi)\wedge(\delta\supset\psi)\vdash_{\mathrm{FV}(\delta)\cup\mathrm{FV}(\theta)}\delta\supset\theta$ is provable in the empty theory, by theorem 9.2, it is valid in every model in any topos. Thus, $\mathcal{U}\models_{\mathcal{E}}\delta\supset\theta$ which means that there is a necessarily unique arrow $N\delta\to N\theta$ in $\mathbb{D}$, making $N\theta = N\phi\times N\psi$. □

**Proposition 10.4.** *For any $\mathcal{U}$, $\mathcal{E}$-model on $\Sigma$ in some topos $\mathcal{E}$,* $\mathrm{LD}(\mathcal{U}) = \langle\mathbb{D},N\rangle$ *is such that $\mathbb{D}$ has finite co-products.*

*Proof.* The sequent $\bot \vdash_{\mathrm{FV}(\phi)} \phi$ is immediately provable for any logical formula $\phi$. By proposition 10.1, there is $N\bot \to N\phi$ in $\mathbb{D}$, which is necessarily unique, making $N\bot$ the initial object of $\mathbb{D}$.

Let $N\phi, N\psi \in \mathrm{Obj}\,\mathbb{D}$ and consider $\theta = \phi \vee \psi$. Since the sequents $\phi \vdash_{\mathrm{FV}(\theta)} \theta$ and $\psi \vdash_{\mathrm{FV}(\theta)} \theta$ are immediately provable, by proposition 10.1, there are $N\phi \to N\theta$ and $N\psi \to N\theta$ in $\mathbb{D}$. Suppose $N\delta \in \mathrm{Obj}\,\mathbb{D}$ is such that $N\phi \to N\delta$ and $N\psi \to N\delta$ are arrows in $\mathbb{D}$. Then $\mathcal{U} \vDash_{\mathcal{E}} \phi \supset \delta$ and $\mathcal{U} \vDash_{\mathcal{E}} \psi \supset \delta$, i.e., the sequents $\phi \vdash_{\mathrm{FV}(\delta)\cup\mathrm{FV}(\theta)} \delta$ and $\psi \vdash_{\mathrm{FV}(\delta)\cup\mathrm{FV}(\theta)} \delta$ are valid in $\mathcal{U}$. Since $(\phi \supset \delta) \wedge (\psi \supset \delta) \vdash_{\mathrm{FV}(\delta)\cup\mathrm{FV}(\theta)} \theta \supset \delta$ is provable in the empty theory, by theorem 9.2, it is valid in every model in any topos. Thus, $\mathcal{U} \vDash_{\mathcal{E}} \theta \supset \delta$ which means that there is a necessarily unique arrow $N\theta \to N\delta$ in $\mathbb{D}$, making $N\theta = N\phi + N\psi$. □

**Proposition 10.5.** *For any $\mathcal{U}$, $\mathcal{E}$-model on $\Sigma$ in some topos $\mathcal{E}$, $\mathrm{LD}(\mathcal{U}) = \langle \mathbb{D}, N\rangle$ is such that $\mathbb{D}$ has exponentiation.*

*Proof.* Let $N\phi, N\psi \in \mathrm{Obj}\,\mathbb{D}$ and consider $\theta = \phi \supset \psi$. Since $\phi \wedge \theta \vdash_{\mathrm{FV}(\theta)} \psi$ is immediately provable in the empty theory, by proposition 10.1, there is $N(\phi \wedge \theta) \to N\psi$ in $\mathbb{D}$. But $N(\phi \wedge \theta) = N\phi \times N\theta$ by proposition 10.3, so there is $N\phi \times N\theta \to N\psi$ in $\mathbb{D}$.

Let $N\delta \in \mathrm{Obj}\,\mathbb{D}$ be such that $N\phi \times N\delta \to N\psi$ is an arrow in $\mathbb{D}$. Then, $\mathcal{U} \vDash_{\mathcal{E}} \phi \wedge \delta \supset \psi$, i.e., the sequent $\phi \wedge \delta \vdash_{\mathrm{FV}(\delta)\cup\mathrm{FV}(\theta)} \psi$ is valid in $\mathcal{U}$. Since $\phi \wedge \delta \supset \phi \vdash_{\mathrm{FV}(\delta)\cup\mathrm{FV}(\theta)} \delta \supset \theta$ is easily provable in the empty theory, by theorem 9.2, it is valid in every model in any topos. Thus, $\mathcal{U} \vDash_{\mathcal{E}} \delta \supset \theta$ which means that there is a necessarily unique arrow $N\delta \to N\theta$ in $\mathbb{D}$, making $N\theta = N\psi^{N\phi}$. □

**Proposition 10.6.** *For any $\mathcal{U}$, $\mathcal{E}$-model on $\Sigma$ in some topos $\mathcal{E}$, $\mathrm{LD}(\mathcal{U}) = \langle \mathbb{D}, N\rangle$ is such that $\mathbb{D}$ is distributive.*

*Proof.* Using the notation of definition 4.1, let $NA, NB, NC \in \mathrm{Obj}\,\mathbb{D}$ and let $\vec{x} \equiv \mathrm{FV}(A) \cup \mathrm{FV}(B) \cup \mathrm{FV}(C)$. Since $A \wedge (B \vee C) \vdash_{\vec{x}} (A \wedge B) \vee (A \wedge C)$ is easily provable in the empty theory, by proposition 10.1, there is $N(A \wedge (B \vee C)) \to N((A \wedge B) \vee (A \wedge C))$ in $\mathbb{D}$. But, by propositions 10.3 and 10.4, $N(A \wedge (B \vee C)) = NA \times (NB + NC)$ and $N((A \wedge B) \vee (A \wedge C)) = (NA \times NB) + (NA \times NC)$, thus $NA \times (NB + BC) \to (NA \times NB) + (NA \times NC)$ is an arrow in $\mathbb{D}$, which is the inverse of $\Delta$ by uniqueness. □

**Proposition 10.7.** *For any $\mathcal{U}$, $\mathcal{E}$-model on $\Sigma$ in some topos $\mathcal{E}$, $\mathrm{LD}(\mathcal{U}) = \langle \mathbb{D}, N\rangle$ is such that all the sub-categories $\mathbb{D}_{\forall x:s.\,A}$ have terminal objects.*

*Proof.* Since $\forall x{:}s.\,A \vdash_{\vec{x}} A[t/x]$, where $\vec{x} \equiv \mathrm{FV}(\forall x{:}s.\,A) \cup \mathrm{FV}(t{:}s)$, is provable in the empty theory for any logical formula $A$ on $\Sigma$ and any $t{:}s \in \mathrm{LTerms}(\Sigma)$, from proposition 10.1, there is $N(\forall x{:}s.\,A) \to NA[t/x] = \Sigma_A(x{:}s)(t{:}s)$ in $\mathbb{D}$, showing that $N(\forall x{:}s.\,A)$ is the vertex of a cone on $\Sigma_A(x{:}s)$, i.e., $N(\forall x{:}s.\,A) \in \mathrm{Obj}\,\mathbb{D}_{\forall x:s.\,A}$.

Let $NB \in \mathrm{Obj}\,\mathbb{D}_{\forall x:s.\,A}$: by definition of the sub-category, $B$ is a logical formula over $\Sigma$ such that $x{:}s \notin \mathrm{FV}(B)$ and, for each $t{:}s \in \mathrm{LTerms}(\Sigma)$, $NB \to NA[t/x]$ is an arrow in $\mathbb{D}$. Then, calling $\vec{y} \equiv \mathrm{FV}(B) \cup \mathrm{FV}(\forall x{:}s.\,A) \cup \mathrm{FV}(t{:}s)$, it holds that $\mathcal{U} \vDash_{\mathcal{E}} B \supset A[t/x]$ for any $t{:}s \in \mathrm{LTerms}(\Sigma)$. In particular, for $x{:}s$, $\mathcal{U} \vDash_{\mathcal{E}} B \supset A$, i.e., $B \vdash_{\vec{y}} A$ is valid in $\mathcal{U}$. Since $x{:}s \notin \mathrm{FV}(B)$, $B \vdash_{\vec{y}} \forall x{:}s.\,A$ holds in $\mathcal{U}$. Hence, there is a unique $NB \to N(\forall x{:}s.\,A)$ in $\mathbb{D}$, that is, $N(\forall x{:}s.\,A)$ is the terminal object of $\mathbb{D}_{\forall x:s.\,A}$. □

**Proposition 10.8.** *For any $\mathcal{U}$, $\mathcal{E}$-model on $\Sigma$ in some topos $\mathcal{E}$, $\mathrm{LD}(\mathcal{U}) = \langle \mathbb{D}, N\rangle$ is such that all the sub-categories $\mathbb{D}_{\exists x:s.\,A}$ have initial objects.*

*Proof.* Since $A[t/x] \vdash_{\vec{x}} \exists x{:}s.\,A$, where $\vec{x} \equiv \mathrm{FV}(\exists x{:}s.\,A) \cup \mathrm{FV}(t{:}s)$, is provable in the empty theory for any logical formula $A$ on $\Sigma$ and any $t{:}s \in \mathrm{LTerms}(\Sigma)$, from proposition 10.1, there is $\Sigma_A(x{:}s)(t{:}s) = NA[t/x] \to N(\exists x{:}s.\,A)$ in $\mathbb{D}$, showing that $N(\exists x{:}s.\,A)$ is the vertex of a co-cone on $\Sigma_A(x{:}s)$, i.e., $N(\exists x{:}s.\,A) \in \mathrm{Obj}\,\mathbb{D}_{\exists x:s.\,A}$.

Let $NB \in \mathrm{Obj}\,\mathbb{D}_{\exists x:s.\,A}$: by definition of the sub-category, $B$ is a logical formula over $\Sigma$ such that $x{:}s \notin \mathrm{FV}(B)$ and, for each $t{:}s \in \mathrm{LTerms}(\Sigma)$, $NA[t/x] \to NB$ is an arrow in $\mathbb{D}$. Then, calling $\vec{y} \equiv \mathrm{FV}(B) \cup \mathrm{FV}(\exists x{:}s.\,A) \cup \mathrm{FV}(t{:}s)$, it holds that $\mathcal{U} \vDash_{\mathcal{E}} A[t/x] \supset B$ for any $t\colon s \in \mathrm{LTerms}(\Sigma)$. In particular, for $x\colon s$, $\mathscr{U} \vDash_{\mathscr{E}} A \supset B$, i.e., $A \vdash_{\vec{y}} B$ is valid in $\mathscr{U}$. Since $x\colon s \notin \mathrm{FV}(B)$, $\exists x\colon s.\, A \vdash_{\vec{y}} B$ holds in $\mathscr{U}$. Hence, there is a unique $N(\exists x\colon s.\, A) \to NB$ in $\mathbb{D}$, that is, $N(\exists x\colon s.\, A)$ is the initial object of $\mathbb{D}_{\exists x\colon s.\, A}$. □

**Proposition 10.9.** *For any* $\mathscr{U}$*,* $\mathscr{E}$*-model on* $\Sigma$ *in some topos* $\mathscr{E}$*,* $\mathrm{LD}(\mathscr{U}) = \langle \mathbb{D}, N\rangle$ *is such that* $N$ *satisfies the condition expressed in point* (6) *of definition* 4.1.

*Proof.* By propositions 10.3 to 10.8, the claim is evident. □

**Proposition 10.10.** *For any* $\mathscr{U}$*,* $\mathscr{E}$*-model on* $\Sigma$ *in some topos* $\mathscr{E}$*,* $\mathrm{LD}(\mathscr{U}) = \langle \mathbb{D}, N\rangle$ *is such that the unique arrow* $N(\exists x\colon s.\, A \times B) \to NA \times N(\exists x\colon s.\, B)$ *in* $\mathbb{D}$*, where* $x\colon s \notin \mathrm{FV}(A)$*, has an inverse.*

*Proof.* Since the sequent $A \wedge \exists x\colon s.\, B \vdash_{\vec{x}} \exists x\colon s.\, A \wedge B$, where $\vec{x} \equiv \mathrm{FV}(A) \cup \mathrm{FV}(\exists x\colon s.\, B)$ and $x\colon s \notin \mathrm{FV}(A)$, holds in the empty theory, by propositions 10.1 and 10.3, there is an arrow $NA \times N(\exists x\colon s.\, B) \to N(\exists x\colon s.\, A \times B)$ in $\mathbb{D}$. By uniqueness of this arrow, it is immediate to see that it is the inverse of $N(\exists x\colon s.\, A \times B) \to NA \times N(\exists x\colon s.\, B)$. □

**Proposition 10.11.** *For any* $\mathscr{U}$*,* $\mathscr{E}$*-model on* $\Sigma$ *in some topos* $\mathscr{E}$*,* $\mathrm{LD}(\mathscr{U}) = \langle \mathbb{D}, N\rangle$ *is a logically distributive category.*

*Proof.* Obvious by propositions 10.2 to 10.10. □

**Proposition 10.12.** *Let* $\mathscr{U}$ *be a* $\mathscr{E}$*-model on* $\Sigma$ *in some topos* $\mathscr{E}$*, and let* $\phi \vdash_{\vec{x}} \psi$ *be a sequent in the language over* $\Sigma$*. Then* $\phi \vdash_{\vec{x}} \psi$ *is valid in* $\mathscr{U}$ *if and only if* $\phi \vdash_{\vec{x}} \psi$ *is valid in* $\mathrm{LD}(\mathscr{U}) = \langle \mathbb{D}, N\rangle$*.*

*Proof.* The sequent $\phi \vdash_{\vec{x}} \psi$ is valid in $\mathscr{U}$ if and only if $\mathscr{U} \vDash_{\mathscr{E}} \phi \supset \psi$, which holds if and only if there is an arrow $N\phi \to N\psi$ in $\mathbb{D}$, that is, exactly when $\phi \vdash_{\vec{x}} \psi$ is valid in $\mathrm{LD}(\mathscr{U})$. □

**Definition 10.2** (Logical functor)**.** A functor $F\colon \mathscr{E} \to \mathscr{F}$ between the topoi $\mathscr{E}$ and $\mathscr{F}$ is said to be a *logical functor* when

(1) $F$ preserves finite limits;
(2) for all $a \in \mathrm{Obj}\,\mathscr{E}$, $F(\wp a)$ with the relation $F(\in_a)$ is a power object for $Fa$.

We remind, see [Gol06, p. 104], that a category $\mathbb{C}$ with products is said to have *power objects* when, for each $a \in \mathrm{Obj}\,\mathbb{C}$, there are $\wp a, E_a \in \mathrm{Obj}\,\mathbb{C}$ and a sub-object $\in_a\colon E_a \rightarrowtail \wp a \times a$ such that, for any $b \in \mathrm{Obj}\,\mathbb{C}$ and $r\colon R \rightarrowtail b \times a$, there is a unique $f_r\colon b \rightarrowtail \wp a$ making the following diagram to be a pullback

$$\begin{array}{ccc} R & \overset{r}{\rightarrowtail} & b \times a \\ \downarrow & & \downarrow{\scriptstyle f_r \times 1_a} \\ E_a & \underset{\in_a}{\rightarrowtail} & \wp a \times a \end{array}$$

Also, any topos $\mathscr{E}$ has power objects, in particular $\wp a \cong \Omega^a$ and $\in_a\colon E_a \rightarrowtail \Omega^a \times a$ is the sub-object of $\Omega^a \times a$ whose character is $\mathrm{ev}_a\colon \Omega^a \times a \to \Omega$, the evaluation map. Moreover, a topos $\mathscr{E}$ **is** a category which is finitely complete and has power objects, thus a logical functor preserves finite limits, finite co-limits, exponentials and the sub-object classifier. Also, a logical functor $F$ has a left adjoint if and only if it has a right adjoint, see [Joh02a, p. 89]. Finally, if the logical functor $F$ has a Cartesian left adjoint, then $F$ is an equivalence [Joh02a, p. 96].

**Definition 10.3** ($\Sigma$-Str($\mathscr{E}$))**.** Fixed a topos $\mathscr{E}$ and a logical signature $\Sigma = \langle S, F, R\rangle$, $\Sigma$-Str($\mathscr{E}$) is the category whose objects are the $\mathscr{E}$-models on $\Sigma$ in $\mathscr{E}$, and whose arrows are the corresponding logical endo-functors which preserve interpretations, that is, $G\colon \mathscr{U}_1 \to \mathscr{U}_2$ in $\Sigma$-Str($\mathscr{E}$) when $G$ is a logical functor $\mathscr{E} \to \mathscr{E}$ such that $M_{\mathscr{U}_2} s = G(M_{\mathscr{U}_1} s)$ for any $s \in S$, $M_{\mathscr{U}_2} f = G(M_{\mathscr{U}_1} f)$ for any $f \in F$, and $M_{\mathscr{U}_2} r = G(M_{\mathscr{U}_1} r)$ for any $r \in R$.

**Fact 10.13.** *If* $F\colon \mathscr{U}_1 \to \mathscr{U}_2$ *is an arrow of* $\Sigma$-Str$(\mathscr{E})$*, then* $M_{\mathscr{U}_2}(t\colon s) = G\left(M_{\mathscr{U}_1}(t\colon s)\right)$ *for any* $t\colon s \in \mathrm{LTerms}(\Sigma)$*, and* $[\![\vec{x}.\phi]\!]_{\mathscr{U}_2} = F\left([\![\vec{x}.\phi]\!]_{\mathscr{U}_1}\right)$ *for any logical formula-in-context* $\vec{x}.\phi$*.*

*Proof.* Evident induction from the definitions. □

**Proposition 10.14.** *For each* $\mathscr{U}$*,* $\mathscr{E}$*-model on* $\Sigma$ *in some topos* $\mathscr{E}$*,* $\mathrm{LD}(\mathscr{E}) \in \Sigma$-Str$(\mathscr{E})$*.*

*Proof.* Evident by definition of $\mathrm{LD}(\mathscr{E})$—see also definition 8.3. □

**Definition 10.4** (LD functor)**.** Given a topos $\mathscr{E}$ and a $\lambda$-signature $\Sigma$, the functor LD extends the map LD, see definition 10.1, as follows:

$$\mathrm{LD}\colon\ \Sigma\text{-Str}(\mathscr{E}) \to \Sigma\text{-Str}$$

and, for each $\mathscr{E}$, $\mathscr{E}$-model on $\Sigma$, $\mathrm{LD}(\mathscr{U})$ is defined by the map LD; for any arrow $F\colon \mathscr{U}_1 \to \mathscr{U}_2$ in $\Sigma$-Str$(\mathscr{E})$, the corresponding $\mathrm{LD}(F)\colon\ \mathrm{LD}(\mathscr{U}_1) = \langle \mathbb{D}_1, N_1\rangle \to \langle \mathbb{D}_2, N_2\rangle = \mathrm{LD}(\mathscr{U}_2)$ is defined by

(1) $\mathrm{LD}(F)\left([\![\mathrm{FV}(\phi).\phi]\!]_{\mathscr{U}_1}\right) = F\left([\![\mathrm{FV}(\phi).\phi]\!]_{\mathscr{U}_1}\right)\left(= [\![\mathrm{FV}(\phi).\phi]\!]_{\mathscr{U}_2}\right)$;
(2) for any arrow $f\colon [\![\mathrm{FV}(\phi).\phi]\!]_{\mathscr{U}_1} \to [\![\mathrm{FV}(\psi).\psi]\!]_{\mathscr{U}_1}$, $\mathrm{LD}(F)(f)$ is the unique arrow such that

$$\Rightarrow_{\mathscr{U}_2} \circ \left(F\left([\![\mathrm{FV}(\phi)\cup\mathrm{FV}(\psi).\phi]\!]_{\mathscr{U}_1}\right), F\left([\![\mathrm{FV}(\phi)\cup\mathrm{FV}(\psi).\psi]\!]_{\mathscr{U}_1}\right)\right) = F\left(\top_{M_{\mathscr{U}_1}(\mathrm{FV}(\phi)\cup\mathrm{FV}(\psi))}\right) .$$

It is evident from the properties of logical functors that LD is, indeed, a functor.

By proposition 10.12, LD preserves and reflects the validity of sequents. Thus, first-order models in elementary topoi are as expressive as logically distributive categories.

It is very unlikely that the LD functor can be ‘inverted’, e.g., by proving that it has a left or a right adjoint. The amount of information present in the topos is much wider than the one coded into the logically distributive category—using Johnstone’s words “every topos thinks to be **Set**, the topos of sets and functions between them”.

## 11. KRIPKE SEMANTICS

The traditional semantics used to interpret first-order intuitionistic logic is Kripke's one. This section is devoted to relate the Kripke semantics with the one based on logically distributive categories.

**Definition 11.1** (Kripke model)**.** Let $P$ be a poset and $\Sigma = \langle S, F, R\rangle$ a logical signature. A *Kripke model on $P$ and* $\Sigma$ is a structure

$$\mathcal{K} = \langle \{\mathcal{U}_p = \langle A_p, F_p, R_p\rangle\}_{p\in \mathrm{Obj}\,P}, \{A_{p,q}\colon A_p \to A_q\}_{p\le q \text{ in } P}\rangle$$

such that

(1) for each $p \in \mathrm{Obj}\,P$, $\mathcal{U}_p$ is a **Set**-model on $\Sigma$ in the topos **Set**;

(2) for $p \le q$, arrow of $P$, $A_{p,q}\colon A_p \to A_q$ is a family $\left\{A^s_{p,q}\colon A^s_p \to A^s_q\right\}_{s\in S}$ such that $A^s_{p,q} \in \mathrm{Hom}_{\mathbf{Set}}(A^s_p, A^s_q)$, and
- for each $t\colon s \in \mathrm{LTerms}(\Sigma)$, $A^s_{p,q}\left(M_p(t\colon s)\right) = M_q(t\colon s)$, where $M_p$ and $M_q$ are the interpretations in $\mathcal{U}_p$ and $\mathcal{U}_q$, respectively;
- for each $r\colon s_1 \times \cdots \times s_n \in R$ and $t_1\colon s_1,\ldots,t_n\colon s_n \in \mathrm{LTerms}(\Sigma)$,

  $$\left(M_p(t_1\colon s_1),\ldots,M_p(t_n\colon s_n)\right)\vec{x} \in M_p r$$

  if and only if

  $$\left(M_q(t_1\colon s_1),\ldots,M_q(t_n\colon s_n)\right)A_{p,q}(\vec{x}) \in M_q r$$

  for every $\vec{x} \in M_p\left(\bigcup_{i=1}^n \mathrm{FV}(t_i\colon s_i)\right)$, and $A_{p,q}(\vec{x}) = \left(A^{s_1}_{p,q} \times \cdots \times A^{s_n}_{p,q}\right)\vec{x}$;

(3) for each $p \in \mathrm{Obj}\,P$ and for each $s \in S$, $A^s_{p,p} = 1_{A^s_p}$;

(4) for every pair $p \le q$, $q \le r$ of arrows in $P$, and for each $s \in S$, $A^s_{p,r} = A^s_{q,r} \circ A^s_{p,q}$.

*Note* 11.1. This definition is slightly different from the standard ones in literature, see, e.g., [BM77, p. 416], in that we deal with a multi-sorted environment. Apart from this and with a preference for a categorical-oriented presentation, the definition given here is easily proved to be equivalent to the various ones found in literature.

**Definition 11.2** (Validity)**.** Given a poset $P$, a logical signature $\Sigma = \langle S, F, R\rangle$, and a Kripke model $\mathcal{K}$ on $P$ and $\Sigma$, $\mathcal{K} \vDash_p \phi$ indicates that the logical formula $\phi$ on $\Sigma$ is *valid in* $p \in \mathrm{Obj}\,P$. The notion $\mathcal{K} \vDash_p \phi$ is inductively defined by

(1) $\phi \equiv \top$: $\mathcal{K} \vDash_p \phi$;

(2) $\phi \equiv \bot$: $\mathcal{K} \nvDash_p \phi$;

(3) $\phi \equiv r(t_1,\ldots,t_n)$, with $r\colon s_1 \times \cdots \times s_n \in R$ and $t_1\colon s_1,\ldots,t_n\colon s_n \in \mathrm{LTerms}(\Sigma)$: $\mathcal{K} \vDash_p \phi$ exactly when

$$\left(M_p(t_1\colon s_1),\ldots,M_p(t_n\colon s_n)\right)\vec{x} \in M_p r$$

for every $\vec{x} \in M_p\left(\bigcup_{i=1}^n \mathrm{FV}(t_i\colon s_i)\right)$, where $M_p$ is the interpretation associated with the $p$-th **Set**-model in $\mathcal{K}$;

(4) $\phi \equiv \psi_1 \wedge \psi_2$: $\mathcal{K} \vDash_p \phi$ if and only if $\mathcal{K} \vDash_p \psi_1$ and $\mathcal{K} \vDash_p \psi_2$;

(5) $\phi \equiv \psi_1 \vee \psi_2$: $\mathcal{K} \vDash_p \phi$ if and only if $\mathcal{K} \vDash_p \psi_1$ or $\mathcal{K} \vDash_p \psi_2$;

(6) $\phi \equiv \psi_1 \supset \psi_2$: $\mathcal{K} \vDash_p \phi$ if and only if, for every $q \in \mathrm{Obj}\,P$ such that $p \le q$ in $P$, $\mathcal{K} \nvDash_p \psi_1$ or $\mathcal{K} \vDash_p \psi_2$;

(7) $\phi \equiv \exists y\colon s.\psi$: $\mathcal{K} \vDash_p \phi$ if and only if there is $a \in A^s_p$, the $s$-th component of the $p$-th **Set**-model of $\mathcal{K}$, such that $\mathcal{K} \vDash_{p,a} \psi$, where $\mathcal{K} \vDash_{p,a} \psi$ is inductively defined exactly as $\mathcal{K} \vDash_p \psi$ except for clause (3), which becomes $\mathcal{K} \vDash_{p,a} r(t_1,\ldots,t_n)$ when

$$\left(M_p(t_1\colon s_1),\ldots,M_p(t_n\colon s_n)\right)(\vec{x},a) \in M_p r$$

for any $\vec{x} \in M_p\left(\bigcup_{i=1}^n \mathrm{FV}(t_i\colon s_i)\right)$ and $y$ is interpreted in $a$;

(8) $\phi \equiv \forall y\colon s.\psi$: $\mathcal{K} \vDash_p \phi$ if and only if, for every $q \in \mathrm{Obj}\,P$ such that $p \le q$ in $P$, and for each $a \in A^s_q$, $\mathcal{K} \vDash_{p,a} \psi$ with the same notation as above.

We say that $\phi$ is *valid in* $\mathcal{K}$, notation $\mathcal{K} \models \phi$, when $\mathcal{K} \models_p \phi$ for any $p \in \operatorname{Obj} P$.

*Note* 11.2. It is possible to analyse Kripke models inside category theory, see, e.g., [Gol06, section 11.6], since definition 11.1 shows that a Kripke model is, in fact, a category. Nevertheless, it is convenient to embed these models in a suitable topos, as we are going to do in the following.

**Definition 11.3** (Kripke functor)**.** Let $P$ be a poset, $\Sigma = \langle S, F, R\rangle$ a logical signature, and $\mathcal{K}$ a Kripke model on $P$ and $\Sigma$. For each $s \in S$, the functor $A^s \colon P \to \mathbf{Set}$ maps each $p \in \operatorname{Obj} P$ in $A^s_p$, the $s$-th component of the $p$-th **Set**-model of $\mathcal{K}$, and each arrow $p \le q$ of $P$ in $A^s_{p,q}$, the $s$-th component of the $(p,q)$-th part of the accessibility relation of $\mathcal{K}$.

*Note* 11.3. By clauses (3) and (4) of definition 11.1, it is immediate to see that $A^s$ as defined is a functor, indeed.

**Definition 11.4** (Kripke model in a topos)**.** Given a poset $P$, a logical signature $\Sigma = \langle S, F, R\rangle$, and a Kripke model $\mathcal{K}$ on $P$ and $\Sigma$ whose parts are denoted as in definition 11.1, we define $\mathcal{K}^* = \langle A^*, F^*, R^*\rangle$ as the $[P,\mathbf{Set}]$-model on $\Sigma$ in the topos of presheaves on $P$ as follows:

(1) $A^* = \{A^s\}_{s\in S}$ where $A^s$ is the functor of definition 11.3;
(2) $F^* = \{f^* \colon A^{s_1} \times \cdots \times A^{s_n} \to A^{s_0}\}_{f \colon s_1\times\cdots\times s_n \to s_0 \in F}$ with $f^*$ the natural transformation having $f^*_p = M_p f$ as $p$-th component;
(3) $R^* = \{r^* \colon A^{s_1} \times \cdots \times A^{s_n} \to \Omega\}_{r \colon s_1\times\cdots\times s_n \in R}$ where $r^*$ is the natural transformation having $r^*_p(a_1,\dots,a_n) = \left\{q \colon p \le q \text{ in } P \text{ and } \left(A^{s_1}_{p,q}(a_1),\dots,A^{s_n}_{p,q}\right) \in M_q r\right\}$ as $p$-th component.

**Proposition 11.1.** *There is a bijective correspondence between Kripke models on $P$ and $\Sigma$ and $[P,\mathbf{Set}]$-models on $\Sigma$.*

*Proof.* See [Gol06, p. 258]. □

**Theorem 11.2.** *For any Kripke model $\mathcal{K}$ on $P$ and $\Sigma$, and associated $[P,\mathbf{Set}]$-model $\mathcal{K}^*$, we have that, for any logical formula $\phi$ on $\Sigma$, $\mathcal{K} \models \phi$ if and only if $\mathcal{K}^* \models_{[P,\mathbf{Set}]} \phi$.*

*Proof.* See [Gol06, p. 264]. □

**Proposition 11.3.** *For each Kripke model $\mathcal{K}$ on $P$ and $\Sigma$, there is a logically distributive category $\langle \mathbb{K}, N\rangle \in \Sigma$-Str such that, for any logical formula $\phi$, $\mathcal{K} \models \phi$ if and only if $\phi$ is valid in $\mathbb{K}$.*

*Proof.* Let $\langle \mathbb{K}, N\rangle = \mathrm{LD}(\mathcal{K}^*)$, with LD the functor of definition 10.4 and $\mathcal{K}^*$ as for definition 11.4. Since the functor LD preserves and reflects validity, the result immediately follows from theorem 11.2. □

## 12. Grothendieck topoi

The main objective of the project funding the research efforts illustrated in this work was to analyse predicative theories and their relation to Grothendieck topoi.

Although the research line originally proposed revealed impossible to achieve, at least for the author, what has been developed and presented till now has the precise aim to fulfil the original promise stated in the research project.

In this section, we will see how allowing 'ideal' elements to find their way in the picture we have drawn so far, is the way to 'invert' the functors mapping, e.g., topos-oriented models to structures in logically distributive categories. This technical result produces a number of consequences: firstly, it will follow that presheaf topoi are sufficient to give a sound and complete meaning to the logical theories we consider; secondly, ideal elements suggest that some sort of impredicative construction is at work. Also, as presheaf topoi are a very special case of Grothendieck topoi, a wider analysis is pursued.

12.1. **Ideal elements.** In this part, we will show how the LD functor of definition 10.4 can be 'inverted' by proving that it has a sort of 'ideal' right adjoint.

**Proposition 12.1.** *Let $F\colon \mathbb{C} \to \mathbb{D}$ be a functor between small categories. Then, the functor $F^*\colon [\mathbb{D}^{\mathrm{op}}, \mathbf{Set}] \to [\mathbb{C}^{\mathrm{op}}, \mathbf{Set}]$ given by $F^*(G) = G \circ F$ has both left and right adjoints $F_! \dashv F^* \dashv F_*$. Moreover, there is a natural isomorphism $F_! \circ \mathbf{Y}_{\mathbb{C}} \cong \mathbf{Y}_{\mathbb{D}} \circ F$, where $\mathbf{Y}$ is the Yoneda embedding:*

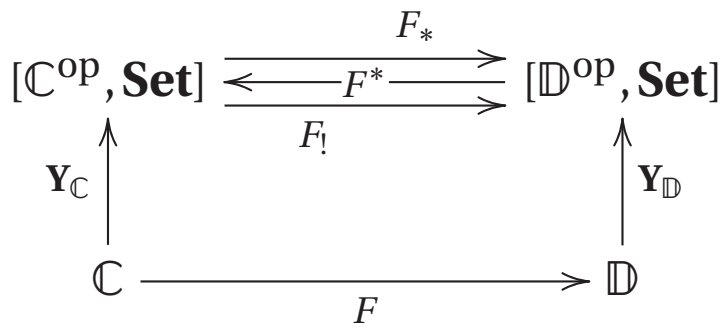

*Proof.* See Corollary 9.17 in [Awo06, p. 239]. □

So, every functor $F\colon \mathbb{C} \to \mathbb{D}$ between small categories has a sort of right adjoint $F^* \circ \mathbf{Y}_{\mathbb{D}}\colon \mathbb{D} \to [\mathbb{C}^{\mathrm{op}}, \mathbf{Set}]$, except that its values are in the 'ideal' elements of the co-completion $[\mathbb{C}^{\mathrm{op}}, \mathbf{Set}]$. Properly speaking, this is a right Kan extension, but we will stick on the adjoint terminology as it conveys a more natural feeling of what is intended.

In particular, the functor LD: $\Sigma$-Str$(\varepsilon) \to \Sigma$-Str, with $\varepsilon$ a topos, see definition 10.4, has an ideal right adjoint in the sense above, DL: $\Sigma$-Str $\to [\Sigma\text{-Str}(\varepsilon)^{\mathrm{op}}, \mathbf{Set}]$, which can be calculated as the right Kan extension composed with the Yoneda embedding, as in proposition 12.1.

**Proposition 12.2.** *Given a logically distributive category $\langle \mathbb{C}, M\rangle$ on a signature $\Sigma$, there is a logically distributive subcategory $\langle \mathbb{D}, N\rangle$, $\mathbb{D} \subseteq \mathbb{C}$ with $N$ surjective such that $M\phi = N\phi$ for every logical formula $\phi$.*

*Proof.* Let $\mathbb{D}$ be the full subcategory of $\mathbb{D}$ whose objects are in the image of $M$, and define $N$ accordingly. The result follows immediately from definition 4.1. □

The meaning of proposition 12.2 is that we can always restrict our attention to the 'working part' of a logically distributive category when interested in logical interpretations: in fact, only this part is relevant for the considered semantics, as shown in the following corollary:

**Corollary 12.3.** *The category $\langle \mathbb{D}, N\rangle$ above is such that $[\![\vec{x}.\, t\colon A]\!]_{N'} = [\![\vec{x}.\, t\colon A]\!]_{M'}$ for every $\vec{x}.\, t\colon A$ term-in-context, and $N'$ and $M'$ $\Sigma$-structures in $\langle \mathbb{D}, N\rangle$ and $\langle \mathbb{C}, M\rangle$, respectively.*

*Proof.* Evident by definition 4.2. □

Obviously, restricting to the class of logically distributive categories $\langle \mathbb{C}, M\rangle$ with $M$ surjective, hereafter *surjective logically distributive categories*, does provide a soundness theorem, like theorem 5.4, essentially with the same proof, and a completeness result, like theorem 6.14 with a similar proof, since the syntactical category $\langle \mathbb{C}_T, M_T\rangle$ of definition 6.2 lies in the class of surjective logically distributive categories. Therefore, all the results obtained in sections 8, 10 and 11 can be rephrased in the restricted framework of surjective logically distributive categories.

In particular, it follows that, modulo isomorphisms, the class of surjective logically distributive categories coincides with the $\Sigma$-Str category of definition 8.4. Thus, in this restricted but not limiting framework, we can say that, fixed any topos $\varepsilon$, it is possible to map each $\Sigma$-structure $\langle \mathbb{D}, M, M_A\rangle$ over $\Sigma$ into a model in a presheaf topos, namely $[\Sigma\text{-Str}(\varepsilon)^{\mathrm{op}}, \mathbf{Set}]$ via the DL functor, the ideal right adjoint of LD: $\Sigma\text{-Str}(\varepsilon) \to \Sigma\text{-Str}$, and this map preserves validity.

The first consequence of this fact is that presheaf topoi and the logical interpretation in them are a sufficiently rich environment to provide a sound and complete semantics for the first-order logical theories based on intuitionistic logic—a fact already clear after Kripke semantics, but obtained through a rather different path.

12.2. **Grothendieck topoi.** A Grothendieck topos is a category $\mathbb{G}$ which is equivalent to the topos $\mathbf{Sh}(\mathbb{C}, J)$ of sheaves on a site $\langle \mathbb{C}, J\rangle$. Since every sheaf is a contravariant functor, there is an evident inclusion $i\colon \mathbf{Sh}(\mathbb{C}, J) \to [\mathbb{C}^{\mathrm{op}}, \mathbf{Set}]$. We may regard this inclusion as a forgetful functor since it 'forgets' the topological structure induced by $J$. Then, as commonly happens with forgetful functors, $i$ has a left adjoint.

**Theorem 12.4.** *Let $\mathbb{C}$ be a small category and $\langle \mathbb{C}, J\rangle$ a site. Then, the inclusion functor $i\colon \mathbf{Sh}(\mathbb{C}, J) \to [\mathbb{C}^{\mathrm{op}}, \mathbf{Set}]$ has a left adjoint* **a***, called the* associated sheaf functor. *Moreover, this functor commutes with finite limits.*

*Proof.* See [MLM92, p. 128]. □

Diagrammatically, the situation can be depicted as follows:

$$\mathbb{G} \underset{G}{\overset{F}{\rightleftarrows}} \mathbf{Sh}(\mathbb{C}, J) \underset{\mathbf{a}}{\overset{i}{\rightleftarrows}} [\mathbb{C}^{\mathrm{op}}, \mathbf{Set}]$$

where $G$ is left adjoint to $F$, and both the unit and the co-unit of the adjunction are natural isomorphisms, and $\mathbf{a}$ is left adjoint to $i$. Thus, the left adjoint of $i \circ F$ exists and it is $G \circ \mathbf{a}$. Moreover, since $\mathbf{a} \circ i$ is naturally isomorphic to the identity $1_{\mathbf{Sh}(\mathbb{C},J)}$, see [MLM92, p. 133], $G \circ \mathbf{a} \circ i \circ F$ is naturally isomorphic to the identity functor $1_{\mathbb{G}}$.

It is not difficult to see that Grothendieck topoi are best related via *geometric morphisms*, i.e., adjunctions where the left adjoint is *left exact*, that is, it preserves finite limits. From theorem 12.4, it follows that the inclusion functor $i$ together with $\mathbf{a}$ is a geometric morphism, a situation that easily generalises to $i \circ F$ above.

Unfortunately, geometric morphisms do not preserve co-limits and have other bad properties when dealing with some of the categorical properties we have used to give a meaning to first-order logical theories. Differently from logical functors, they do not preserve the topos structure, in the elementary sense, but rather the topological structure, which is not preserved by logical functors. In fact, a specific logic has been invented to support and to capture this structure, *geometric logic*, an infinitary first-order system, extensively described in [Joh02b, part D] and [MLM92, p. 532–541].

To some extent, every first order theory $T$ can be transformed into a geometric theory $T_G$ in such a way that $T_G$ is 'equivalent' to $T$: precisely, the two theories validate the same models in **Set**, and, more in general, in every Boolean topos. This transformations is called *Morleyisation* and it is described and analysed in [Joh02b, D.1.5.13].

Since every Grothendieck topos is also an elementary topos, the soundness theorem holds: if $T \vdash \phi$ then, for every model $M$ in any Grothendieck topos $\mathbb{G}$ such that $M \vDash_{\mathbb{G}} T$, $M \vDash_{\mathbb{G}} \phi$. Moreover, as every presheaf topos $[\mathbb{C}^{\text{op}}, \mathbf{Set}]$ is a Grothendieck topos in a trivial way, a completeness result holds as well: if for each model $M$ in a Grothendieck topos $\mathbb{G}$ such that $M \vDash_{\mathbb{G}} T$, it holds that $T \vdash \phi$. The proof is simple: the hypothesis applies to all the presheaf topoi, and we already have shown a completeness result on them.

But, in our approach, we relate models in a topos via logical functors, see definition 10.2, which are not geometric in the general case, thus we loose the characterising aspect of Grothendieck topoi, namely the fact of exploiting the topological properties of the sites they are built from. In short, this fact prevented the development of the project beside this work as originally stated and planned, and there are good reasons to believe that the situation cannot be naturally recovered—but, since this part is essentially negative and does not provide useful insights, it will not appear in this report.

Nevertheless, something more can be said about Grothendieck topoi in our analysis of first-order logic.

Let $\langle \mathbb{C}, M \rangle$ be a surjective logically distributive category. The Yoneda embedding $\mathbf{Y}_{\mathbb{C}} \colon \mathbb{C} \to [\mathbb{C}^{\text{op}}, \mathbf{Set}]$ is full and faithful, see, e.g., [Bor94a, p. 19]. Hence, the category $\langle [\mathbb{C}^{\text{op}}, \mathbf{Set}], M' \rangle$, where $M'\phi = \text{Hom}_{\mathbb{C}}(-, M\phi)$, $\text{Hom}_{\mathbb{C}}$ being the usual contravariant representation functor, is a logically distributive category, not necessarily surjective, as it is immediate to verify from definition 4.1 remembering that $[\mathbb{C}^{\text{op}}, \mathbf{Set}]$ is an elementary topos and that every topos is a distributive category.

Since $\mathbf{Y}_{\mathbb{C}}$ is full and faithful, $M'$ identifies a model in $[\mathbb{C}^{\text{op}}, \mathbf{Set}]$ which is isomorphic to the model $M$ in $\mathbb{C}$. This fact shows that the presheaf topoi are powerful enough to provide a sound and complete semantics for first-order theories, and it does so in another way complementing those already introduced.

In addition, whenever a sub-canonical Grothendieck topology $J$, i.e., when all the representable functors are sheaves, is imposed over $\mathbb{C}$, we get that $\langle \mathbf{Sh}(\mathbb{C}, J), M'' \rangle$ with $M''$ defined as $M'$ but restricted to the topos of sheaves, is a logically distributive category whose model $M''$ is isomorphic to $M$. In a diagram:

$$\begin{array}{ccc} \langle \mathbf{Sh}(\mathbb{C},J), M''\rangle & \underset{\mathbf{a}}{\overset{i}{\rightleftarrows}} & \langle [\mathbb{C}^{\text{op}},\mathbf{Set}], M\rangle \\ & {}_{\mathbf{Y}_{\mathbb{C},J}}\nwarrow \qquad \nearrow_{\mathbf{Y}_{\mathbb{C}}} & \\ & \langle \mathbb{C}, M\rangle & \end{array}$$

Thus, via the same reasoning done more than once, the class of sheaf topoi on sub-canonical topologies provides enough room for a sound and complete semantics of first-order theories (and the related $\lambda$-calculi).

We want to conclude our analysis of Grothendieck topoi with a suggestion for a possible line of research that may be pursued as a consequence of the achievements obtained so far. When we consider two Grothendieck topoi $\mathbb{G}$ and $\mathbb{H}$ which are *Morita-equivalent*, i.e., they share the same classifying topos, modulo isomorphisms, see [MLM92, chapters 8 and 10] and [Joh02b, D.1.4.9], and $\mathbb{G}$ and $\mathbb{H}$ arise from the logical theories $T_{\mathbb{G}}$ and $T_{\mathbb{H}}$ expressed in the language of geometric logic, validity is ‘transported’ by the equivalence, that is, $T_{\mathbb{G}} \vdash \phi$ if and only if $T_{\mathbb{H}} \vdash \theta$, where $\theta$ is uniquely identified given $\phi$ and the Morita-equivalence, and vice versa. This reasoning technique has been introduced by O. Caramello in [Car09] and extensively studied in her subsequent works [Car].

Now, supposing $\mathbb{G} \simeq \mathbf{Sh}(\mathbb{C}, J)$ where $\langle \mathbb{C}, M \rangle$ is a (surjective) logically distributive category, there is $M_{\mathbb{G}}$ such that $\langle \mathbb{G}, M_{\mathbb{G}} \rangle$ is logically distributive and the model $M_{\mathbb{G}}$ is isomorphic to $M$. Making the same assumption for $\mathbb{H}$, i.e., $\langle \mathbb{H}, M_{\mathbb{H}} \rangle$ is both a Grothendieck topos and a logically distributive category, when $\mathbb{G}$ is Morita-equivalent to $\mathbb{H}$, there is a strict link between $M_{\mathbb{G}}$ and $M_{\mathbb{H}}$, and rarely a trivial one. We conjecture that a technique analogous to Caramello’s can be developed to show deep ‘entanglements’ among first-order theories arising from Morita-equivalent, logically distributive Grothendieck topoi.

## 13. IMPREDICATIVE THEORIES

The objective of the research project this report provides an account of, was to analyse predicative logical theories in the framework of categorical logic.

In the first place, it must be said that the main reason to consider first-order theories oppositely to higher-order ones, is that the latter are intrinsically impredicative, thus making impossible to analyse them in the proposed perspective.

But first-order theories may be impredicative too. Usually, predicativity is associated to the possibility to assign a computational meaning to proofs in such a way that no oracles are needed: in other terms, whenever a statement $\Gamma \vdash A$ gets proved, there is an algorithm whose preconditions are the elements of $\Gamma$ and whose postcondition is $A$. This fact is trivially true in the framework of logically distributive categories as any proof is mapped into a suitable $\lambda$-terms $\theta$. But $\theta$ is computable only when each axiom of the theory can be turned into a 'program' whose postcondition is the axiom. Evidently, this is not always the case: for example, when the principle of excluded middle is part of the theory, making it classical, we are asking for a (collection of) procedure(s) such that, for every formula $\phi$, either an output satisfying $\phi$ or a counterexample to $\phi$ is produced on demand. From Gödel's incompleteness theorem this is not always possible, so, in general, the theory would be computable modulo a number of non-effective oracles, the ones deciding $\phi \vee \neg\phi$.

Moreover, even assuming all axioms to be associated with computable oracles is not enough. In fact, circular definitions are common in the mathematical practice and one has to show that they are well-founded. These definitions introduce new elements in the language which need appropriate decision procedures, i.e., oracles, and nothing in principle prevents them to be uncomputable.

A predicative theory does not admit the usage of these circular definitions and, thus, it shows a better computational character, as its effectiveness is defined exactly by the oracles associated with its axioms.

13.1. **An illustrating example.** Although it is possible to give a completely abstract presentation, it is easier to study a concrete and real case of impredicative definition to understand how the phenomenon impacts on the computational interpretation of a real theory, and the kind of structures impredicativity generates.

The Zermelo-Frænkel set theory, ZF for short, is usually presented as being based on the signature $ZF = \langle\{\mathbf{Set}\}, \emptyset, \{\in : \mathbf{Set} \times \mathbf{Set}\}\rangle$, i.e., a single sort, the one of collections, and a single relation symbol, membership. Nevertheless, the theory is almost immediately extended by a number of definitions, even to introduce its axioms. In the following, we mainly follow the presentation as in the classical text [Kun80].

Among the axioms in ZF, we have the *comprehension scheme*: for each formula $\phi$

(1) $$\forall z. \exists y. y = \{x \in z : \phi\}, \quad \text{with } y \not\in \mathrm{FV}(\phi).$$

It is evident that (1) does not belong to the language of ZF. In fact, a number of definitions is needed to write the axiom scheme, each one introducing a new symbol. In the first place, we introduce equality, $= : \mathbf{Set} \times \mathbf{Set}$ by means of the definition:

(2) $$\begin{aligned} \forall x, y.\ \big(x = y \supset \big(\forall z.\ (z \in x \supset z \in y) \wedge (z \in y \supset z \in x)\big)\big) \wedge \\ \wedge \big(\big(\forall z.\ (z \in x \supset z \in y) \wedge (z \in y \supset z \in x)\big) \supset x = y\big) \quad . \end{aligned}$$

In the following, we will use bi-implication, notation $A \equiv B$, which stands for $(A \supset B) \wedge (B \supset A)$. This is not a new symbol in the logical language, but just an abbreviation. With this notation, we can rewrite (2) as:

(3) $$\forall x, y. x = y \equiv \big(\forall z. z \in x \equiv z \in y\big) \quad .$$

The pattern of (2) is typical:

$$\forall x_1, \ldots, x_n. p(x_1, \ldots, x_n) \equiv \phi$$

where $p$ is the new relation symbol we are introducing, and $\phi$ is the formula in the language deprived of $p$ that the new symbol wants to abbreviate. As $p(t_1,\dots,t_n)$ is logically equivalent to $\phi[t_1/x_1,\dots,t_n/x_n]$, the extension is harmless as whatever we can prove in the new language on the signature $ZF_= = \langle\{\mathbf{Set}\},\emptyset,\{\in,=\colon \mathbf{Set}\times\mathbf{Set}\}\rangle$ can be reduced to a formula in $ZF$ thanks to the axiom (2). It is worth remarking that (2) is indeed an axiom whose purpose is to fix the only possible interpretation of the = symbol.

We leave to the reader to check that, fixed a logically distributive category $\langle\mathbb{C},M\rangle$ and an interpretation $M_{\mathrm{Ax}}$ for the axioms of ZF over the signature $ZF$, we can always find a new $ZF_=$-structure $\langle\mathbb{C}',M',M'_{\mathrm{Ax}}\rangle$ in a canonical way so that $M'$ and $M'_{\mathrm{Ax}}$ restricted to $ZF$ reduce to $M$ and $M_{\mathrm{Ax}}$, respectively, and $M'$ interprets $t_1=t_2$ as an object isomorphic to $M(\forall z.\, z\in t_1\equiv z\in t_2)$, for each $t_1,t_2\in \mathrm{LTerms}(ZF_=)$, with the proviso that $z$ does not occur in $t_1$ and $t_2$. Of course, the converse holds, too: given $\langle\mathbb{C}',M',M'_{\mathrm{Ax}}\rangle$, a $ZF_=$-structure, we can always find a $ZF$-structure $\langle\mathbb{C},M,M_{\mathrm{Ax}}\rangle$ for the theory in the $ZF$-language, simply by restricting $M'$ and $M'_{\mathrm{Ax}}$ to the $ZF$ signature. Evidently, in both cases, valid formulae in the ZF theory are still valid in the extended interpretation, and vice versa.

From a computational point of view, the oracle associated to each new atomic formula of $ZF_=$ is the same as the one associated to the *definiens* part of (2): so, if we prove $t_1=t_2$, the associated algorithm is obtained from an appropriate composition, given by the proof, of the oracles associated with the axioms of ZF, since the role of (2) is literally to convert each instance of equality into a formula in $ZF$.

Although all the previous remarks are well-known and they are generally considered as trivial, it is important to understand them with respect to the computational interpretation of proofs and formulae: essentially, as one expects, equality in $ZF_=$ is just an abbreviation so it does not change in any way the computational character of the theory.

It is important to notice that (2) is a predicative axiom since no new set is 'constructed': if $x$ and $y$ do exist, all the ingredients to check whether $x=y$ are already in place. On the contrary, (1) is impredicative: given an already constructed $z$, we say that there is a newly constructed set $y$ satisfying the rule implicit in $\{x\in x\colon \phi\}$. But $z$ ranges over all sets, even the value of $y$, so a circular application of (1) is possible, in principle.

Since it is easy to prove that equality is transitive, it immediately follows that the $y$ in (1) is unique when $\{x\in z\colon \phi\}$ is a term. And, of course, we need an axiom to define these terms in a language extension.

The pattern to introduce a new function symbol $f$ is

$$\forall x_1,\dots,x_n.\,(\exists! y.\,\psi(x_1,\dots,x_n,y))\wedge\psi(x_1,\dots,x_n,f(x_1,\dots,x_n))\ .$$

It is clear that such a pattern is impredicative as the universal quantification ranges over the newly constructed element **before** the construction takes place. In a multi-sorted environment this may be avoided when constructing a $y$ which has a different sort from those of $x_1,\dots,x_n$. Also, in a recursive definition, the $y$ is constructed in stages so, although formally impredicative, the definition becomes effective. But, in the general case, impredicativity generates undesired computational effects, and this is the case in the comprehension axiom scheme.

In fact, we want to extend the $ZF_=$ signature to

$$ZF^+ = \langle\{\mathbf{Set}\},\{\{\alpha\in_\colon\phi\}\colon \mathbf{Set}\to\mathbf{Set}\}_{\phi\in\mathrm{LForms}(ZF^\alpha_=)},\{\in,=\colon\mathbf{Set}\times\mathbf{Set}\}\rangle\ ,$$

that is, we want to add a new function symbol for each formula in $ZF_=$, with a distinct variable $\alpha$ which may occur in it. The function symbol uses a fancy syntax, $\{\alpha\in_\colon\phi\}$.

Posing $f(z)=\{\alpha\in z\colon\phi\}$ and $\psi(z,y)\equiv\forall x.\,x\in y\equiv(x\in z\wedge\phi[x/\alpha])$ in the function definition pattern, one gets

$$\forall z.\,(\exists! y.\,\forall x.\,x\in y\equiv x\in z\wedge\phi[x/\alpha])\wedge(\forall x.\,x\in\{\alpha\in z\colon\phi\}\equiv x\in z\wedge\phi[x/\alpha])\ .$$

Splitting the conjunction, with a few manipulations, one gets a pair of axioms: (1) and

(4) $$\forall x,z.\,x\in\{\alpha\in z\colon\phi\}\equiv x\in z\wedge\phi[x/\alpha]\ .$$

The idea behind the definition of a new function symbol is that, whenever we can prove $\Gamma \vdash \theta$, where $\theta \equiv \theta'[\{\alpha \in t : \phi\}/x]$, with $x$ new in the proof, we may immediately prove $\Gamma \vdash \exists x. \theta' \wedge x = \{\alpha \in t : \phi\}$ and such a value of $x$ is unique.

Summarising, when an atomic formula on $ZF^+$ contains an occurrence of $\{\alpha \in t_2 : \phi\}$, we may reduce that formula to an equivalent one in the language of the ZF theory by means of the following rules, justified by the axioms:

(5) $$t_1 \in \{\alpha \in t_2 : \phi\} \rightsquigarrow t_1 \in t_2 \wedge \phi[t_1/\alpha];$$

(6) $$t_1 = \{\alpha \in t_2 : \phi\} \rightsquigarrow \forall z. z \in t_1 \equiv z \in t_2 \wedge \phi[z/\alpha];$$

(7) $$\{\alpha \in t_2 : \phi\} \in t_1 \rightsquigarrow \exists w. w \in t_1 \wedge w = \{\alpha \in t_2 : \phi\}.$$

where $z$ and $w$ are variables.

Apparently, we have been able to eliminate all the occurrences of $\{\alpha \in t_2 : \phi\}$ from the atomic formulae. More in general, as all the reduction arrows are logical equivalences, we may eliminate occurrences of the newly defined function symbol from any formula in $ZF^+$ by means of a recursive procedure which yields an equivalent formula in the base language over $ZF$.

Of course, we may get an alternative version of (1), restricted to the language of ZF:

(8) $$\forall z. \exists y. x \in y \equiv x \in z \wedge \phi[x/\alpha] \quad , \quad y \notin \mathrm{FV}(\phi) \ .$$

As a remark, the above listed reduction rules (5), (6) and (7) are able to eliminate occurrences of the $\{\alpha \in t : \phi\}$ terms even in the case $\phi$ is a formula of the language over $ZF^+$, thus justifying a further extension of the language, as is normally done in set theory. For our purposes, such an extension is useless, so we stick on the language generated by $ZF^+$ as it suffices to show the effect of having an impredicative definition or axiom in the theory.

When we consider a $ZF^+$-structure $\langle \mathbb{C}, M', M'_{\mathrm{Ax}} \rangle$ which is a model for the extended ZF theory, we immediately obtain a $ZF$-model $\langle \mathbb{C}, M, M_{\mathrm{Ax}} \rangle$ for ZF simply by restricting $M'$ and $M'_{\mathrm{Ax}}$ to the $ZF$ signature, thanks to the equivalences justifying the reduction (5), (6) and (7) above. Moreover, it long but easy to prove that, whenever $\Gamma \vdash \psi$ in the extended ZF theory, then $\Gamma^r \vdash \psi^r$ in ZF, where $\Gamma^r$ is the reduction of $\Gamma$ and $\psi^r$ the reduction of $\psi$ according to (5), (6) and (7). In fact, this is the property which allows us to say that the extended ZF theory is a *definitional* extension of ZF: apart denotation, the deductive power of the extended theory is exactly the same as that of ZF.

In the case of (2) the converse also holds: from any model of ZF, we can build an essentially equivalent model over $ZF_=$ in a canonical way.

But when we interpret $\{\alpha \in t_2 : \phi\} \in t_1$, we get a surprise due to the impredicative nature of the definition of comprehension. In fact, $M'\left(\{\alpha \in t_2 : \phi\} \in t_1\right)$ is isomorphic to $M'\left(\exists w. w \in t_1 \wedge w = \{\alpha \in t_2 : \phi\}\right)$, which is the initial object of the corresponding existential subcategory, thus it generates the co-cone:

$$\begin{array}{c} M'\left(\exists w. w \in t_1 \wedge w = \{\alpha \in t_2 : \phi\}\right) \\ \uparrow \cdots \uparrow \\ \left\{M'(s \in t_1) \times M'\left(s = \{\alpha \in t_2 : \phi\}\right)\right\}_{s \in \mathrm{LTerms}(ZF^+)} \end{array}$$

in any logically distributive category over $ZF^+$. So,

$$M'\left(\{\alpha \in t_2 : \phi\} \in t_1\right) \times M'\left(\{\alpha \in t_2 : \phi\} = \{\alpha \in t_2 : \phi\}\right)$$

is a possible instance, when $s = \{\alpha \in t_2 : \phi\}$, a perfectly legitimate term. Thus, there is an evident circularity if we want to canonically extend a $ZF$-model $\langle \mathbb{C}, M, M_{\mathrm{Ax}} \rangle$ to $ZF^+$ by reversing the reductions (5), (6) and (7), i.e., by posing $M'\left(\{\alpha \in t_2 : \phi\} \in t_1\right) \cong M'\left(\exists w. w \in t_1 \wedge w = \{\alpha \in t_2 : \phi\}\right)$. Moreover, we may semantically justify $\{\alpha \in t_2 : \phi\} \in t_1$ either by conjunction elimination or by exists introduction from the very same instance.

Nevertheless, in a $ZF$-model $\langle \mathbb{C}, M, M_{\mathrm{Ax}}\rangle$, $M\left(\left(\{\alpha \in t_2 : \phi\} \in t_1\right)^r\right)$ has a meaning, that is, it is semantically justified by being the initial element of the corresponding existential subcategory, so it generates the co-cone:

$$\begin{array}{c} M\left(\exists w.\,(w \in t_1)^r \wedge \left(w = \{\alpha \in t_2 : \phi\}\right)^r\right) \\ \uparrow \cdots \uparrow \\ \left\{M\left((s \in t_1)^r\right) \times M\left(\left(s = \{\alpha \in t_2 : \phi\}\right)^r\right)\right\}_{s \in \mathrm{LTerms}(ZF)} \end{array}$$

Here, the main difference from the above co-cone is not the presence of reductions, but the range of possible instances, as no value of $s$ justifies the top formula by conjunction elimination. This fact shows how the impredicative nature of the definition of the new function symbols introduces a number of ideal elements which modify the structure of the semantics in a subtle way.

The ideal elements are the $\{\alpha \in t_2 : \phi\}$ terms above in a pointwise approach, and the $\{\alpha \in t_2 : \phi\} \in t_1$ atomic formulae in the point-free approach of logically distributive categories. In a pointwise setting, the interpretation of these ideal elements poses a problem: how can one say that these sets do exist? Precisely, since their construction depends on a circular definition, are we allowed to say this is a legitimate construction?

In a point-free setting, the interpretation of the corresponding ideal elements is not a question of existence, as terms do not have a universe where they are interpreted. It raises a question of structure, which is twofold: either we interpret these ideal elements as proper language extensions, or we interpret them as convenient shorthands to write complex formulae—a meta-linguistic feature introduced to second the purposes of conciseness and human readability. In the former case, our analysis has shown that new structure arises in the interpretation, allowing for detours in the justifications which are not reducible to normalisation of formulae, as one can easily see by considering just canonical models as for definition 6.2 with respect to note 3.1. In the latter case, impredicativity does not directly rise problems as it is reduced to 'syntactical sugar' which has no meaning in the formal sense. Better, $\{\alpha \in t_2 : \phi\}$ is **not** a set, but just a super-structure with no role in the semantics: it is just a convenient notation to provide the illusion to speak about an entity which is completely non-representable in the language in any direct way, and thus has no need to be constructed—and whose existence is not justified by the mathematical content of the theory, as we capture it essentially by substitution, and it cannot be used to this purpose.

From a computational point of view, that is, looking at the logical distributive categories as models for the $\lambda$-calculus of section 3, it is clear that the former alternative above introduces an asymmetry: restricting our view to the classifying categories, see proposition 6.13, the one on the extended language admits reductions, computational paths or proofs, if you prefer, that are redundant and which have no distinct correspondent in the model of the non-extended language. In other words, the map of arrows from the classifying model of the extended language to the classifying model of the base language is not injective—which is rather expected— but it reduces structurally different proofs, differing not just for the presence of the new terms, to the very same proof—which is rather unexpected and a direct consequence of impredicativity. Rather obviously, the latter alternative does not change anything in the computational interpretation, as the reference model remains unchanged; but the notation has no more direct correspondence in the computational rendering of a theory, thus no computational property of the new 'elements' can be canonically inferred.

In conclusion, impredicativity manifests itself in a point-free semantics like the one based on logically distributive categories as a phenomenon that introduces unnecessary detours in the semantical explanation of formulae in the extended language, and, thus, gratuitous inefficiencies in the computational interpretation, beyond those already present in the non-extended language.

Although this fact has been shown in a specific example, it is completely general and it is easy for the reader to construct other examples in different theories. An abstract treatment it possible but not yet developed in the full generality: we preferred a presentation based on single and concrete example as the effort required to understand the general and abstract treatment is overwhelming, at least for the author, without bringing much more than what we have shown.

## 14. Some philosophical remarks

In the previous sections, the technical and mathematical aspects of logically distributive categories have been introduced, studied and discussed. Along with well-known topics, to make the presentation self-contained, a number of novel concepts and results have been presented: a point-free categorical semantics for first-order intuitionistic theories has been defined and proved to be sound and complete, both as logical and computational systems, in the form of $\lambda$-calculi. The relation of this semantics with other interpretations for the same class of theories has been discussed. Finally, the role of (im)predicative theories has been analysed.

All these results together provide a somewhat complete picture of the *view* behind the technical definition of logically distributive categories and the semantics on them.

In this section, we want to discuss the above mentioned view in some aspects that may be of interest for the readers having an attention to the philosophy of mathematics.

This section does not aim at being exhaustive or at attaining any degree of depth, rather it wants to illustrate some of the ideas and thoughts that guided the author in the development of the previously illustrated results, and it aims at remarking some points which are not immediately evident in the technical exposition. Because of its purpose and nature, this section is discursive and, to some degree, it reflects a debatable personal view.

14.1. **A point-free approach.** The approach to semantics of logically distributive categories is point-free, as technically explained in section 13. This choice for providing a meaning to a theory may appear as bizarre and counter-intuitive.

In fact, the historical and traditional way of interpreting intuitionistic theories is the Kripke semantics, which models the idea of a class of possible worlds, each one providing a 'local' explanation to formulae in the classical Tarski's sense. These worlds are accessible one to the others in a way that allows to lift the 'local' meanings to the global level of all possible worlds.

Kripke's approach is 'point-wise' since each world provides a universe where terms get interpreted as elements (points). This raises a problem: we have to assume those universes do exist. In classical mathematics, this assumption is not perceived as an obstacle as soon as these universes can be described as sets, which is the case in Kripke's semantics. On the contrary, in constructive mathematics, existence is a more delicate concept: an 'object' exists when it can be constructed, which normally means that there is an effective building procedure; in the case of universes, we should require at least the ability to test membership, to list elements, and to check whether two elements are equal. As soon as we require these procedures to be computable, most usual 'universes', e.g., real numbers, do not exist in this stronger sense.

Traditional semantics based on category theory, as Heyting categories or internal models in topoi, circumvent the problem by abstracting over universes: assuming that sorts can be interpreted as objects in a suitable category, terms and formulae are rendered by arrows and sub-objects, respectively. Sorts are not required to have any internal structure; instead, their properties are derived from their external structure, i.e., how they are related to each other via arrows. So, a term is a *link* among sorts, carrying in its construction the amount of information that gives it a meaning. Similarly, a formula is a sub-object, an arrow in the end, which is interpreted as the sub-part of a collection of sorts that verifies the formula, intuitively. But, again, it is just the formal structure of the formula, and how it gets encoded by the semantics that really convey a meaning. Nevertheless, each object in those categories **is** a 'collection' of arrows representing terms, since this is the use the interpretation makes of them.

On this line of thought, it is possible to say that the traditional categorical models are point-free, as we do not assume the existence of universes as collections of elements to interpret terms in; nevertheless, we assume the existence of universes as abstract entities, having a few structural properties, expressed as 'links' with other, similar entities, which say that they behave *as if* they are collections of elements in a weak sense.

Logically distributive categories take a step further: recognising that the only role of variables is that of being substituted by terms, the objects of the category are used to interpret formulae. The structure of the category is responsible for linking together simple formulae so to form interpretations of complex ones. The distinctive element in the structure of the category is a direct notion of substitution, rendered via the $\Sigma$ functors, which allows to interpret quantifiers. In this sense, terms and variables are out of the scene: terms appear just as collections of indexes, as an effect of the application of $\Sigma$ functors; in particular, terms and variables have no interpretation in the models. Hence, no concrete neither abstract universe is required to exist in order to interpret terms: more simply, terms are **not** interpreted, instead the substitution operation is, in the way universal and existential sub-categories are constructed.

Although it may be disturbing for the working mathematician, our results prove that existence of universes is not necessary, thus it is not a problem of mathematical logic. In a sense, logically distributive categories are moderately agnostic in their epistemological requirements: they do not force to assume the existence of universes — but they allow to, if one wants — and they are less demanding in the power of meta-theory with respect to traditional approaches based on category theory. In fact, the categorical constructions they rely on are finite products and co-products, exponentiation, and the possibility to select subsets of objects on a syntactic basis via the substitution functors. So they do not require finite limits and co-limits, which reduces to say that they do not require the existence of equalisers and co-equalisers, and, more important, they do not require sub-object classifiers, i.e., a strong notion of sub-object. Nevertheless, as one may expect, the ability to select subsets, even in a limited way, is not for free: nothing prevents, in principle, these subsets to be non-constructible, thus the semantics is not 'constructive' in a strict sense.

The choice of developing a point-free semantics has a consequence that may be non immediately evident: we can restrict our attention to models having a denumerable cardinality of objects when dealing with effective theories. In fact, the classifying models have a size which is given by the cardinality of formulae, even when the 'intended universe' is much larger. It should be remarked that this fact is in line with the classical Löwenheim-Skolem theorems, so it is not unexpected. But, in the case of effective theories, this fact means that the 'power' required to capture the needed subsets is limited by the $2^{\aleph_0}$ cardinal — result that puts the meta-theory beyond the constructive threshold, but not too far.

14.2. **Inconsistent theories.** In the classical categorical models, see sections 7, 9, 11, and 12, an inconsistent theory has no model. On the contrary, logically distributives categories provides an interpretation also to inconsistent theories.

At a first glance, this behaviour of the semantics is disturbing because it breaks the link between syntax and semantics which says that a theory is syntactically free from contradictions exactly when it has a model. However, this symmetry is easily reestablished by noticing that a theory is syntactically consistent exactly when it admits at least one non-trivial model, i.e., a model which interprets $\top$ and $\bot$ into two distinct and non-isomorphic objects. In particular, a theory is consistent exactly when the classifying model is non-trivial, as it is immediate to prove.

From a different point of view, it is perfectly natural that inconsistent theories have a model. In fact, when the model provides an account to the computational content of a theory, as is the case with logically distributive categories since they allow to interpret both the logical theory and the associated $\lambda$-calculus, we should recognise the fact that inconsistent theories have a non-trivial computational content, which is expressed in a 'trivial' model.

Although out of the scope of the present work, one may regard an inconsistent theory as the model of the 'world' of a faulty software project, where the fault lies in the overall design rather than in the algorithms or in their implementation. When specifications are contradictory but the designer is not aware of, it may be of interest to model that world, so to understand what actions can be performed to correct the problem, minimising the impact on the already developed code. In this respect, traditional approaches to semantics provide no clues since the 'world' of the program does not exist; on the contrary, a point-free semantics models the faulty 'world', capturing its computational content, and so the actual code, at least in an abstract sense. A possible line of investigation of potential interest to computer scientists would be to study how these inconsistent models can be transformed into consistent ones preserving part of their computational content, which is the mathematical way to rephrase the corrective action one wishes to perform in that situation.

14.3. **Equality.** In the classical, set-theoretic approach to semantics, equality is treated as a special case: instead of being a generic atomic relation, its interpretation is forced to be the diagonal relation in the universe, i.e., two terms are equal exactly when they denote the very same element. Although this interpretation is perfectly reasonable since it captures the intended idea of 'being equal', it induces also a number of problems.

In fact, the diagonal relation is not exactly described by the first-order axioms about equality: in a first-order theory, the equality axioms are satisfied by any congruence relation, not only the diagonal relation. This fact has been explored in the case of Kripke models, showing that the two interpretations are not equivalent: the congruence interpretation is evidently more general, and it allows to prove completeness results for some intermediate logics where the stricter diagonal interpretation fails to.

From a computational point of view, the diagonal relation is a problematic interpretation for the equality symbol: for example, in the theory of real numbers, it is easy to write an algorithm to decide whether two numbers are different — simply considering the usual decimal expansion, if the $n$-th digit is different in the to inputs, we can conclude that the numbers are different as well; of course, if the two inputs are equal, the algorithm does not terminate in finite time, since it requires to check an infinite amount of digits. In our terms, the oracle associated to equality is non-effective (not recursive) in the general case, undermining the possibility to develop a really constructive and predicative theory.

For these reasons, the point-free interpretation based on logically distributive categories relies on axioms and does not reserve a special interpretation to equality. It should be remarked that this is a choice: in fact, it is possible to force the diagonal interpretation via a suitable quotient of terms, which reflects into the categorical structure via the $\Sigma$ functors. Such a quotient reveals how much force one has to exercise to induce the stricter interpretation. It should be clear that the 'design' choice of a free interpretation of equality has the consequence of a better control over the computational properties of logical theories, at the price of a strictly constructive attitude which results somehow 'built-in' the semantics.

14.4. **Ideal elements.** It has been shown in section 13 how impredicative definitions introduce ideal elements. Precisely, the point-free approach clarifies how the notation a definition introduces is interpreted in a redundant computational structure which is justified by a reduction rule to a simpler structure. Of course, this modus operandi poses a number of problems, mainly related to circularity, which have been illustrated before.

Ideal elements are 'fictitious' in the sense that they are understood via a set of suitable reductions to the base language arising from their definitions. When considering ideal elements as first-class components of the language, that amounts to relate the classifying models over the base and the extended language, we loose the direct correspondence between these models, as the ideal elements induce detours and structural redundancies, i.e., arrows, which have no correspondent in the base model.

Therefore, if our aim in giving a definition is to single out an element of interest in a theory, we have to admit that, in general, we are unable to achieve the desired result since the new element so defined is not completely described by its definition. This fact is due to our way of interpreting: since our semantics does not provide a meaning for terms, 'elements' if you prefer, but it uses terms as a glue to link together propositions, definitions capture or constrain the interpretation by means of the way terms relate formulae to each other via substitution. When this 'external' description of terms is not enough to fix the meaning of a term, i.e., how it determines the computational behaviour of formulae and proofs, in a way that avoids the introduction of new pieces of information in the semantics, we have an asymmetry in the interpretation and, thus, an ideal element.

Not every definition introduces ideal elements: in fact, predicative definitions do not. They preserve the computational structure of classifying models, thus providing a strict correspondence between base and extended languages that allows to transfer both logical and computational content between the two classifying models in both directions. This is the reason why predicative definitions are mostly relevant and interesting: they do preserve computational content and, we claim, this is their distinctive character.

So, if we want to draw a conclusion from this research effort, we can say that is the interplay between logic and computation that guides the degree of *constructiveness* of a theory, and this fact is best understood and studied in a point-free perspective.

DIPARTIMENTO DI SCIENZE TEORICHE E APPLICATE, UNIVERSITÀ DEGLI STUDI DELL'INSUBRIA, VIA MAZZINI 5, I-21100 VARESE (VA), ITALY

*E-mail address*: `marco.benini@uninsubria.it`

*URL*: `http://marcobenini.wordpress.com/`